\documentclass[12pt,leqno]{article}
\title{ A classical approach on cyclotomic fields and Fermat-Wiles theorem - Version 1.2}
\author{Roland Qu\^eme}
\usepackage{amsmath,amsthm,amstext,amsbsy,amsfonts}
\newtheorem{thm}{Theorem}[section]
\newtheorem{cor}[thm]{Corollary}
\newtheorem{prop}{Proposition}
\newtheorem{lem}[thm]{Lemma}
\newcommand{\N}{\mathbb{N}}
\newcommand{\R}{\mathbb{R}}
\newcommand{\C}{\mathbb{C}}
\newcommand{\Q}{\mathbb{Q}}
\newcommand{\Z}{\mathbb{Z}}
\newcommand{\s}{\mathbf{s}}
\newcommand{\modu}{\ \mbox{mod}\ }
\newcommand{\be}{\begin{equation}}
\newcommand{\ee}{\end{equation}}
\date{2002 october 31}
\begin{document}
\tableofcontents
\clearpage
%
%
\abstract $ $
Roland Qu\^eme

13 avenue du ch\^ateau d'eau

31490 Brax

France

2002 october 31

mailto: roland.queme@free.fr

home page: http://roland.queme.free.fr/index.html

**********************************************
\clearpage
\begin{itemize}
\item
This  monograph extends two  previous works entitled:
\begin{itemize}
\item
{ \it A generalization of  Eichler's criterion for Fermat's Last Theorem }
\item
{ \it A classical approach on Fermat-Wiles theorem}
\end{itemize}
\item
This monograph contains  a first part on \underline{\bf cyclotomic fields} dealing with class group, unit group and Hilbert class field of cyclotomic fields, independantly of Fermat-Wiles Theorem.

\underline{\bf Class group of cyclotomic field $\Q(\zeta)$}
\begin{itemize}
\item
This topic is studied in sections \ref{s09121} p. \pageref{s09121}.
Let $p\in\N,\quad p>5$ be a prime. Let $\zeta$ be a root of the equation $X^{p-1}+X^{p-2}+\dots+X+1=0$.
Let $\Q(\zeta)$ be the $p$-cyclotomic field and $\Z[\zeta]$ be the ring of integers of $\Q(\zeta)$.
Let $C_p, C_p^+,C_p^-$ be respectively the $p$-class group of $\Q(\zeta)$, the $p$-class group of
$\Q(\zeta+\zeta^{-1})$ and the relative $p$-class group $C_p^-=C_p/C_p^+$.
Let $C_p^1\subset C_p$ be the subgroup of the $p$-class group of $\Q(\zeta)$ whose elements are of order $1$ or $p$. Let $r_p$ be the rank of $C_p^1$. Let us note $C_p^1=\oplus_{i=1}^{r_p} \Gamma_i$, where $\Gamma_i$ is a cyclic group of order $p$. Let $u\in\N$ be a primitive root $\modu p$. Let $\sigma:\zeta\rightarrow\zeta^u$ be a $\Q$-isomorphism of $\Q(\zeta)$. We know that $\sigma$ generates $G=Gal(\Q(\zeta)/\Q)$.
Let $\mathbf b_i,\quad i=1,\dots,r_p$, be  a not principal ideal of $\Q(\zeta)$ whose class belongs to  group $\Gamma_i$. Observe at first that $\sigma(\mathbf b_i)\simeq \mathbf b_i^{\mu_i}$ where $\simeq$ is notation for class equivalence and $\mu_i\in{\bf F}_p^*$ and  where ${\bf F}_p^*$ is the set of no null elements of  the finite field of cardinal $p$.
Let the ideal $\mathbf b=\prod_{i=1}^{r_p}\mathbf b_i$.  Let $d\in\N,\quad p-1\equiv 0\modu d$. Let us define the minimal polynomial $P_{r_d}(U)$ of degree $r_d$ in the indeterminate $U$ such that
$P_{r_d}(\sigma^d)$ annihilates the ideal class of $\mathbf b$ written also
$\mathbf b^{P_{r_d}(\sigma^d)}\simeq \Z[\zeta]$. We show that the polynomial $P_{r_d}(U)$ verifies a relation of form
$P_{r_d}(U)=\prod_{i=1}^{r_d}(U-\mu_i^d),\quad \mu_i\in{\bf F}_p^*$,
and then we obtain the following results:
\begin{enumerate}
\item
For $d=1$ then  $r_1$ is the index of irregularity of $\Q(\zeta)$ (the number of even Bernoulli Numbers $B_{p-1-2m}\equiv 0 \modu p$ for $1\leq m\leq \frac{p-3}{2})$.
\item
For $Gcd(d_1,d_2)=1$ then $r_{d_1}\times r_{d_2}\geq r_1$.
\item
Let the ideal $\pi=(\zeta-1)\Z[\zeta]$.
There exists $B_i\in\Z[\zeta],\quad i=1,\dots,r_p$, with $B_i\Z[\zeta]=\mathbf b_i^{n_i p},\quad n_i
\in{\bf F}_p^*,\quad \pi|(B_i-1)$ and with
 $\sigma(B_i)=B_i^{\mu_i}\times\alpha_i^p,\quad \alpha_i\in\Q(\zeta),\quad \mu_i\in{\bf F}_p^*$.
For $\mu_i=u^{m_i},\quad m_i\in\N,\quad 1\leq m_i\leq p-2$, then $\pi^{m_i}| {B_i-1}$.
\item
Let $r_p^+$ be the rank of the $p$-class group of $\Q(\zeta+\zeta^{-1})$.
We can precise the previous result: with a certain reordering of indexing of $B_i,\quad i=1,\dots,r_p$,
\begin{itemize}
\item
For $i=1,\dots,r_p^+$  then the $B_i$ are singular primary, so
$\pi^p |(B_i-1)$,
\item
for $i=r_p^++1,\dots,r_p$ then $\pi^{m_i}\| (B_i-1)$.
\end{itemize}
\item
$\mu_i=u^{2m_i+1}$ with $1\leq m_i\leq \frac{p-3}{2}$  corresponds to an ideal $\mathbf b_i$ whose class belongs to $C_p^-$, relative $p$-class group of $\Q(\zeta)$. In that case define
$C_i=\frac{B_i}{\overline{B}_i}$, so with $C_i\in\Q(\zeta)$. If $2m_i+1>\frac{p-1}{2}$ then it is possible to prove the explicit very simple  formula
for $C_i$:
\begin{displaymath}
C_i\equiv 1-\frac{\gamma_i}{1-\mu_i} \times
(\zeta+\mu_i^{-1}\zeta^u+\dots+\mu_i^{-(p-2)}\zeta^{u^{p-2}})\modu \pi^{p-1},\quad \gamma_i\in{\bf F}_p.
\end{displaymath}
\end{enumerate}

\underline
{\bf  unit group $F=\Z[\zeta+\zeta^{-1}]^*/(\Z[\zeta+\zeta^{-1}]^*)^p$.}
\item
This topic is studied in section \ref{s210163} p. \pageref{s210163}.
We apply in following results to unit group $\Z[\zeta+\zeta^{-1}]^*$ the method applied to $p$-class group in previous results:
\begin{enumerate}
\item
There exists a fundamental system of units $\eta_i,\quad i=1,\dots,\frac{p-3}{2}$, of the group
$F=\Z[\zeta+\zeta^{-1}]^*/(\Z[\zeta+\zeta^{-1}]^*)^p$ verifying the relations:
\begin{equation}\label{e201274}
\begin{split}
& \eta_i\in\Z[\zeta+\zeta^{-1}]^*,\quad i=1,\dots,\frac{p-3}{2},\\
& \sigma(\eta_i)=\eta_i^{\mu_i}\times\varepsilon_i^p,\quad\mu_i\in{\bf F}_p^*,
\quad \varepsilon_i\in\Z[\zeta+\zeta^{-1}]^*,
\quad i=1,\dots,\frac{p-3}{2},\\
& \eta_i\equiv 1\modu \pi^{2m_i},\quad \mu_i=u^{2m_i},\quad i=1,\dots,
\frac{p-3}{2},\\
& \sigma(\eta_i)\equiv \eta_i^{\mu_i}\modu\pi^{p+1},\quad i=1,\dots,\frac{p-3}{2}.
\end{split}
\end{equation}
\item
With a certain reordering of indexing of $i=1,\dots,\frac{p-3}{2}$, then
\begin{itemize}
\item
for $i=1,\dots,r_p$, then $\eta_i$ are primary, so $\pi^p| (\eta_i-1)$,
\item
for $i=r_p+1,\dots,\frac{p-3}{2}$, then $\pi^{2m_i}\|(\eta_i-1)$.
\end{itemize}
\item
If $2m_i>\frac{p-1}{2}$ then it is possible to prove the very simple  explicit formula
for $\eta_i$:
\begin{displaymath}
\eta_i\equiv 1-\frac{\gamma_i}{1-\mu_i} \times
(\zeta+\mu_i^{-1}\zeta^u+\dots+\mu_i^{-(p-2)}\zeta^{u^{p-2}})\modu \pi^{p-1},\quad \gamma_i\in{\bf F}_p.
\end{displaymath}
\end{enumerate}

\underline{\bf $p$-elementary Hilbert class field of $\Q(\zeta)$.}
\item
Section \ref{s110021} p. \pageref{s110021} deals of Hilbert class field :  we recall some results on Hilbert class field and on Furtwangler and Hecke Theorem of interest in study of Fermat's equation.
In section \ref{s207272} p. \pageref{s207272}  we study some relations between Hilbert class field of $\Q(\zeta)$ and of subfield $\Q(\zeta+\zeta^{-1})$.
\item
We call $p$-elementary Hilbert class field $H$ of $\Q(\zeta)$ the subfield  of the Hilbert class field of
$\Q(\zeta)$ verifying by Artin map $Gal(H/\Q(\zeta))\cong C_p^1$, where the group $C_p^1$ has been defined above.
The following result is representative of our  results on Hilbert class field:
\item
With a certain reordering  of index $i=1,\dots,r_p$, and with the same meaning of $B_i,\quad i=1,\dots,r_p$, than in previous results,
the $p$-elementary Hilbert class field verifies the  structure:
\begin{displaymath}
H=\Q(\zeta,\omega_1,\dots,\omega_{r_p}),\quad \omega_i^p\in\Q(\zeta),\quad \omega_i\not\in\Q(\zeta),
\quad i=1,\dots,r_p,
\end{displaymath}
\begin{itemize}
\item
where, for $i=1,\dots,r_p$, the $\omega_i^p$ are primary,
\item
where, for $i=1,\dots,r_p^+$,
\begin{displaymath}
\begin{split}
& \omega_i^p= B_i,\quad B_i\in\Z[\zeta],\quad B_i\Z[\zeta]
=\mathbf b_i^p,\quad Cl(\mathbf b_i)\in C_p^-,\\
&  B_i\overline{B}_i=\beta_i^p,\quad \beta_i\in\Z[\zeta+\zeta^{-1}]-\Z[\zeta+\zeta^{-1}]^*,\\
\end{split}
\end{displaymath}
\item
where, for  $i=r_p^++1,\dots,r_p$, the $\omega_i^p \in \Z[\zeta+\zeta^{-1}]^*$.
\end{itemize}
\underline{\bf Connection of Mirimanoff polynomials with  $p$-class group of $\Q(\zeta)$.}

Recall that, for $m\in\N,\quad 1\leq m\leq\frac{p-3}{2}$, the polynomials
$\phi_{2m+1}(T)=\sum_{i=1}^{p-1} i^{2m}\times T^i$ are the odd Mirimanoff polynomials of the indeterminate $T$.
Suppose that $p$ divides  the class number $h$ of $\Q(\zeta)$.
Let $r_p^->0$ be the relative $p$-class group $C_p^-$.
We have proved that,
for $i=1,\dots,r_p^-$, there exist   different natural integers $m_i\in\N,\quad 1\leq m_i\leq\frac{p-3}{2}$, and ideals $\mathbf b_i$ of $\Z[\zeta],\quad Cl(\mathbf b_i)\in C_p^-,\quad
\sigma(\mathbf b_i)\simeq\mathbf b_i^{\mu_i},\quad \mu_i=u^{2m_i+1},$
and   natural integers $b_i, \quad 1<b_i\leq p-1$,  verifying Mirimanoff polynomial congruences following
\begin{displaymath}
\phi_{2m_i+1}(b_i)\times\phi_{p-2m_i-1}(b_i)\equiv 0\modu p,\quad i=1,\dots,r_p^-.
\end{displaymath}
\end{itemize}
\item
This monograph contains a second part  dealing with \underline{\bf Fermat-Wiles} theorem on a classical approach.
The Wiles's proof of Fermat's Last Theorem lets open the question to know if there exists also a proof of first case of Fermat's Last Theorem resting on Classical Algebraic Number Theory. In this context, this monograph deals with several researchs on Fermat's equation  lying on:
\begin{itemize}
\item
some elementary  properties of the ring of integers of the $p$-cyclotomic number field $\Q(\zeta)$,
\item
on the representations of the Galois group $G$ of $\Q(\zeta)$ over the finite field  ${\bf F}_p$, with the  $p$-class group of $\Q(\zeta)$ seen as  a ${\bf F}_p[G]$ module,
\item
on  the Hilbert class field of $\Q(\zeta)$.
\end{itemize}
In this second part:
\begin{itemize}
\item
Let $p>5$ be a prime. Let $\zeta$ be a primitive $p$-root of unity. Let $\Q(\zeta)$ be the $p$-cyclotomic number field. Let $h$ be the class number of $\Q(\zeta)$. Let $e_p\in\N$ defined by $h=p^{e_p}\times h_2,\quad h_2\in \N, \quad h_2\not\equiv 0 \modu p$.
Let $r_p$ be the $p$-rank of the class group of $\Q(\zeta)$.
Let $x,y,z\in\Z-\{0\}$ be mutually co-prime.
Let us assume that $x,y,z$ verify the Fermat's equation $x^p+y^p+z^p=0$ with $xyz\not\equiv 0 \modu p$ (first case).
Let $t\in\Z,\quad t\equiv-\frac{x}{y}\modu p$.
Let  $\phi_{2m+1}(T)=T+2^{2m}\times T^2+\dots+(p-1)^{2m}\times T^{p-1},\quad m=1,\dots,\frac{p-3}{2},$ be the odd Mirimanoff polynomials of the indeterminate $T$.
In section \ref{s4} p. \pageref{s4}, we give   an {\bf elementary} proof that, if $r_p<\frac{p-1}{2}$, then
\begin{displaymath}
Card\{m\quad |\quad 1\leq m\leq \frac{p-3}{2},\quad \phi_{2m+1}(t)\equiv 0\modu p\}
\geq \frac{p-3}{2}-r_p.
\end{displaymath}
\item
In sections \ref{s108252} p. \pageref{s108252} and \ref{s14101} p. \pageref{s14101}, we study Mirimanoff's congruences in intermediate fields between $\Q$ and $\Q(\zeta)$ and obtain strong generalizations on Mirimanoff polynomials. As an example, let $p-1= f\times g$ with $f$ odd. Let $K$ be the intermediate field
$\Q\subset K\subset \Q(\zeta)$ and $[K:\Q]=g$. Suppose that $p$ does not divide the class number of $K/\Q$. Then we have the Mirimanoff's polynomials congruences
$\modu p$
\begin{displaymath}
\phi_{f(2n+1)}(t)\equiv 0\modu p,\quad n=0,1,\dots,\frac{g-2}{2}.
\end{displaymath}
\item
Section \ref{s18121} p. \pageref{s18121} deals, in a strictly elementary appoach, of a set of explicit polynomials congruences
$P(t)\equiv 0\modu p,\quad P(T)\in {\bf F}_p[T],
\quad t\equiv -\frac{x}{y}$; we explain  why it is legitimous to conjecture that all these congruences are not simultaneously possible.
As an example, we prove the particular case of Fermat (first and second case): $p\|(x-y)\times(y-z)\times (z-x)$ implies that
$x^p+y^p+z^p\not=0$.
\item
Section \ref{s103112} p. \pageref{s103112} generalizes to polynomial congruences
$\phi_m(t)\equiv 0\modu p,\quad 1\leq m\leq p-1$, some results obtained in previous sections for congruences
$\phi_{2m+1}(t)\equiv 0 \modu p,\quad 1\leq m \leq \frac{p-3}{2}$.
\item
In section \ref{seckum} p. \pageref{seckum}, we  give new {\it summation criteria} depending only on $p$ for the first case of FLT with an extension of methods of nineteenth century resting on Kummer's result on Jacobi resolvents and on Stickelberger theorem.
\begin{itemize}
\item
the first known example seen in the litterature being the Cauchy criterion : if   $\sum_{j=1}^{(p-1)/2} \frac{1}{j^3}\not\equiv 0 \modu p$, then the first case of FLT holds for $p$). Different other summation criteria seen in the literature are  Vandiver, Schwindt, Yamada, Ribenboim, Skula, Ciranek, Dilcher-Skula and others.
\item
In this monograph, we give several new summation criteria. As an example of results obtained, we give the $(p-2)(p-3)$ new summation criteria:
for each $d_1,d_2\in\N, \quad 1\leq d_1<d_2 \leq p-2$,
if $\sum_{l=1}^{d_1}\sum_{m=1}^{d_2}
\sum_{j\in[\frac{lp}{d_1+1},\frac{lp}{d_1}]\cap[\frac{mp}{d_2+1},\frac{mp}{d_2}]}
\frac{1}{j^3}\not\equiv 0 \modu p$,
then the first case of FLT holds for $p$.
\end{itemize}
\item
In section \ref{s110022} p. \pageref{s110022} dedicated to second case, we prove, as an example, the second case of FLT when $p\parallel y$, with  a strictly  elementary method in $\Z$ resting on Barlow-Abel formulas proved in first part on nineteenth century. We give also {\it canonical} forms of the Fermat's equation resting on classical theory of $\Q(\zeta)$  when the class number $h^+$ of the field $\Q(\zeta+\zeta^{-1})$ is divisible or not divisible by $p$.
\item
In this monograph, it is also important to note, as explained in section \ref{s08081} p. \pageref{s08081}, that almost our results on the First case of FLT  can be applied to a wider set of diophantine equations than the Fermat equation, in fact for a prime $p\in \N, \quad p>5$,
\begin{itemize}
\item
all the equations of the form:
\begin{equation}\label{e05031a}
\frac{x^p+y^p}{x+y}=t_1^p,\quad x,y, t_1\in\Z, \quad xy(x^2-y^2)\not\equiv 0 \modu p.
\end{equation}
Observe that
$xy(x^2-y^2)\not\equiv 0\modu p $ is hypothesis of Terjanian conjecture for equation (\ref{e05031a})  to have only a trivial Kummer's system of congruences.
\item
All the equations of the form:
\begin{equation}\label{e05032a}
x^p+y^p+c\times z^p=0, \quad x,y,z,c\in\Z,\quad
xy (x^2-y^2)\not\equiv 0 \modu p,
\end{equation}
where all the prime factors $q$ of $c$  verify $q\not\equiv 1 \modu p$, belong to this category.
\end{itemize}
\item
For a first quick outlook of first case of Fermat-Wiles Theorem , see at first the {\bf Table Of Content} of this monograph and the {\bf  Comparative survey with First Case bibliography} section \ref{s19102}
on page \pageref{s19102}, where we compare some of our results with bibliography, with, for each result, the ident number  and page number of the theorem.
\end{itemize}
\end{itemize}
\begin{itemize}
\item
The important results (subjective) are tagged with the symbol ***.  Our main interrogation is to know those of them that Professionnal Number Theorists  feel correct and, among them, those that they feel  not  known in literature.
\item
The {\it Table of Content} completes this abstract to see more precisely  the methods and perimeter of our study of cyclotomic fields and of Fermat-Wiles theorem.
\item
{\bf dependances between sections and reading order}:
\begin{itemize}
\item
sequence:  $1 \rightarrow 2 \rightarrow 3\rightarrow 4\rightarrow 5\rightarrow 6$, cyclotomic fields
\item
sequence:  $1 \rightarrow 7 \rightarrow 8\rightarrow 9\rightarrow 10\rightarrow 11$, FLT
\item
sequence:  $1 \rightarrow 7 \rightarrow 8\rightarrow 12$, FLT
\item
sequence:  $1 \rightarrow 2 \rightarrow 3\rightarrow 7\rightarrow 8\rightarrow 13\rightarrow 14$, FLT
\item
sequence:  $1 \rightarrow 2 \rightarrow 3\rightarrow 4\rightarrow 15$, FLT
\item
sequence:  $1 \rightarrow 7 \rightarrow 8\rightarrow 16$, FLT
\end{itemize}
\end{itemize}
\endabstract
\maketitle
%

%
%
\section*{PART ONE: ON CYCLOTOMIC FIELDS}
This first part of the monograph deals of cyclotomic field $\Q(\zeta_p)$ for $p$ prime, without any connection with Fermat-Wiles Theorem. It contains:
\begin{itemize}
\item
A study of the structure of the $p$-class group of the field $\Q(\zeta)$.
\item
A study of the structure of the unit group $\Z[\zeta+\zeta^{-1}]^*$.
\item
A study of the structure of the Hilbert class field of $\Q(\zeta)$.
\item
A study of the structure of the Hilbert class field of $\Q(\zeta+\zeta^{-1})$.
\end{itemize}
%
%
\clearpage
\section{Some definitions}
In this section, we fix notations used in all this monograph.
\begin{itemize}
\item
For $a\in\R^+$, we note $[a]$ the integer part of $a$ or the integer immediately below $a$.
\item
We denote $[a,b],\quad a,b\in \R$, the closed interval bounded by $a,b$.
\item
Let us  denote $< a > $ the cyclic group generated by the element $a$.
\item
Let $p\in\N,\quad p > 5,$ be a prime.
\item
Let $u\in\N$ be a primitive root $\modu p$.
\item
For $i\in\N$, let us denote $u_i\equiv u^i \modu p,\quad 1\leq u_i\leq p-1$.
For $i\in\Z,\quad i<0$, this is to be understood as $u_i u^{-i}\equiv 1 \modu p$.
This notation follows the convention adopted in Ribenboim \cite{rib}, last paragraph of page 118.
This notation is largely used in the sequel of this monograph.
\item
Let $\Q(\zeta)$ be the $p$-cyclotomic number field.
\item
Let $\Z[\zeta]$ be the ring of integers of $\Q(\zeta)$.
\item
Let $\Z[\zeta]^*$ be the group of units of $\Z[\zeta]$.
\item
Let $\Q(\zeta+\zeta^{-1})$ be the maximal real subfield of $\Q(\zeta)$, with
$[\Q(\zeta):\Q(\zeta+\zeta^{-1})]=2$. The ring of integers of $\Q(\zeta+\zeta^{-1})$ is $\Z[\zeta+\zeta^{-1}]$. Let $\Z[\zeta+\zeta^{-1}]^*$ be the group of units of
$\Z[\zeta+\zeta^{-1}]$.
\item
Let ${\bf F}_p$ be the finite field with $p$ elements.
Let ${\bf F}_p^*={\bf F}_p-\{0\}$.
\item
Let us denote $\mathbf a$ the integral ideals of $\Z[\zeta]$. Let us note
$\mathbf a\simeq \mathbf b$ when the two ideals $\mathbf a$ and $\mathbf b$ are in the same class of the class group of $\Q(\zeta)$. $\mathbf a\simeq\Z[\zeta]$ means that  the ideal $\mathbf a$ is principal.
\item
Let us  note $Cl(\mathbf a)$ the class of the ideal $\mathbf a$ in the class group of $\Q(\zeta)$.
$<Cl(\mathbf a)>$ is the finite group generated by the class $Cl(\mathbf a)$.
\item
If $a\in\Z[\zeta]$, we note $a\Z[\zeta]$ a principal integral ideal of $\Z[\zeta]$.
\item
We have $p\Z[\zeta]=\pi^{p-1}$ where $\pi$ is the principal prime ideal $(1-\zeta)\Z[\zeta]$.
\item
Let us denote $\lambda= \zeta-1$.
\item
Let us denote $\tau : \Q(\zeta) : \rightarrow \Q(\zeta)$, the
$\Q$-isomorphism of $\Q(\zeta)$ defined
by $\tau(\zeta)=\zeta^{-1}$, (complex conjugation). For $a\in \Q(\zeta)$,
we note also $\overline{a}=\tau(a)$.
\item
For $i=1,\dots,p-1$, let us denote $\tau_i:\Q(\zeta) \rightarrow \Q(\zeta)$, the
$\Q$-isomorphism of $\Q(\zeta)$ defined
by $\tau_i(\zeta)=\zeta^{i} $.
\item
Let $G=Gal(\Q(\zeta/\Q)$ be the Galois group of the field $\Q(\zeta)$.
\item
Let $\sigma :\Q(\zeta)\rightarrow \Q(\zeta)$ be the $\Q(\zeta)$-isomorphism defined by $\sigma(\zeta)=\zeta^u$ where $u$ is a primitive root $\modu p$. The $\Q$-isomorphism $\sigma$ generates the Galois group $G=Gal(\Q(\zeta)/\Q)$.
\item
For $d\in\N,\quad p-1\equiv 0\modu d$, let $G_d$ be the cyclic subgroup of $G$ generated by $\sigma^d$, so with $G_1=G$.
\item Let $C_p$ be the $p$-class group of the field $\Q(\zeta)$, the subgroup of the class group whose elements are of $p$-power order.
\item
Let $C_p^+$ be the $p$-class group of the field $\Q(\zeta+\zeta^{-1})$.
\item
Let $C_p^-$ be the relative class group defined by
$C_p^-=C_p/C_p^+$.
\item
Let  $h$ be the class number of $\Q(\zeta)$. The class number $h$ verifies the formula $h=h^-\times h^+$, where $h^+$ is the class number of the maximal real field $\Q(\zeta+\zeta^{-1})$, so called also second factor, and $h^-$  is the relative class number, so called first factor.
\item
Let us define $e_p$ by $h=p^{e_p}\times h_2,\quad h_2\not \equiv 0 \modu p$.
Let us define $e_p^-$ by $h^-=p^{e_p^-}\times h_2^-,\quad h_2^-
\not \equiv 0 \modu p$.
\item
Let $r_p, r_p^+, r_p^-,$ be respectively  the $p$-rank of the class group of $\Q(\zeta)$, of the class group of $\Q(\zeta+\zeta^{-1})$ and of the relative class group.
\item
The abelian group $C_p$ of order $p^{e_p}$ is $C_p=\oplus_{i=1}^{r_p} C_i$ where $C_i,\quad i=1,\dots,r_p$, are cyclic  groups of order $p^{e_i}$ with $e_p=e_1+\dots+e_{r_p}$. Let $C_p^1$ be the subgroup of $C_p$ whose elements are of order $1$ or $p$. The group $C_p^1$ is a group of order $p^{r_p}$ with
$C_p^1=\oplus_{i=1}^{r_p} C_i^\prime$ where $C_i^\prime$ are cyclic group of order $p$. The group $C_p^1$ plays a central role in this FLT study.
\item
Let $C_p^{1+}$ be the subgroup of $C_p^+$ whose elements are of order $1$ or $p$.
\item
Let $C_p^{1-}$ be the subgroup of $C_p^-$ whose elements are of order $1$ or $p$.
\item
Let us note $\phi_m(T)=\sum_{i=1}^{p-1}i^{m-1}\times T^i$
the Mirimanoff polynomials of the indeterminate $T$.
\end{itemize}
%
\clearpage
\section{On structure of $p$-subgroup of the class group of $\Q(\zeta)$}\label{s09121}
The three first subsections \ref{s20111} p.\pageref{s20111}, \ref{s111171} p. \pageref{s111171} and \ref{s20112} p.\pageref{s20112} give some definitions and general properties of the $p$-class group of the extension $\Q(\zeta)/\Q$.
The subsection \ref{s108031} p. \pageref{s108031} gives some results obtained in the structure of the $p$-class group $C_p$ of $\Q(\zeta)$ and on class number $h$ of $\Q(\zeta)$.

%
\subsection{Some definitions and notations}\label{s20111}
In this subsection, we fix or recall some notations used in all this section.
\subsection{Representations
$\rho_d : G_d \rightarrow {\bf GL}({\bf F}_p.V_d).$}\label{s111171}
\begin{itemize}
\item
Let  $c_i= Cl(\sigma^i(\mathbf b)),\quad i=0,\dots,p-2,$  be the class of $\sigma^i(\mathbf b)$ in the $p$-class group of $\Q(\zeta)$.
\item
With these notations, we have $c_{i d}=Cl(\sigma^{i d}(\mathbf b)),
\quad i=0,\dots,\frac{p-1}{d}-1$.
\item
The groups $<c_{i d}>,\quad i=0,\dots,\frac{p-1}{d}-1$, are cyclic groups of order $p$.
\item
Let, for $d$ given, $C_{p,d}$ be the subgroup of $C_p^1$ generated by the set of
 ideal classes
\begin{displaymath}
E_d=\{ c_{i d}\quad | \quad i=0,\dots,\frac{p-1}{d}-1\}.
\end{displaymath}
Let, for $d$ given,   $r_d$ be the $p$-rank of $C_{p,d}$. Clearly, from these definitions, it results immediatly that
\begin{equation}\label{e105121}
1 \leq  r_d \leq \frac{p-1}{d}.
\end{equation}
Show that the group $C_{p,d}$ is also generated by the set of classes of ideals
\begin{equation}\label{e011251}
E_d(r_d)=\{ c_{i d}\quad | \quad i=0,\dots,r_d-1\}:
\end{equation}
\begin{itemize}
\item
It is clear when $r_d=\frac{p-1}{d}$.
\item
Suppose that $r_d<\frac{p-1}{d}$ and that the group $C(E_d(r_d))$ generated by the set of classes $E_d(r_d)$ verifies $C(E_d(r_d))\subset C_{p,d}$ and $C(E_d(r_d))\not= C_{p,d}$ and search for a contradiction:

\begin{itemize}
\item Suppose at first that $c_{id},\quad i=0,\dots,r_d$, are not linearly independant:  from linear relation between the $c_{i d},\quad i=0,\dots, r_d$, by conjugation $\sigma^d$, we would get the same relation between $c_{i d},\quad i=1,\dots, r_d+1$ and then, by elimination, a linear relation between $c_{i d},\quad i=0,\dots, r_d-1$ and $c_{(r_d+1)d}$, and in the same way, between $c_{i d},\quad i=0,\dots, r_d-1$ and $c_{(r_d+2) d}$, and finally $C(E_d(r_d))=C_{p,d}$, contradiction;
\item
then $c_{i d},\quad i=0,\dots, r_d$, should be linearly independant and therefore the $p$-rank of $C(E_d(r_d))$ should be $\geq r_d+1 > r_d$, which is impossible.
\end{itemize}
\item
Elsewhere,  the group $C_{p,d}$ is not generated  by the set of class of ideals
\begin{displaymath}
E_d(l)=\{ c_{i d}\quad | \quad i=0,\dots,l-1, \quad l<r_d, \}
\end{displaymath}
if not the $p$-rank would be smaller than $r_d$.
\end{itemize}
\item
Note that $C_{p,d}\subset C_p^1$,  without assuming any other hypothesis: $C_{p,d}=C_p^1$ or
$C_{p,d}\not= C_p^1$ and therefore
\begin{displaymath}
r_d\leq r_p,
\end{displaymath}
where $r_p$ is the $p$-rank of the class group of $\Q(\zeta)$.
\item
Observe that the property
$\prod_{i=0}^{p-2}\sigma^{i }(\mathbf b)$  is a principal ideal implies that
$r_1<p-1$.
\item
$F(r_d) = \oplus_{i=0}^{r_d-1} <c_{i d}>$  can be seen as a $\bf F_p$-vector space of dimension $r_d$ that we denote ${\bf F}_p.V_d$ in the sequel.
\item
With this  definition, the set $\{c_{i d}\ |\ i=0,\dots,r_d-1 \}$ is a basis  of the vector space ${\bf F}_p.V_d$,
\item
With this  definition of the vector space ${\bf F}_p.V_d$, corresponds {\bf one} representation
$\rho_d : G_d \rightarrow  {\bf GL}({\bf F}_p.V_d)$ of degree $r_d$ defined by:
\begin{displaymath}
\begin{split}
& \rho_d(\sigma^d)(c_{i d})= c_{(i+1)d}, \quad i=0,\dots, r_d-2,\\
& \rho_d(\sigma^d)(c_{(r_d-1)d})=-\sum_{i=0}^{r_d-1}\lambda_{i, d} c_{i d},\\
& \lambda_{i,d}\in {\bf F}_p,
\quad \lambda_{0,d}\not=0, \lambda_{r_d,d}=1,
\end{split}
\end{displaymath}
the second relation meaning that
\begin{displaymath}
\s^{\lambda_{0,d}}\times\sigma^d(\s)^{\lambda_{1,d}}\times\dots\times
\sigma^{d(r_d-1)}(\s)^{\lambda_{r_d-1,d}}\times\sigma^{d r_d}(\s)\simeq \Z[\zeta],
\end{displaymath}
(with notation $\simeq\Z[\zeta]$ for principal).
\item
Observe that $\rho_d(\sigma^d)( c_{(r_d-1)d})=c_{r_d d}$,
therefore, we have the linear dependance relation
\begin{equation}\label{e01071a}
L_{d}: \quad  c_{r_d d}+\sum_{i=0}^{r_d-1} \lambda_{i, d} c_{i d} = 0.
\end{equation}
The vector $c_{(r_d+1)d}$ can be computed by
\begin{equation}\label{e01071a}
c_{(r_d+1) d}+\sum_{i=0}^{r_d-1} \lambda_{i, d} c_{(i+1) d} = 0.
\end{equation}
The vectors $c_{(i+1) d},\quad i= r_d,\dots, \frac{p-1}{d}-2$, can be explicitly computed successively in function of $\lambda_{i, d},\quad i=0,\dots, r_d-1$ in this way.
\item
Observe that if $d=1$ then $r_1 <p-1$, and so
$\rho_1 : G \rightarrow {\bf GL}({\bf F}_p.V_1)$ corresponds to this definition.
\item
${\bf F}_p.V_d$ can be seen as the sub-${\bf F}_p[G_d]$-module of $C_p^1$ generated by action of the group $G_d$ on group $<c_d> =<Cl(\sigma^d(\mathbf b))>$.
\end{itemize}
\begin{itemize}
\item
Recall that $Cl(\mathbf b)$ is the class of the ideal $\mathbf b$ of $\Z[\zeta]$.
Observe that exponential notations $\mathbf b^\sigma$ can be used indifferently in the sequel.
With this notation, we have
\begin{itemize}
\item
$\mathbf b^{\sigma^d}=\sigma^d(\mathbf b)$.
\item
For $\lambda\in {\bf F}_p$, we have
$ \mathbf b^{\sigma+\lambda} = \mathbf b^\lambda\times  \sigma(\mathbf b).$
\item
Let $P(\sigma)=
\sigma^m+\lambda_{m-1}\sigma^{m-1}+\dots+\lambda_1\sigma+\lambda_0
 \in {\bf F}_p[\sigma]$; then
$ \mathbf b^{P(\sigma)} =
\sigma^m(\mathbf b)\times \sigma^{m-1}(\mathbf b)
^{\lambda_{m-1}}\times\dots\times\sigma(\mathbf b)^{\lambda_1}
\times \mathbf b^{\lambda_0}$.
\item
Let us note $ \mathbf b^{P(\sigma)} \simeq \Z[\zeta]$,  if the ideal
$\sigma^m(\mathbf b)\times \sigma^{m-1}(\mathbf b)^{\lambda_{m-1}}\dots \sigma(\mathbf b)^{\lambda_1}\times \mathbf b^{\lambda_0}$ is principal.
\item
Let $P(\sigma), Q(\sigma)\in {\bf F}_p[\sigma]$;
if $\mathbf b^{P(\sigma)}\simeq \Z[\zeta]$,
then $ \mathbf b^{(Q(\sigma)\times P(\sigma))}\simeq \Z[\zeta]$.
\item
Observe that trivially $ \mathbf b^{(\sigma^{p-1}-1)} \simeq \Z[\zeta].$
\end{itemize}
\item
To relation $L_d$  defined in relation (\ref{e01071a}) p.\pageref{e01071a}, corresponds the  monic minimal  polynomial $P_{r_d}(V)\in {\bf F}_p[V]$, polynomial ring of the indeterminate $V$; the polynomial  $P_{r_d}(V)$  verifies the relation, for $V=\sigma^d$:
\begin{equation}\label{e30101}
 \mathbf b^{P_{r_d}(\sigma^d)}
=  \mathbf b^{(\sigma^{d r_d}+\sum_{i=0}^{r_d-1} \lambda_{i, d} \sigma^{d i})}
\simeq \Z[\zeta].
\end{equation}
It is minimal because  for all polynomials $R(V)\in {\bf F}_p(V), \quad R(V)\not=0,\quad
deg(R(V))<deg(P_{r_d}(V))$, we have $  \mathbf b^{ R(\sigma^d)}\not\simeq \Z[\zeta]$.
It means, with an other formulation in term of ideals, that
$\prod_{i=0}^{r_d} \sigma^{i d}(\mathbf b)^{\lambda_{i,d}}$ is a principal ideal and
that $\prod_{i=0}^{\alpha} \sigma^{i d}(\mathbf b)^{\beta_i}$ is not principal when $\alpha<r_d$ and $\beta_i,\quad i=0,\dots,\alpha$, are not all simultaneously null.
\end{itemize}
%
\subsection{Representations of Galois group $Gal(\Q(\zeta)/\Q)$ in characteristic $p$.}\label{s20112}
In this subsection we give some general properties of representations of $G=Gal(\Q(\zeta)/\Q)$ in characteristic $p$ and obtain some results on the structure of  the $p$-class group of $\Q(\zeta)$.
Observe that we never use characters theory.
%
\subsubsection {Some preliminary results} \label{s07061a}
\begin{lem} \label{l04101}
Let $d\in\N, \quad p-1\equiv 0\modu d$, with the representation
\begin{displaymath}
\rho_d : G_d\rightarrow {\bf GL}({\bf F}_p.V_d)
\end{displaymath}
defined in previous subsection \ref{s20111} p. \pageref{s20111}.
Let $V$ be an indeterminate. Then the minimal polynomial $P_{r_d}(V)$ of the
representation $\rho_d$ verifies the factorization

\begin{displaymath}
P_{r_d}(V)=\prod_{i=1}^{r_d} (V-\mu_{i, d}),\quad \mu_{i,d}\in{\bf F}_p,\quad
 i_1\not= i_2 \Rightarrow \mu_{i_1}\not= \mu_{i_2}.
\end{displaymath}
\begin{proof}
Let us consider the polynomials $A(V)=V^{p-1}-1$ and $P_{r_d}(V)\in {\bf F}_p[V]$. It is possible to divide the polynomial  $A(V)$ by $P_{r_d}(V)$ in the polynomial ring
${\bf F}_p[V]$ to obtain
\begin{displaymath}
\begin{split}
& A(V)=P_{r_d}(V)\times Q(V)+R(V),\quad Q(V),R(V)\in{\bf F}_p[V],\\
& d_R=deg_V(R(V))<  r_d = deg_V(P_{r_d}(V)).
\end{split}
\end{displaymath}
For $V=\sigma^d$, we get
 $ \mathbf b^{(\sigma^{d(p-1)}-1)} \simeq \Z[\zeta]$ and
$  \mathbf b^{P_{r_d}(\sigma^d)} \simeq \Z[\zeta]$, so  $  \s^{R(\sigma^d)} \simeq \Z[\zeta]$.
Suppose that $R(V)=\sum_{i=0}^{d_R} R_i V^i,\quad R_i\in {\bf F}_p$, is not identically null; then, it leads to the relation
\begin{displaymath}
\mathbf b ^{(\sum_{i=0}^{d_R} R_i \sigma ^{d i})}\simeq \Z[\zeta],
\end{displaymath}
where the $R_i$ are not all zero,
with $d_R <r_d$, which contradicts the minimality of the relation $L_d$, relation (\ref{e01071a}) p.\pageref{e01071a}, and of the minimal polynomial
$P_{r_d}(V)$.
Therefore, $R(V)$ is identically null and we have
\begin{displaymath}
V^{p-1}-1= P_{r_d}(V)\times Q(V).
\end{displaymath}
The factorization of $V^{p-1}-1$ in ${\bf F}_p[V]$ is
$V^{p-1}-1
=\prod_{i=1}^{p-1} (V-i)$.
The factorization is unique in the euclidean ring ${\bf F}_p[V]$ and so
$P_{r_d}(V)=\prod_{i=1}^{r_d} (V-\mu_{i,d}),\quad \mu_{i,d}\in{\bf F}_p,
\quad  i_1\not= i_2 \Rightarrow \mu_{i_1}\not= \mu_{i_2}$, which achieves the proof.
\end{proof}
\end{lem}
%
\begin{lem}\label{l22101}
Let  the  representation
$\rho_1 : G\rightarrow {\bf GL}({\bf F}_p.V_1)$ of degree $r_1$ with minimal polynomial $P_{r_1}(U)$ of the indeterminate $U$.
Let $d\in\N,\quad p-1\equiv 0\modu d$, be such that
$\rho_d : G_d\rightarrow {\bf GL}({\bf F}_p.V_d)$ is a   representation of degree  $r_d$ of minimal polynomial $P_{r_d}(W)$ of the indeterminate $W$, as defined in previous subsection \ref{s20111} p. \pageref{s20111}.
Then
\begin{enumerate}
\item
$ P_{r_1}(U)=\prod_{i=1}^{r_1} (U-\mu_i),\quad \mu_i\in{\bf F}_p.$
\item
$ P_{r_d}(U^d)=\prod_{i=1}^{r_d} (U^d-\mu_i^d)
= P_{r_1}(U)\times Q_d(U),\quad r_d\leq r_1,
\quad Q_d(U)\in{\bf F}_p[U]$.
\item
The $p$-ranks $r_1$ and $r_d$ verify the inequalities
\begin{equation}\label{e011181}
r_d\times d\geq  r_1\geq r_d.
\end{equation}
\item
Let $K_d$ be the intermediate field $\Q\subset K_d\subset\Q(\zeta),
\quad [K_d:\Q]=d$. Suppose that $p$ does not divide the class number of $K_d/\Q$;
then $\mu_i^d\not=1$ for $i=1,\dots,r_d$. In particular   $\mu_i\not=1$ for $i=1,\dots,r_1$.
\end{enumerate}
\begin{proof}$ $
\begin{itemize}
\item
Observe, at first, that $deg_U(P_{r_d}(U^d))=d\times r_d\geq r_1$: if not, for the polynomial $P_{r_d}(U^d)$ in the indeterminate $U$, whe should have
$deg_U(P_{r_d}(U^d))<r_1$ and
$P_{r_d}(\sigma^d) \circ \mathbf b \simeq \Z[\zeta]$ and, as previously, the polynomial  $P_{r_d}(U^d)$ of the indeterminate $U$ should be identically null.
\item
We  apply euclidean algorithm in the polynomial ring ${\bf F}_p[U]$ of the indeterminate $U$.
Therefore,
\begin{displaymath}
\begin{split}
& P_{r_d} (U^d)=P_{r_1}(U)\times Q(U)+R(U),\quad Q(U),R(U)\in{\bf F}_p[U],\\
& deg(R(U))<deg(P_{r_1}(U)).
\end{split}
\end{displaymath}
But we have $\mathbf b^{P_{r_d}(\sigma^d)} \simeq\Z[\zeta],
\quad  \mathbf b^{P_{r_1}(\sigma )} \simeq \Z[\zeta]$, therefore
$\mathbf b^{R( \sigma)}\simeq \Z[\zeta]$.
Then, similarly to proof of lemma \ref{l04101} p.\pageref{l04101}, $R(U)$ is identically null and
$P_{r_d}(U^d)= P_{r_1}(U) \times Q(U)$.
\item
Applying lemma \ref{l04101} p.\pageref{l04101},  we obtain
\begin{displaymath}
\begin{split}
&P_{r_1}(U)=\prod_{i=1}^{r_1} (U-\mu_{i}),\quad \mu_{i}\in{\bf F}_p,\\
&P_{r_d}(U^d)=\prod_{i=1}^{r_d} (U^d-\mu_{i,d}),\quad \mu_{i,d}\in{\bf F}_p.
\end{split}
\end{displaymath}
\item
Then, we get
\begin{displaymath}
P_{r_d}(U^d)=\prod_{i=1}^{r_d} (U^d-\mu_{i,d})=\prod_{i=1}^{r_1} (U-\mu_i)\times Q(U).
\end{displaymath}
There exists at least one $i, \quad 1\leq i \leq r_d$, such that
$(U^d-\mu_{i,d})=(U-\mu_1)\times Q_1(U)$:
if not, for all $i=1,\dots,r_d$, we should have
$U^d- \mu_{i,d}\equiv R_i\modu (U-\mu_1),\quad R_i\in {\bf F}_p-\{0\}$, a contradiction because $\prod_{i=1}^{r_d} R_i\not=0$.
We have $\mu_{i,d}=\mu_1^d$: if not $U-\mu_1$ should divide $U^d-\mu_{i,d}$ and $U^d-\mu_1^d$ and also $U-\mu_1$ should divide
$(\mu_{i,d}-\mu_1^d)\in {\bf F}_p-\{0\}$, a contradiction.
Therefore, there exists at least one $i, \quad 1 \leq i\leq r_d$, such that $\mu_{i,d}=\mu_1^d$ and $U^d-\mu_{i,d}= U^d-\mu_1^d= (U-\mu_1)\times Q_1(U)$.
Then, generalizing to $\mu_{i,d}$ for all $i=1,\dots,r_d$, we get with a certain reordering of index $i$
\begin{displaymath}
P_{r_d}(U^d)=\prod_{i=1}^{r_d} (U^d-\mu_i^d)
=\prod_{i=1}^{r_1} (U-\mu_i)\times Q(U).
\end{displaymath}
\item
We have
\begin{displaymath}
P_{r_d}(U^d)=\prod_{i=1}^{r_d} (U^d-\mu_i^d).
\end{displaymath}
This relation leads to
\begin{displaymath}
P_{r_d}(U^d)=\prod_{i=1}^{r_d}\prod_{j=1}^d(U-\mu_i\mu_d^j),
\end{displaymath}
where $\mu_d\in {\bf F}_p,\quad \mu_d^d=1$.
We have shown that $P_{r_d}(U^d)=P_{r_1}(U)\times Q_d(U)$ and so
$deg_U(P_{r_d}(U))=d\times r_d\geq r_1$;  thus $d\times r_d\geq r_1$.
\item
We finish by the proof of enumeration 4): suppose that, for some $i,\quad 1\leq i\leq r_d$, we have $\mu_i^d=1$ and search for a contradiction:
there exists, for the indeterminate V,  a polynomial $P_1(V)\in {\bf F}_p(V)$ such that
$P_{r_d}(V)=(V-\mu_i^d)\times P_1(V)=(V-1)\times P_1(V)$.
But for $V= \sigma^d$, we have $ \mathbf b ^{P_{r_d}(\sigma^d) }\simeq \Z[\zeta]$,
so $ \mathbf b ^{(\sigma^d P_{1}(\sigma^d)- P_1(\sigma^d))}\simeq\Z[\zeta]$. So, $\mathbf b^{P_1(\sigma^d)} $ is the class of an  ideal $\mathbf c$ of $\Z[\zeta]$ with $Cl(\sigma^d(\mathbf c))=Cl(\mathbf c)$;
then  $Cl(\sigma^{2d}(\mathbf c))=Cl(\sigma^d(\mathbf c))
=Cl(\sigma^d(\mathbf c))=Cl(\mathbf c)$.
Then $Cl(\sigma^d(\mathbf c)\times\sigma^{2d}(\mathbf c)\times\dots
\times\sigma^{(p-1)d/d}(\mathbf c))
=Cl(\mathbf c^{(p-1)/d})$.
Let $\tau=\sigma^{d}$; then
$Cl(\tau(\mathbf c)\times\tau^{2 }(\mathbf c)\times\dots\times
\tau^{(p-1)/d}\mathbf c))=Cl(\mathbf c^{(p-1)/d})$;
Then we deduce that
$Cl(N_{\Q(\zeta)/K_d} (\mathbf c))=Cl(\mathbf c^{(p-1)/d})$ and thus
$\mathbf c$ is a principal ideal because the ideal  $N_{\Q(\zeta)/K_d} (\mathbf c)$
of $K_d$ is principal,  (from hypothesis, $p$ does not divide $h(K_d/\Q)$); so $ \mathbf b^{P_1(\sigma^d)} \simeq\Z[\zeta]$, which contradicts the minimality of
the minimal polynomial equation $ \mathbf b ^{P_{r_d}(\sigma^d)}\simeq\Z[\zeta]$ because, for the indeterminate $V$, we would have  $deg(P_{1}(V))<deg(P_{r_d}(V))$, which achieves the proof.
\end{itemize}
\end{proof}
\end{lem}
%
\subsubsection {Preliminary  on representations}\label{s07062a}
We summarize some results obtained in:
\begin{lem}\label{t31101a}
Let $\mathbf b$ be an ideal of $\Z[\zeta],\quad \mathbf b^p\simeq \Z[\zeta],\quad
\mathbf b\not \simeq \Z[\zeta]$.
Let $\rho_1:G\rightarrow {\bf GL}({\bf F}_p.V_1)$ be the   representation  corresponding to $\sigma$,
see definition in subsection \ref{s20111} p. \pageref{s20111}. Let $r_1$ be the degree of $\rho_1$.
Let $d\in\N,\quad p-1\equiv 0\modu d$.
Let $\sigma^d :\zeta \rightarrow \zeta^{u_d}$ be the $\Q$-isomorphism of $\Q(\zeta)$ corresponding to $d$.
Let $\rho_d : G_d \rightarrow {\bf GL}({\bf F}_p.V_d)$ be the representation corresponding to $\sigma^d$.
Let $r_d$ be the degree of the representation $\rho_d$.
Let, for the indeterminate $W$,
\begin{displaymath}
P_{r_d}(W)=W^{r_d}+\lambda_{r_d-1,d} W^{r_d-1}+\dots+\lambda_{1,d}W+\lambda_{0,d}
\end{displaymath}
be the minimal polynomial corresponding to the representation $\rho_d$.
Then there exists $\mu_1,\mu_2,\dots,\mu_{r_1}\in{\bf F}_p$,
with $i\not=i^\prime\Rightarrow \mu_i\not=\mu_i^\prime$,  such that, for the indeterminate $U$,
\begin{itemize}
\item
the minimal polynomials $P_{r_1}(U)$ and $P_{r_d}(U^d)$ are respectively given by
\begin{displaymath}
\begin{split}
&P_{r_1}(U)=\prod_{i=1}^{r_1} (U-\mu_i),\\
&P_{r_d}(U^d)=\prod_{i=1}^{r_d} (U^d-\mu_i^d),\quad r_d\leq r_1,\\
&P_{r_1}(U) \ | \ P_{r_d}(U^d).
\end{split}
\end{displaymath}
\item
The coefficients of $P_{r_d}(U^d)$ are explicitly computable by
\begin{displaymath}
\begin{split}
& P_{r_d}(U^d)=\\
& U^{d r_d} -S_1(d)\times U^{d(r_d-1)}+S_2(d)\times U^{d(r_d-2)}
+\dots+(-1)^{{r_d}-1} S_{{r_d}-1}(d)\times U^d
+(-1)^{r_d} S_{r_d}(d),\\
& S_0(d)=\lambda_{r_d,d}=1,\\
& S_1(d)=-\lambda_{{r_d}-1,d}=\sum_{i=1,\dots, r_d} \mu_i^d,\\
& S_2(d)=\lambda_{{r_d}-2,d}=
\sum_{1\leq i_1<i_2\leq r_d} \mu_{i_1}^d\mu_{i_2}^d,\\
&\vdots\\
&S_{r_d}(d)=(-1)^{r_d} \lambda_{0,d}=\mu_1^d\mu_2^d\dots\mu_{r_d}^d.
\end{split}
\end{displaymath}
\item
Then the ideal
\begin{equation}\label{e04112}
\prod_{i=0}^{r_d} \sigma^{di}(\mathbf b)^{(-1)^{r_d-i}\times S_{r_d-i}(d)}
=\prod_{i=0}^{r_d}\sigma^{di}(\mathbf b)^{\lambda_{i,d}} =\mathbf b^{P_{r_d}(\sigma^d)}
\end{equation}
is a principal ideal.
\end{itemize}
\end{lem}
{\bf Remark:} For other annihilation methods of $Cl(\Q(\zeta)/\Q)$ more involved, see for instance Kummer, in Ribenboim \cite{rib} p 119, (2C) and (2D) and Stickelberger in Washington \cite{was} p 94 and 332.
%
\subsection{On the structure of the $p$-class group of subfields of $\Q(\zeta)$}
\label{s108031}
In this subsection on $p$-class group of cyclotomic fields, strictly independant of FLT equation,  we deal of:
\begin{itemize}
\item
a formulation, with our notations,  of a Ribet's  result on irregularity index.
\item
the action of $Gal(\Q(\zeta)/\Q)$ on  $p$-class groups  of subfields of $\Q(\zeta)$,
\item
an inequality involving the $p$-rank of the group generated by the action of  groups $G,
G_d, G_g, \quad d|p-1,\quad g|p-1,\quad gcd(d,g)=1$ on ideals $\mathbf b$ of $\Q(\zeta)$,
\item
some relations obtained from Stickelberger theorem in $\Q(\zeta)$,
\item
some relations obtained from Stickelberger theorem in subfields of $\Q(\zeta)$,
\item
some congruences $\modu \pi^{p+1}$ connected to $p$-class group $C_p$ of $\Q(\zeta)$ with consequences on structure of relative $p$-class group $C_p^-$,
\item
some congruences on all prime factors of class number $h$ of $\Q(\zeta)$.
\item
an explicit formula for the computation of relative class number $h^-$.
\end{itemize}
%
\subsubsection{Some definitions and notations}
\begin{itemize}
\item
Recall that:
\begin{itemize}
\item
$C_p$ is the $p$-class group of $\Q(\zeta)$, subgroup of the class group $C$ whose elements are of
$p$-power order,
\item
$r_p$ is the $p$-rank of the class group of $\Q(\zeta)$,
\item
$C_p^1$ is the subgroup of $C_p$ whose elements are of order $1$ or $p$.
\end{itemize}
\item
Let $G=Gal(\Q(\zeta)/\Q)$ be the Galois group of the field $\Q(\zeta)$. Let,
for $d|p-1$,  $G_d$ be the subgroup of $d$ powers $\sigma^{di}$ of elements $\sigma^i$ of $G$.
\item
Suppose that $r_p>0$.
There always exists an ideal class with a representant $\mathbf b\subset\Z[\zeta]$,
with $\mathbf b^p\simeq \Z[\zeta],\quad \mathbf b\not\simeq\Z[\zeta]$, which verifies, in term of representations, for some  ideals  $\mathbf b_i$ of $\Z[\zeta],\quad i=1,\dots,r_p$,
\begin{equation}\label{e107162}
\begin{split}
& \mathbf b\simeq\prod_{i=1}^{r_p} \mathbf  b_i,\\
& \mathbf b_i^p\simeq \Z[\zeta],\quad \mathbf b_i\not\simeq \Z[\zeta],\quad i=1,\dots,r_p,\\
& \sigma(\mathbf b_i)\simeq\mathbf b_i^{\mu_i},
\quad \mu_i\in{\bf F}_p,\quad \mathbf b_i+\pi=\Z[\zeta],\quad i=1,\dots,r_p,\\
& C_p^1\cong \oplus_{i=1}^{r_p} <Cl(\mathbf b_i)>,\\
& P_{r_1}(U)=\prod_{i=1}^{r_1} (U-\mu_i),\quad \mathbf b^{P_{r_1}(\sigma)}
\simeq \Z[\zeta],\quad 1\leq r_1\leq r_p,
\end{split}
\end{equation}
where $P_{r_1}(U)$ is the minimal polynomial in the indeterminate $U$ for the action of $G$ on an ideal  $\mathbf b$, see theorem \ref{t31101a} p \pageref{t31101a}.
Recall that it is possible to  encounter the case  $\mu_i=\mu_j$ in the set
$\{\mu_1,\dots,\mu_{r_p}\}$; by opposite if $U-\mu_i$ and $U-\mu_j$ divide the minimal polynomial $P_{r_1}(U)$ then $\mu_i\not=\mu_j$. $r_1$ is therefore the degree of the minimal polynomial $P_{r_1}(U)$.
\item
Let us denote $M_{r_1}=\{\mu_i\quad|\quad i=1,\dots,r_1\}$.
\item
Let $d\in\N,\quad d|p-1,\quad 2\leq d\leq\frac{p-1}{2}$. Let $K_d$ be the field
$\Q\subset K_d\subset\Q(\zeta),\quad [K_d:\Q]=d$.
\item
Let $P_{r_d}(V)$ be the minimal polynomial of the action of the group $G_d$ on the ideal class group $<\mathbf b>$ of order $p$.
Let $r_d$ be the degree of $P_{r_d}(V)$.
\end{itemize}
%
\subsubsection{ On the irregularity index}
Recall that $r_p$ is the $p$-rank of the group $C_p$.
The irregularity index is the number
\begin{displaymath}
i_p=Card\{B_{p-1-2m}\ | \ B_{p-1-2m}\equiv 0\modu p,\quad 1\leq m\leq\frac{p-3}{2}\},
\end{displaymath}
where $B_{p-1-2m}$ are even Bernoulli Numbers.
The next theorem connects irregularity index and $p$-rank $r_1$ defined in relation
(\ref{e107162}) p. \pageref{e107162}.
%
\begin{thm}{ *** }\label{t204251}
With meaning of $r_1$ defined in relation (\ref{e107162}) p. \pageref{e107162} then $i_p=r_1$. Moreover, if
$h^+\not\equiv 0\modu p$ then $i_p=r_1=r_p$.
\begin{proof}
Let us consider in relation (\ref{e107162}) the set of ideals  $\{\mathbf b_i\ |\ i=1,\dots,r_p\}$.
The result of Ribet using theory of modular forms \cite{rie} mentionned in Ribenboim \cite{rib} (8C) p 190 can be formulated, with our notations,
\begin{equation}\label{e204251}
B_{p-1-2m}\equiv 0\modu p \Leftrightarrow \exists i,\quad 1 \leq i\leq r_p,\quad
\mathbf b_i^{\sigma-u_{2m+1}}\simeq \Z[\zeta].
\end{equation}
There exists at least one such $i$, but it is possible for $i\not=i^\prime$ that
$\mathbf b_i^{\sigma-u_{2m+1}}\simeq \mathbf b_{i^\prime}^{\sigma-u_{2m+1}}\simeq \Z[\zeta]$.
The relation (\ref{e204251}) implies that $i_p=r_1$.
When $h^+\not\equiv 0$ then $i_p=r_p$, see for instance Washington \cite{was} in proof of theorem 6.18 p 102.
\end{proof}
\end{thm}
{\bf Remark:} when $h^+\not\equiv 0$ the relation $r_p=r_1$ is a reformulation of theorem
\ref{t011241} p. \pageref{t011241} proved  in FLT context in an elementary way.
%
%

%
\subsubsection{Inequalities involving  degrees $r_1, r_d, r_g $ of minimal polynomials $P_{r_1}(V),P_{r_d}(V),P_{r_g}(V)$ annihilating ideal $\mathbf b$.}
In this subsection, we always assume that $\mathbf b$ is defined by relation  (\ref{e107162}) p. \pageref{e107162}.

Let $p>3$ be a prime.
Let $d, g\in \N$, with $gcd(d,g)=1$ and $ d\times g\ |\ p-1$.
Recall that are   the degree of the minimal polynomials
$P_{r_1}(V), P_{r_d}(V), P_{r_g}(V)$ of the indeterminate $V$ with
$\mathbf b^{P_{r_1}(\sigma^d)}\simeq \mathbf b^{P_{r_d}(\sigma^d)}
\simeq \mathbf b^{P_{r_g}(\sigma^g)}\simeq \Z[\zeta]$.
%
\begin{lem}\label{l210021}
Let $d\in\N,\quad 1 < d <p-1,\quad p-1\equiv 0\modu d$.
If $d\times r_1\geq p-1$ then $r_d<r_1$.
\begin{proof}
For each $i=1,\dots,r_1$, there are $d$ values $\mu_{j i},\quad j=1,\dots,d$,
such that $\mu_{i j}^d=\mu_i^d$.
Therefore there are $d\times r_1$ values $\mu_{i,j}$ which cannot be all different, thus there exists at least two pairs $(i_1,j_1)$ and $(i_1,j_2)$ such that $\mu_{i_1,j_1}=\mu_{i_1,j_2}$ and so $r_d<r_1$.
\end{proof}
\end{lem}
%
The next theorem is a relation between the three degree  $r_1, r_d$ and $r_g$.
\begin{thm}{ *** }\label{l108081}
Let $d, g\in\N, \quad gcd(d,g)=1,\quad d\times g\ |\ p-1$.
Then
\begin{equation}\label{e108081}
r_d\times r_g\geq r_1.
\end{equation}
and if $r_d=1$ then $r_g=r_1$.
\begin{proof} $ $
\begin{itemize}
\item
Let us consider the minimal polynomials $P_{r_d}(U^d)=\prod_{i=1}^{r_d}
(U^d-\mu_i^d)$ and $P_{r_g}(U^g)=\prod_{i=1}^{r_g}(U^g-\nu_j^g)$ of the indeterminate $U$ with $\mathbf b^{P_{r_d}(\sigma^d)}\simeq \Z[\zeta]$ and
$\mathbf b^{P_{r_g}(\sigma^g)}\simeq \Z[\zeta]$.
\item
From lemma \ref{l22101} p.\pageref{l22101},
we have seen that $P_{r_1}(U)|P_{r_d}(U^d)$ and that similarly
$P_{r_1}(U)|P_{r_g}(U^g)$, thus $P_{r_1}(U)|gcd(P_{r_d}(U^d), P_{r_g}(U^g))$.
\item
Let $M_{r_1}=\{\mu_i\quad|\quad i=1,\dots,r_1\}$.
Let us define
\begin{displaymath}
C_1(\mu_i)=\{\mu_i\times \alpha_j\quad|\quad\alpha_j^d=1,\quad j=1,\dots,d\}\cap M_{r_1}, \quad i=1,\dots,r_d.
\end{displaymath}
Let us define in the same way
\begin{displaymath}
C_2(\nu_i)=\{\nu_i\times \beta_j\quad|\quad\beta_j^g=1,\quad j=1,\dots,g\}\cap M_{r_1}, \quad i=1,\dots,r_g.
\end{displaymath}
\item
We have proved in lemma \ref{l22101} p.\pageref{l22101} that $P_{r_1}(U)|P_{r_d}(U^d)$. Therefore the sets
$C_1(\mu_i),\quad i=1,\dots,r_d$, are a partition of $M_{r_1}$ and
$r_1=\sum_{i=1}^{r_d} Card(C_1(\mu_i))$.
\item
In the same way $P_{r_1}(U)|P_{r_g}(U^g)$. Therefore the sets
$C_2(\nu_i),\quad i=1,\dots,r_g$, are a partition of $M_{r_1}$ and
$r_1=\sum_{i=1}^{r_g} Card(C_2(\nu_i))$.
\item
There exists at least one $i\in\N,\quad 1\leq i\leq r_d$, such that
$Card(C_1(\mu_i))\geq \frac{r_1}{r_d}$. For this $i$, let
$\nu_1 =\mu_i\times \alpha_1,\quad \alpha_1^d=1,\quad \nu_1\in M_{r_1}$ and, in the same way, let
$\nu_2 =\mu_i\times \alpha_2,\quad \alpha_2^d=1,\quad \nu_2\in M_{r_1},\quad \nu_2\not=\nu_1$.
We have $\nu_1^g\not=\nu_2^g$ : if not we should simultaneously have
$\alpha_1^d=\alpha_2^d$ and $\alpha_1^g=\alpha_2^g$, which should imply, from
$gcd(d,g)=1$, that $\alpha_1=\alpha_2$, contradicting $\nu_1\not=\nu_2$ and therefore we get $C_2(\nu_1)\not=C_2(\nu_2)$.
\item
Therefore, extending the same reasoning to all elements of $C_1(\mu_i)$, we get  $\frac{r_1}{r_d}\leq Card(C_1(\mu_i))\leq r_g$, which leads to the result.
\item
If $r_d=1$ then $r_g\geq r_1$ and in an other part $r_g\leq r_1$ and so $r_g=r_1$.
\end{itemize}
\end{proof}
\end{thm}
%
{\bf Remark:} Observe that $r_d=0$ implies directly that $r_1=0$.
%
%
%
\subsubsection {On Stickelberger's ideal in field $\Q(\zeta)$}
In this subsection, we give a result   resting on the annihilation of class group of $\Q(\zeta)$ by Stickelberger's ideal.
\begin{itemize}
\item
Let us denote $\mathbf a\simeq \mathbf c$ when the two ideals $\mathbf a$ and $\mathbf c$ of $\Q(\zeta)$ are in the same ideal class.
\item
Let $G=Gal(\Q(\zeta)/\Q)$.
\item
Let $\tau_a :\zeta \rightarrow \zeta^a,\quad a=1,\dots,p-1$, be the $p-1$ $\Q$-isomorphisms of
the field $\Q(\zeta)/\Q$.
\item
Recall that $u$ is a primitive root $\modu p$,  and that

$\sigma : \zeta \rightarrow \zeta^u$
is a $\Q$-isomorphism of the field $\Q(\zeta)$
which generates $G$.  Recall that, for $i\in \N$, then we denote $u_i$ for
$u^i \modu p$ and $1\leq u_i\leq p-1$.
\item
Let $\mathbf b$ be the not principal ideal defined in relation (\ref{e107162}) p.\pageref{e107162}.
Let $P_{r_1}(\sigma)\in {\bf F}_p[G]$ be the polynomial of minimal degree such that
$P_{r_1}(\sigma)$ annihilates $\mathbf b$, so such that
$\mathbf b^{P_{r_1}(\sigma)}$ is principal ideal, see lemma \ref{l04101} p.\pageref{l04101}, and so
\begin{displaymath}
P_{r_1}(\sigma)=\prod_{i=1}^{r_1} (\sigma-\mu_i),\quad \mu_i\in {\bf F}_p,
\quad i\not=i^\prime \Rightarrow \mu_i\not=\mu_i^\prime.
\end{displaymath}
\end{itemize}
In the next result we shall explicitly use the annihilation of class group of $\Q(\zeta)$ by the Stickelberger's ideal.
%
\begin{lem}\label{l103021}
Let $P_{r_1}(U)=\prod_{i=1}^{r_1}(U-\mu_i)$ be the polynomial of the indeterminate $U$, of minimal degree,  such that $\mathbf b^{P_{r_1}(\sigma)}$ is principal.  Then $\mu_i\not= u,\quad i=1,\dots,r_1$.
\begin{proof}
Let $i\in\N, \quad 1\leq i\leq r_1$.
From relation (\ref{e107162}) p.\pageref{e107162}, there exists   ideals $\mathbf b_i\in \Z[\zeta],\quad i=1,\dots,r_p$, not principal and such that $\mathbf b=\prod_{i=1}^{r_p} \mathbf b_i$,  with $\mathbf b_i^{\sigma-\mu_i}$  principal.
Suppose that $\mu_i=u$, and search for a contradiction.
Let us consider  $\theta=\sum_{a=1}^{p-1}\frac{a}{p}\times\tau_a^{-1}\in\Z[G]$.
$p\theta\in\Z[G]$ and the ideal
$\s^{p\theta}$ is principal from Stickelberger's theorem, see for instance Washington \cite{was}, theorem 6.10 p 94.
We can set  $a=u^{m},\quad a=1,\dots,p-1$, and $m$ going through all the set $\{0,1,\dots,p-2\}$, because $u$ is a primitive root $\modu p$.
Then $\tau_a :\zeta\rightarrow \zeta^a$ and so
$\tau_a^{-1} :\zeta \rightarrow \zeta^{(a^{-1})}=\zeta^{((u^m)^{-1})}
=\zeta^{(u^{-m})}
=\zeta^{(u^{p-1-m})}=\sigma^{p-1-m}=\sigma^{-m}$.
Therefore,
$\theta=\sum_{m=0}^{p-2} u^m\sigma^{-m}$.
The element $\sigma-\mu_i=\sigma-u$ annihilates the class of $\mathbf b_i$ and also the element
$u\times\sigma^{-1}-1$ annihilates the class of $\mathbf b_i$. Therefore
$u^{m}\sigma^{-m}-1,\quad m=0,\dots,p-2,$ annihilates the class of $\mathbf b_i$ and finally $p-1$ annihilates the class of $\mathbf b_i$, so $\mathbf b_i^{p-1}$ is principal, but $\mathbf b_i^p$ is also principal, and finally $\mathbf b_i$ is principal which contradicts our hypothesis and achieves the proof.
\end{proof}
\end{lem}
%
%

\subsubsection{$\pi$-adic congruences  connected  to $p$-relative class group $C_p^-$}\label{s108311}
In this subsection , we shall describe  some $\pi$-adic congruences  connected  to $p$-relative class group $C_p^-$.
\begin{itemize}
\item
Let $C_p$ be the $p$-class group of $\Q(\zeta)$, $C_p^1$ be the subgroup of $C_p$ whose elements are of order $1$ or $p$.
\item
Let $r_p$ be the $p$-rank of $C_p$, let  $r_p^+$ be the $p$-rank of $C_p^+$ and $r_p^-$ be the relative $p$-rank of $C_p^-$.
Let us consider the ideals $\mathbf b$ defined in relation (\ref{e107162}) p.\pageref{e107162}:
\begin{equation}\label{e203031}
\begin{split}
& \mathbf b=\mathbf b_1\times\dots\times\mathbf b_{r_p^-}
\times \mathbf b_{r_p^-+1}\times\dots\times \mathbf b_{r_p},\\
& C_p^1=\oplus_{i=1}^{r_p} <Cl(\mathbf b_i)>,\\
& \mathbf b_i^p\simeq\Z[\zeta],\quad \mathbf b_i\not\simeq \Z[\zeta],
\quad i=1,\dots,r_p,\\
&\sigma(\mathbf b_i)\simeq \mathbf b_i^{\mu_i},\quad\mu_i\in{\bf F}_p^*,
\quad i=1,\dots,r_p,\\
& Cl(\mathbf b_i)\in C_p^-,\quad i=1,\dots,r_p^-,\\
& Cl(\mathbf b_i)\in C_p^+,\quad i=r_p^-+1,\dots,r_p,\\
& \mathbf b^{P_{r_1}(\sigma)}\simeq\Z[\zeta],\\
& (\frac{\mathbf b}{\overline{\mathbf b}})^{P_{r_1^-}(\sigma)}\simeq\Z[\zeta].
\end{split}
\end{equation}
Recall that $P_{r_1}(\sigma)\in {\bf F}_p[G]$ is the minimal polynomial
such that $\mathbf b^{P_{r_1}(\sigma)}\simeq\Z[\zeta]$ with $r_1\leq r_p$.
Recall that $P_{r_1^-}(\sigma)\in {\bf F}_p[G]$ is the minimal polynomial
such that
$(\frac{\mathbf b}{ \overline{\mathbf b}})^{P_{r_1}(\sigma)}\simeq\Z[\zeta]$
 with $r_1^-\leq r_p^-$.
\item
We say that $C\in \Q(\zeta)$ is  singular if $C\Z[\zeta]=\mathbf c^p$ for
some ideal $\mathbf c$ of $\Q(\zeta)$. We say that $C$ is singular primary
if $C$ is singular and $C\equiv c^p\modu \pi^p,\quad c\in \Z,\quad c\not\equiv 0\modu p$.
\end{itemize}
%
\begin{lem}\label{l108161}
There exists $B_i\in \Z[\zeta],\quad i=1,\dots,r_p^-$, such that
\begin{equation}
\begin{split}
& B_i\Z[\zeta]=\mathbf b_i^p,\\
& \sigma(\frac{B_i}{\overline{B_i}})
\times(\frac{B_i}{\overline{B_i}})^{-\mu_i}
=(\frac{\alpha_i}{\overline{\alpha_i}})^p,\quad\alpha_i\in\Q(\zeta),
\quad \alpha_i\Z[\zeta]+\pi=\Z[\zeta],\\
& \sigma(\frac{B_i}{\overline{B_i}})\equiv
(\frac{B_i}{\overline{B_i}})^{\mu_i}\modu \pi^{p+1},\quad \mu_i=u_{2m_i+1}.
\end{split}
\end{equation}
\begin{proof}
Observe that we can neglect the $\mu_i=u_{2m_i}$ such that $\sigma-\mu_i$ annihilates ideal classes  $\in C_p^+$, because we consider only
quotients $\frac{B_i}{\overline{B_i}}$, with ideal classes  $Cl(\mathbf b_i)$ in $C_p^-,\quad i=1,\dots,r_p^-$.
The ideal $\mathbf b_i^p$ is principal. So let one $\beta_i\in\Z[\zeta]$ with
$\beta_i\Z[\zeta]=\mathbf b_i^p$.
We have seen in relation (\ref{e107162}) p.\pageref{e107162} that   $\sigma(\mathbf b_i)\simeq \mathbf b_i^{\mu_i}$, therefore there exists
$\alpha_i\in \Q(\zeta)$ such that
$\frac{\sigma(\mathbf b_i)}{\mathbf b_i^{\mu_i}}=\alpha_i\Z[\zeta]$, also
$\frac{\sigma(\mathbf \beta_i)}{\mathbf \beta_i^{\mu_i}}= \varepsilon_i\times \alpha_i^p,\quad \varepsilon_i\in\Z[\zeta]^*$.
Let $B_i=\delta_i^{-1}\times \beta_i,\quad \delta_i\in \Z[\zeta]^*$,  for a choice of
the unit $\delta_i$ that whe shall explicit in the next lines.
We have
\begin{displaymath}
\sigma(\delta_i\times B_i)=\alpha_i^p\times (\delta_i\times B_i)^{\mu_i}\times\varepsilon_i.
\end{displaymath}
Therefore
\begin{equation}\label{e108162}
\sigma(B_i)=
\alpha_i^p\times B_i^{\mu_i} \times
(\sigma(\delta_i^{-1})\times\delta_i^{\mu_i}\times\varepsilon_i).
\end{equation}
From Kummer's lemma on units, we can write
\begin{displaymath}
\begin{split}
& \delta_i=\zeta^{v_1}\times\eta_1,\quad v_1\in \Z,
\quad \eta_1\in\Z[\zeta+\zeta^{-1}]^*,\\
& \varepsilon_i=\zeta^{v_2}\times\eta_2,\quad v_2\in \Z,
\quad \eta_2\in\Z[\zeta+\zeta^{-1}]^*.
\end{split}
\end{displaymath}
Therefore
\begin{displaymath}
\sigma(\delta_i^{-1})\times \delta_i^{\mu_i}\times\varepsilon_i=
\zeta^{-v_1 u+v_1\mu_i+v_2}\times\eta,\quad \eta\in\Z[\zeta+\zeta^{-1}]^*.
\end{displaymath}
From lemma \ref{l103021} p.\pageref{l103021}, we deduce that $\mu_i\not=u$, therefore there exists one
$v_1$ with $-v_1 u+v_1\mu_i+v_2\equiv 0\modu p$.
Therefore, chosing this value $v_1$ for the unit $\delta_i$,
\begin{equation}\label{e112212}
\begin{split}
& \sigma(B_i)=\alpha_i^p\times B_i^{\mu_i}\times\eta,
\quad \alpha_i\Z[\zeta]+\pi=\Z[\zeta],\quad \eta\in\Z[\zeta+\zeta^{-1}]^*,\\
& \sigma(\overline{B_i})=\overline{\alpha_i}^p\times\overline{B_i}^{\mu_i}
\times\eta.
\end{split}
\end{equation}
We have $\alpha_i\equiv \overline{\alpha}_i\modu \pi$ and we have proved in proposition  \ref{p2} p.\pageref{p2} that $\alpha_i^p\equiv\overline{\alpha}_i^p\modu \pi^{p+1}$,  which leads to the result.
\end{proof}
\end{lem}
%
\begin{lem}\label{l108171}
For each $i=1,\dots,r_p^-$, there exists  $B_i\in \Z[\zeta]$, such that
\begin{equation}\label{e108191}
\begin{split}
& \mu_i=u_{2m_i+1},\quad i=1,\dots,r_p^-,\\
& B_i\Z[\zeta]=\mathbf b_i^p,\\
& C_i=\frac{B_i}{\overline{B_i}}\equiv 1
\modu \pi^{2m_i+1}.
\end{split}
\end{equation}
\begin{proof} $ $
\begin{itemize}
\item
Fix, without loss of generality, the index $i,\quad 1\leq i\leq r_p$, in this proof.
We have $B_i\Z[\zeta]+\pi=\Z[\zeta]$.
Let $C(\zeta)=\frac{B_i}{\overline{B_i}}$.
From the definition of $C(\zeta)$ we see that the valuation $v_\pi(C(\zeta))=0.$
We can write $C(\zeta)$ in the form
\begin{equation}\label{e108201}
\begin{split}
& C(\zeta)=\frac{B_i}{\overline{B_i}}=c_0+c_1\times \zeta+\dots
+c_{p-2}\times \zeta^{p-2},\\
& c_i\in \Q, \quad v_p(c_i)\geq 0.
\end{split}
\end{equation}
Observe at first that, from the definitions,   $C(\zeta)\equiv 1\modu\pi$.
\item
Let $X$ be an indeterminate and $m\in\{0,1\}$.
Let us denote $C(X, m)=\sum_{i=0}^{p-2} c_i\times X^{i u^m}$.
\item
We have $C(\zeta,0)=\sum_{i=0}^{p-2} c_i\times \zeta^i$ and
$C(\zeta,1)=\sum_{i=0}^{p-2} c_i\times\zeta^{ u i}=\sigma(C(\zeta,0))$.
We have proved in lemma \ref{l108161} p.\pageref{l108161} that
\begin{equation}\label{e108202}
C(\zeta,1)=\sigma(C(\zeta,0))\equiv C(\zeta,0)^{\mu}\modu \pi^{p+1}.
\end{equation}
Let us denote $C^\prime(\zeta,m)$ the derivative on $X$ for $X=\zeta$.
Using exactly the same technics of derivation on the indeterminate $X$ that in subsection \ref{s12112} p.\pageref{s12112} and foundation theorem \ref{s12113} p.\pageref{s12113}, we get
\begin{equation}\label{e108181}
C^\prime(\zeta,1)\equiv \mu\times C^\prime(\zeta,0)\times C(\zeta,0)^{\mu-1}\modu \pi^{p-2}.
\end{equation}
\item
We have $C(X,1)=\sum_{i=0}^{p-2} c_i \times X^{i u}$ and so
$C^\prime(X,1)=\sum_{i=1}^{p-2} c_i\times  i\times  u\times  X^{i u-1}$.
We have $C(X,0)=\sum_{i=0}^{p-2} c_i \times X^{i }$ and so
$C^\prime(X,0)=\sum_{i=1}^{p-2} c_i \times i  \times X^{i -1}$. We have $C(\zeta,0)\equiv 1\modu \pi$. Thus we obtain
\begin{displaymath}
u\times(\sum_{i=1}^{p-2} c_i\times  i \times \zeta^{i u-1})\equiv
\mu\times(\sum_{i=1}^{p-2} c_i\times  i\times\zeta^{i  -1})\modu \pi.
\end{displaymath}
A first congruence is $(u-\mu)\times (\sum_{i=1}^{p-2} c_i\times  i)\equiv 0\modu \pi$. We have seen in lemma \ref{l103021} p.\pageref{l103021} that $u\not=\mu$,
so
\begin{equation}\label{e108182}
C^\prime(1,0)=\sum_{i=1}^{p-2} c_i\times  i\equiv 0\modu p.
\end{equation}
\item
Let us define $C_1(X,m)=C^\prime(X,m)\times X$.
We have, multiplying relation (\ref{e108181}) p.\pageref{e108181} by $\zeta$
\begin{displaymath}
C_1(\zeta,1)\equiv \mu\times C_1(\zeta,0)\times C(\zeta,0)^{\mu-1}\modu\pi^{p-2}.
\end{displaymath}
With a derivation on $X$ for the value $X=\zeta$, we get
\begin{equation}\label{e108183}
\begin{split}
& C^\prime_1(\zeta,1)\equiv \mu\times C^\prime_1(\zeta,0)\times
C(\zeta,0)^{\mu-1}\\
&+\mu\times (\mu-1)\times C_1(\zeta,0)\times C(\zeta,0)^{\mu-2}
\times C^\prime(\zeta,0)\modu\pi^{p-3}.
\end{split}
\end{equation}
Therefore
\begin{displaymath}
\begin{split}
& C^\prime_1(1,1)\equiv \mu\times C^\prime_1(1,0)\times
C(1,0)^{\mu-1}\\
&+\mu\times (\mu-1)\times C_1(1,0)\times C(1,0)^{\mu-2}
\times C^\prime(1,0)\modu p.
\end{split}
\end{displaymath}
We have $C_1(X,1)=u\times (\sum_{i=1}^{p-2} c_i\times i\times X^{i u})$ and so
$C^\prime_1(X,1)=u^2\times (\sum_{i=1}^{p-2} c_i\times i^2\times X^{i u-1})$.
We have $C_1(X,0)=(\sum_{i=1}^{p-2} c_i\times i\times X^{i})$, and so
$C^\prime_1(X,1)=(\sum_{i=1}^{p-2} c_i\times i^2\times X^{i-1 })$.
We have seen that $C(1,0)\equiv 1\modu p$.
We have proved in relation (\ref{e108182}) p.\pageref{e108182}   that $C_1(1,0)\equiv
C^\prime(1,0)\equiv 0\modu p$.
Gathering these results, we obtain
\begin{displaymath}
u^2\times(\sum_{i=1}^{p-2} c_i\times i^2)
\equiv \mu\times(\sum_{i=1}^{p-2} c_i\times i^2) \modu p,
\end{displaymath}
which leads to
\begin{displaymath}
(u_2-\mu)\times(\sum_{i=1}^{p-2} c_i\times i^2)\equiv 0\modu p.
\end{displaymath}
We have seen in lemma \ref{l108161} p. \pageref{l108161} that $\mu=u_{2m+1}$ so that $u_2-\mu\not=0$, therefore
\begin{equation}\label{e108186}
\sum_{i=1}^{p-2} c_i\times i^2\equiv 0\modu p.
\end{equation}
\item
Pursuing this induction, let us define
\begin{displaymath}
\begin{split}
& C_k(X,1)=C_{k-1}^\prime(X,1)\times X,\quad k=2,\dots,p-2,\\
& C_k(X,0)=C_{k-1}^\prime(X,0)\times X,\quad k=2,\dots,p-2.
\end{split}
\end{displaymath}
For $k=2$ we get, multiplying by $\zeta$ the relation (\ref{e108183}) p.\pageref{e108183},
\begin{equation}
\begin{split}
& C_2(\zeta,1)\equiv \mu\times C_2(\zeta,0)\times
C(\zeta,0)^{\mu-1}\\
&+\mu\times (\mu-1)\times C_1(\zeta,0)^2\times C(\zeta,0)^{\mu-2}
\modu\pi^{p-3}.
\end{split}
\end{equation}
As an other example, we have for $k=3$,
\begin{displaymath}
\begin{split}
& C_3(\zeta,1)=\mu C_3(\zeta,0)C(\zeta,0)^{\mu-1}\\
& +3\mu(\mu-1)C_2(\zeta,0)C_1(\zeta,0)C(\zeta,0)^{\mu-2}
+\mu(\mu-1)(\mu-2)C_1(\zeta,0)^3 C(\zeta,0)^{\mu-3}.
\end{split}
\end{displaymath}
Pursuing, we have for $k>1$ relations of general form,
\begin{equation}\label{e108185}
\begin{split}
& C_k(\zeta,1)\equiv \mu\times C_k(\zeta,0)\times
C(\zeta,0)^{\mu-1}\\
&+\sum_l \gamma_l\times C(\zeta,0)^\alpha\times C_{l_1}(\zeta,0)^{\alpha_1}\times \dots\times C_{l_\nu}(\zeta,0)^{\alpha_\nu}
\modu \pi^{p-1-k},\\
&\gamma_l\in\Q(\zeta),\quad v_\pi(\gamma_l)\geq 0,\quad l_1,\dots,l_\nu\in\N-\{0\},\quad l_1 \alpha_1+\dots+l_\nu\alpha_\nu=k,\quad\nu>1.
\end{split}
\end{equation}
By induction for $k<2m+1$, we have $C_l(\zeta,0)\equiv 0\modu \pi,\quad  l=1,\dots,k$, and
$C_l(1,0)\equiv 0\modu p,\quad l=1,\dots,k$:
this results of the induction formula
\begin{displaymath}
\begin{split}
& (u_k-\mu)\times C_k(1,0)=(u_k-\mu)\times(\sum_{i=1}^{p-2} c_i\times i^k)\equiv \\
&\sum_l \gamma_l \times C_{l_1}(1,0)^{\alpha_1}\times \dots\times C_{l_\nu}(1,0)^{\alpha_\nu}\equiv 0
\modu p,
\quad \nu>1,\quad l_1\alpha_1+\dots+l_\nu\alpha_\nu=k.
\end{split}
\end{displaymath}
\item
Pursuing
for all $k<2m+1$, we see that $u_k-\mu\not=0$
and so
\begin{displaymath}
C_k(1,0)=\sum_{i=1}^{p-2} c_i\times i^k\equiv 0\modu p.
\end{displaymath}
Therefore
\begin{displaymath}
C_k(1,0)\equiv 0\modu p,\quad k=1,\dots,2m.
\end{displaymath}
Similarly to relation (\ref{e22021}) p.\pageref{e22021} in foundation theorem reciprocal proof, we deduce
\begin{displaymath}
C^{(k)}(1,0)\equiv 0\modu p,\quad k=1,\dots, 2m,
\end{displaymath}
observing that $C^\prime_{k-1}(X,0)\times X= C_k(X,0)$ and noting $C^{(k)}(X,0)$ for the $k$-eme derivate on $X$  of $C(X,0)$.
We deduce from Taylor formula, with $\lambda=\zeta-1$, that
\begin{displaymath}
C(\zeta,0)\equiv C(1,0)+\lambda \times C^\prime(1,0)+\dots+\frac{\lambda^{2m}}{(2m)!}\times  C^{(2m)}(1,0)
\modu \pi^{2m+1}
\end{displaymath}
which gives  \begin{displaymath}
C(\zeta,0)\equiv 1\modu \pi^{2m+1},
\end{displaymath}
which is the relation (\ref{e108191}) p.\pageref{e108191}.
\item
Observe that the process stops at $k=2m$, because we have $u_{2m+1}-\mu=0$ and we can say nothing on $C_{2m+1}(1,0)$ in this way.
\end{itemize}
\end{proof}
\end{lem}
%
{\bf Remarks:}
\begin{enumerate}
\item
The previous proof allows to say more: else
$\pi^{2m+1}\|C(\zeta,0)-1$ and $C(\zeta,0)$ is not singular primary else $C(\zeta,0)\equiv 1\modu\pi^p$ and $C(\zeta,0)$ is singular primary: if $C_{2m+1}(1,0)\equiv 0\modu p$ then relation (\ref{e108185}) p. \pageref{e108185} shows that $C_{2m+2}(1,0)\equiv \dots\equiv C_{2(2m+1)}(1,0)\equiv 0\modu p$ and so on; then $C_1(1,0)\equiv\dots\equiv C_{p-2}(1,0)\equiv 0\modu p$ and so $C(\zeta,0)\equiv 1\modu \pi^{p-1}$.
We derive that $C(\zeta,0)\equiv 1\modu \pi^p$ from the complementary property that
$C(\zeta,1)-\mu C(\zeta,0)\equiv 0\modu \pi^{p+1}$: we have seen that
$C(\zeta,0)=1+pV,\quad V\in\Q(\zeta),\quad v_\pi(V)\geq 0$; then $C(\zeta,1)=
1+p\sigma(V)$ and $C^\mu\equiv 1+p\mu V\modu \pi^{p}$; $C(\zeta,1)-C(\zeta,0)^\mu
\equiv p(\sigma(V)-\mu V)\modu \pi^{p}$, so $\sigma(V)-\mu V\equiv 0\modu \pi$, also $(1-\mu)V\equiv 0\modu \pi$, and finally $V\equiv 0\modu \pi$.
%
\item
In previous lemma \ref{l108171} p. \pageref{l108171} we have seen that
\begin{displaymath}
C(\zeta,0)\equiv 1+\frac{C_{2m+1}(1,0)}{(2m+1)!}\lambda^{2m+1}\modu \pi^{2m+2},
\end{displaymath}
where $C_1(1,0)\equiv\dots\equiv C_{2m}(1,0)\equiv 0\modu p$ implies that $C_{2m+1}(1,0)\equiv C^{(2m+1)}(1,0)\modu \pi^{p-1}$.
The choice of the determination ideal $\mathbf b_i,\quad \sigma(\mathbf b_i)\simeq \mathbf b_i^{\mu_i}$
among the $p-1$  not principal classes of $<Cl(\mathbf b_i)>$ is arbitrary and if $C(\zeta,0)$ is not singular primary, equivalent to $C_{2m+1}(1,0)\not\equiv 0\modu p$,
it is always possible to choose the determination $\mathbf b_i$ such that
$C_{2m+1}(1,0)\equiv 1\modu p$, determination chosen in the sequel.
\end{enumerate}
%
\begin{thm} {*** }\label{t112311}
Let $C=C_{1}^{\alpha_1}\dots C_{m}^{\alpha_m}$ with $\alpha_i\in{\bf F}_p^*,\quad i=1,\dots,m$.
Then $C$ is singular primary if and only if all the $C_i,\quad i=1,\dots,m$, are all singular primary.
\begin{proof}$ $
\begin{itemize}
\item
If $C_i,\quad i=1,\dots,m$, are all singular primary, then $C$ is clearly singular primary.
\item
Suppose that $C_i,\quad i=1,\dots, l$, are not singular primary and that $C_i,\quad i=l+1,\dots,m$, are singular primary. Then, from lemma \ref{l108171} p.\pageref{l108171} and remark following it, $\pi^{2m_i+1}\| C_i-1,\quad i=1,\dots,l$, where we suppose, without loss of generality, that
$1<2m_1+1<\dots < 2m_l+1$. Then $\pi^{2m_1+1}\| C-1$ and so $C$ is not singular primary.
\end{itemize}
\end{proof}
\end{thm}

%

\begin{lem}\label{l202231}
Let $C(\zeta)$ defined in relation (\ref{e108201}) p. \pageref{e108201}. If $C(\zeta)$ is not singular primary, then
\begin{displaymath}
C(\zeta)\equiv 1+ V(\mu)\modu \pi^{p-1},\quad \mu=u_{2m+1}
\quad V(\mu)\in\Z[\zeta],
\end{displaymath}
where $V(\mu)\modu p$ depends only on $\mu$ with $\pi^{2m+1}\| V(\mu)$.
\begin{proof}
If $C(\zeta)$ is not singular primary then, from relation (\ref{e108185}) p. \pageref{e108185}, $C_{2m+2}(1,0),\dots, C_{p-2}(1,0)$ and so
$C^{(2m+2)}(1,0),\dots, C^{(p-2)}(1,0)$ can be derived of $C_{2m+1}(1,0)$ with relations depending only of $\mu$, which achieves the proof.
\end{proof}
\end{lem}
%
\begin{thm}{ *** }\label{l203031}
Let $C_1, C_2$ defined by relation (\ref{e108201}) p. \pageref{e108201}.
If $\mu_1=\mu_2$ then $C_1$ and $C_2$ are singular primary.
\begin{proof}
Let $\mu_1=\mu_2=\mu=u_{2m+1}$.
From previous lemma \ref{l202231} p. \pageref{l202231} we get
\begin{displaymath}
\begin{split}
& C_1=1+V(\mu)+p W_1,\quad W_1\in\Q(\zeta), \quad v_\pi(V(\mu))\geq 2m+1,
\quad v_\pi(W_1)\geq 0,\\
& C_2=1+V(\mu)+p W_2,\quad W_2\in\Q(\zeta), \quad v_\pi(V(\mu))\geq 2m+1,
\quad v_\pi(W_2)\geq 0.
\end{split}
\end{displaymath}
Elsewhere, $C_1, C_2$ verify
\begin{displaymath}
\begin{split}
& \sigma(C_1)\equiv C_1^\mu\modu\pi^{p+1},\\
& \sigma(C_2)\equiv C_2^\mu\modu\pi^{p+1},
\end{split}
\end{displaymath}
which leads to
\begin{displaymath}
\begin{split}
& 1+\sigma(V(\mu))+p \sigma(W_1)\equiv
1+ A(\mu)+ p \mu W_1 \modu \pi^{p+1},\\
& 1+\sigma(V(\mu))+p \sigma(W_2)\equiv
1+ A(\mu)+ p \mu  W_2 \modu \pi^{p+1},
\end{split}
\end{displaymath}
where $A(\mu)\in \Q(\zeta),\quad v_\pi(A(\mu)\geq 0$ depends only on $\mu$.
By difference, we get
\begin{displaymath}
p(\sigma(W_1-W_2))\equiv p  \mu (W_1- W_2)\modu \pi^{p+1},
\end{displaymath}
which implies that
\begin{displaymath}
\sigma(W_1-W_2)\equiv   \mu (W_1- W_2)\modu \pi^{2}.
\end{displaymath}
Let $W_1-W_2= a\lambda+b,\quad a,b\in\Z,\quad \lambda=\zeta-1$.
The previous relation implies that $b(1-\mu)\equiv 0\modu p$ and so
that $b\equiv 0\modu p$.
Thus $W_1-W_2\equiv 0\modu \pi$
and finally $C_1\equiv C_2\modu \pi^p$ and also
$C_1 C_2^{-1}\equiv 1\modu \pi^p$ and $C_1 C_2^{-1}$ is singular primary, which shall imply that $C_1$ and $C_2$ are singular primary and shall achieve the proof:
\begin{enumerate}
\item
$C_1$ not singular primary and $C_2$ singular primary should imply that $C_1 C_2^{-1}$ is not singular primary.
\item
$C_2$ not singular primary and $C_1$ singular primary should imply that $C_1 C_2^{-1}$ is not singular primary.
\item
$C_1$ not singular primary and $C_2$ not singular primary should imply that
$C_1 C_2^{-1}$ is not singular primary from theorem \ref{t112311} p. \pageref{t112311}.
\end{enumerate}
\end{proof}
\end{thm}
%
\begin{thm}{ *** On structure of $p$-class group $C_p$. }\label{t203021}

Let $\mathbf b_i$ be the ideals defined in relation (\ref{e203031}) p. \pageref{e203031}.
Let $C_p^{1-}=\oplus _{i=1}^{r_p^-} <Cl(\mathbf b_i)>$.
Let $C_i=\frac{B_i}{\overline{B}_i},\quad B_i\Z[\zeta]=\mathbf b_i^p,
\quad Cl(\mathbf b_i)\in C_p^-,
\quad i=1,\dots,r_p^-$,  where
$B_i$ is defined in relation (\ref{e112212}) p. \pageref{e112212}.
Let $r_1^-$ be the degree of the minimal polynomial $P_{r_1^-}(\sigma)$ defined in relation (\ref{e203031}) p. \pageref{e203031}.
With a certain ordering of $C_i,\quad i=1,\dots,r_p^-$,
\begin{enumerate}
\item
$C_i$ are  singular primary for $i=1,\dots, r_p^+$, and $C_i$ are not singular primary for $i=r_p^+,\dots, r_p^-$.
\item
If $\mu_i=\mu_j$ then $1\leq i< j\leq r_p^+\leq r_p^-$ and   $C_i, C_j$ are singular primary.
\item
If $h^+\equiv 0\modu p$ then $r_p^--r_p^++1\leq r_1^-\leq r_p^-$.
\item
If $h^+\not \equiv 0\modu p$ then $r_1^-=r_p^-$.
\end{enumerate}
\begin{proof}$ $
\begin{enumerate}
\item
It is an application of Hilbert class field theorem \ref{t201013} \pageref{t201013} which is a forward independant reference : there are $r_p^+$ singular primary $C_i$ with
$Cl(\mathbf b_i)\in C_p^-$.
\item
lemma \ref{l203031} p. \pageref{l203031}
\item
Reformulation of previous items, observing that if $r_p^+>0$ then $\mu_i,\quad i=1,\dots,r_p^+$ take at least one value $\mu$ and at most $r_p^+$ different values $\mu$.
\item
If $h^+\not\equiv 0\modu p$ there are no $C_i$ singular primary.
\end{enumerate}
\end{proof}
\end{thm}
%
%
\begin{lem}\label{l203171}
Let, for $i=1,\dots,r_p^-$,  the ideals $\mathbf b_i$ with
$Cl(\mathbf b_i)\in C_p^-$ defined in relation (\ref{e203031}) p. \pageref{e203031}.
Let $C_i=\frac{B_i}{\overline{B_i}},\quad B_i\in\Z[\zeta],\quad i=1,\dots,r_p^-$, defined in relation (\ref{e108191}) p. \pageref{e108191}. There exists
$B_{i,1}= \frac{B_i^2}{\eta_i},\quad\eta_i\in\Z[\zeta+\zeta^{-1}]^*$, such that
\begin{equation}\label{e203271}
\begin{split}
& \mu_i=u_{2m_i+1},\quad i=1,\dots,r_p^-,\\
& \sigma(B_{i,1})=B_{i,1}^{\mu_i}\times \alpha_i^p,\quad \alpha_i\in\Q(\zeta),\\
& B_{i,1}\Z[\zeta]=\mathbf b_i^{2p},\quad B_{i,1}\in\Z[\zeta],\\
& B_{i,1}^{p-1}\equiv 1\modu \pi^{2m_i+1}.
\end{split}
\end{equation}
\begin{proof}
We have $C_i\Z[\zeta]=\mathbf b_i^p$ where the ideal $\mathbf b_i$ verifies
$\sigma(\mathbf b_i)\simeq \mathbf b_i^{\mu_i}$ and
$Cl(\mathbf b_i)\in C_p^-$. Let us note $C,B,\mu$ for $C_i, B_i,\mu_i$ to simplify notations.
From  relation (\ref{e112212}) p. \pageref{e112212} and from
$Cl(\mathbf b_i)\in C_p^-$, we can choose $B$ such that
\begin{displaymath}
\begin{split}
& C=\frac{B}{\overline{B}},\quad B\in\Z[\zeta], \\
& B\overline{B}=\eta\times\gamma^p,\quad \eta\in\Z[\zeta+\zeta^{-1}]^*,
\quad \gamma\in\Q(\zeta),\quad v_\pi(\gamma)=0,\\
& \sigma(B) = B^\mu\times\alpha^p\times \varepsilon,
\quad \mu=u_{2m+1},\quad \alpha\in\Q(\zeta),\quad v_\pi(\alpha)=0,
\quad \varepsilon\in\Z[\zeta+\zeta^{-1}]^*.
\end{split}
\end{displaymath}
We derive  that
\begin{displaymath}
\sigma(B\overline{B})= (B\overline{B})^\mu\times (\alpha\overline{\alpha})^p\times \varepsilon^2,
\end{displaymath}
and so
\begin{displaymath}
\sigma(\eta)=\eta^\mu\times\varepsilon^2\times\varepsilon_1^p,\quad
\varepsilon_1\in \Z[\zeta+\zeta^{-1}]^*.
\end{displaymath}
We have seen that
\begin{displaymath}
\sigma(B^2)=B^{2\mu}\times  \alpha^{2p}\times \varepsilon^2,
\end{displaymath}
and so
\begin{displaymath}
\sigma(B^{2})=B^{2\mu}\times\alpha^{2p}\times (\sigma(\eta)\eta^{-\mu}\varepsilon_1^{-p})
\end{displaymath}
which leads to
\begin{displaymath}
\sigma(\frac{B^2}{\eta})=(\frac{B^{2}}{\eta})^\mu\times \alpha_2^p,\quad
\alpha_2\in\Q(\zeta),\quad v_\pi(\alpha_2)=0.
\end{displaymath}
Let us note $B_1=\frac{B^2}{\eta},\quad B_1\in\Z(\zeta),\quad v_\pi(B_1)=0$.
We get
\begin{equation}\label{e203092}
\sigma(B_1)=B_1^\mu\times \alpha_2^p.
\end{equation}
This relation (\ref{e203092}) is similar to hypothesis of lemma  \ref{l108171} p. \pageref{l108171} which leads to $B_1\equiv d_1^p\modu \pi^{2m+1},\quad d_1\not\equiv 0\modu p$; Therefore $B_1^{p-1}\equiv 1\modu \pi^{2m+1}$,
 which achieves  the proof.
\end{proof}
\end{lem}
%
%
\begin{lem}\label{l203272}
Let the ideals $\mathbf b_i,\quad i=r_p^-+1,\dots,r_p$, such that
$Cl(\mathbf b_i)\in C_p^+$ defined in relation (\ref{e203031}) p. \pageref{e203031}. There exists
$B_i\in\Z(\zeta)$ such that:
\begin{equation}\label{e205041}
\begin{split}
& \mu_i=u_{2m_i},\quad i=r_p^-+1,\dots,r_p,\\
& \sigma(\mathbf b_i)\simeq \mathbf b_i^{\mu_i},\\
& B_i\Z[\zeta]=\mathbf b_i^p,\\
& \sigma(B_i)=B_i^{\mu_i}\times\alpha_i^p,\quad\alpha_i\in\Q(\zeta),\\
& B_i\equiv 1\modu \pi^{2m_i}.\\
\end{split}
\end{equation}
\begin{proof}
Similarly to relation (\ref{e112212}) p. \pageref{e112212}, there exists $B_i$ with $B_i\Z[\zeta]=\mathbf b_i^p$ such that
\begin{displaymath}
\sigma(B_i)=B_i^{\mu_i}\times\alpha^p\times \eta,
\quad \alpha\in\Q(\zeta),\quad \eta\in\Z[\zeta+\zeta^{-1}]^*.
\end{displaymath}
From relation (\ref{e201274}) p. \pageref{e201274}, independant forward reference in section dealing of
unit group $\Z[\zeta+\zeta^{-1}]^*$, we can write
\begin{displaymath}
\begin{split}
& \eta= (\prod_{j=1}^N\eta_j^{\lambda_i})\times (\prod_{j=N+1}^{(p-3)/2} \eta_{j}^{\lambda^j}),
\quad \lambda_j\in {\bf F}_p,\quad 1\leq  N<\frac{p-3}{2},\\
& \sigma(\eta_j)=\eta_j^{\mu_i}\times\beta_j^p,
\quad \eta_j,\beta_j\in\Z[\zeta+\zeta^{-1}]^*,\quad j=1,\dots,N,\\
& \sigma(\eta_{j})=\eta_{j}^{\nu_{j}}\times\beta_j^p,
\quad \eta_j,\beta_j\in\Z[\zeta+\zeta^{-1}]^*, \quad j=N+1,\dots,\frac{p-3}{2},\\
& \nu_j\not=\mu_i,\quad j=N+1,\dots,\frac{p-3}{2}.\\
\end{split}
\end{displaymath}
Let us note
\begin{displaymath}
E=\prod_{j=1}^N \eta_j^{\lambda_j},\quad
U=\prod_{j=N+1}^{(p-3)/2} \eta_j^{\lambda_j}.
\end{displaymath}
Show that there exists $V\in\Z[\zeta+\zeta^{-1}]^*$ such that
$\sigma(V)\times V^{-\mu_i}=U\times \varepsilon^p,\quad \varepsilon\in\Z[\zeta+\zeta^{-1}]^*$:
Let us set $V=\prod_{j=N+1}^{(p-3)/2} \eta_j^{\rho_i}$.
It suffices that
\begin{displaymath}
\eta_j^{\rho_j \nu_j}\times \eta_j^{-\rho_j\mu_i}
=\eta_j^{\lambda_j}\times \varepsilon_j^p,
\quad \varepsilon_j\in\Z[\zeta+\zeta^{-1}]^*,\quad j=N+1,\dots,\frac{p-3}{2}.
\end{displaymath}
It suffices that
\begin{displaymath}
\rho_j\equiv \frac{\lambda_j}{\nu_j-\mu_i}\modu p,\quad j=N+1,\dots,\frac{p-3}{2},
\end{displaymath}
which is possible, because $\nu_j\not\equiv \mu_i,\quad j=N+1,\dots, \frac{p-3}{2}$.
Therefore, for $B_i^\prime=B_i\times V$, we get simultaneously
\begin{displaymath}
\begin{split}
& \sigma(B_i^\prime)=(B_i^\prime)^{\mu_i}\times \alpha_i^p\times E,\quad \alpha_i\in\Q(\zeta),\\
& \sigma(E)=E^{\mu_i}\times \varepsilon_1^p,\quad \varepsilon_1\in\Z[\zeta+\zeta^{-1}]^*.\\
\end{split}
\end{displaymath}
If $E\in (\Z[\zeta+\zeta^{-1}]^*)^p\Leftrightarrow \lambda_j=0,\quad j=1,\dots,N$, then we get
\begin{equation}\label{e203281}
\sigma(B_i^\prime)=(B_i^\prime)^{\mu_i}\times (\alpha_i^\prime)^p.
\end{equation}
If $E\not\in(\Z[\zeta+\zeta^{-1}]^*)^p$ then  by conjugation,
\begin{displaymath}
\begin{split}
& \sigma(B_i^\prime)=(B_i^\prime)^{\mu_i}\times \alpha_i^p\times E,\quad a_i\in\Q(\zeta),\\
& \sigma^2(B_i^\prime)=\sigma(B_i^\prime)^{\mu_i}\times  E^{\mu_i}\times b_i^p,
\quad b_i\in\Q(\zeta),\\
\end{split}
\end{displaymath}
and so
\begin{displaymath}
c^p\times \sigma(B_i^\prime)^{\mu_i}(B_i^\prime)
^{-\mu_i^2}=\sigma^2(B_i^\prime)\sigma(B_i^\prime)^{-\mu_i},
\quad c\in\Q(\zeta)
\end{displaymath}
which leads to
\begin{displaymath}
(B_i^\prime)^{(\sigma-\mu_i)^2}=c^p,\quad c\in\Q(\zeta).
\end{displaymath}
Elsewhere $(B_i^\prime)^{\sigma^{p-1}-1}=1$, so
\begin{displaymath}
(B_i^\prime)^{gcd((\sigma^{p-1}-1,
(\sigma-\mu_i)^2})= (B_i^\prime)^{\sigma-\mu_i}=\alpha_i^p,
\quad \alpha_i\in\Q(\zeta),
\end{displaymath}
and so $\sigma(B_i^\prime)=(B_i^\prime)^{\mu_i}\times\alpha_i^p$. The end of proof is similar to
proof of previous lemma \ref{l203171} p. \pageref{l203171}.
\end{proof}
\end{lem}
%
%
The next important theorem summarize the two previous lemmas and give explicit congruences $\modu p$ connected to $p$-class group of $\Q(\zeta)$.
\begin{thm} { *** }{ $\pi$-adic structure of $p$-class group $C_p$}\label{t203281}

Let the ideals $\mathbf b_i,\quad i=1,\dots,r_p$, such that
$Cl(\mathbf b_i)\in C_p$ and defined by relation (\ref{e203031}) p. \pageref{e203031}.
Then, there exists $B_i\in\Z(\zeta),\quad i=1,\dots,r_p$, such that
\begin{equation}\label{e205042}
\begin{split}
& \mu_i=u_{m_i},\quad i=1,\dots,r_p,\\
& \sigma(\mathbf b_i)\simeq \mathbf b_i^{\mu_i},\\
& B_i\Z[\zeta]=\mathbf b_i^p,\\
& \sigma(B_i)=B_i^{\mu_i}\times\alpha_i^p,\quad \alpha_i\in\Q(\zeta),\\
& B_i\equiv 1\modu \pi^{m_i}.\\
\end{split}
\end{equation}
Moreover, among them:
\begin{enumerate}
\item
the $r_p^+$ algebraic integers $B_i,\quad i=1,\dots,r_p^+$, corresponding to $\mathbf b_i\in C_p^-$
and $C_i=\frac{B_i}{\overline{B}_i}$ singular primary, verify
$B_i\equiv 1\modu \pi^p$.
\item
the $r_p-r_p^+$ other algebraic integers $B_i,\quad i=r_p^++1,\dots,r_p$ verify
$\pi^{m_i}\| (B_i-1)$.
\end{enumerate}
\begin{proof}
Apply lemmas \ref{l203171} p. \pageref{l203171} and \ref{l203272} p. \pageref{l203272} and theorem \ref{t203021} p. \pageref{t203021}.
\end{proof}
\end{thm}
%
%
\subsubsection{The case $\mu=u_{2m+1},\quad 2m+1>\frac{p-1}{2}$}\label{s201111}
In the next lemma  we shall investigate more deeply the consequences of the congruence $C\equiv 1\modu \pi^{2m+1}$ of lemma \ref{l108171} p. \pageref{l108171} when $2m+1> \frac{p-1}{2}$.
%
\begin{lem}\label{l201091}
Let $C$ with $\mu=u_{2m+1},\quad 2m+1>\frac{p-1}{2}$.
\begin{displaymath}
\begin{split}
& C=1+\gamma+\gamma_0\zeta+\gamma_1 \zeta^u
+\dots + \gamma_{p-3}\zeta^{u_{p-3}},\\
& \gamma\in\Q,\quad v_p(\gamma)\geq 0,\quad  \gamma_i\in\Q,\quad
v_p(\gamma_i)\geq 0,\quad i=0, \dots, p-3,\\
& \gamma+\gamma_0\zeta+\gamma_1\zeta^u+\dots + \gamma_{p-3}\zeta^{u_{p-3}}\equiv 0\modu \pi^{2m+1},
\quad 2m+1> \frac{p-1}{2}.
\end{split}
\end{displaymath}
Then $C$ verifies  the congruences
\begin{displaymath}
\begin{split}
& \gamma\equiv -\frac{\gamma_{p-3}}{\mu-1}\modu p,\\
& \gamma_0\equiv -\mu^{-1} \times \gamma_{p-3}\modu p,\\
& \gamma_1\equiv -(\mu^{-2}+\mu^{-1})\times \gamma_{p-3}\modu p,\\
&\vdots\\
& \gamma_{p-4}\equiv -(\mu^{-(p-3)}+\dots+\mu^{-1})\times \gamma_{p-3}\modu p.
\end{split}
\end{displaymath}
\begin{proof}
We have seen in lemma \ref{l108161} p. \pageref{l108161} that $\sigma(C)\equiv C^{\mu}\modu\pi^{p+1}$.
From $2m+1 > \frac{p-1}{2}$ we derive that
\begin{displaymath}
C^\mu\equiv 1+\mu\times(\gamma+\gamma_0\zeta+\gamma_1\zeta^u
+\dots +\gamma_{p-3} \zeta^{u_{p-3}})\modu\pi^{p-1}.
\end{displaymath}
Elsewhere, we get by conjugation
\begin{equation}\label{e201111}
\sigma(C)=1+\gamma+\gamma_0\zeta^u+\gamma_1\zeta^{u_2}
+\dots +\gamma_{p-3} \zeta^{u_{p-2}}.
\end{equation}
We have the identity
\begin{displaymath}
\gamma_{p-3} \zeta^{u_{p-2}}=-\gamma_{p-3}-\gamma_{p-3}\zeta
-\dots -\gamma_{p-3}\zeta^{u_{p-3}}.
\end{displaymath}
This leads to
\begin{displaymath}
\sigma(C)=1+\gamma-\gamma_{p-3}-\gamma_{p-3}
\zeta+(\gamma_0-\gamma_{p-3})\zeta^u
+\dots+(\gamma_{p-4}-\gamma_{p-3})\zeta^{u_{p-3}}.
\end{displaymath}
Therefore, from the congruence $\sigma(C)\equiv C^{\mu}\modu \pi^{p+1}$ we get
the congruences in the basis $1,\zeta,\zeta^{u},\dots,\zeta^{u_{p-3}}$,
\begin{displaymath}
\begin{split}
& 1+\mu \gamma\equiv 1+\gamma-\gamma_{p-3}\modu p,\\
& \mu \gamma_0\equiv -\gamma_{p-3}\modu p,\\
& \mu \gamma_1\equiv \gamma_0-\gamma_{p-3}\modu p,\\
& \mu \gamma_2\equiv \gamma_1-\gamma_{p-3}\modu p,\\
& \vdots\\
& \mu \gamma_{p-4}\equiv \gamma_{p-5}-\gamma_{p-3}\modu p,\\
& \mu \gamma_{p-3}\equiv \gamma_{p-4}-\gamma_{p-3}\modu p.
\end{split}
\end{displaymath}
From these congruences, we get $\gamma\equiv -\frac{\gamma_{p-3}}{\mu-1}\modu p$ and
$\gamma_0\equiv -\mu^{-1} \gamma_{p-3}\modu p$ and then
$\gamma_1\equiv \mu^{-1}(\gamma_0-\gamma_{p-3})\equiv \mu^{-1}(-\mu^{-1}
\gamma_{p-3}-\gamma_{p-3})
\equiv -(\mu^{-2}+\mu^{-1})\gamma_{p-3}\modu p$
and $\gamma_2
\equiv\mu^{-1}(\gamma_1-\gamma_{p-3})
\equiv \mu^{-1}(-(\mu^{-2}+\mu^{-1})\gamma_{p-3}-\gamma_{p-3})
\equiv -(\mu^{-3}+\mu^{-2}+\mu^{-1})\gamma_{p-3}\modu p$ and so on.
\end{proof}
\end{lem}
%
%
%
The next theorem gives an explicit  important formulation  of $C$ when $2m+1>\frac{p-1}{2}$.
\begin{thm}{ *** }\label{l202211}
Let $\mu=u_{2m+1},\quad p-2\geq 2m+1>\frac{p-1}{2}$, corresponding to $C$ defined in lemma \ref{l108171} p. \pageref{l108171}, so $\sigma(C)\equiv C^\mu\modu\pi^{p+1}$.
Then $C$ verifies the formula:
\begin{equation}\label{e202211}
C\equiv 1-\frac{\gamma_{p-3}}{\mu-1}\times
(\zeta+\mu^{-1}\zeta^u+\dots+\mu^{-(p-2)}\zeta^{u_{p-2}}) \modu \pi^{p-1}.
\end{equation}
\begin{proof}
From definition of $C$, setting $C=1+V$, we get :
\begin{displaymath}
\begin{split}
& C=1+V,\\
& V=\gamma+\gamma_0\zeta+\gamma_1\zeta^u+\dots+\gamma_{p-3}\zeta^{u_{p-3}},\\
& \sigma(V)\equiv \mu\times V\modu \pi^{p+1}.
\end{split}
\end{displaymath}
Then, from lemma \ref{l201091} p. \pageref{l201091},  we obtain the relations
\begin{displaymath}
\begin{split}
& \mu=u_{2m+1},\\
& \gamma \equiv -\frac{\gamma_{p-3}}{\mu-1}\modu p,\\
& \gamma_0\equiv -\mu^{-1} \times \gamma_{p-3}\modu  p,\\
& \gamma_1\equiv -(\mu^{-2}+\mu^{-1})\times \gamma_{p-3}\modu p,\\
&\vdots\\
& \gamma_{p-4}\equiv -(\mu^{-(p-3)}+\dots+\mu^{-1})\times \gamma_{p-3}\modu  p,\\
& \gamma_{p-3}\equiv -(\mu^{-(p-2)}+\dots+\mu^{-1})\times \gamma_{p-3}\modu  p.
\end{split}
\end{displaymath}
From these relations we get
\begin{displaymath}
V\equiv  -\gamma_{p-3}\times
(\frac{1}{\mu-1}+\mu^{-1}\zeta+(\mu^{-2}+\mu^{-1})\zeta^{u}+\dots+
(\mu^{-(p-2)}+\dots+\mu^{-1})\zeta^{u_{p-3}})\modu p.
\end{displaymath}
Then
\begin{displaymath}
V\equiv-\gamma_{p-3}\times (\frac{1}{\mu-1}+
\mu^{-1}(\zeta+(\mu^{-1}+1)\zeta^u+\dots+(\mu^{-(p-3)}+\dots+1)
\zeta^{u_{p-3}}))\modu p.
\end{displaymath}
Then
\begin{displaymath}
V\equiv -\gamma_{p-3}\times (\frac{1}{\mu-1}+
\mu^{-1}(
\frac{(\mu^{-1}-1)\zeta+(\mu^{-2}-1)\zeta^u+\dots
+(\mu^{-(p-2)}-1)\zeta^{u_{p-3}}}{\mu^{-1}-1}))
\modu p.
\end{displaymath}
Then
\begin{displaymath}
V\equiv -\gamma_{p-3}\times (\frac{1}{\mu-1}+
\mu^{-1}(
\frac{\mu^{-1}\zeta+\mu^{-2}\zeta^u+\dots
+\mu^{-(p-2)}\zeta^{u_{p-3}}
-\zeta-\zeta^u-\dots -\zeta^{u_{p-3}}}{\mu^{-1}-1}))
\modu p.
\end{displaymath}
Then
\begin{displaymath}
V\equiv -\gamma_{p-3}\times (\frac{1}{\mu-1}+
\mu^{-1}(
\frac{1+\mu^{-1}\zeta+\mu^{-2}\zeta^u+\dots
+\mu^{-(p-2)}\zeta^{u_{p-3}}+\mu^{-(p-1)}\zeta^{u_{p-2}}}{\mu^{-1}-1}))
\modu p.
\end{displaymath}
Then
\begin{displaymath}
V\equiv -\gamma_{p-3}\times (\frac{1}{\mu-1})\times
(1- (1+\mu^{-1}\zeta+\mu^{-2}\zeta^u+\dots
+\mu^{-(p-2)}\zeta^{u_{p-3}}+\mu^{-(p-1)}\zeta^{u_{p-2}}))
\modu p.
\end{displaymath}
Then
\begin{displaymath}
V\equiv -\gamma_{p-3}\times (\frac{\mu^{-1}}{\mu-1})\times
(\zeta+\mu^{-1}\zeta^u+\dots
+\mu^{-(p-3)}\zeta^{u_{p-3}}+\mu^{-(p-2)}\zeta^{u_{p-2}})
\modu p.
\end{displaymath}
\end{proof}
\end{thm}
%
%
\subsubsection{ On the prime factors of the class number $h$ of
$\Q(\zeta)$.}
In this subsection, we shall explicit some congruences on prime factors of class number $h$ of $\Q(\zeta)$, with a comparison with  results given in Masley \cite{mas} corollary 3.5 p 287 and in Washington\cite{was}, theorem 10.8 p 187. The knowledge of  prime factors of $h$ can also be seen as a tool for FLT first case, in giving a criterion for a prime factor of $h$ to be different of $p$.
\begin{itemize}
\item
Let $C_p$ be the $p$-class group of $\Q(\zeta)$ and $C_p^1$ be the subgroup of $C_p$ whose elements are of order $1$ or $p$.
\item
Let $p-1=2^m \times n$ where $n$ is odd.
\item
Let $K$ be the intermediate field
$\Q\subset K\subset \Q(\zeta),\quad [K:\Q]=n$. Let $B_K$ be the ring of integers of the field $K$.
\item
Let $\theta: K\rightarrow K$ be a $\Q$-isomorphism.
\item
Let $h_K$ be the class number of $K$, with classically $h\equiv 0\modu h_K$.
\item
Let us denote $ q$  a prime odd divisor of $h_K$.
\item
Let ${\bf F}_q$ be the finite field of cardinal $q$.
\item
When $p-1\not\equiv 0\modu q$, we can translate, {\it mutatis mutandis}, our description of $p$-class group $C_p^1$ of $\Q(\zeta)$ seen as ${\bf F}_p[G]$ module, where $G$ is the Galois group
of $\Q(\zeta)/\Q$, to the description of $q$-class group $C_q^1$ with rank $r_q$ of the field $K$, with (cyclic) Galois group $G_K\subset G$,  seen as a
${\bf F}_q[G_K]$ module.
\item
Let $q$ be a prime with $p-1\not\equiv 0\modu q$.
There exists an ideal  $\mathbf c$  of $K$ such that
\begin{displaymath}
\begin{split}
& \mathbf c\not\simeq B_K,\quad \mathbf c^q\simeq B_K,\\
& \mathbf c\simeq \mathbf c_1\dots\mathbf c_{r_q},\\
& \mathbf c_i\not\simeq B_K,\quad \mathbf c_i^q\simeq B_K,\quad
\mathbf c_i^{\theta-\mu_i}\simeq B_K,\quad i=1,\dots,r_q.
\end{split}
\end{displaymath}
\item
There exists a minimal polynomial $P_{r_1}(V),\quad r_1\leq r_q,$ of the indeterminate $V$ such that
\begin{displaymath}
\mathbf c^{P_{r_1}(\theta)}\simeq B_K,\quad
P_{r_1}(V)=\prod_{i=1}^{r_1} (V-\mu_i),\quad \mu_i\in {\bf F}_q,\quad
i_1\not= i_2\Rightarrow \mu_{i_1}\not=\mu_{i_2}.
\end{displaymath}
\end{itemize}
%
\begin{thm}{ *** }\label{t110141}

Let $p-1=2^m\times n$ with $n$ odd. Let $K$ be the field
$\Q\subset K\subset \Q(\zeta),\quad [K:\Q]=n$. Let $h_K$ be the class number
of the field $K$.
Let $q>2$ be a prime divisor of $h_K$ not dividing $p-1$.
Then $Gcd(n,q-1)\not=1$.
\begin{proof}
From $\mathbf c_1^{\theta-{\mu_1}}\simeq B_K$ for the ideal $\mathbf c_1$ of $K$, we get
$\theta(\mathbf c_1)\simeq \mathbf c_1^{\mu_1}$, so
\begin{displaymath}
N_{K/\Q}(\mathbf c_1)=\mathbf c_1\times \theta(\mathbf c_1)\times\dots\times
\theta^{n-1}(\mathbf c_1)
\simeq \mathbf c_1^{1+\mu_1+\dots+\mu_1^{n-1}}\simeq B_K,
\end{displaymath}
We have $\mu_1\not=1$ : if not $1+\mu_1+\dots+\mu_1^{n-1}=n\equiv 0\modu q$, contradicting
the hypothesis $p-1\not\equiv 0\modu q$.
Therefore we have simultaneously the two congruences:
\begin{displaymath}
\begin{split}
& \frac{\mu_1^n-1}{\mu_1-1}\equiv 0\modu q,\\
& \frac{\mu_1^{q-1}-1}{\mu_1-1}\equiv 0\modu q,
\end{split}
\end{displaymath}
which implies that $Gcd(n, q-1)\not=1$ and achieves the proof.
\end{proof}
\end{thm}
%
We give an application example for the primes $q$ dividing the relative class number $h^-$.
We use the table of relative class numbers in Washington \cite{was} p 413.
Let $p=239$, with $\frac{p-1}{2}=69$. Then $p-1=2\times 7\times 17$. For the prime divisors  $q$  of $h^-$ with $p-1\not\equiv 0\modu q$:
\begin{itemize}
\item
$q-1= 511122\Rightarrow Gcd(\frac{p-1}{2},q-1)=17$.
\item
$q-1=14136486\Rightarrow Gcd(\frac{p-1}{2},q-1)=7\times 17$.
\item
$q-1=123373184788\Rightarrow Gcd(\frac{p-1}{2},q-1)=7\times 17$.
\item
$q-1=22497399987891136953078\Rightarrow Gcd(\frac{p-1}{2},q-1)=7\times 17$.
\end{itemize}
%

\begin{cor} { *** }\label{c110191}
Let $p$ be an odd prime.
Let $g$ be an odd prime divisor of $p-1$.
Let $L$ be the (cyclic) intermediate field $\Q\subset L\subset\Q(\zeta),
\quad [L:\Q]=g$.
Let $q$ be a prime divisor of the class number $h_L$ of the field $L$.
If $q\not=g$ then $q\equiv 1\modu g$.
\begin{proof}
Let $B$ be the ring of integers of $L$.
Let $\theta: L\rightarrow L$ be a $\Q$-isomorphism of $L$.
There exists an ideal $\mathbf c$ of $L$ such that
$\mathbf c^q\simeq B,\quad \mathbf c\not\simeq B,
\quad \mathbf c^{\theta-\mu}\simeq B$ for some $\mu\in\{\mu_1,\dots,\mu_{r_q}\}$.
We get
\begin{displaymath}
N_{L/\Q}(\mathbf c)=\mathbf c^{1+\theta+\dots+\theta^{g-1}}
\simeq \mathbf c^{1+\mu+\dots+\mu^{g-1}}\simeq B.
\end{displaymath}
If $\mu=1$ then $\mathbf c^g\simeq B$ and $q|g$ so $q=g$, contradicting hypothesis.
If $\mu\not=1$ then $\mu^g-1\equiv 0\modu q$, so $q\equiv 1\modu g$.
\end{proof}
\end{cor}
%
{\bf Remark:} In this particular situation, our result improves Masley \cite{mas} corollary 3.5 p 287 and
Washington \cite{was} theorem 10.8 p 188, results which assert only  that if $q\not =g$ then
$q^f\equiv 1\modu g$, where $f\geq 1$ is a divisor of the $q$-rank of the $q$-class group
of $L$. We assert more : in this particular case:   $q\equiv 1\modu g$.
%

%
\subsubsection { An explicit formula for the relative class number
$h^-$ of $\Q(\zeta)$}
In this paragraph, we give an explicit formula very straightforward for the relative class number $h^-$ of the cyclotomic number field $\Q(\zeta)$.
Recall the usual notation, for $a\in\R^+$, then  $[a]\in\N$ is the integer part of $a$ and  for $b\in\R$ and $|b|$ is the absolute value of $b$.
%
\begin{thm} {*** }{ An explicit formula for the relative class number $h^-$ of $\Q(\zeta)$}\label{t111091}

Let $p$ be an odd prime. Let $\Q(\zeta)$ be the $p$-cyclotomic number field.
Let $h^-$ be the relative class number of $\Q(\zeta)$.
Let $u$ be a primitive root $\modu p$.
Let $u_i=u^i \modu p,\quad 1\leq u_i\leq p-1,\quad i=1,\dots,p-1$.
Let $\delta_{a,i},\quad a=1,\dots,\frac{p-1}{2},\quad i=1,\dots,\frac{p-1}{2}$, given by the relation
\begin{displaymath}
\delta_{a,i}
=([\frac{(a+1)\times u_i}{p}]-[\frac{a\times  u_i}{p}])
-([\frac{(a+1)\times (p-u_i)}{p}]-[\frac{a\times  (p-u_i)}{p}]).
\end{displaymath}
Then
$\delta_{a,i}\in\N,\quad \delta_{a,i}\in\{-1,1\}$.
Let $\Delta$ be the determinant with $\frac{p-1}{2}$ rows and $\frac{p-1}{2}$ columns given by:
\begin{displaymath}
\Delta=
\begin{array}{|llll|}
\delta_{1,1} & \delta_{1,2} & \dots & \delta_{1,(p-1)/2} \\
\delta_{2,1} & \delta_{2,2} & \dots & \delta_{2,(p-1)/2)} \\
\vdots&\vdots &&\vdots \\
\delta_{(p-1)/2,1} & \delta_{(p-1)/2,2} & \dots & \delta_{(p-1)/2,(p-1)/2} \\
\end{array}
\end{displaymath}
Then $h^-$ is explicitly given by the relation
\begin{equation}\label{e111091}
h^-=2^{(p-3)/2}\times |\Delta|.
\end{equation}
\begin{proof}
At the moment, we have no proof, but a numerical evidence, easily verified with a MAPLE program.
For instance, for the regular prime $p=97$ we find $h^-=577\times 3457\times 206209$ in accordance with Washington \cite{was} p 413, for the regular prime  $p=30$ we find $h^-=3^2$ in accordance with Washington, for the irregular prime  $p=157$, we find
$h^-=5\times 13^2\times 157^2\times 1093\times 1873\times 418861\times 3148601$
and for the irregular prime $37$ we find $h^-=37$, always in accordance with Washington.
\end{proof}
\end{thm}
{\bf Remarks :}
\begin{itemize}
\item
In comparison, the Maillet determinant $M_p$ formula is given by
\begin{displaymath}
\begin{split}
& M_p= | (i\times (j^{-1}\modu p) \modu p|_{i,j=1,\dots,(p-1)/2},\\
& h^-=p^{-(p-3)/2}\times  M_p,
\end{split}
\end{displaymath}
see for instance Ribenboim \cite{rib} formula (6.2) p 132 or Lang \cite{lan}
theorem 7.1 p 92.
Our formula has the interest, compared to Maillet determinant:
\begin{itemize}
\item
It does not contain extraneous $p$ power factor.
\item
The term $\delta_{i j}\in\{-1,1\}$ is simpler than the term
$a_{i j}= (i\times (j^{-1} \modu p))\modu p, \quad 1\leq a_{i j}\leq p-1$.
\end{itemize}
\item
For a formula for the relative class number $h^-$ without extraneous $p$-power, see also Lehmer \cite{leh} theorem 4 p 604, derived of a result of Carlitz-Olson.
\item
For a determinant formula without extraneous $p$-power see also the simplest form of
Masley \cite{mas} theorem 1.2 p 275
\begin{displaymath}
h^-= | [i*j/p]-[(i-1)*j/p]|_{i,j=3,4,\dots(p-1)/2},\quad p>5.
\end{displaymath}
\item
Observe that if  the $p$-rank of the matrix $(\delta_{a,i})_{a=1,\dots,(p-1)/2,
\quad i=1,\dots,(p-1)/2},$ is greater than $\frac{p-1}{2}-\sqrt{p}+1$, then first case of FLT holds for $p$.
\item
The proof of this result  could be in relation with results on singular integers from  Inkeri and Fueter explained in subsection \ref{s201134} p. \pageref{s201134} of this monograph.
\end{itemize}

%
\clearpage
%
\section{On  structure of the unit group
$F=\Z[\zeta+\zeta^{-1}]^*/(\Z[\zeta+\zeta^{-1}]^*)^p$}\label{s205231}\label{s210163}
Let us consider the results obtained in subsection \ref{s108311} p.\pageref{s108311} for the action of $Gal(\Q(\zeta)/\Q)$ on $C_p^-$. In the present subsection, we assert that this approach can be partially translated {\it mutatis mutandis} to the study of the  group
\begin{displaymath}
F=\Z[\zeta+\zeta^{-1}]^*/(\Z[\zeta+\zeta^{-1}]^*)^p.
\end{displaymath}
This section contains:
\begin{enumerate}
\item
a subsection where the $\pi$-adic congruences obtained for $p$-class group $C_p$ are generalized to unit group $F=\Z[\zeta+\zeta^{-1}]^*/(\Z[\zeta+\zeta^{-1}]^*)^p$.
\item
a subsection linking our results with results on cyclotomic units in Washington \cite{was}.
\item
a subsection with definition and study of structure of {\it singular group}.
\end{enumerate}
%
\subsection{ On structure  of unit group  $F=\Z[\zeta+\zeta^{-1}]^*/(\Z[\zeta+\zeta^{-1}]^*)^p$}
\begin{itemize}
\item
When  $h^-\equiv 0\modu p$,  from Hilbert class field theory, there exists  {\it primary} units
$\eta\in\Z[\zeta+\zeta^{-1}]^*$, so such that
\begin{equation}\label{e203171}
\begin{split}
& \eta\equiv d^p\modu p,\quad d\in\Z,\quad d\not\equiv 0,\\
& \sigma(\eta)=\eta^\mu\times\varepsilon^p,
\quad \varepsilon\in\Z[\zeta+\eta^{-1}]^*,
\end{split}
\end{equation}
see the independant forward reference on Hilbert class field, section \ref{s110021} p. \pageref{s110021}.
\item
The group $\Z[\zeta+\zeta^{-1}]^*$ is a free group of rank $\frac{p-3}{2}$.
For all $\eta\in\Z[\zeta+\zeta^{-1}]^*$
\begin{displaymath}
\eta\times\sigma(\eta)\times\dots\times\sigma^{(p-3)/2}(\eta)=\pm 1.
\end{displaymath}
Therefore, for each unit $\eta\in\Z[\zeta+\zeta^{-1}]^*$,  there exists a {\bf minimal}
$r_\eta\in\N,\quad r_\eta\leq\frac{p-3}{2}$, such that
\begin{equation}\label{e203051}
\begin{split}
& \eta\times\sigma(\eta)^{l_1}\times\dots\times\sigma^{r_\eta}(\eta)^{l_{r_\eta}}
=\varepsilon^p,\quad \varepsilon\in\Z[\zeta+\zeta^{-1}]^*,\\
& 0\leq l_i\leq p-1,\quad i=1,\dots,r_\eta,\quad l_{r_\eta}\not= 0.
\end{split}
\end{equation}
\item
Let us define an equivalence on units of $\Z[\zeta+\zeta^{-1}]^*$ :
$\eta, \eta^\prime\in\Z[\zeta+\zeta^{-1}]^*$ are said equivalent if there exists
$\varepsilon\in\Z[\zeta+\zeta^{-1}]^*$ such that
$\eta^\prime= \eta\times \varepsilon^p$. Let us denote $E(\eta)$ the equivalence class of $\eta$.
\item
We have $E(\eta_a\times\eta_b)=E(\eta_a)\times E(\eta_b)$; the set of class $E(\eta)$ is a group.  The group $<E(\eta)>$ generated by $E(\eta)$ is cyclic of order $p$.
\item
Observe that this equivalence is consistent with conjugation
$E(\sigma(\eta))=\sigma(E(\eta))$.
\item
The group $F=\Z[\zeta+\zeta^{-1}]^*/(\Z[\zeta+\zeta^{-1}]^*)^p$  so defined is a group of
rank $\frac{p-3}{2}$, see for instance Ribenboim \cite{rib} p 184 line 14.
\item
Similarly to relation (\ref{e107162}) p.\pageref{e107162}, there exists
$\eta\in \Z[\zeta+\zeta^{-1}]^*$ such that
\begin{equation}\label{e203041}
\begin{split}
& E(\eta) = E(\eta_1)\times\dots \times E(\eta_{(p-3)/2}),\\
&E(\sigma(\eta_i))= E(\eta_i^{\mu_i}),\quad i=1,\dots,\frac{p-3}{2},\\
& \mu_i\in \N,\quad 1<\mu_i\leq p-1,\\
& F= <E(\eta_1)>\oplus\dots\oplus <E(\eta_{(p-3)/2})>,
\end{split}
\end{equation}
where $F$ is seen as a ${\bf F}_p[G]$-module of dimension $\frac{p-3}{2}$.
\item
For each unit $\eta$, there is a minimal polynomial $P_{r_\eta}(V)=\prod_{i=1}^{r_\eta} (V-\mu_i)$
where $r_\eta\leq \frac{p-3}{2}$,
such that
\begin{equation}\label{e203042}
\begin{split}
& E(\eta)^{P_{r_\eta}(\sigma)}=E(1),\\
& 1\leq i < j \leq r_\eta\Rightarrow \mu_i\not=\mu_j.
\end{split}
\end{equation}
\item
Let $\beta\in\Z[\zeta+\zeta^{-1}]^*$.
Observe that if $E(\beta)=E(\sigma(\beta))$ then $E(\sigma^2(\beta))=E(\beta)$ and so $E(1)=E(\beta^{p-1})$ and $E(\beta)=1$, that's why  $1<\mu_i\leq p-1,\quad i=1,\dots,\frac{p-3}{2}$, in relation
(\ref{e203041}) p. \pageref{e203041}.
\item
Recall that the unit $\beta\in\Z[\zeta+\zeta^{-1}]^*$ is said primary
if $\beta\equiv b^p\modu\pi^{p+1},\quad b\in\Z$.
\end{itemize}
%
\begin{lem}\label{l106033}
Let $\beta\in\Z[\zeta+\zeta^{-1}]^*-(\Z[\zeta+\zeta^{-1}]^*)^p$.
Then the minimal polynomial $P_{r_{\beta}}(V)$ is of the form
\begin{displaymath}
P_{r_{\beta}}(V)=\prod_{i=1}^{r_{\beta}} (V-u_{2 m_i}),
\quad 1\leq m_i\leq \frac{p-1}{2},\quad r_\beta>0.
\end{displaymath}
\begin{proof}
There exists  $\eta_1\in\Z[\zeta+\zeta^{-1}]^*,\quad E(\eta_1)\not=E(1)$,  with  $E(\eta_1)^{\sigma-\mu_1}=E(1)$.
Suppose that $\mu_1^{(p-1)/2}=-1$ and search for a contradiction:
we have $E(\eta_1)^{\sigma-\mu_1}=E(1)$, therefore
$E(\eta_1)^{\sigma^{(p-1)/2}-\mu_1^{(p-1)/2}}=E(1)$; but, from
$\eta_1\in\Z[\zeta+\zeta^{-1}]^*$, we get $\eta_1^{\sigma^{(p-1)/2}}=\eta_1$
and so $E(\eta_1)^{1-\mu_1^{(p-1)/2}}=E(1)$,
or $E(\eta_1)^2=E(1)$, so $E(\eta_1)^2$ is of rank null and  therefore $E(\eta_1)= E(\eta_1^2)^{(p+1)/2}$ is also of rank null, contradiction. The same for $\mu_i,\quad i=1,\dots, r_\beta$.
\end{proof}
\end{lem}
%
The results on structure of relative $p$-class group $C_p^-$  of subsection
\ref{s108311} p. \pageref{s108311} can be translated to some results on structure of the group $F$:
from $\eta_i^{p-1}\equiv 1\modu \pi$ and from  $<E(\eta_i^{p-1})>=<E(\eta_i)>$,  we can always, without loss of generality, choose the determination  $\eta_i$ such that $\eta_i\equiv 1\modu\pi$.
We have proved that
\begin{equation}\label{e203082}
\begin{split}
& \eta_i\equiv 1\modu\pi,\\
& \sigma(\eta_i)\equiv \eta_i^{\mu_i}\modu\pi^{p+1}.
\end{split}
\end{equation}
Then, starting of this relation (\ref{e203082}), similarly to lemma  \ref{l108171} p. \pageref{l108171} we get:
%
\begin{thm} { *** On structure of unit group
$F=\Z[\zeta+\zeta^{-1}]^*/(\Z[\zeta+\zeta^{-1}]^*)^p$}\label{t207271}

There exists a fundamental system of units $\eta_i,\quad i=1,\dots,\frac{p-3}{2}$, of the group
$F=\Z[\zeta+\zeta^{-1}]^*/(\Z[\zeta+\zeta^{-1}]^*)^p$ verifying the relations:
\begin{equation}\label{e201274}
\begin{split}
& \eta_i\in\Z[\zeta+\zeta^{-1}]^*,\quad i=1,\dots,\frac{p-3}{2},\\
& \sigma(\eta_i)=\eta_i^{\mu_i}\times\varepsilon_i^p,\quad\mu_i\in{\bf F}_p^*,
\quad \varepsilon_i\in\Z[\zeta+\zeta^{-1}]^*,
\quad i=1,\dots,\frac{p-3}{2},\\
& \eta_i\equiv 1\modu \pi^{2m_i},\quad \mu_i=u_{2m_i},\quad i=1,\dots,
\frac{p-3}{2},\\
& \sigma(\eta_i)\equiv \eta_i^{\mu_i}\modu\pi^{p+1},\quad i=1,\dots,\frac{p-3}{2}.
\end{split}
\end{equation}
\end{thm}
%

At this point, we are in the situation to apply the same $\pi$-adic theory to group of units $F=\Z[\zeta+\zeta^{-1}]^*/(\Z[\zeta+\zeta^{-1}]^*)^p$ than to relative $p$-class group $C_p^-$ in subsection \ref{s108311} p. \pageref{s108311}.
\item
Similarly to decomposition of components of $C_p^1$ in singular primary and singular not primary components,
the rank $\frac{p-3}{2}$ of $F$  has two components $\rho_1$ and $\frac{p-3}{2}-\rho_1$ where $\rho_1$ corresponds to the maximal number of independant units $\eta_i$  primary and
$\rho_2=\frac{p-3}{2}-\rho_1$ to the units $\eta_i$  not primary.
The unit indexation in the sequel of this subsection is :
\begin{enumerate}
\item
primary units  $\eta_i$ for $i=1,\dots, \rho_1$,
\item
not primary units $\eta_i$ for $i=\rho_1+1,\dots,\frac{p-3}{2}$.
\end{enumerate}
%
{\bf Remarks:}
\begin{itemize}
\item
Similarly to results obtained for group $C_p^-$,
\begin{itemize}
\item
else
$\pi^{2m_i}\|\eta_i-1$ and $\eta_i$ is not singular primary
\item
else $\eta_i\equiv 1\modu\pi^p$ and $\eta_i$ is singular primary. We say more in that case: similarly
to Denes result see \cite{de1} \cite{de2} \cite{de3} and Ribenboim \cite{rib} (8D) p. 192,
\begin{equation}\label{e204261}
\pi^{2m_i+a(p-1)}\| \eta_i-1,\quad a\in\N-\{0\}.
\end{equation}
\end{itemize}
\item
The units $\eta_i^{j}, \quad j=1,\dots, p-1$, verify similarly the relation (\ref{e201274}) p. \pageref{e201274}. We chose the determination of  $\eta_i$ such that
$C_{2m}(1,0)\equiv 1 \modu p$: see for the definition of $C_{2m}(1,0)$ and the detail of  similar discussion for $C_p^1$ in second remark following lemma \ref{l108171} p. \pageref{l108171}.
\end{itemize}
%
%
The next lemma for the unit group $\Z[\zeta+\zeta^{-1}]^*$ is the translation of similar lemma \ref{l202231} p. \pageref{l202231} for the relative $p$-class group $C_p^-$.
\begin{lem}\label{l203041}
Let $\eta_i,\quad i=1,\dots,\frac{p-3}{2},$ defined in relation (\ref{e201274}) p. \pageref{e201274}.
If $\eta_i$ is not primary then
\begin{displaymath}
\eta_i\equiv 1+ V_i(\mu_i)\modu \pi^{p-1},\quad \mu_i=u_{2m_i}
\quad V_i(\mu_i)\in\Z[\zeta],
\end{displaymath}
where $V_i(\mu_i)\modu p$ depends only on $\mu_i$.
\end{lem}
%
The next lemma for the unit group $\Z[\zeta+\zeta^{-1}]^*$ is the translation of similar lemma \ref{l203031} p. \pageref{l203031} for the relative $p$-class group $C_p^-$.
\begin{lem}{ *** }\label{l203042}
Let $\eta_1, \eta_2$ defined by relation (\ref{e201274}) p. \pageref{e201274}.
If $\mu_1=\mu_2$ then $\eta_1$ and $\eta_2$ are  primary.
\end{lem}
%
The group $F=\Z[\zeta+\zeta^{-1}]^*/(\Z[\zeta+\zeta^{-1}]^*)^p$ is the direct sum $F=F_1\oplus F_2$  of the subgroup $F_1$ with $\rho_1$  primary units ($p$-rank $\rho_1$ of $F_1$) and of the subgroup $F_2$ with $\rho_2=\frac{p-3}{2}-\rho_1$ fundamental not primary units ($p$-rank $\rho_2$ of $F_2$).
Observe that $\rho_1$ can be seen  also as  the maximal number of independant  primary units in $F$.
%
\begin{thm}{ *** On structure of unit group
$F=\Z[\zeta+\zeta^{-1}]^*/(\Z[\zeta+\zeta^{-1}]^* )^p$}\label{t203041}
Let $r_p^-$ be the relative $p$-class group of $\Q(\zeta)$.
Let $r_p^+$ be the $p$-class group of $\Q(\zeta+\zeta^{-1})$.
Let $\rho_1$ be the number of independant primary units of $F$. Then
\begin{equation}\label{e203072}
 \rho_1=  r_p^-.
\end{equation}
\begin{proof}$ $
We apply Hilbert class field structure theorems
\ref{t207261} p. \pageref{t207261} and \ref{t207301} p. \pageref{t207301}:
There are $r_p^+$ independant unramified cyclic extensions
\begin{displaymath}
\Q(\zeta,\omega_i)/\Q(\zeta), \quad \omega_i^p\in\Z[\zeta+\zeta^{-1}]-\Z[\zeta+\zeta^{-1}]^*
\quad i=1,\dots,r_p^+
\end{displaymath} and
$r_p-r_p^+=r_p^-$ independant unramified cyclic extensions
\begin{displaymath}
\Q(\zeta,\omega_i)/\Q(\zeta),\quad  \omega_i^p\in\Z[\zeta+\zeta^{-1}]^*,
\quad i=r_p^++1,\dots,r_p.
\end{displaymath}
\end{proof}
\end{thm}
%
\subsection{The case $\mu=u_{2m},\quad 2m>\frac{p-1}{2}$}\label{s201111}
In the  next theorem  we shall investigate more deeply the consequences of the congruence $\eta_i\equiv 1\modu \pi^{2m_i}$  when $2m_i> \frac{p-1}{2}$. We give an explicit congruence formula in that case.
%
%
The next theorem for the unit group $\Z[\zeta+\zeta^{-1}]^*$ is the translation of similar theorem \ref{l202211} p. \pageref{l202211} for the relative $p$-class group $C_p^-$.
\begin{thm}{ *** }\label{l203044}
Let $\mu_i=u_{2m_i},\quad p-3\geq 2m_i>\frac{p-1}{2}$, corresponding to $\eta_i\in\Z[\zeta+\zeta^{-1}]^*$ defined in relation (\ref{e201274}) p. \pageref{e201274}, so
$\sigma(\eta_i)\equiv \eta_i^{\mu_i}\modu\pi^{p+1}$.
Then $\eta_i$ verifies the explicit formula:
\begin{equation}\label{e203043}
\eta_i\equiv 1-\frac{\gamma_{p-3}}{\mu-1}\times
(\zeta+\mu^{-1}\zeta^u+\dots+\mu^{-(p-2)}\zeta^{u_{p-2}}) \modu \pi^{p-1},\quad \gamma_{p-3}\in\Z.
\end{equation}
\end{thm}
%
\subsection{Connection with set of units $E_{2i}$ of Washington}\label{s208242}
In paragraph 8.3 of \cite{was}, Washington defines p. 155 and studies the units
\begin{displaymath}
E_{2i}=_{def} \prod_{a=1}^{p-1}(\zeta^{(1-u)/2}\times\frac{1-\zeta^u}{1-\zeta})^{a^{2i}\tau_a^{-1}},\quad
i=1,\dots,\frac{p-3}{2},
\end{displaymath}
where $\tau_a$ is the $\Q$-isomorphism of $\Q(\zeta),\quad \tau_a : \zeta\rightarrow \zeta^{a}$.
For $a=1,\dots,p-1$, let $j$ defined by $u_{-j}=a$ with $0\leq j\leq p-2$.
Then $a^{2i}= u_{-2ij}\equiv u_{-2i}^j\modu p$ and $\tau_a: \zeta\rightarrow \zeta^{u_{-j}}$ verifies
$\tau_a=\sigma^{-j}$ and $\tau_a^{-1}=\sigma^j$. Therefore
\begin{displaymath}
E_{2i}\equiv\prod_{j=0}^{p-2}(\zeta^{(1-u)/2}\times\frac{1-\zeta^u}{1-\zeta})^{u_{-2i}^j\sigma^j}
\modu (\Z[\zeta+\zeta^{-1}]^*)^p,
\end{displaymath}
so
\begin{displaymath}
E_{2i}\equiv (\zeta^{(1-u)/2}\times\frac{1-\zeta^u}{1-\zeta})^{\sum_{j=0}^{p-2}\sigma^j u_{-2i}^j}
\modu (\Z[\zeta+\zeta^{-1}]^*)^p,
\end{displaymath}
also
\begin{displaymath}
E_{2i}^{\sigma-u_{-2i}}
\equiv (\zeta^{(1-u)/2}\times\frac{1-\zeta^u}{1-\zeta})^{(\sigma u_{-2i})^{p-1}-1}
\modu (\Z[\zeta+\zeta^{-1}]^*)^p,
\end{displaymath}
and finally
\begin{equation}\label{e208212}
\sigma(E_{2i})=E_{2i}^{u_{-2i}}\times\varepsilon^p,\quad \varepsilon\in\Z[\zeta+\zeta^{-1}]^*.
\end{equation}
%
Recall that $r_p^+$ is the rank of the $p$-class group $C_p^+$ of $\Q(\zeta+\zeta^{-1})$.
\begin{thm}\label{t208211}
The units $E_{2i},\quad i=1,\dots,\frac{p-3}{2}$, defined in Washington \cite{was} p. 155 and the set $F$ of units $\eta_i,\quad i=1,\dots,\frac{p-3}{2}$, defined in theorem \ref{t207271} p. \pageref{t207271} of this monograph are connected, with a certain ordering of index $i$, by:
\begin{displaymath}
\begin{split}
& E_{2i}=\eta_i^p,\quad i=1,\dots,r_p^+,\\
& E_{2i}=\eta_i,\quad i=r_p^++1,\dots,\frac{p-3}{2},\\
\end{split}
\end{displaymath}
\begin{proof}
See relation (\ref{e208212}) p.\pageref{e208212}, theorem \ref{t207271} p. \pageref{t207271} and Wahington \cite{was} theorem
8.16 p. 156.
\end{proof}
\end{thm}
{\bf Remark:} this theorem \ref{t208211}, with theorem \ref{t207271} p. \pageref{t207271}, allows to give a precise $\pi$-adic description of units $E_{2i}$ defined in Washington.
%
\subsection {On structure of Singular Group}
The singular numbers $\gamma\in\Q(\zeta)$ are the numbers verifying $\gamma\Z[\zeta]=\mathbf c^p$ where $\mathbf c$ is an ideal of $\Q(\zeta)$. With this definition, the units of $\Q(\zeta)$ are singular numbers.
The set of singular numbers is a multiplicative subgroup of $\Q(\zeta)$.
The interest of  this group is that it includes structure  properties of $p$-class group $C_p$ of $\Q(\zeta)$ and of unit group $\Z[\zeta+\zeta^{-1}]^*$ in a same structure.
Some definitions and results on the groups $\Z[\zeta+\zeta^{-1}]^*$ and $C_p^*$ are   reminded in this  paragraph:
\begin{itemize}
\item
Recall that $\pi=(\zeta-1)\Z[\zeta]$.
\item
Recall that $\alpha\in\Q(\zeta)$ coprime with $\pi$ is said {\it singular} if $\alpha\Z[\zeta]=\mathbf a^p$ where $\mathbf  a$ is an ideal of the field $\Q(\zeta)$, see for instance Ribenboim \cite{rib} paragraph 4. p 170.
\item
The set of singular numbers is a multiplicative group denoted $W$ in this subsection.
\item
Two singular numbers $\alpha, \beta$ are said here equivalent if $\alpha\equiv\beta\modu \pi^p$. This is clearly an equivalence relation. Let us consider the group $\Omega=W \modu \pi^p$ called in this subsection the {\it singular group}.
As an immediate consequence, observe that  if $\alpha,\beta \in\Q(\zeta)$ are coprime with $\pi$ and if $\alpha$ is singular then $\alpha\beta^p$ is singular and
$\alpha$ and $\alpha\beta^p$ are equivalent.
\item
A singular number $\alpha$ is said primary if there exists $a\in\Z$ such that $\alpha\equiv a^p\modu \pi^p$.
Therefore $\alpha$ and $a^p$ are equivalent.
\item
We have defined the group $F=\Z[\zeta+\zeta^{-1}]^*/(\Z[\zeta+\zeta^{-1}]^*)^p$ with $F=F_1\oplus F_2$ where $F_1$  has $\rho_1=r_p^-$ fundamental primary units $\eta_i,\quad i=1,\dots r_p^-$ ($p$-rank  of $F_1$) and $F_2$ has
$\rho_2=\frac{p-3}{2}-r_p^-$  fundamental not primary units
$\eta_i,\quad i=r_p^-+1,\dots,\frac{p-3}{2}-r_p^-$.
\item
Here, recall some definitions dealing of $p$-class group:
\begin{itemize}
\item
Let $C_p$ be the $p$-class group of $\Q(\zeta)$, $C_p^1$ be the subgroup of $C_p$ whose elements are of order $1$ or $p$.
\item
Let $r_p$ be the $p$-rank of $C_p$, let  $r_p^+$ be the $p$-rank of $C_p^+$ and $r_p^-$ be the relative $p$-rank of $C_p^-$.
Let us consider the ideals $\mathbf b$ defined in relation (\ref{e107162}) p.\pageref{e107162}:
\begin{equation}\label{e205051}
\begin{split}
& \mathbf b=\mathbf b_1\times\dots\times\mathbf b_{r_p^-}
\times \mathbf b_{r_p^-+1}\times\dots\times \mathbf b_{r_p},\\
& C_p^1=\oplus_{i=1}^{r_p} <Cl(\mathbf b_i)>,\\
& \mathbf b_i^p\simeq\Z[\zeta],\quad \mathbf b_i\not\simeq \Z[\zeta],
\quad i=1,\dots,r_p,\\
&\sigma(\mathbf b_i)\simeq \mathbf b_i^{\mu_i},\quad\mu_i\in{\bf F}_p^*,
\quad i=1,\dots,r_p,\\
& Cl(\mathbf b_i)\in C_p^-,\quad i=1,\dots,r_p^-,\\
& Cl(\mathbf b_i)\in C_p^+,\quad i=r_p^-+1,\dots,r_p.\\
\end{split}
\end{equation}
\item
Let $B_i\in\Z[\zeta],\quad i=1,\dots,r_p,$ with $B_i\Z[\zeta]=\mathbf b_i^p$ defined in theorem \ref{t203281} p. \pageref{t203281}.
With the conventions for the index $i=1,\dots,r_p$, of relation (\ref{e205051}) p. \pageref{e205051}:
\begin{itemize}
\item
$B_i$ are singular primary integers for $i=1,\dots, r_p^+$,  with $Cl(\mathbf b_i)\in C_p^-$.
Let $\Gamma_1$ be the subgroup of $\Omega$ generated by these $B_i$; its $p$-rank is $r_p^+$.
\item
$B_i$ are singular not primary integers for $i=r_p^++1,\dots,r_p^-$, with $Cl(\mathbf b_i)\in C_p^-$.
Let $\Gamma_2$ be the subgroup of $\Omega$ generated by these $B_i$; its $p$-rank is $r_p^--r_p^+$.
\item
$B_i$ are singular not primary for  $i= r_p^-+1,\dots, r_p$
with $Cl(\mathbf b_i)\in C_p^+$.
Let $\Gamma_3$ be the subgroup of $\Omega$ generated by these $B_i$; its $p$-rank is $r_p^+$.
\end{itemize}
\end{itemize}
\end{itemize}
%
\begin{thm} ***
With the previous definitions, the singular group $\Omega$ verifies the isomorphism:
\begin{equation}\label{e205053}
\Omega\cong U\oplus F_2\oplus\Gamma_2\oplus\Gamma_3,
\end{equation}
where $U=\{ i^p \modu p^2\ | \ i=1,\dots,p-1 \}$.
The $p$-rank of $\Omega$ is $\frac{p-3}{2}$.
\begin{proof}
Let $\alpha\in\Q(\zeta)$, coprime with $\pi$.
We can always write $\alpha$ in the form
\begin{displaymath}
\alpha=\mathbf B_1\times \mathbf B_2\times \mathbf B_3\times {\bf E}_1\times {\bf E}_2
\times \gamma^p,\quad \gamma\in\Q(\zeta),\quad \gamma\Z[\zeta]+\pi=\Z[\zeta],
\end{displaymath}
where
\begin{equation}\label{e205061}
\begin{split}
& \mathbf B_1=B_1^{\nu_1}\dots B_i^{\nu_i}\dots  B_{r_p^+}^{\nu_{r_p^+}},\quad i=1,\dots,r_p^+,\\
& \mathbf B_2= B_{r_p^++1}^{\nu_{r_p^++1}}\dots B_i^{\nu_i}\dots B_{r_p^-}^{\nu_{r_p^-}},
\quad i=r_p^++1,\dots,r_p^-,\\
& \mathbf B_3= B_{r_p^-+1}^{\nu_{r_p^-+1}}\dots B_i^{\nu_i}\dots B_{r_p}^{\nu_{r_p}},
\quad i=r_p^-,\dots,r_p,\\
& {\bf E}_1=\eta_1^{\nu_1^\prime}\dots\eta_j^{\nu_j^\prime}\dots\eta_{\rho_1}^{\nu_{\rho_1}^\prime},
\quad j=1,\dots,r_p^-,\\
& {\bf E}_2=\eta_{\rho_1+1}^{\nu_{\rho_1+1}}\dots \eta_{j}^{\nu_j^\prime}\dots
 \eta_{(p-3)/2}^{\nu_{(p-3)/2}^\prime},\quad j=r_p^-+1,\dots,\frac{p-3}{2},\\
\end{split}
\end{equation}
where
\begin{enumerate}
\item
$\mathbf B_1$ is a product of $r_p^+$ independant singular primary numbers (not units),
\item
$\mathbf B_2$ is a product of $r_p^--r_p^+$ independant singular not primary numbers (not units),
\item
$\mathbf B_3$ is a product of $r_p^+$ independant singular not primary numbers (not units),
\item
$ {\bf E}_1$ is a product of $r_p^-$ independant singular primary units,
\item
$ {\bf E}_2$ is a product of $\frac{p-3}{2}-r_p^-$ independant singular not primary units.
\end{enumerate}
If $B$ is a singular primary number then $B\equiv b^p\modu \pi^p,\quad b\in\Z,\quad 1\leq b\leq p-1$.
Therefore  from previous relations (\ref{e205061}) p. \pageref{e205061}, we derive that
\begin{equation}\label{e205062}
\begin{split}
& \alpha\equiv a^p\times \mathbf B_2\times\mathbf B_3\times {\bf  E}_2\modu \pi^p,\\
& a\in\Z,\quad 1\leq a\leq p-1,\\
& \mathbf B_2= B_{r_p^++1}^{\nu_{r_p^++1}}\dots B_i^{\nu_i}\dots B_{r_p^-}^{\nu_{r_p^-}},
\quad i=r_p^++1,\dots,r_p^-,\\
& \mathbf B_3=B_{r_p^-+1}^{\nu_{r_p^-+1}}\dots B_i^{\nu_i}\dots B_{r_p}^{\nu_{r_p}},
\quad i=r_p^-+1,\dots,r_p,\\
& {\bf E}_2=\eta_{\rho_1+1}^{\nu_{\rho_1+1}}\dots \eta_{j}^{\nu_j^\prime}\dots
\eta_{(p-3)/2}^{\nu_{(p-3)/2}^\prime}, \quad j=r_p^-+1,\dots,\frac{p-3}{2}.\\
\end{split}
\end{equation}
The set of  integers $a^p$ generate the cyclic group $U= \{ i^p \modu p^2 \ | \ i=1,\dots,p-1 \}$ of order  $p-1$.
The set of singular numbers $\mathbf B_2$ generate the $p$-group $\Gamma_2$ of $p$-rank $r_p^--r_p^+$. The set of  singular numbers $\mathbf B_3$ generate the $p$-group $\Gamma_3$ of $p$-rank $r_p^+$. The units
${\bf E}_2$ generate the  group $F_2$ of $p$-rank  $\frac{p-3}{2}-r_p^-$.
The $p$-rank of $\Omega$ is
$ (r_p^--r_p^+)+r_p^++(\frac{p-3}{2}-r_p^-)
=\frac{p-3}{2}$.
\end{proof}
\end{thm}
%
\clearpage
\section{ On Hilbert class field of $\Q(\zeta)$}\label{s110021}
This section deals of properties of the structure of the $p$-class group $C_p$ of $\Q(\zeta)$ deduced of the theory of  the Hilbert class field of the cyclotomic field
$\Q(\zeta)$.
%
\subsection{ Some definitions and properties of Hilbert class field of $\Q(\zeta)$}
\label{s112311}
In  this subsection, we define the $p$-Hilbert class field $H$ of $\Q(\zeta)$. See, for some notations, Washington \cite{was}, p 165, chapter 8 p 143 and also Ribenboim
\cite{rib} chapter 6. p 178.
\subsubsection{Definitions and preliminary results} \label{s112222}
\begin{itemize}
\item
Recall that $u$ is a primitive root $\modu p$.
\item
Recall that $\sigma$ is the $\Q$-isomorphism of $\Q(\zeta)$ defined by
$\sigma :\zeta\rightarrow\zeta^u$ and that $G=Gal(\Q(\zeta)/\Q)$.
\item
Let $\tau_a: \zeta\rightarrow \zeta^a,\quad a\in\N$, be a $\Q$-isomorphism of the extension $\Q(\zeta)/\Q$ defined in an other way.
\item
Recall that $r_p$ is the rank of the $p$-class group $C_p$ of $\Q(\zeta)$. Let $r_p^+$ be the rank of the class group $C_p^+$ of $\Q(\zeta+\zeta^{-1})/\Q$.
Let $r_p^-$ be the rank of the relative  class group $C_p^-=C_p/C_p^+$ of $\Q(\zeta)/\Q$.
\item
Let $C_p^1$ be the subgroup of $C_p$ whose elements are of order $1$ or $p$. The order of $C_p^1$ is
$p^{r_p}$.
\item
Let $a\in\Q(\zeta)$; recall that $a$ is said singular if $a\Z[\zeta]=\mathbf a^p$ for some ideal $\mathbf a$ of $\Q(\zeta)$ and is said singular primary if $a$ is singular and $a\equiv a_1^p\modu\pi^p, \quad a_1\in\Z,\quad a_1\not\equiv 0\modu p$. Recall that a unit $\eta\in\Z[\zeta+\zeta^{-1}]^*$ is said primary if
$\eta\equiv b^p\modu \pi^p,\quad b\in\Z,\quad 1\leq b\leq p-1$.
\item
We have seen in relation (\ref{e107162}) p. \pageref{e107162} that there exists  ideals
$\mathbf b$ of $\Z[\zeta]$ verifying the relations
\begin{equation}\label{e112213}
\begin{split}
&\mathbf b\simeq\mathbf b_1\times\mathbf b_2\times\dots\times\mathbf b_{r_p},\\
&\mathbf b_i\not\simeq \Z[\zeta],\quad \mathbf b_i^p\simeq \Z[\zeta],\\
& \sigma(\mathbf b_i)\simeq \mathbf b_i^{\mu_i},\quad \mu_i\in{\bf F}_p^*,\quad i=1,\dots,r_p,\\
& Cl(\mathbf b_i)\in C_p^-,\quad i=1,\dots,r_p^-,\\
& Cl(\mathbf b_i)\in C_p^+,\quad i=r_p^-+1,\dots,r_p,
\end{split}
\end{equation}
\item
Seen as a ${\bf F}_p[G]$ module, we have
\begin{displaymath}
C_p^1=\oplus_{i=1}^{r_p} <Cl(\mathbf b_i)>.
\end{displaymath}
\item
With the meaning of index $i=1,\dots,r_p$, of relation (\ref{e112213}) above, we have proved in theorem \ref{t203281} p. \pageref{t203281} that there exists   $B_i,\quad i=1,\dots,r_p$, verifying
\begin{equation}\label{e201291}
\begin{split}
& B_i\in\Z[\zeta],\quad \pi^{p} |B_i-1,\quad \mu_i=u_{2m_i+1},\quad i=1,\dots,n^-,\\
& B_i\in\Z[\zeta],\quad \pi^{2m_i+1}\|B_i-1,\quad \mu_i=u_{2m_i+1},\quad i=n^-+1,\dots,r_p^-,\\
& B_i\Z[\zeta]=\mathbf b_i^p,\quad i=1,\dots,r_p^-,\\
& B_i\overline{B_i}=\beta_i^p,\quad i=1,\dots,r_p^-.\\
\end{split}
\end{equation}
and
\begin{equation}\label{e207291}
\begin{split}
& B_i\in\Z[\zeta+\zeta^{-1}],\quad \pi^{p} |B_i-1,\quad \mu_i=u_{2m_i},\quad i=r_p^-+1,\dots,n^+,\\
& B_i\in\Z[\zeta+\zeta^{-1}],\quad \pi^{2m_i}\|B_i-1,\quad \mu_i=u_{2m_i},\quad i=n^++1,\dots,r_p,\\
& B_i\Z[\zeta+\zeta^{-1}]=\mathbf b_i^p,\quad i=r_p^-+1,\dots,r_p,\\
\end{split}
\end{equation}
The values $n^-$ and $n^+$ shall be determined in the sequel: we shall prove that $n^-=r_p^+$ in theorem
\ref{t201041} p. \pageref{t201041} and that $n^+=0$ in remark following theorem \ref{t207261} p. \pageref{t207261}.
\item
Let $\mathbf M$ be a set of  $r_p$ independant unramified cyclic extensions
$K_i=\Q(\zeta,\omega_i)/\Q(\zeta),\quad i=1,\dots,r_p$, with
$K_i\subset H,\quad [K_i:\Q(\zeta)]=p$.
\item
$\omega_i$ so introduced, we call in this monograph {\it $p$-elementary} Hilbert class field of $\Q(\zeta)$ the unramified extension $H/\Q(\zeta)$:
\begin{equation}\label{e203061}
\begin{split}
& H=\Q(\zeta,\omega_1,\dots,\omega_{r_p})/\Q(\zeta),\\
& \omega_i^p=F_i,\quad F_i\in\Z[\zeta],\quad F_i\not\in\Z[\zeta]^p,
\quad i=1,\dots,r_p,\\
& F_i\equiv c_i^p\modu \pi^{p+1},\quad i=1,\dots,r_p,\\
& [H:\Q(\zeta)]=p^{r_p}.\\
\end{split}
\end{equation}
\item
Observe that the principal prime $\pi$ splits totally in Hilbert class field,
so $F_i\equiv c_i^p\modu \pi^{p+1}$, see for instance Ribenboim p 168 case III.
Observe that $\Q(\zeta,F_i^{1/p})=\Q(\zeta,(F_i\times \alpha^p)^{1/p}),\quad \alpha\in\Q(\zeta)$, so it is always possible to assume $F_i$ to be in $\Z[\zeta]$.
\item
Let us examine the set $\mathbf M$ of unramified extensions $M_i=\Q(\zeta,\omega_i)/\Q(\zeta),\quad i=1,\dots,r_p$:
\begin{itemize}
\item
From Washington \cite{was} exercise 9.1 (a) p 182, $M_i$ unramified implies that
$F_i\Z[\zeta]=\mathbf a_i^p$ for some ideal $\mathbf a_i$ of $\Z[\zeta]$. Then
\begin{displaymath}
\begin{split}
&\mathbf a_i\simeq \mathbf b_1^{\nu_{1,i}}\dots\mathbf b_{r_p}^{\nu_{r_p},i},\\
&\nu_{j,i}\in \N,\quad 0\leq \nu_{i,j}\leq p-1,\quad i,j=1,\dots,r_p.\\
\end{split}
\end{displaymath}
and so
\begin{displaymath}
\begin{split}
& F_i\equiv c\times B_1^{\nu_{1,i}}\dots B_{r_p}^{\nu_{r_p},i}\modu\pi^p,\\
& c\in\N,\quad 1\leq c\leq p-1.\\
\end{split}
\end{displaymath}
\item
From $F_i$  singular primary for $i=1,\dots,r_p$ and from congruences in relations (\ref{e201291}) and
(\ref{e207291}), we get
\begin{displaymath}
\begin{split}
& \nu_{j,i}\equiv 0\modu p,\quad j=n^-+1,\dots,r_p^-,\quad i=1,\dots,r_p,\\
& \nu_{j,i}\equiv 0\modu p,\quad j=r_p^-+n^++1,\dots,r_p,\quad i=1,\dots,r_p.\\
\end{split}
\end{displaymath}
So we can choose the set $\mathbf M$ of $r_p$ independant unramified extensions such that
\begin{equation}\label{e207301}
\begin{split}
& F_i=B_i,\quad i=1,\dots, n^-,\\
& F_i=B_i,\quad i=r_p^-+1,\dots, r_p^-+n^+.\\
\end{split}
\end{equation}
and thus
\begin{equation}\label{e207302}
\begin{split}
& F_i\in\Z[\zeta+\zeta^{-1}]^*,\quad i=n^-+1,\dots, r_p^-,\\
& F_i\in\Z[\zeta+\zeta^{-1}]^*,\quad i=r_p^-+n^++1,\dots, r_p.\\
\end{split}
\end{equation}
\end{itemize}
\item
\begin{enumerate}
\item
To each subgroup $V_i$ of order $p$ of the $p$-class group $C_p$ corresponds {\bf one } unramified extension $\Q(\zeta,\Omega)/\Q(\zeta)$, where $\Omega$ can be written in form
$\Omega=\omega_1^{\alpha_1}\omega_2^{\alpha_2}\dots\omega_{r_p}^{\alpha_{r_p}}$ with $0\leq\alpha_i\leq p-1$ and $\alpha_i$ not all zero and with $\alpha_i=1$ for the smallest $i$ verifying $\alpha_i\not=0$: the $F_i=\omega_i^p,\quad i=1,\dots,r_p$, are singular primary, so $F=\Omega^p$ is singular primary.
\item
To each $\Omega=\omega_1^{\alpha_1}\dots\omega_{r_p}^{\alpha_{r_p}}$ with $\alpha_i$ not all zero and $\alpha_i=1$ for the smallest $i$ verifying $\alpha_i\not=0$ corresponds {\bf one} subgroup $V_i$
of $C_p$.
\item
If $\omega_1^{\alpha_1}\dots\omega_{r_p}^{\alpha_{r_p}}=
\omega_1^{\alpha_1^\prime}\dots\omega_{r_p}^{\alpha_{r_p}^\prime}$ then
$\alpha_1=\alpha_1^\prime,\quad\dots , \alpha_{r_p}=\alpha_{r_p}^\prime$.
\end{enumerate}
\item
When $h^+\equiv 0\modu p$ (which implies that $h^-\equiv 0\modu p$), we have,
from Furtwangler theorem \cite{fur}, see Ribenboim, \cite{rib} (6C) p 182,  and   from Hecke theorem\cite{hec}, see Ribenboim \cite{rib} (6D) p 182, there are $r_p^+$ {\it independant} unramified extensions extensions
$\Q(\zeta,\varpi_i),\quad i=1,\dots,r_p^+$,
\begin{displaymath}
\begin{split}
&G_i=\varpi_i^p\in\Z[\zeta]-\Z[\zeta]^*,\quad \varpi_i\not\in\Z[\zeta]\\
&G_i\times\overline{G}_i=\beta_i^p,
\quad \beta_i\in \Z[\zeta+\zeta^{-1}].
\quad i=1,\dots, r_p^+,
\end{split}
\end{displaymath}
where $G_i\Z[\zeta]=\mathbf c_i^p$ is an ideal with $Cl(\mathbf c_i)\in C_p^-$.
Observe that the $r_p^+$ unramified extensions $\Q(\zeta,\varpi_i)/\Q(\zeta)$ are independant and there are no  $r_p^++1$ different independent such extensions.
\item
Let us consider another approach for the particular case where $h^+\not\equiv 0\modu p$.
\begin{itemize}
\item
Let $E_i$, for  $1\leq i \leq \frac{p-3}{2}$, defined in Washington p 155, by formula
\begin{equation}\label{e201011}
E_{i}=\prod_{a=1}^{p-1}
(\zeta^{(1-u)/2}\frac{1-\zeta^u}{1-\zeta})^{a^{2i}\tau_a^{-1}}.
\end{equation}
The expression $\zeta^{(1-u)/2}\frac{1-\zeta^u}{1-\zeta}\in\Z[\zeta+\zeta^{-1}]^*$ and so $E_i\in \Z[\zeta+\zeta^{-1}]^*$ is a unit, and then we can chose the determination $\omega_i=E_i^{1/p}\in\R$.
\item
Let $r_i$ be  the index of {\it i}rregularity,
\begin{displaymath}
r_i=Card(\{B_{2j}\quad |\quad B_{2 j}\equiv 0\modu p,\quad j=1,2,\dots,\frac{p-3}{2}\},
\end{displaymath}
where $B_{2j}$ are the even Bernoulli numbers.
\item
Let $i_j, \quad j=1,\dots, r_i,\quad 1\leq i_j\leq \frac{p-3}{2}$, be the irregular indices  of $\Q(\zeta)$.
\item
The $p$-elementary Hilbert class field can be explicitly described, see Washington \cite{was}, exercise 8.9 p 165:
\begin{equation}\label{e109072}
H=\Q(\zeta,E_{i_1}^{1/p},\dots,E_{i_{r_p^-}}^{1/p}).
\end{equation}
and $r_i=r_p^-$.
\end{itemize}
\end{itemize}
%
\subsubsection{Some results on structure of  $p$-Hilbert class field of $\Q(\zeta)$}\label{s112221}
In this subsection, we give several general results  on structure of $p$-elementary Hilbert class field.
%
\begin{thm}{ *** }\label{t109071}
Let $\mathbf b_i,\quad i=1,\dots,r_p$, be the ideals defined in relation
(\ref{e112213}) p. \pageref{e112213}.
Let $H$ be the $p$-elementary Hilbert class field of $\Q(\zeta)$.
Then $H$ verifies the structure's properties:
\begin{equation}\label{e109156}
\begin{split}
& H=\Q(\zeta,\omega_1,\dots,\omega_{r_p})/\Q(\zeta),\\
& \omega_i^p\in\Z[\zeta],\quad \omega_i\not\in \Z[\zeta],\quad i=1,\dots,r_p,\\
&\omega_i^p\equiv c_i^p \modu \pi^{p+1},\quad c_i\in\Z,
\quad c_i\not\equiv 0\modu p,
\quad i=1,\dots,r_p,\\
& \omega_i^p\Z[\zeta]=\mathbf b_i^{\delta_i p},\quad \delta_i\in\{0,1\},
\quad i=1,\dots,r_p.
\end{split}
\end{equation}
\begin{proof}
Apply
results obtained in definition subsection \ref{s112222} p. \pageref{s112222}.
\end{proof}
\end{thm}
%
{\bf Some definitions:}
\begin{itemize}
\item
Recall that $r_p^-$ is the $p$-rank of the relative $p$-class group $C_p^-$.
Recall that we have $\sigma(\mathbf b_i)\simeq \mathbf b_i^{\mu_i},\quad
\mu_i=u_{2m_i+1}$.
\item
Let us define the set $E=\{B_1,\dots,B_{r_p^-}\}$. The set $E$ is partitioned in the  two sets $E_1=\{B_i\in E\ |\ B_i\mbox{\ singular primary} \}$ and
$E_2=\{B_i\in E\ |\ B_i\mbox{\ not singular primary}\}$
\item
Recall that, when $B_i\in E_1$ then $B_i\equiv 1\modu \pi^p$, so the extension
$\Q(\zeta, B_i^{1/p})/\Q(\zeta)$ is unramified and we can take for $B_i$ defined in relations
(\ref{e201291}) p.\pageref{e201291} and (\ref{e207291}) p.\pageref{e207291} the value $F_i=B_i$.
\end{itemize}
%
\begin{thm} { *** }\label{t201041}
Let us consider the $B_i,\quad i=1,\dots,r_p^-$, defined in relations (\ref{e201291}) p.\pageref{e201291} and (\ref{e207291}) p.\pageref{e207291}.
There are $r_p^+$ singular primary numbers $B_i,\quad i=1,\dots,r_p^+$, with
\begin{equation}\label{e207292}
\begin{split}
& B_i\in\Z[\zeta+\zeta^{-1}]-\Z[\zeta+\zeta^{-1}]^*,\quad i=1,\dots,r_p^+,\\
& B_i\Z[\zeta]= \mathbf b_i^{p},\quad  B_i\overline{B}_i=\beta_i^p,\quad \beta_i\in\Z[\zeta+\zeta^{-1}],\\
& Cl(\mathbf b_i)\in C_p^-\quad i=1,\dots,r_p^+. \\
\end{split}
\end{equation}
So $Card(E_1)=r_p^+$, or $n^-=r_p^+$ in relation (\ref{e201291}).
\begin{proof}
From Furtwangler, see Ribenboim \cite{rib} (6C) p 182
and from Hecke, see Ribenboim \cite{rib} (6D) p 182,
there are exactly $r_p^+$ unramified {\bf independant} extensions verifying
\begin{displaymath}
\begin{split}
& L_i=\Q(\zeta,\varpi_i)/\Q(\zeta),\quad i=1,\dots,r_p^+,\\
& \varpi_i^p\not\in\Z[\zeta]^*,\quad i=1,\dots,r_p^+,\\
& \varpi_i^p\overline{\varpi}_i^p=\alpha_i^p,
\quad \alpha_i\in\Q(\zeta+\zeta^{-1}),\\
& \varpi_i^p=G_i\equiv c_i^p\modu \pi^{p+1},
\quad c_i\in\Z,\quad c_i\not\equiv 0\modu p,\\
& G_i\Z[\zeta]=\mathbf c_i^p,\quad Cl(\mathbf c_i)\in C_p^-.
\end{split}
\end{displaymath}
Let us consider the $B_i,\quad i=1,\dots,r_p^-$. Let $n^-$ be the number of
the $B_i$ singular primary among them.
We have seen that $F_i= B_i\not\in\Z[\zeta+\zeta^{-1}]^*$ for $i=1,\dots,n^-$ and that
$F_i\in\Z[\zeta+\zeta^{-1}]^*$ for $i=n^-+1,\dots,r_p$.
From $\mathbf c_i\in C_p^-$ in relation (\ref{e207292}), we can write
\begin{displaymath}
G_i=B_1^{\nu_{1,i}}\times\dots\times B_{n^-}^{\nu_{n^{-},i}}\times
 F_{n^{-}+1}^{\nu_{n^{-}+1,i}}\times\dots\times F_{r_p}^{\nu_{r_p,i}}.
\end{displaymath}
Observing that $G_i\overline{G}_i=\alpha_i^p$ and that $B_i\overline{B}_i=\beta_i^p$, we get
\begin{displaymath}
\alpha_i^p= \beta_1^p\times\dots\times\beta_{n^-}^p\times
F_{n^-+1}^{2\nu_{n^-+1,i}}\times\dots\times F_{r_p}^{2\nu_{r_p,i}},
\quad i=1,\dots,r_p.
\end{displaymath}
The unramified extensions $\Q(\zeta,F_i^{1/p})/\Q(\zeta),\quad i=n^-+1,\dots,r_p$, are independant, therefore $\nu_{i,n^-+1}=\dots =\nu_{i,r_p^-}\equiv 0\modu p$, and so
\begin{displaymath}
G_i=B_1^{\nu_{1,i}}\times\dots\times B_{n^-}^{\nu_{n^{-},i}},\quad i=1,\dots,r_p^+.
\end{displaymath}
The $r_p^+$ extensions $\Q(\zeta,G_i^{1/p})/\Q(\zeta),\quad i=1,\dots,r_p^+$, are independant, and from Hecke theorem there is no $r_p+1$ such independant extensions.
Therefore $n^-=r_p^+$, which achieves the proof.
\end{proof}
\end{thm}
%
%
%
%
In the next theorem, we refine the Hilbert class field theorem \ref{e109156} p. \pageref{e109156}.
\begin{thm} { *** }\label{t201013}
With the same global meaning of index $i$ in all the theorem,
the Hilbert $p$-class field verifies the  structure:
\begin{equation}\label{e201014}
\begin{split}
& H=\Q(\zeta,\omega_1,\dots,\omega_{r_p}),\\
& \omega_i^p= B_i,\quad B_i\in\Z[\zeta],\quad B_i\Z[\zeta]
=\mathbf b_i^p,\\
&  B_i\overline{B}_i=\beta_i^p,\quad \beta_i\in\Z[\zeta+\zeta^{-1}]-\Z[\zeta+\zeta^{-1}]^*,
\quad i=1,\dots,r_p^+,\\
& Cl(\mathbf b_i)\in C_p^-, \quad i=1,\dots,r_p^+,\\
& \omega_i^p \in \Z[\zeta+\zeta^{-1}]^*,\quad i=r_p^++1,\dots,r_p^-,\\
& \omega_i^p \Z[\zeta+\zeta^{-1}]= \mathbf b_i^{\delta_i},\quad
Cl(\mathbf b_i)\in C_p^+,\quad \delta_i\in\{0,1\},
\quad i=r_p^-+1,\dots,r_p,
\end{split}
\end{equation}
where $B_i,\quad i=1,\dots,r_p^+$, are singular primary.
\begin{proof}
Apply theorem \ref{t109071} p. \pageref{t109071}, theorem \ref{t201041} p. \pageref{t201041}.
\end{proof}
\end{thm}
%
\begin{thm}\label{t201061}
Let us consider the set of   singular primary elements
$B_i,\quad i\in I_1=\{1,\dots,r_p^+\}\cup\{r_p^-+1,\dots,r_p^-+n^+\}$ .
Let us consider the set of  not singular primary elements
$B_i,\quad i\in I_2 =\{r_p^++1,\dots, r_p^-\}\cup\{r_p^-+n^++1,\dots, r_p\}$.
Let $\Gamma\in\Q(\zeta)$
\begin{displaymath}
\Gamma=c^p\times \prod_{i\in I_1} B_i^{\nu_i}\times\prod_{i\in I_2} B_i^{\nu_i},
\quad c\in \Q(\zeta).
\end{displaymath}
If $\ \Gamma$ is singular primary then $\ \nu_{i}\equiv 0\modu p$ for all $i\in I_2$.
\begin{proof}
Apply  theorem \ref{t201041} and relations (\ref{e201291}) p.\pageref{e201291}
and (\ref{e207291}) p. \pageref{e207291}.
\end{proof}
\end{thm}
%
%
%
%
\subsubsection{ On Galois group of extension $H/\Q$}\label{s20627}.
\begin{itemize}
\item
Let $H=\Q(\zeta,F_1^{1/p},\dots,F_{r_p}^{1/p})/\Q(\zeta)$ be the $p$-elementary Hilbert's extension.
\item
Let $G_H$ be the Galois group of $H/\Q(\zeta)$.
Let us consider the $\Q$-isomorphisms $\theta$ of $H/\Q(\zeta)$ defined by:
\begin{displaymath}
\begin{split}
&\theta_i^{j_i}: \omega_i\rightarrow \omega_i\times\zeta^{j_i},
\quad i=1,\dots,r_p,\quad j_i=0,\dots,p-1,\\
&\theta_{j_1,\dots,j_{r_p}}=\theta_1^{j_1}\circ\dots\circ
\theta_{r_p}^{j_{r_p}}.
\end{split}
\end{displaymath}
\item
Let $B$ be the ring of integers of $H$ and $B_H^*$ the group of units of $B_H$.
Let $\theta$ be the  $\Q(\zeta)$-isomorphisms of $H$ previously defined.
\item
Let $\sigma :\zeta \rightarrow \zeta^u$ be the $\Q$-isomorphism of $\Q(\zeta)$.
Observe that $\sigma$ can be extended to $H$ unambiguously because
$F_j,\quad j=1,\dots,r_p$, singular primary implies that
$\sigma(F_j)$ is also singular primary,   therefore
\begin{displaymath}
\sigma(\omega_{j})=_{def}(\sigma(F_j))^{1/p}\in H,
\end{displaymath}
if not the unramified extension $\Q(\zeta,\sigma(\omega_{j}))\not\subset H$;
Then, we get
\begin{displaymath}
\sigma(\omega_j)\in H,\quad j=1,\dots,r_p.
\end{displaymath}
Therefore
\begin{equation}\label{e102091}
\sigma(\omega_j)=P_j(\omega_1,\dots,\omega_{r_p}),\quad j=1,\dots,r_p,
\end{equation}
where $P_j(\omega_1,\dots,\omega_{r_p})$ is a polynomial of
$\omega_1,\dots,\omega_{r_p}$, with coefficients in $\Q(\zeta)$.
We have in the same way
\begin{displaymath}
\sigma\circ\sigma(\omega_j)=
\sigma(P_j(\omega_1,\dots,\omega_{r_p}))
=P_j(\sigma(\omega_1,\dots,\sigma(\omega_{r_p}))\in H,
\end{displaymath}
and so on,
with the choice of determinations of
$\sigma^k(\omega_j)=_{def}(\sigma^k(F_j))^{1/p},\quad j=1,\dots,r_p$, such that
\begin{equation}\label{e102092}
\sigma^k(\omega_j)=
\sigma\circ\sigma^{k-1}(\omega_j)
=\sigma\circ(\sigma^{k-1}(F_j))^{1/p}=_{def}(\sigma^k(F_j))^{1/p},
\quad k=1,\dots,p-1.
\end{equation}
\end{itemize}
Therefore we have proved :
%
\begin{lem}\label{l101022}
Let the $p$-elementary Hilbert class field $H$ of $\Q(\zeta)$, defined by $H=\Q(\zeta,F_1^{1/p},\dots,F_{r_p}^{1/p})$.
Then $H/\Q$ is a Galois extension, (not abelian) defined by the $\Q$-isomorphisms
$\sigma^i\circ \theta,\quad i=0,\dots,p-2,$ where $\sigma$ is the extension of the $\Q$-isomorphisms $\sigma$ of $\Q(\zeta)$ to $H$ defined above, and where $\theta$ are the $\Q(\zeta)$-isomorphisms of $H$ and
\begin{equation} \label{e101023}
[H:\Q]=(p-1)\times p^{r_p}.
\end{equation}
\end{lem}
%
\begin{thm}{ *** }\label{l207301}
Let $C_1$ be a subgroup of order $p$ of $C_p$ such that $\sigma(C_1)=C_1$.
Let $M=\Q(\zeta,\Omega)/\Q(\zeta)$ be a cyclic unramified extension of $\Q(\zeta)$ corresponding by Artin map to $C_1$. Then $M/\Q$ is Galois.
\end{thm}
%
\subsection{ Congruences $\modu \pi$ in $p$-elementary Hilbert class field $H$}
\begin{lem}\label{l205291}
Let $H$ be the $p$-elementary Hilbert's class field of the field $\Q(\zeta)$.
Let $\pi_H$ be a prime of $H$ above $\pi=(\zeta-1)\Z[\zeta]$.
Let $\alpha,\beta\in H$

with $\alpha\not\equiv 0 \mod \pi_H$
and $\alpha\equiv\beta\modu \pi_H$.
Then $\alpha^p\equiv \beta^p  \modu \pi_H^{p+1}$.
\begin{proof}
The absolute norm of the ideal $\pi_H$ verifies $N_{H/\Q}(\pi_H)=p$.
Let $\lambda=(\zeta-1)$.
We have
$\alpha-\beta\equiv 0 \modu \pi_H
\Rightarrow
(\alpha-\zeta^k\beta)\equiv 0 \modu \pi_H$ for $k=0,1,\dots,p-1$.
Therefore, for all
$k, \quad 0\leq k \leq p-1$, there exists $a_k\in \N,\quad 0\leq a_k \leq p-1$, such that
$(\alpha-\zeta^k\beta)\equiv \lambda a_k \modu \pi_H^2$.
For another value $l,\quad 0\leq l\leq p-1$, we have, in the same way,
$(\alpha-\zeta^l\beta)\equiv \lambda a_l \modu \pi_H^2$,
hence $(\zeta^k-\zeta^l)\beta\equiv \lambda(a_k-a_l) \modu \pi_H^2$.
For $k\not=l$ we get $a_k\not= a_l$, because $\pi_H \| (\zeta^k-\zeta^l)$ and
because hypothesis $\alpha\not\equiv 0 \modu \pi_H$ implies that $\beta\not\equiv 0 \modu \pi_H$.
Therefore, there exists one and only one $k$ such that
$(\alpha-\zeta^k\beta)\equiv 0 \modu \pi_H^2$. Then, we have
$ \prod_{j=0}^{p-1}(\alpha-\zeta^j\beta)
= (\alpha^p-\beta^p)\equiv 0 \modu \pi_H^{p+1}$.
\end{proof}
\end{lem}

%
\begin{lem}\label{l012261}
Let $H/\Q(\zeta)$ be the $p$-elementary Hilbert class field of $\Q(\zeta)$. Let $B_H$ be the ring of integers of $H$.
Suppose that
\begin{displaymath}
\beta\equiv \alpha \modu \pi,\quad \alpha,\beta\in B_H
\end{displaymath}
Then
\begin{displaymath}
\beta^p-\alpha^p\equiv 0 \modu \pi^{p+1}.
\end{displaymath}
\begin{proof}
The principal ideal $\pi$ of $\Q(\zeta)$ splits totally in $p$-elementary Hilbert's class field
$H$ of $\Q(\zeta)$. For  each prime  ideal $\pi_i$ of $H$ above $\pi$, which is  of residual degree 1,   we can apply  lemma  \ref{l205291} p. \pageref{l205291}
to get $\beta^p-\alpha^p\equiv 0\modu \pi_i^{p+1}$, which  leads to the proof.
\end{proof}
\end{lem}
%
%
Let $u$ be a primitive root $\modu p$.
Let $\sigma:\zeta \rightarrow \zeta^u$ be a $\Q$-isomorphism of $\Q(\zeta)$.
Let us consider the $p$-Hilbert class field $H=\Q(\zeta,\omega_1,\dots,\omega_{r_p})/\Q(\zeta)$.
We know
that $H/\Q$ is a Galois extension and so $\sigma\in Gal(\Q(\zeta)/\Q)$ can be extended (see lemma \ref{l101022} p. \pageref{l101022}),
to $\sigma_H\in G_H$ where  $G_H=Gal(H/\Q)$, where  $\sigma_H(\omega_i)\in H,\quad i=1,\dots,r_p$,
 and $\sigma_H$ chosen with $\sigma_H^{p-1}(\omega_i)=\omega_i,\quad i=1,\dots,r_p$, from definitions of subsection \ref{s20627} p. \pageref{s20627}. Therefore, for this extension $\sigma_H$, the map
$\sigma_H^{p-1}$ is the identity map. Let $\pi_i,\quad i=1,\dots,p^{r_p}$, be the prime ideals of $H$ above the prime ideal $\pi$ of $\Q(\zeta)$. The notation $\sigma_H^{(p-1)/2}(\pi_i),\quad i=1,\dots,r_p$ is meaningful and $\sigma_H^{(p-1)/2}(\pi_i)$ is a prime ideal of $H$ above $\pi$.
The next  theorem describes the action of conjugation $\sigma_H^{(p-1)/2}$ on prime ideals of $H$ above $\pi$.
%
\begin{thm}{ *** }\label{l205311}
There exists at least one  prime ideal $\pi_H$ of $H$ above $\pi$ such that
$\sigma_H^{(p-1)/2}(\pi_H)=\pi_H$.
\begin{proof}$ $
\begin{itemize}
\item
Let us note here $\tau=\sigma_H^{(p-1)/2}$ for $\sigma_H\in G_H$.
Let $H=\Q(\zeta,\omega_1,\dots,\omega_{r_p})/\Q(\zeta)$ be the Hilbert $p$-class field of $\Q(\zeta)$.
\item
Let us consider the unramified extension $M_1=\Q(\zeta,\omega_1)/\Q(\zeta)$.
Let $\pi_1,\dots,\pi_p$ be the $p$ prime ideals of $K_1$ above $\pi$. Let $B_1$ be the ring of integers of
$M_1$. Then $\pi B_1=\pi_1\dots\pi_p$. By conjugation we get
$\pi B=\tau(\pi) B= \pi_1\dots\pi_p=\tau(\pi_1)\dots\tau(\pi_p)$, where we note that
$\tau(\pi_1),\dots,\tau(\pi_p)$ are also prime ideals of $M_1$ above $\pi$ because $M_1/\Q$ is Galois.
Therefore $\tau(\pi_1)=\pi_k$ and $\pi_1=\tau^2(\pi_1)=\tau(\pi_k)$. But there is an odd number $p$ of primes $\pi_i,\quad i=1,\dots,p$; therefore there exists at least one $l\in\N,\quad 1\leq l\leq p$ such that $\pi_l=\tau(\pi_l)$.
\item
Let us consider the unramified extension $M_2=\Q(\zeta,\omega_1,\omega_2)/\Q(\zeta)$.
Let $\pi_{1,1},\dots,\pi_{1,p}$ be the $p$ prime ideals of $M_2$ under $\pi_1$. Let $B_2$ be the ring of integers of
$M_2$. Then $\pi_1 B_2=\pi_{1,1}\dots\pi_{1,p}$. By conjugation we get
$\pi_1 B=\tau(\pi_1) B= \pi_{1,1}\dots\pi_{1,p}=\tau(\pi_{1,1})\dots\tau(\pi_{1,p})$, where we note that
$\tau(\pi_{1,1}),\dots,\tau(\pi_{1,p})$ are prime ideals of $K_2$ under $\pi_1$.
Therefore $\tau(\pi_{1,1})=\pi_{1,k}$ and $\pi_{1,1}=\tau^2(\pi_{1,1})=\tau(\pi_{1,k})$. But there is an odd number $p$ of primes $\pi_{1,i},\quad i=1,\dots,p$; therefore there exists at least one
$l\in\N,\quad 1\leq l\leq p$ such that $\pi_{1,l}=\tau(\pi_{1,l})$.
\item
And so on up to $K_{r_p}=H=\Q(\zeta,\omega_1,\dots,\omega_{r_p})$ : finally there exists at least one prime ideal  $\pi_H$ of $H$ verifying $\tau(\pi_H)=\pi_H$.
\end{itemize}
\end{proof}
\end{thm}
%
\begin{thm}\label{l205271}
Let $H/\Q(\zeta)$ be the $p$-elementary Hilbert class field of $\Q(\zeta)$.  Let $\pi=(\zeta-1)\Z[\zeta]$.
Let $\sigma :\zeta \rightarrow \zeta^u$ be a $\Q$-isomorphism of $\Q(\zeta)$ . Let $\sigma_H$ be the extension of $\sigma$ to a $\Q$-isomorphism of $H$ with $\sigma_H^{p-1}$ identity map.
Let $\tau=\sigma_H^{(p-1)/2}$.
Let $\alpha\in H$ where $\alpha$ is coprime with $\pi$.
Then there exists a prime ideal $\pi_H$ of $H$ above $\pi$ such that
\begin{equation}\label{e205271}
\frac{\alpha}{\tau(\alpha)}\equiv 1\modu \pi_H.
\end{equation}
\begin{proof}$ $
Let us consider $\alpha\in H$ with $\alpha B+\pi B= B$. Let $\pi_H$ be a prime of $H$ above $\pi$ with

$\pi_H=\tau(\pi_H)$.
$N_{H/\Q}(\pi_H)=p$ implies that there exists $a\in\N,\quad 1\leq a\leq p-1$, such that
$\alpha\equiv a\modu \pi_H$. By conjugation, we get $\tau(\alpha)\equiv a\modu \tau(\pi_H)=\pi_H$.
and so $\frac{\alpha}{\tau(\alpha)}\equiv 1\modu\pi_H$.
\end{proof}
\end{thm}
%
%
\clearpage
\section{On  Hilbert class field of  $\Q(\zeta+\zeta^{-1})$}\label{s207272}
Let $r_p$ be the rank of the $p$-class group $C_p$. Let $r_p^+$ be the rank of the $p$-class group $C_p^+$.
In theorem \ref{t201013} p.\pageref{t201013}, we have seen, from a result on Hilbert class field of Hecke and Futwangler that $r_p\geq 2 r_p^+$. Starting of this result, in this subsection, we investigate some relations between the Hilbert class field of $\Q(\zeta)$ and of $\Q(\zeta+\zeta^{-1})$.
%
\subsection{Some definitions}
\label{s112311}
\begin{itemize}
\item
Let $p>3$ be a prime.
\item
Let $\zeta\in \C$ be  a root of the equation $X^{p-1}+X^{p-2}+\dots+X+1=0$.
\item
Let $K=\Q(\zeta)$ be the $p$-cyclotomic field.
\item
Let $u$ be a primitive root $\modu p$. Let $\sigma: \zeta\rightarrow\zeta^u$ be a $\Q$-isomorphism generating $Gal(K/\Q)$.
\item
Let $h$ be the class number of $K$.
\item
Let $H$ be the Hilbert class field of $K$.
\item
Let $K^+=\Q(\zeta+\zeta^{-1})$ be the maximal totally real subfield of $K$,  verifying $[K:K^+]=2$.
\item
Let $h^+$ be the class number of the field $K^+$.
\item
Let $H^+$ be the Hilbert class field of $K^+$.
\item
Let $C_p$ be the $p$-class group of $K$.
\item
Let $C_p^+$ be the $p$-class group of $K^+$, so with $C_p^+\subset C_p$.
\item
Let $C_p^-=C_p/C_p^+$ be the relative $p$-class group.
\item
Let the ideals $\mathbf b_i,\quad i=1,\dots,r_p$, defined in relation \ref{e112213} p. \pageref{e112213}, where
$Cl(\mathbf b_i)\in C_p^-,\quad i=1,\dots, r_p^-$ and $Cl(\mathbf b_i)\in C_p^+,\quad i=r_p^-+1,\dots,r_p$,
where $r_p=r_p^-+r_p^+$.
Suppose that $h^+\equiv 0\modu p$, we have seen that:
\begin{itemize}
\item
there exists  at least one cyclic group $C_1^+=<Cl(\mathbf b_+)>\ \subset C_p^+$   of order $p$, where $\mathbf b_+$ is an ideal not principal of $K$ with $\mathbf b_+^p\simeq \Z[\zeta]$ and $\sigma(\mathbf b_+)\simeq \mathbf b_+^{\mu_+},\quad \mu_+\in {\bf F}_p^*,\quad \mu_+^{(p-1)/2}=1$.
\item
From theorem \ref{t201013} p. \pageref{t201013} derived of Hecke and Furtwangler, there exists  a corresponding cyclic group $C_1^-=<Cl(\mathbf b_-)>\ \subset C_p^-$   of order $p$, where $\mathbf b_-$ is an ideal not principal of $\Q(\zeta)$ with $\mathbf b_-^p\simeq \Z[\zeta]$ and $\sigma(\mathbf b_-)\simeq \mathbf b_-^{\mu^-},\quad \mu^-\in {\bf F}_p^*,\quad \mu_-^{(p-1)/2}=-1$.
\end{itemize}
\item
From Hilbert class field theorem,
$h^+\equiv 0\modu p$  implies that there exists at least one field,  cyclic unramified extension $L/K^+$, with
$K^+\subset L\subset H^+,\quad [L:K^+]=p$.
\end{itemize}
%
\subsection{Relations between Hilbert class fields $H$ and $H^+$ }
\begin{lem}\label{l207181}
Suppose that $h^+\equiv 0\modu p$ and let $L$ be the field defined immediatly above.
The extension $L/K^+\subset H^+/K^+$ is of form
$L=K^+(\Omega_1+\frac{b}{\Omega_1})$ where $\Omega_1^p=F\in K-K^+$ where $b\in K^+$ and $F\overline{F}=b^p$.
\begin{proof}$ $
\begin{itemize}
\item
 The field $L$ can be written in form
\begin{equation}\label{e207051}
\begin{split}
& L=K^+(\Omega),\\
&P(\Omega)=\Omega^p+a_{p-1}\Omega^{p-1}+\dots+a_1\Omega+a_0=0,\\
&a_i\in K^+,\quad i=0,\dots,p-1,\quad a_0\not=0,\\
\end{split}
\end{equation}
where $P(\Omega)$ is the minimal polynomial of $\Omega\in\C$. From $\zeta\in\C$ and $a_i\in K^+$, we derive that $a_i\in \R$. The coefficients of $P(\Omega)$ are in $\R$ and the degree $p$ of $P(\Omega)$ is odd,  therefore there exists at least one root $\Omega\in\R$. In the sequel, we choose this root. Therefore, if we  note $\overline{\Omega}$ the complex conjugate of $\Omega$, then $\Omega=\overline{\Omega}$.
\item
Let us denote $M=K(\Omega)$.  Then $[M:K]\times[K:K^+]=[M:L]\times[L:K^+]$ and so
$[M:K]\times 2=[M:L]\times p$, which implies that $[M:K]=p$ and [M:L]=2.
\item
Show that $M/K$ is a cyclic extension: $L/K^+$ is a cyclic extension. Let $\theta_L$ be a $K^+$-isomorphism generating the Galois group $Gal(L/K^+)$. We define a $K$-isomorphim $\theta_M$ generating the Galois group
$Gal(M/K)$ by $\theta_M(a)=\theta_L(a)$ for $a\in L\subset M$ and $\theta_M(\zeta)=\zeta$ for $\zeta\in M,\quad \zeta\not\in L$.
\item
From $\zeta\in K$ and from $[M:K]=p$ and $M/K$ cyclic, we derive that $M/K$ is a Kummer cyclic extension. Therefore there exists $\Omega_1\in M,\quad \Omega_1\not\in K$ such that
$M=K(\Omega_1)$ where $\Omega_1\in \C$ is one of the $p$ roots of the equation $\Omega_1^p= F$ with $F\in K$.  A $K$-isomorphim
$\theta_M$ generating $Gal(M/K)$ is defined by
$\theta_M : \Omega_1 \rightarrow \Omega_1\times\zeta,\quad \zeta \rightarrow \zeta$.
\item
Show that the extension $M/K$ is  unramified:
\begin{itemize}
\item
Let $\pi_{K^+}$ be the prime ideal of $K^+$ above $p$.
Let $\pi_{K}=(1-\zeta)\Z[\zeta]$ be the prime ideal of $K$ above $\pi_{K^+}$.
Let $\pi_M$ be a prime of $M$ above $\pi_K$.
Let $\pi_L$ be a prime of $L$ above $\pi_{K^+}$.
The extension $L/K^+$ is unramified.
The  multiplicativity of ramification indexes in the extensions $(M/K, K/K^+)$ and $(M/L, L/K^+)$ implies that $\pi_K$ does not ramify in $M/K$:
\begin{displaymath}
\begin{split}
& e(\pi_M|\pi_{K^+})=e(\pi_M|\pi_K)\times e(\pi_{K}|\pi_{K^+})=e(\pi_M|\pi_K)\times 2\geq 2,\\
& e(\pi_M|\pi_{K^+})=e(\pi_M|\pi_L)\times e(\pi_{L}|\pi_{K^+})=e(\pi_M|\pi_L)\leq 2=[M:L],\\
\end{split}
\end{displaymath}
which implies that $e(\pi_M|\pi_K)=1$.
\item
Let $q\not=p$ be a prime.
Let $\mathbf q_{K^+}$ be a prime ideal of $K^+$ above $q$.
Let $\mathbf q_{K}$ be a prime ideal of $K$ above $\mathbf q_{K^+}$.
Let $\mathbf q_M$ be a prime of $M$ above $\mathbf q_K$.
Let $\mathbf q_L$ be a prime of $L$ above $\mathbf q_{K^+}$.
The extension $L/K^+$ is unramified.
The  multiplicativity of ramification indexes in the extensions $(M/K, K/K^+)$ and $(M/L, L/K^+)$ implies that $\mathbf q_K$ does not ramify in the extension $M/K$:
\begin{displaymath}
\begin{split}
& e(\mathbf q_M|\mathbf q_{K^+})=e(\mathbf q_M|\mathbf q_K)\times e(\mathbf q_{K}|\mathbf q_{K^+})
=  e(\mathbf q_M|\mathbf q_K),\\
& e(\mathbf q_M|\mathbf q_{K^+})=e(\mathbf q_M|\mathbf q_L)\times e(\mathbf q_{L}|\mathbf q_{K^+})
=e(\mathbf q_M|\mathbf q_L).\\
\end{split}
\end{displaymath}
But $M=L(\zeta),\quad [M:L]=2$ and $q\not = p$ implies that $e(\mathbf q_M|\mathbf q_L)=1$,
which leads to  $e(\mathbf q_M|\mathbf q_K)=1$.
\item
Therefore, $M/K$ is unramified and cyclic, therefore $M\subset H$.
\end{itemize}
\item
Let us denote $\overline{\Omega}_1$ the complex conjugate of $\Omega_1$.
Let $\Omega_2=\Omega_1\overline{\Omega}_1$. From the definitions  $\Omega_2\in\R$. We get $\Omega_2^p= F\overline{F}$, where $F\overline{F}\in  K^+\subset \R$.
\item
In the sequel, $(F\overline{F})^{1/p}$ is to be understood as the determination corresponding to the real  $p^{th}$ root of $F\overline{F}\in K^+\subset\R$.
\item
We  show that  the case $(F\overline{F})^{1/p}\not\in K^+$  lead to a contradiction.
From $M=K(\Omega)=K(\Omega_1)$,  we derive that $\Omega_1=P_1(\Omega)$ where
$P_1(\Omega)=\sum_{i=0}^{p-1} b_i\Omega^i,\quad b_i\in K$. By complex conjugation, observing
that $\Omega=\overline{\Omega}$, we get $\overline{\Omega}_1=\sum_{i=0}^{p-1} \overline{b}_i \Omega^i$.
Thus
\begin{displaymath}
\Omega_2=\Omega_1\overline{\Omega}_1=
(\sum_{i=0}^{p-1} b_i\Omega^i)(\sum_{i=0}^{p-1}\overline{b}_i\Omega^i)
\in K^+(\Omega)=L,
\end{displaymath}
and so $K^+(\Omega_2)\subset L$.
From $\Omega_2^p= F\overline{F}\in K^+$ and $(F\overline{F})^{1/p}\not\in K^+$, we derive that
$[K^+(\Omega_2):K^+]= p$ and so $L=K^+(\Omega_2)$. The extension $L/K^+$ is cyclic and so the
roots $\Omega_2\zeta^i,\quad i=0,\dots, p-1$, of the equation $X^p-F\overline{F}=0$ are in $L$ and so
$\zeta\in L$, contradiction because $M\not=L$.
\item
Suppose now that $b=(F\overline{F})^{1/p}\in \Z[\zeta+\zeta^{-1}]$:
\begin{itemize}
\item
From $\Omega_1^p\overline{\Omega}_1^p=F\overline{F}=b^p$, we derive that, without ambiguity in the determinations,  $\Omega_1\overline{\Omega}_1=b$ where $\Omega_1\overline{\Omega}_1\in \R,\quad b\in\R$.
Then let $\Omega_3=\Omega_1+\overline{\Omega}_1=\Omega_1+\frac{b}{\Omega_1}$.
At first note that $\Omega_3\not \in K^+$, if not we should have simultaneously
$\Omega_1\overline{\Omega}_1\in K^+$ and $\Omega_1+\overline{\Omega}_1\in K^+$ and $\Omega_1$ and
$\overline{\Omega}_1$ should be the two roots of an equation $X^2+a_1 X+a_2=0,\quad a_1,a_2\in K^+$, which contradicts $\Omega_1^p=F,\quad F\in K,\quad F^{1/p}\not\in K$.
This implies that $\Omega_3=\sum_{i=0}^{p-1} (b_i+\overline{b}_i)\Omega^i\in K^+(\Omega)$.
Therefore from $K^+(\Omega_3)\subset K^+(\Omega)$ and from $\Omega_3\not\in K^+$, we derive that
$L=K^+(\Omega_3)=K^+(\Omega)$.
\item
The $K$-isomorphims of $M/K$ are defined by:
\begin{displaymath}
\theta_M^i : \Omega_1\rightarrow \Omega_1\times\zeta^i,\quad \zeta\rightarrow \zeta.
\end{displaymath}
\item
Observe that $\Omega_1\not=\overline{\Omega}_1$, if not we should have $\Omega_1=\frac{b}{\Omega}_1$ and so
$\Omega_1^2=b,\quad b\in K^+$, which clearly contradicts $\Omega_1^p=F,\quad F\in K,\quad F^{1/p}\not\in K$.
The $K^+$-isomorphisms $\theta_L^i$ of $L/K^+$, restrictions of $K$-isomorphisms of $M/K$ are therefore  defined by:
\begin{displaymath}
\theta_L^i : \Omega_1+\frac{b}{\Omega_1}\rightarrow
\Omega_1 \zeta^i+\frac{b}{\Omega_1\zeta^i},\quad \zeta+\zeta^{-1}\rightarrow \zeta+\zeta^{-1},
\end{displaymath}
because
\begin{displaymath}
\begin{split}
& \Omega_1\zeta^i+\frac{b}{\Omega_1\zeta^i}=\Omega_1\zeta^i+\overline{\Omega}_1\zeta^{-i}
=\sum_{j=0}^{p-1} b_j\zeta^i\Omega^j+\sum_{j=0}^{p-1}\overline{b}_j\zeta^{-i}\Omega^j\\
&=\sum_{j=0}^{p-1} (b_j\zeta^i+\overline{b_j\zeta^i})\Omega^j\in K^+(\Omega)=L.\\
\end{split}
\end{displaymath}
\end{itemize}
\end{itemize}
\end{proof}
\end{lem}
%
\begin{lem}\label{l207182}
The field $M=K(\Omega_1)$ is a CM-field.
\begin{proof}$ $
\begin{itemize}
\item
$[M:L]=2$ and $M=L(\zeta)$. It is sufficient to prove that $L=K^+(\Omega_1+\overline{\Omega}_1)$ is totally real.
\item
From $\Omega_1^p=F$ then $\Omega_1=F^{1/p}$ for one choice of the root $F^{1/p}$.
Recall that $\overline{\Omega}_1$ is the complex conjugate of $\Omega_1$. Then $\Omega_1^p\overline{\Omega}_1^p=F\overline{F}$ and $\Omega_1\overline{\Omega}_1=(F\overline{F})^{1/p}=$ where $(F\overline{F})^{1:p}=b\in K^+$ is the real root defined unambiguously.
\item
Let $\sigma\in Gal(K/\Q)$. Then $\sigma(\Omega_1)\times\sigma(\overline{\Omega}_1)=\sigma(b)$.
In an other part $\sigma(\Omega_1)^p\overline{\sigma(\Omega_1)}^p=\sigma(F\overline{F})$, which gives the unambiguous definition
$\sigma(\Omega_1)\overline{\sigma(\Omega_1)}=(\sigma(F\overline{F}))^{1/p}=\sigma(b)$ for $\sigma(b)\in\R$ defined unambiguously  and so
$\sigma(\overline{\Omega}_1)=\overline{\sigma(\Omega_1)}$.
\item
Then $\sigma(\Omega_1+\overline{\Omega}_1)=\sigma(\Omega_1)+\sigma(\overline{\Omega}_1)=
\sigma(\Omega_1)+\overline{\sigma(\Omega_1)}\in \R$.
In the same way
$\sigma(\Omega_1\zeta+\overline{\Omega}_1\zeta^{-1})=\sigma(\Omega_1)\sigma(\zeta)
+\sigma(\overline{\Omega}_1)\sigma(\zeta^{-1})=
\sigma(\Omega_1)\zeta+\overline{\sigma(\Omega_1)\zeta}\in \R$.
\item
Therefore, all the conjugates $\sigma^i\circ\theta^j(\Omega_1+\overline{\Omega}_1),\quad i=0,\dots,p-2,
\quad j=0,\dots,p-1$ are real and $L/\Q$ is totally real, which achieves the proof.
\end{itemize}
\end{proof}
\end{lem}
%
\subsection{Artin map for Hilbert class fields $H$ and $H^+$}
The next results aims to explicit Artin map of Hilbert class field of $\Q(\zeta)$ and
of $\Q(\zeta+\zeta^{-1})$.
\begin{lem} \label{l207051}$ $
\begin{itemize}
\item
In Artin map $C_p\longleftrightarrow Gal(H/K)$:
\begin{itemize}
\item
$<Cl(\mathbf b_+ )>\ \subset C_p^+\subset C_p$ corresponds to unramified cyclic $M/K=K(\Omega_1)/K$ with:
\begin{itemize}
\item
$\Omega_1^p=F\in K-K^+,\quad F\overline{F}=b^p,\quad b\in K$.
\item
$F\Z[\zeta]=\mathbf b_-^p,\quad Cl(\mathbf b_-)\in C_p^-$.
\end{itemize}
\item
$<Cl(\mathbf b_-)>\ \subset C_p^-$ corresponds to unramified cyclic extension $N=K(\omega)/K$ with:
\begin{itemize}
\item
$\omega^p=F^\prime\in \Z[\zeta+\zeta^{-1}]^*$.
\item
$F\Z[\zeta]=\mathbf b_-^p,\quad Cl(\mathbf b_-)\in C_p^-$.
\end{itemize}
\end{itemize}
\item
In Artin map $C_p^+\longleftrightarrow Gal(H^+/K^+)$:
$<Cl(\mathbf b_+)>$ corresponds to unramified cyclic extension $L/K^+=K^+(\Omega_1+\overline{\Omega_1})/K^+=K^+(\Omega_1+\frac{b}{\Omega_1})/K^+$  with:
\begin{itemize}
\item
$\Omega_1^p=F\in K-K^+,\quad F\overline{F}=b^p,\quad b\in K$.
\item
$F\Z[\zeta]=\mathbf b_-^p,\quad Cl(\mathbf b_-)\in C_p^-$.
\end{itemize}
\end{itemize}
\begin{proof}$ $
\begin{itemize}
\item
Let $C_1^+= <Cl(\mathbf b_+)>\ \subset C_p^+\subset C_p$ be the cyclic group of order $p$ corresponding by Artin map to $Gal(L/K^+)$.
From
\begin{displaymath}
\begin{split}
& K^+\subset K,\quad L\subset M,\quad L\subset H^+,\quad M\subset H,\\
& [L:K^+]=[M:K]=p,\quad [M:L]=[K:K^+]=2,\\
\end{split}
\end{displaymath}
then by Artin map, we get (see for instance Washington \cite{was} p. 400):
\begin{itemize}
\item
$C_1^+= <Cl(\mathbf b_+\Z[\zeta])> \cong Gal(M/K)$ by Artin map,  where $\mathbf b_+$ can be chosen ideal of
$K^+$ with $\mathbf b_+^p$ principal.
\item
$ <Cl(N_{K/K^+}(\mathbf b_+))>=<Cl(\mathbf b_+^2)>\cong Gal(L/K^+)$ by Artin map.
\item
Let $\theta_L$ generating $Gal(L/K^+)$ and  $\theta_M$ generating $Gal(M/K)$. Then $\theta_L$ is the restriction of $\theta_M$ to $L\subset M$.
\end{itemize}
\item
With another point of view for the same result:
Let $i \in \N,\quad 0\leq i\leq p-1$. Then by Artin map
\begin{itemize}
\item
$Cl(\mathbf b_+^i) \leftrightarrow
\theta_L^i: \Omega_1+\overline{\Omega}_1\rightarrow \Omega_1\zeta^i+\overline{\Omega}_1\zeta^{-i}, \quad \theta_L^i\in Gal(L/K^+),
\quad L/K^+\subset H^+/K^+$, by Artin map in $H^+$.
\item
$\theta_L^i\in Gal(L/K^+)\leftrightarrow \theta_M^i :\Omega_1\rightarrow \Omega_1\zeta^i,\quad \theta_M^i\in Gal(M/K)$,
\item
$\theta_M^i\in Gal(M/K)\leftrightarrow Cl( \mathbf b_+^i)$ by Artin map in $H$.
\end{itemize}
\end{itemize}
\end{proof}
\end{lem}
%
The next theorem summarizes Artin map in global context of $p$-elementary Hilbert class fields of $\Q(\zeta)$ and of $\Q(\zeta+\zeta^{-1})$.
\begin{itemize}
\item
In relations (\ref{e207301}) p. \pageref{e207301} and (\ref{e207302}) p. \pageref{e207302}, we have defined the set $\mathbf M$ of $r_p$ independant unramified extensions
$M_i=\Q(\zeta,\omega_i)/\Q(\zeta),\quad i=1,\dots,r_p$, with
\begin{displaymath}
\begin{split}
& \omega_i^p =B_i,\quad i=1,\dots, r_p^+,\\
& \omega_i^p \in \Z[\zeta+\zeta^{-1}]^*,\quad i=r_p^++1,\dots,r_p^-,\\
& \omega_i^p= B_i,\quad i=r_p^-+1,\dots, r_p^-+n^+,\\
& \omega_i^p \in\Z[\zeta+\zeta^{-1}]^*,\quad i=r_p^-+n^++1,\dots,r_p.\\
\end{split}
\end{displaymath}
\item
For $i=r_p^++1,\dots,r_p^-$ and for $i=r_p^-+n^++1,\dots,r_p$, we choose $\omega_i$ such that
the extension $\Q(\zeta,\omega_i)/\Q(\zeta)$  corresponds by Artin map to the group
$<Cl(\mathbf b_i)>$ of order $p$ for the same index $i$. Thus the set $\mathbf M$ of $r_p$ independant unramified extensions is completely defined. This definition of set $\mathbf M$  is used in the sequel.
\end{itemize}
%
\begin{thm} {***}\label{t207261}
Let $\mathbf b_i,\quad i=1,\dots,r_p$, be the ideals defined in relation (\ref{e112213}) p. \pageref{e112213}. With the meaning of index $i$ of theorem
\ref{t201013} p. \pageref{t201013},
the Artin map  $C_p\leftrightarrow H$ and  $C_p^+\leftrightarrow H^+$
verify:
\begin{itemize}
\item
In Artin map $C_p\leftrightarrow Gal(H/K)$ for $i=1,\dots,r_p^+$ then
\begin{itemize}
\item
$<Cl(\mathbf b_{r_p^-+i} )>\ \subset C_p^+$ corresponds to unramified cyclic $M_i/K=K(\omega_i)/K$  with:
\begin{itemize}
\item
$\omega_i^p=F_i\in K-K^+,\quad F_i\overline{F_i}=b_i^p,\quad b_i\in K$.
\item
$F_i\Z[\zeta]=\mathbf b_i^p,\quad Cl(\mathbf b_i)\in C_p^-$.
\end{itemize}
\item
$<Cl(\mathbf b_i)>\ \subset C_p^-$ corresponds to unramified cyclic extension $M_{r_p^-+i}/K
=K(\omega_{r_p^-+i})/K$  with:
\begin{itemize}
\item
 $\omega_{r_p^-+i}^p=F_i^\prime\in \Z[\zeta+\zeta^{-1}]^*$.
\item
$F_i\Z[\zeta]=\mathbf b_i^p,\quad Cl(\mathbf b_i)\in C_p^-$.
\end{itemize}
\end{itemize}
\item
In Artin map $C_p\rightarrow Gal(H/K)$ for $i=r_p^++1,\dots,r_p^-$ then
$<Cl(\mathbf b_i)>\ \subset C_p^-$ corresponds to unramified cyclic extension $M_i/K
=K(\omega_i)/K$  with:
\begin{itemize}
\item
 $\omega_{i}^p=F_i^\prime\in \Z[\zeta+\zeta^{-1}]^*$.
\item
$F_i\Z[\zeta]=\mathbf b_i^p,\quad Cl(\mathbf b_i)\in C_p^-$.
\end{itemize}
\item
In Artin map $C_p^+\rightarrow Gal(H^+/K^+)$ for $i=1,\dots,r_p^+$ then
$<Cl(\mathbf b_{r_p^-+i})>\ \subset C_p^+$ corresponds to unramified cyclic extension
$L_i/K^+=K^+(\omega_i+\overline{\omega_i})/K^+=K^+(\omega_i+\frac{b_i}{\omega_i})/K^+$  with:
\begin{itemize}
\item
$\omega_i^p=F_i\in K-K^+,\quad F_i\overline{F_i}=b_i^p,\quad b_i\in K^+$.
\item
$F_i\Z[\zeta]=\mathbf b_i^p,\quad Cl(\mathbf b_i)\in C_p^-$.
\end{itemize}
\end{itemize}
\end{thm}
%
%
{\bf Remark:} from this result on Artin map, we see that $n^+=0$, which allows to simplify in a strong way formulation of the Hilbert class field theorem \ref{e109156} p. \pageref{e109156}
and \ref{t201013} p. \pageref{t201013}.
%
\begin{thm} { *** }\label{t207301}
With the same global meaning of index $i$ in all the theorem,
the Hilbert $p$-class field verifies the  structure:
\begin{equation}\label{e201014}
\begin{split}
& H=\Q(\zeta,\omega_1,\dots,\omega_{r_p}),\\
& \omega_i^p= B_i,\quad B_i\in\Z[\zeta],\quad B_i\Z[\zeta]
=\mathbf b_i^p,\quad Cl(\mathbf b_i)\in C_p^-,\quad i=1,\dots,r_p^+,\\
&  B_i\overline{B}_i=\beta_i^p,\quad \beta_i\in\Z[\zeta+\zeta^{-1}]-\Z[\zeta+\zeta^{-1}]^*,\quad i=1,\dots,r_p^+,\\
& \omega_i^p \in \Z[\zeta+\zeta^{-1}]^*,\quad i=r_p^++1,\dots,r_p,\\
\end{split}
\end{equation}
where $B_i,\quad i=1,\dots,r_p^+$, are singular primary.
\end{thm}
%
%
\begin{thm}{ *** }\label{t207311}
Let $C_1=<Cl(\mathbf b)>$ be a subgroup of order $p$ of $C_p$  with
$\sigma(\mathbf b)\simeq \mathbf b^\mu,
\quad \mu\in{\bf F}_p^*$. Let $M=\Q(\zeta,\Omega)$ be the cyclic unramified extension with
$Gal(M/\Q(\zeta))$ corresponding to $C_1$ by Artin map.
Then $\sigma(\Omega)=\Omega^{\nu}\times\gamma,\quad\nu\in{\bf F}_p^*,\quad \gamma\in\Q(\zeta)$.
Then $\nu\equiv u\times\mu^{-1}\modu p$.
\begin{proof}$ $
\begin{itemize}
\item
By Artin map, $Cl(\mathbf b)$ corresponds to the $\Q(\zeta)$-isomorphism
$\theta: \Omega\leftrightarrow\Omega\zeta$ of $M$.
\item
By Artin map $Cl(\mathbf b^\mu)$ corresponds to the $\Q(\zeta)$-isomorphism
$\theta^\mu: \Omega\leftrightarrow\Omega\zeta^\mu$ of $M$.
\item
Show that, by Artin map, $Cl(\sigma(\mathbf b))=Cl(\mathbf b^\mu)$ corresponds to the $\Q(\zeta)$-isomorphism
$\sigma(\Omega)\leftrightarrow\sigma(\Omega\zeta)$ of $M$:
\begin{itemize}
\item
Let $\mathbf q$ be a prime of $\Q(\zeta)$ with $\mathbf q\simeq\mathbf b$. Let $B_M$ be the ring of integers of $M$. Let $\mathbf q_M$ be a prime of $M$ above $\mathbf q$.
\item
There exists   a  Frobenius automorphism $\Phi_1$ generating  Galois group of $B_M/\mathbf q_M$ over $\Z[\zeta]/\mathbf q$ and there exists
a   Frobenius automorphism $\Phi_2$ generating  Galois group of
$B_M/\mathbf \sigma(\mathbf q_M)$ over $\Z[\zeta]/\mathbf \sigma(\mathbf q)$ with :
\begin{equation}\label{e207311}
\begin{split}
& \Phi_1(\Omega)=\Omega\zeta\equiv \Omega^{q}\modu \mathbf q B_M,\\
& \Phi_2(\sigma(\Omega))\equiv \sigma(\Omega)^{q}\modu \mathbf \sigma(\mathbf q) B_M,\\
& q=N_{\Q(\zeta)/\Q}(\mathbf q)=N_{\Q(\zeta)/\Q}(\mathbf \sigma(\mathbf q)),\\
\end{split}
\end{equation}
(see for instance Marcus, \cite{mar} theorem 32 p. 109).
\item
This implies that $\Phi_2(\sigma(\Omega))$, of form $\Phi_2(\sigma(\Omega))=\sigma(\Omega)\zeta^v
\quad v\in\N,\quad v\not\equiv 0\modu p$,
verifies $\Phi_2(\sigma(\Omega))=\sigma(\Omega)\zeta^v\equiv \sigma(\Omega)^q
\modu \sigma(\mathbf q)B_M$. But by conjugation, we get also
$\sigma(\Omega\zeta)=\sigma(\Omega)\zeta^u\equiv \sigma(\Omega)^q \modu\sigma(\mathbf q)B_M$.
Therefore $v=u$ and $\Phi_2(\sigma(\Omega))=\sigma(\Omega)\zeta^u=\sigma(\Omega\zeta)$.
\end{itemize}
\item
We have shown that, by Artin map, $Cl(\mathbf b^\mu)$ corresponds to the $\Q(\zeta)$-isomorphism
$\sigma(\Omega)\rightarrow \sigma(\Omega)\zeta^u$, thus $\Omega^\nu\leftrightarrow\Omega^\nu\zeta^u$ of $M$.
\item
By Artin map, $Cl(\mathbf b^\mu)$ corresponds to the $\Q(\zeta)$-isomorphism
$\Omega^\nu\leftrightarrow\Omega^\nu\zeta^{\mu\nu}$ of $M$.
Therefore $\mu\nu\equiv u\modu p$, which achieves the proof.
\end{itemize}
\end{proof}
\end{thm}
{\bf Remark:} Suppose that $h^+\equiv 0\modu p$. Let $\mathbf b_+$ be a not principal ideal
with $\mathbf b_+^p\simeq\Z[\zeta],\quad Cl(\mathbf b^+)\in C_p^+$. Let $M/K$ be the cyclic unramified
extension of $K$ with $<Cl(\mathbf b_+)>\ \leftrightarrow Gal(M/K)$ by Artin map. In that particular case
we have $\mu^{(p-1)/2}=1$ and so $\nu^{(p-1)/2}=-1$, which corresponds well to the fact that, in that case,
$\Omega\in K-K^+$.

%
%
\clearpage
\section{ Connection between roots $\modu p$ of Mirimanoff polynomials and $p$-class group of $\Q(\zeta)$}
%
Let $b\in\N,\quad 1 < b< p-1$. Let $m\in\N$. Let us note $\phi_m(b)$ the Mirimanoff polynomial value
$\phi_m(b)=\sum_{i=1}^{p-1} i^{m-1}\times  b^i$.
In this subsection we give some  arithmetic congruences $\modu p$  of Mirimanoff polynomials $\phi_m(b)$, strongly connected to class group structure of the field $\Q(\zeta)$.
%
\begin{lem}\label{l209053}
Let $b\in\N,\quad 1<b<p-1$.
Let $m\in\N,\quad 1\leq m \leq \frac{p-3}{2}$.
Let $D(b,2m+1)=\prod_{i=0}^{p-2}(\frac{b-\zeta^{u_i}}{b-\zeta^{-u_i}})^{u_{-(2m+1)i}}$.
\begin{enumerate}
\item
If, and only if,  $\phi_{2m+1}(b)\not\equiv 0\modu p$ then $D(b,2m+1)$ verifies:
\begin{displaymath}
\begin{split}
& D(b,2m+1)\not\in\Q(\zeta)^p,\\
& \sigma(D(b,2m+1))=D(b,2m+1)^{u_{2m+1}}\times\gamma^p,\quad\gamma\in\Q(\zeta),\\
& \sigma(D(b,2m+1))\equiv D(b,2m+1)^{u_{2m+1}}\modu \pi^{p+1},\\
& \pi^{2m+1}\|(D(b,2m+1)-1).\\
\end{split}
\end{displaymath}
\item
If, and only if, $\phi_{2m+1}(b)\equiv 0\modu p$ then $D(b,2m+1)\equiv 1\modu p$.
\end{enumerate}
\begin{proof}
Let us note in the proof $D$ for $D(b,2m+1)$.
Observe at first that
$D\not=1$ and that $D\equiv 1\modu \pi$.
\begin{itemize}
\item
Show at first that $D\not\in\Q(\zeta)^p$. Suppose that $D\in\Q(\zeta)^p$ and search for a contradiction:
$D\in\Q(\zeta)^p$ should imply that $\prod_{i=0}^{p-2}(\frac{b-\zeta^{u_i}}{b-\zeta^{-u_i}})^{u_{-(2m+1)i}}\equiv 1\modu p$. Similarly to foundation theorem \ref{s12113} p. \pageref{s12113} this should imply that
$\phi_{2n+1}(b)\times\sum_{i=0}^{p-2} u_{(2n+1-(2m+1))i}\equiv 0\modu p$ for $n=1,\dots,\frac{p-3}{2}$,
and thus $\phi_{2m+1}(b)\equiv 0\modu p$, contradicting hypothesis assumed.
\item
We get
\begin{displaymath}
\sigma(D)=\prod_{i=0}^{p-2}(\frac{b-\zeta^{u_{i+1}}}{b-\zeta^{-u_{i+1}}})^{u_{-(2m+1)i}}
\end{displaymath}
Let $i^\prime=i+1$. Then
\begin{displaymath}
\sigma(D)=
\prod_{i^\prime=1}^{p-1}(\frac{b-\zeta^{u_{i^\prime}}}{b-\zeta^{-u_{i^\prime}}})^{u_{-(2m+1)(i^\prime-1)}}.
\end{displaymath}
But $u_{-(2m+1)(i^\prime-1)}\equiv u_{-(2m+1)i^\prime}\times u_{2m+1}\modu p$, so
\begin{displaymath}
\begin{split}
& \sigma(D)=
\prod_{i^\prime=1}^{p-1}(\frac{b-\zeta^{u_{i^\prime}}}{b-\zeta^{-u_{i^\prime}}})
^{u_{-(2m+1)i^\prime}\times u_{2m+1}}\times \gamma^p,\quad\gamma\in\Q(\zeta),\\
& \sigma(D)\equiv
\prod_{i^\prime=1}^{p-1}(\frac{b-\zeta^{u_{i^\prime}}}{b-\zeta^{-u_{i^\prime}}})
^{u_{-(2m+1)i^\prime}\times u_{2m+1}}\modu\pi^{p+1}.\\
\end{split}
\end{displaymath}
Then
\begin{displaymath}
\begin{split}
&i^\prime=p-1\Rightarrow (\frac{b-\zeta^{u_{i^\prime}}}{b-\zeta^{-u_{i^\prime}}})
^{u_{-(2m+1)(i^\prime-1)}}
=(\frac{b-\zeta}{b-\zeta^{-1}})^{u_{2m+1}}.\\
&i=0\Rightarrow(\frac{b-\zeta^{u_{i}}}{b-\zeta^{-u_{i}}})
^{u_{-(2m+1)i}}
=\frac{b-\zeta}{b-\zeta^{-1}}.\\
\end{split}
\end{displaymath}
Therefore
\begin{displaymath}
\begin{split}
& \sigma(D)= D^{u_{2m+1}}\times\gamma^p,\\
& \sigma(D)=\prod_{i=1}^{p-2}(\frac{b-\zeta^{u_{i}}}{b-\zeta^{-u_{i}}})
^{u_{-(2m+1)i}\times u_{2m+1}}\times(\frac{b-\zeta}{b-\zeta^{-1}})^{u_{2m+1}}
\equiv D^{u_{2m+1}}\modu \pi^p.\\
\end{split}
\end{displaymath}
$D\equiv 1 \modu p$ should imply that $\prod_{i=0}^{p-2}(\frac{b-\zeta^{u_i}}{b-\zeta^{-u_i}})^{u_{-(2m+1)i}}\equiv 1\modu p$. Similarly to foundation theorem \ref{s12113} p. \pageref{s12113} this should imply that
$\phi_{2n+1}(b)\times\sum_{i=0}^{p-2} u_{(2n+1-(2m+1))i}\equiv 0\modu p$ for $n=1,\dots,\frac{p-3}{2}$,
and thus $\phi_{2m+1}(b)\equiv 0\modu p$, contradicting hypothesis assumed.
Similarly to proof of  lemma  \ref{l108171} p. \pageref{l108171} and remark following this lemma,  it   implies that $\pi^{2m+1}\| (D(b,2m+1)-1)$.
\end{itemize}
\end{proof}
\end{lem}
{\bf Remark:} This result generalizes the results obtained in the  references:
\begin{enumerate}
\item
the description of the $\pi$-adic structure of the
$p$-class group $C_p$, seen  in subsection \ref{s108311}  p. \pageref{s108311},
\item
the structure of the unit group
$\Z[\zeta+\zeta^{-1}]^*$ seen in subsection \ref{s205231} p. \pageref{s205231}, particularly theorem
\ref{t207271} p. \pageref{t207271}.
\end{enumerate}

%
Here a generalization to odd and even Mirimanoff polynomial of the previous lemma dealing with odd Mirimanoff polynomials:
\begin{lem}\label{l211021}
Let $b\in\N,\quad 1<b<p-1$.
Let $m\in\N,\quad 1< m \leq p-2$.
Let $D(b,m)=\prod_{i=0}^{p-2}(b-\zeta^{u_i})^{u_{-m i}}$.
\begin{enumerate}
\item
If, and only if, $\phi_{m}(b)\not\equiv 0\modu p$ then $D(b,m)$ verifies:
\begin{displaymath}
\begin{split}
& D(b,m)\not\in\Q(\zeta)^p,\\
& \sigma(D(b,m))=D(b,m)^{u_m}\times\gamma^p,\quad\gamma\in\Q(\zeta),\\
& \sigma(D(b,m))\equiv D(b,m)^{u_m}\modu \pi^{p+1},\\
& \pi^{m}\|(D(b,m)-d),\quad d\in\Z,\quad d\not\equiv 0\modu p.\\
\end{split}
\end{displaymath}
\item
If, and only if $\phi_{m}(b)\equiv 0\modu p$ then $D(b,m)\equiv d\modu p,\quad d\in\Z$.
\end{enumerate}
\begin{proof}
Same proof as lemma \ref{l209053} p. \pageref{l209053}.
\end{proof}
\end{lem}
%
The next theorem connects the existence of some roots $b$ of Mirimanoff polynomials  $\phi_n(b)\equiv 0\modu p$ to the structure of the $p$-class group $C_p$ of the field $\Q(\zeta)$.
\begin{thm}{ *** }\label{t211021}
Suppose that $p$ divides  the class number $h$ of $\Q(\zeta)$.
Let $r_p^-$ be the relative $p$-class group $C_p^-$.
For $i=1,\dots,r_p^-$, there exist   different natural integers $m_i\in\N,\quad 1\leq m_i\leq\frac{p-3}{2}$, and ideals $\mathbf b_i$ of $\Z[\zeta],\quad Cl(\mathbf b_i)\in C_p^-,\quad
\sigma(\mathbf b_i)\simeq\mathbf b_i^{\mu_i},\quad \mu_i=u_{2m_i+1},$
and   natural integers $b_i, \quad 1<b_i\leq p-1$,  verifying Mirimanoff polynomial congruences
\begin{equation}\label{e211021}
\phi_{2m_i+1}(b_i)\times\phi_{p-2m_i-1}(b_i)\equiv 0\modu p,\quad i=1,\dots,r_p^-.
\end{equation}
\begin{proof}$ $
\begin{itemize}
\item
The existence of the $\mathbf b_i,\quad i=1,\dots,r_p^-$, is shown in relation \ref{e203031} p. \pageref{e203031}.
From lemma \ref{l108161} p. \pageref{l108161} there exist, for $i=1,\dots,r_p^-$, singular
$C_i\in\Q(\zeta)$ with $C_i\equiv 1\modu \pi, \quad \sigma(C_i)= C_i^{\mu_i}\times\alpha_i^p,
\quad \alpha_i\in\Q(\zeta)$.
\item
If $C_i$ is singular primary, then $C_i\equiv 1\modu\pi^p$ and from lemma \ref{l209053} p. \pageref{l209053} there exists at least one $b_i\in\N,\quad 1<b_i\leq p-1$, with $\phi_{2m_i+1}(b_i)\equiv 0\modu p$:
If $\phi_{2m_i+1}(b_i)\not\equiv 0\modu p$ for all $b_i,\quad 1<b_i\leq p-1$, then it should imply that
 $\pi^{2m_i+1}\|(D(b_i,2m_i+1)-1)$ for all $b_i,\quad 1<b_i\leq p-1$,  and so,
from $C_i\equiv D(b_i,2m_i+1)\modu\pi^{p-1}$, we should obtain
$\pi^{2m_i+1}\|(C_i-1)$, contradiction.
\item
If $C_i$ is singular not primary then from theorem \ref{t207311} p. \pageref{t207311} on class field theory, there exists a primary unit $F_i=\omega_i^p\in\Z[\zeta+\zeta^{-1}]^*$ with
$\sigma(F_i)=F_i^{p-2m_i-1}\times\varepsilon_i^p,\quad\varepsilon_i\in\Z[\zeta+\zeta^{-1}]^*$.
Then $F_i\equiv c_i^p,\quad c_i\in\Z$.
If $\phi_{p-2m_i-1}(b_i)\not\equiv 0\modu p$ for all $b_i,\quad 1<b_i\leq p-1$, then, from lemma
\ref{l211021} p. \pageref{l211021},  it should imply that
 $\pi^{p-2m_i-1}\|(D(b_i,p-2m_i-1)-1)$ for all $b_i,\quad 1<b_i\leq p-1$, and so
from $F_i\equiv d_i\times D(b_i,p-2m_i-1)\modu\pi^{p-1},\quad d_i\in\Z$, that
$\pi^{p-2m_i-1}\|(F_i-1)$, contradiction.
\end{itemize}
\end{proof}
\end{thm}
%
\clearpage
\section*{PART TWO: ON FERMAT-WILES THEOREM}
This section contains a study of Fermat-Wiles equation in first and second case with an elementary approach and with use of results on cyclotomic fields obtained in the the first part of the monography. It contains:
\begin{itemize}
\item a study of FLT first case:
\begin{itemize}
\item
Several congruences $\modu p$ of Mirimanoff polynomials connected to FLT.
\item
Several congruences $\modu p$ of Mirimanoff polynomials connected to FLT in intermediate fields
$K,\quad \Q\subset K\subset \Q(\zeta)$.
\item
Several strictly elementary congruence criterions.
\item
An approach with the point of view of representations of $Gal(\Q(\zeta)/\Q)$.
\item
An approach with the point of view of Hilbert class field of $\Q(\zeta)$.
\item
A formulation of explicit criterions  with Jacobi resolvents method
\item
A comparative survey of our results   with bibliography
\end{itemize}
\item
A study on FLT second case:
\begin{itemize}
\item
 An elementary approach.
\item
A  cyclotomic field  approach.
\item
An  Hilbert class field approach.
\end{itemize}
\end{itemize}
%
\clearpage
\section{Some preliminaries on FLT}\label{s210161}
This section contains some definitions on Fermat-Wiles theorem, some general preliminary results and the computation of a rough  upper bound of the  class number $h$ of $\Q(\zeta)$ used in the sequel.
%
\subsection{Some definitions on FLT}\label{s210151}
\begin{itemize}
\item
For $m=1,\dots,\frac{p-3}{2}$, the polynomials
\begin{displaymath}
\phi_{2m+1}(T)= 1^{2m}\times T+2^{2m}\times T^2+\dots+(p-1)^{2m}\times T^{p-1}
\end{displaymath}
are the odd Mirimanoff polynomials of the indeterminate $T$.
\item
Let us denote $x^p+y^p+z^p=0$ the Fermat equation, where, as usual,
$x\in\Z-\{0\}, \quad y\in\Z-\{0\},\quad z\in\Z-\{0\}$, and $x,y, z$ are pairwise coprime. In all the sequel of this monograph , $x, y, z$ are the solutions of the Fermat equation. Let $\zeta$ be a root of the equation $\frac{X^p-1}{X-1}=0$.
Classically, we have
\begin{equation}\label{e1}
\begin{split}
&(x+\zeta^{i} y)\Z[\zeta] = \pi^\delta\s_i^p,\quad i=1,\dots,p-1,\\
&\delta=0 \mbox{ if } z\not\equiv 0 \modu p,\\
&\delta=1 \mbox{ if } z\equiv 0 \modu p.
\end{split}
\end{equation}
where $\s_i$ is an integral ideal of $\Z[\zeta]$. Recall that, when $\s_1$ is a principal ideal, there is a classical proof of the first case of Fermat's Last Theorem.
\item
In the sequel, we shall use the classical acronym FLT for Fermat's Last Theorem or now  Fermat-Wiles Theorem.
\item
The Fermat equation can be written
\begin{displaymath}\label{e2}
\prod_{i=0}^{p-1}(x+\zeta^{i} y) +z^p=0.
\end{displaymath}
It is always possible to choose the pair
$\{x,y\}\subset \{x, y, z\}$ such that $x\not\equiv y \modu p$.
\begin{itemize}
\item
It is clear if $z\equiv 0 \modu p$.
\item
If $xyz\not\equiv 0 \modu p$, from $p>5$ assumed,  it is always possible to choose a pair $(x,y)$ such that $x\not\equiv y \modu p$, because $x\equiv y\equiv z \modu p$ would imply that $p=3$.
\end{itemize}
It is  important to note that in the sequel of this article,we suppose that
$x\not\equiv y \modu p$.
\item
In all the sequel of this monograph, let us denote $t\in\N$ for
$t\equiv -\frac{x}{y}\modu p,\quad 0\leq t\leq p-1$.
In first case of FLT, we can assume that $t$ does not belong to the set of three values
$\{0,1,p-1\}$.
\item
$t=0$ would imply $x\equiv 0 \modu p$.

\item
$t=p-1$ would imply $y\equiv x\modu p$ excluded by previous bullet item.
\item
$t=1$ would imply $x+y\equiv 0 \modu p$ and so $z\equiv 0 \modu p$.
\end{itemize}
{\bf Remark:} In almost all the results of this article, we assume only relation (\ref{e1}) corresponding to
the diophantine equation in $\Z$, for $x,y\in\Z-\{0\}$ mutually co-prime,
\begin{equation}\label{e14113}
\frac{x^p+y^p}{x+y}=p^\delta t_1^p,\quad t_1\in\Z, \quad \delta\in\{0,1\}.
\end{equation}
The Fermat equation $x^p+y^p+z^p=0$ implies diophantine equation (\ref{e14113}), but  the converse is false. In that meaning, almost our results deal with a wider set of diophantine equations than the Fermat equation.
  This remark is detailed in the section \ref{s08081} p. \pageref{s08081} dealing on the first case of the Generalized Fermat-Wiles equation $x^p+y^p+c\times z^p=0,\quad xy\times (x^2-y^2)\not\equiv 0 \modu p$ and the partial Barlow equation $\frac{x^p+y^p}{x+y}=t_1^p,\quad xy(x^2-y^2)\not\equiv 0 \modu p$.
\subsection{Preliminary results}\label{s1}
This section contains some elementary congruences $\modu p$ in $\Z[\zeta]$ and an upper bound estimate of the class number of $\Q(\zeta)$.
%
\subsubsection{Elementary congruences $\modu p$ in $\Z[\zeta]$ }\label{s2}
\begin{prop}\label{p2}
Let $\alpha,\beta\in\Z(\zeta)$
with $\alpha\not\equiv 0 \mod \pi$
and $\alpha\equiv\beta\modu \pi$.
Then $\alpha^p\equiv \beta^p  \modu \pi^{p+1}$.
\begin{proof}
Let $\lambda=(\zeta-1)$.
We have
$\alpha-\beta\equiv 0 \modu \pi
\Rightarrow
(\alpha-\zeta^k\beta)\equiv 0 \modu \pi$ for $k=0,1,\dots,p-1$.
Therefore, for all
$k, \quad 0\leq k \leq p-1$, there exists $a_k\in \N,\quad 0\leq a_k \leq p-1$, such that
$(\alpha-\zeta^k\beta)\equiv \lambda a_k \modu \pi^2$.
For another value $l,\quad 0\leq l\leq p-1$, we have, in the same way,
$(\alpha-\zeta^l\beta)\equiv \lambda a_l \modu \pi^2$,
hence $(\zeta^k-\zeta^l)\beta\equiv \lambda(a_k-a_l) \modu \pi^2$.
For $k\not=l$ we get $a_k\not= a_l$, because $\pi \| (\zeta^k-\zeta^l)$ and
because hypothesis $\alpha\not\equiv 0 \modu \pi$ implies that $\beta\not\equiv 0 \modu \pi$.
Therefore, there exists one and only one $k$ such that
$(\alpha-\zeta^k\beta)\equiv 0 \modu \pi^2$. Then, we have
$ \prod_{j=0}^{p-1}(\alpha-\zeta^j\beta)
= (\alpha^p-\beta^p)\equiv 0 \modu \pi^{p+1}$.
\end{proof}
\end{prop}
{\bf Remark:}
The slight improvement in $\pi^2$, comparing with the result $\alpha^p\equiv \beta^p\modu pZ[\zeta]=\pi^{p-1}$, quoted in the literature, will be used for the proof of second case of Fermat's Last Theorem when $p\parallel y$,
see thm \ref{t14102} p. \pageref{t14102}.
%
\begin{cor} \label{p1}
Let $\alpha\in\Z[\zeta],\quad \alpha\not\equiv 0 \modu \pi$. Then, there exists $b\in\Z$ such that  $\alpha^p\equiv b^p \mod \pi^{p+1}$ where $b\in\Z$.
\begin{proof}
Immediate consequence of proposition \ref{p2} p. \pageref{p2} where there exists $b\in\Z$ with $b\equiv \beta \modu \pi$.
\end{proof}
\end{cor}
{\bf Remark:} This result is a slight improvement in $\pi^2$ of the result quoted in the literature :
$a^p\equiv b \mod p$, see for instance Washington \cite{was}, lem. 1.8 p 5.
In a former version, we have shown $\alpha^p\equiv b\modu \pi^p$. Professor G. Terjanian mentionned us  the improvement up to corollary \ref{p1} p.\pageref{p1}.
%
\subsubsection{An upper bound estimate of the class number of the cyclotomic field}
In the next theorem, we compute an explicit rough upper bound of the class number of $\Q(\zeta)$. This one shall be used in the sequel of this article.
\begin{lem}\label{thmhb}
Let $h$ be the class number of $\Q(\zeta)$. Then
$h <p^{p-2}$.
\begin{proof}
Recall that in this paper, we assume $p>5$.
Let $h=h^*\times h^+$ be the class number of $\Q(\zeta)$ where $h^*$ is the first factor and $h^+$ is the second factor.
From Ribenboim(\cite{rib}) p 132, derived from Lepisto and Metsankyla, formula 6.7
\begin{displaymath}
h^*<2p\times (\frac{p}{24})^{(p-1)/4}.
\end{displaymath}
From (\cite{rib}) formula 1.24 p 97,
\begin{displaymath}
h^+=\frac{2^{(p-3)/2}}{R}\times
\prod_{k=1}^{(p-3)/2}|\sum_{j=0}^{(p-3)/2} \eta^{2kj}log|1-\zeta^{g^j}||,
\end{displaymath}
where $\eta\in\C$ is a $(p-1)$-primitive root of $1,\quad g\in\N$ is a primitive root $\modu p$ and $R\in\R_+$ is the Regulator of $\Q(\zeta)$.
Zimmert  has shown in \cite{zim} corollary p 375, that for a number field $K$ of signature $(r_1,r_2)$,
of Regulator $R$,  which contains $w$ roots of unity, we have the inequality
$2\frac{R}{w}\geq 0.04\times exp(0.46r_1+0.1r_2)$. Then, for the cyclotomic field $\Q(\zeta)$, we have $w=2p$, $r_1=0$, $r_2=\frac{p-1}{2}$  and
$R>0.04\times \frac{w}{2}\times exp(0.1\frac{(p-1)}{2})$;
thus
$R>0.04\times p\times exp(0.1\times\frac{(p-1)}{2})$.
We have the rough estimate
$|\eta^{2kj}\times log|1-\zeta^{g^j}||=|log|1-\zeta^{g^j}||
< log(p)$ because $p>5$.
Therefore,
\begin{displaymath}
h^+<  \frac{2^{(p-3)/2}}{(0.04\times p\times exp(0.1\times\frac{(p-1)}{2}))}
 \times(\frac{(p-1)}{2}log(p))^{(p-3)/2}.
\end{displaymath}
From these upper bounds of $h^*$ and of $h^+$ we get
\begin{displaymath}
\begin{split}
& h<  Q(p)=\{2p\times(\frac{p}{24})^{(p-1)/4}\}\times\\
& \frac{2^{(p-3)/2}}{(0.04\times p\times exp(0.1\times\frac{(p-1)}{2}))}
 \times(\frac{(p-1)}{2}log(p))^{(p-3)/2}. \\
\end{split}
\end{displaymath}
A rough calculation shows that $Q(p)<p^{p-2}$ for all primes $p>5$.
\end{proof}
\end{lem}
%
\clearpage
\section{FLT  first case: Some congruences $\modu p$ }\label {s3}
In this section, we study the first case : we suppose that $xyz\not\equiv 0 \modu p$. Our method is  similar to Eichler's method in (\cite{eic}) and (\cite{rib}). Implicitly, we suppose that $e_p<p-2$, a result that we have proved in theorem \ref{thmhb} p. \pageref{thmhb}.
The next subsection gives  principle of the method to obtain congruences $\modu p$ connected to FLT.
\subsection{Principle of the method}\label{s12111}
This lemma is purely technical and does not take part in the method.
\begin{lem}\label{lem1}
Let $g\in\N,\quad 1\leq g\leq p-1$.
Then $(x+\zeta^g y)\Z[\zeta]$ has at least one prime divisor different of $\pi$.
\begin{proof}
From $xyz\not\equiv 0 \modu p$, we deduce that $x+y\not\equiv 0 \modu \pi$.
It is enough to prove that $(x+\zeta^g y)\not\in\Z[\zeta]^*$.
Suppose that $(x+\zeta^g y)\in\Z[\zeta]^*$ and search for a contradiction:
then, we have $\frac{x^p+y^p}{x+y}=\pm 1$, hence $(x^p+y^p)=\pm (x+y)$.
In an other part, we have $x^p+y^p\equiv x+y\modu p$.
From these two relations, we obtain $x^p+y^p=x+y$.
If $\frac{x}{y}>0$, this relation is clearly impossible because we not have $|x|=|y|=1$.
Then, suppose that $\frac{x}{y}<0$.
Suppose, without loss of generality, that $-1 <\frac{x}{y} <0$ and $x>0$:
consider the function $f(u)=x^u-(-y)^u-x-y$ of the variable $u\in\R,\quad u\geq 1$.
We have $f(1)=0$ and $f_u^\prime(u)=u(x^{u-1}-(-y)^{u-1})>0$ for $1<u<\infty$ and therefore $f(u)\not=0$
and thus $x^p+y^p\not= x+y$.
\end{proof}
\end{lem}
%
The next proposition gives the elementary principle of the method.
\begin{prop}\label{p4}
Let $h=p^{e_p}\times h_2,\quad h_2\not\equiv 0 \modu p$, be the class number of $\Q(\zeta)/\Q$.
Let  $f\in\N, \quad p-1\geq f \geq e_p+2$.
Let  $I_f$ be {\bf any} set of $f$ natural distinct numbers
\begin{equation}\label{ea30}
I_f =\{g_i \quad | \quad g_i\in\N,\quad i=1,\dots,f,\quad 1\leq g_i \leq p-1  \}
\end{equation}
such that the two sets $I_f$
 and $I_f^\prime=\{p-g_i\quad |\quad i=1,\dots,f \}$ are not equal.

Then, there exists at least one set of $f$ natural numbers $L(I_f)$, not all simultaneously null,
\begin{equation}\label{e4}
L(I_f)=\{\ l_i \ | \quad l_i \in\N,\quad 1\leq i\leq f,
\quad 0\leq l_i\leq p-1\ \}
\end{equation}
such that simultaneously,
\begin{equation}\label{e5}
\begin{split}
&\prod_{i=1}^f (x+\zeta^{g_i} y)^{l_i}=\eta\times
\gamma^p,\quad \gamma\in\Z[\zeta],\quad \eta\in\Z[\zeta+\zeta^{-1}]^*,\\
&\prod_{i=1}^f (\frac{(x+\zeta^{g_i}y)}
{(x+\zeta^{-g_i}y)})^{l_i}=(\frac{\gamma}{\overline{\gamma}})^p,\\
&\prod_{i=1}^f (x+\zeta^{g_i} y)^{l_i} - \prod_{i=1}^f(x+\zeta^{-g_i}  y)^{l_i} \not= 0, \\
&\prod_{i=1}^f (x+\zeta^{g_i} y)^{l_i} - \prod_{i=1}^f(x+\zeta^{-g_i}  y)^{l_i}
\equiv 0 \modu \pi^{p+1}.\\
\end{split}
\end{equation}
\begin{proof} $ $
From theorem \ref{thmhb} p. \pageref{thmhb}, we have $e_p<p-2$.
Let  $I_f$ be a set verifying relation \ref{ea30} p. \pageref{ea30}.
Consider the set $P$ of algebraic numbers defined by
\begin{displaymath}\label{ea50}
 P=\{\ \prod_{i=1}^f (x+\zeta^{g_i} y)^{m_i}\quad |\quad 0\leq m_i\leq p-1\quad,
m_i\mbox{ not all } 0  \}.
\end{displaymath}
We have $(x+\zeta^{g_i}y)\Z[\zeta]=\s_{g_i}^p$ where
$\s_{g_i},\quad i=1,\dots,f$, is an integral ideal of $\Q(\zeta)$.
Consider the set
$E=\{ \prod_{i=1}^f \s_{g_i}^{m_i}\quad | \quad 0\leq m_i\leq p-1,
\quad m_i \mbox{ not all } 0 \}$.
We have $Card(E)=p^f-1 \geq p^{e_p+2}-1 > p^{e_p+1}$.
Therefore, from Dirichlet boxes principle, there exists   a subset $E_0$ of $E$ with at least $p+1$
ideals  $\mathbf a_j,\quad j=1,\dots,p+1$, having the  same $p-$component in the class group $Cl(\Q(\zeta))$.
For each pair
$(\mathbf a_1,\mathbf a_j)\subset E_0, \quad j=2,\dots,p+1$,
we can write , from classical factorization of Fermat equation in the cyclotomic field,
\begin{equation}\label{e6}
\begin{split}
&P_1 = \prod_{i=1}^f (x+\zeta^{g_i} y)^{m_i}\Z[\zeta]= \mathbf a_1^p,\quad m_i\mbox{ not all }0,\\
&P_j = \prod_{i=1}^f (x+\zeta^{g_{i,j}} y)^{n_{i,j}}\Z[\zeta]
= \mathbf a_j^p,\quad n_{i,j}\mbox{ not all }0,\quad j=2,\dots,p+1,
\end{split}
\end{equation}
where the two integral ideals $\mathbf a_1$ and $\mathbf a_j,\quad j=2,\dots,p+1,$ are in the same class
of the class group of $\Q(\zeta)$, so
$\frac{\mathbf a_1}{\mathbf a_j}= \gamma_j\Z[\zeta],\quad \gamma_j\in\Q(\zeta)$.
Therefore, we have
\begin{displaymath}\label{e7}
\prod_{i=1}^f \frac{(x+\zeta^{g_i} y)^{n_{i,j}}}{(x+\zeta^{g_i} y)^{m_i}}
= \varepsilon_j\times\gamma_j^p,\quad \varepsilon_j\in\Z[\zeta]^*.
\end{displaymath}
From Kummer's lemma, we have $\varepsilon_j=\zeta^{v_j}\times\eta_j$, where $v_j\in \N,\quad 0\leq v_j\leq p-1$,
and where
$\eta_j\in\Z[\zeta+\zeta^{-1}]^*$.
Therefore, we deduce that
\begin{displaymath}\label{e8}
\prod_{i=1}^f \frac{(x+\zeta^{g_i} y)^{n_{i,j}}}{(x+\zeta^{g_i} y)^{m_i}}=
\zeta^{v_j}\times\eta_j\times\gamma_j^p.
\end{displaymath}
Because there are more than $p$ different ideals $\mathbf a$ defined by relation (\ref{e6}) p. \pageref{e6} in the same class, there are at least two different pairs of ideals
$(\mathbf a_1,\mathbf a_j),\quad (\mathbf a_1,\mathbf a_{j^\prime})$ with the same $v$.
Therefore, we have
\begin{displaymath}\label{e8f}
\prod_{i=1}^f \frac{(x+\zeta^{g_i} y)^{n_{i,j^\prime}}}{(x+\zeta^{g_i} y)^{m_i}}=
\zeta^{v_j}\times\eta_{j^\prime}\times\gamma_{j^\prime}^p,\quad \eta_{j^\prime}\in\Z[\zeta+\zeta^{-1}]^*,\quad\gamma_{j^\prime}\in\Q(\zeta),
\end{displaymath}
hence
\begin{displaymath}
\begin{split}
&\prod_{i=1}^f (x+\zeta^{g_i} y)^{n_{i,j}-n_{i,j^\prime}} =\eta_1\times\gamma_1^p,
\quad n_{i,j}-n_{i,j^\prime} \mbox{ not all }0,\\
&\eta_1\in\Z[\zeta+\zeta^{-1}]^*, \quad \gamma_1\in\Z(\zeta).
\end{split}
\end{displaymath}
Therefore, there exists a set
$L(I_f)=\{ l_i \quad | \quad i\in I_f,\quad 0\leq l_i \leq p-1\quad \}$,
where the $l_i$ are not all simultaneously null, such that

\begin{displaymath}\label{ea80}
\prod_{i=1}^f (x+\zeta^{g_i} y )^{l_i}=\eta_1\times\gamma_1^p,\quad
 \eta_1\in\Z[\zeta+\zeta^{-1}]^*,\quad\gamma_1\in\Z[\zeta].
\end{displaymath}
Then, we have
\begin{displaymath}
\begin{split}
&\prod_{i=1}^f (x+\zeta^{-g_i}  y )^{l_i}=\eta_1\times \overline{\gamma_1^p},\\
&\prod_{i=1}^f (\frac{(x+\zeta^{g_i}y)}
{(x+\zeta^{-g_i}y)})^{l_i}=(\frac{\gamma_1}{\overline{\gamma_1}})^p.
\end{split}
\end{displaymath}
From proposition \ref{p2} p.\pageref{p2} we deduce that
\begin{displaymath}\label{ea81}
\prod_{i=1}^f (x+\zeta^{g_i} y )^{l_i}-\prod_{i=1}^f (x+\zeta^{-g_i}  y )^{l_i}
\equiv 0 \modu \pi^{p+1}.
\end{displaymath}
Show that we have :
\begin{displaymath}\label{ea82}
A=\prod_{i=1}^f (x+\zeta^{g_i} y )^{l_i}-\prod_{i=1}^f (x+\zeta^{-g_i}  y )^{l_i}\not= 0.
\end{displaymath}
equivalent to
\begin{displaymath}\label{e09081}
\prod_{i=1}^f (x+\zeta^{g_i} y )^{l_i}-\prod_{i=1}^f (x+\zeta^{p-g_i}  y )^{l_i}\not= 0.
\end{displaymath}
Suppose that $A=0$ and search for a contradiction: from hypothesis on $I_f$, we  have
$\{g_1,g_2,\dots, g_f \} \not= \{p-g_1,p-g_2,\dots,p-g_f \}$; therefore there exists
$i,\quad 1\leq i \leq f$ with
$g_i\not\in \{p-g_1,\dots,p-g_f\}$. For this $g_i$, from lemma(\ref{lem1}),
there exists $\theta\in \{p-g_j\quad| \quad j=1,\dots,f\},\quad g_i\not=\theta$
and a prime ideal $\mathbf q\not= \pi $ of $\Z[\zeta]$
such that
$(x+\zeta^{g_i} y)\equiv (x+\zeta^{\theta}  y)\equiv 0 \modu \mathbf q$,
so
$ (\zeta^{g_i}-\zeta^{\theta} )y\equiv 0 \modu \mathbf q$,
hence
$ y\equiv 0 \modu \mathbf q$,
so
$ x\equiv 0 \modu \mathbf q$,
contradiction, because $x$ and $y$ are coprime.
\end{proof}
\end{prop}
%
\begin{thm}{ *** }\label{t2009}
Let $h=p^{e_p} \times h_2,\quad h_2\not\equiv 0 \modu p$, be the class number of $\Q(\zeta)/\Q$.
Let  $f\in\N, \quad p-1\geq f\geq \mbox{ min}(e_p+2,\frac{p-1}{2})$.
Let $I_f$ be {\bf any} set of $f$ distinct natural numbers
\begin{equation}\label{ea30bis}
I_f =\{g_i \quad | \quad g_i\in\N,\quad i=1,\dots,f,\quad 1\leq g_i \leq p-1  \}.
\end{equation}
such that the sets $I_f$ and $I_f^\prime=\{p-g_i\quad |\quad  i=1,\dots,f\}$ are not equal.
Then, there exists at least one set of $f$ natural numbers $L(I_f)$ not all simultaneously null
\begin{equation}\label{e4bis}
L(I_f)=\{\ l_i \ | \quad l_i \in\N,\quad 1\leq i\leq f,\quad 0\leq l_i\leq p-1\ \}
\end{equation}
such that simultaneously,
\begin{displaymath}\label{e5bis}
\begin{split}
&\prod_{i=1}^f (\frac{(x+\zeta^{g_i}y)}{(x+\zeta^{-g_i}y)})^{l_i}=\beta^p,
\quad\beta\in\Q(\zeta),\quad\beta\not=1,\\
&\prod_{i=1}^f (x+\zeta^{g_i} y)^{l_i}
- \prod_{i=1}^f(x+\zeta^{-g_i}  y)^{l_i} \equiv 0 \modu \pi^{p+1}.
\end{split}
\end{displaymath}
\begin{proof}
If $f\leq\frac{p-1}{2}$ in the proposition \ref{p4} p. \pageref{p4}, the theorem is proved.
If $f>\frac{p-1}{2}$ in the relation (\ref{e5}) p. \pageref{e5} of proposition \ref{p4} p. \pageref{p4}, then
$\prod_{i=1}^f (x+\zeta^{g_i} y)^{l_i}$ has some factors of the form

$(x+\zeta^{g_i}y)(x+\zeta^{-g_i}y)$ and also
$\prod_{i=1}^f (x+\zeta^{-g_i} y)^{l_i}$ has the same factors of the form
$(x+\zeta^{g_i}y)(x+\zeta^{-g_i}y)$.
It is then  possible to simplify all such factors in relation (\ref{e5}) p. \pageref{e5} of proposition \ref{p4} such that we have, with
$f^\prime\leq\frac{p-1}{2}$,
\begin{displaymath}
\begin{split}
&\prod_{i=1}^{f^\prime} (x+\zeta^{g_i} y)^{l_i}
-\prod_{i=1}^{f^\prime} (x+\zeta^{-g_i}y)^{l_i} \not= 0, \\
&\prod_{i=1}^{f^\prime} (\frac{(x+\zeta^{g_i}y)}{(x+\zeta^{-g_i}y)})^{l_i}=\beta^p,
\quad\beta\in\Q(\zeta),\quad\beta\not=1,\\
&\prod_{i=1}^{f^\prime} (x+\zeta^{g_i} y)^{l_i}
-\prod_{i=1}^{f^\prime} (x+\zeta^{-g_i}y)^{l_i} \equiv 0 \modu \pi^{p+1}.
\end{split}
\end{displaymath}
\end{proof}
\end{thm}
%
{\bf Remarks:}$ $
\begin{itemize}
\item
Note that this definition of $f$, $I_f$ and $L(I_f)$ is suitable in all the sequel  of this article.
\item
It is possible to find "many" sets $I_f$ verifying the hypothesis of proposition \ref{p4} p. \pageref{p4}.
For instance, all the  sets $I_f$ defined by
\begin{displaymath}
\{g_1,\dots,g_{f-1},g_{f}=p-1 \quad|\quad  2\leq g_1<g_2<\dots<g_{f-1}< g_f=p-1\}
\end{displaymath}
verify the hypothesis : if not, $I_f=I_f^\prime$ with
$g_{f}=p-1$ would imply $ p-g_{f}=1\in I_f$ which contradicts $1\not\in I_f=\{g_1,g_2,\dots,g_{f}\}$.
\item
Observe that it is possible to replace $e_p$ by $e_p^-\leq e_p$ in theorem \ref{t2009} p.\pageref{t2009}, because components of the $p$-class group $C_p^+$ does not intervene in fractions $\frac{x+\zeta^{g_i} y}{x+\zeta^{-g_i}y}$.
\end{itemize}
%
\subsection{Definitions and notations}\label{s3s1}
Here, we fix some notations used in the sequel of this paper.
\begin{displaymath}\label{ea85}
\begin{split}
&f\in\N,\quad p-1\geq f\geq min(e_p+2,\frac{p-1}{2}),\\
& I_f= \{g_i\quad |\quad g_i\in \N,\quad 1\leq g_i\leq p-1,\quad i=1,\dots,f\},
\quad Card(I_f)=f, \\
&I_f^\prime=\{p-g_i\quad|\quad g_i\in I_f,\quad i=1,\dots,f\},\quad I_f\not= I_f^\prime \quad \mbox{assumed}\\
& L(I_f) \quad \mbox {defined in proposition \ref{p4} p.\pageref{p4}}.
\end{split}
\end{displaymath}
Let us define, for the indeterminate $X$, the product $P_L(X)=\prod_{i=1}^f (x+X^{g_i} y)^{l_i}$.
Then, we have $P_L(X^{-1})=\prod_{i=1}^f (x+X^{-g_i} y)^{l_i}$.
We have
\begin{displaymath}
\begin{split}
& P_L(\zeta)=\prod_{i=1}^f (x+\zeta^{g_i} y)^{l_i},\\
& P_L(\zeta^{-1})=\overline{P_L(\zeta)}=\prod_{i=1}^f (x+\zeta^{-g_i}  y)^{l_i}.
\end{split}
\end{displaymath}
Let us denote also
\begin{displaymath}
\begin{split}
& Q_L(\zeta)=\prod_{i=1}^f (y+\zeta^{g_i} x)^{l_i},\\
& Q_L(\zeta^{-1})=\overline{Q_L(\zeta)}=\prod_{i=1}^f (y+\zeta^{-g_i}  x)^{l_i}.
\end{split}
\end{displaymath}
where $f, I_f, L(I_f)$ are  defined in proposition \ref{p4} p.\pageref{p4} and theorem \ref{t2009} p. \pageref{t2009}.
For this set $L(I_f)$, let $S_j(I_f),\quad j=1,\dots,p-1, \quad S_j(I_f) \in \N, \quad 0\leq S_j(I_f)\leq p-1$,
be the  functions defined by
\begin{equation}\label{e18}
\begin{split}
&  S_1(I_f)=\sum_{i=1}^f l_i\times g_i \modu p, \\
&  S_2(I_f)=\sum_{i=1}^f l_i\times (g_i)^2\modu p, \\
&  S_3(I_f)=\sum_{i=1}^f l_i\times (g_i)^3\modu p, \\
& \vdots
\end{split}
\end{equation}
Let $\lambda=\zeta-1$. We have $\zeta^{g_i}=(\lambda+1)^{g_i}
\Rightarrow
\zeta^{g_i}\equiv (1+\lambda\times  g_i) \modu \pi^2$,
$\zeta^{g_i}\equiv (1+\lambda\times g_i+\lambda^2 \times\frac{g_i(g_i-1)}{2}) \modu \pi^3,\dots$.
In the same way, we have
$\zeta^{-g_i} =\zeta^{p-g_i}=(\lambda+1)^{p-g_i}$,
hence
$\zeta^{-g_i} \equiv (1+\lambda (p-g_i)) \modu \pi^2$,
also
$\zeta^{-g_i} \equiv (1-\lambda\times g_i) \modu \pi^2$
and
$\zeta^{-g_i} \equiv (1-\lambda \times g_i +\lambda^2\times \frac{g_i(g_i+1)}{2})
\modu \pi^3,\dots$.
%
\subsection{Some general  congruences $\modu \pi^\nu$ in $\Z[\zeta]$.}
\label{s12112}
Let $F(X)\in\Z[X]$ where $X$ is an indeterminate. In this section, we note $F_X^\prime(X)$ the derivative polynomial of $F(X)$ for the indeterminate $X$. Observe that this section is a purely technical intermediate result.
%
\begin{lem}\label{p5}
Let $\mu,\quad \mu\in \N,\quad  0 < \mu \leq p-1$. Let
$F[X]\in\Z[X]$, where $X$ is an indeterminate.
If $F(\zeta) \equiv 0\modu \pi^\mu $,
then $F_X^\prime(\zeta)\equiv 0 \modu \pi^{\mu-1}$. Moreover, if $0<\mu < p-1$ and $\pi^\mu\| F(\zeta)$ then $\pi^{\mu-1}\| F^\prime_X(\zeta)$.
\begin{proof}
We have
$F(\zeta)\equiv 0\modu \pi^\mu$,
so

$ (1-\zeta)^{p-1-\mu} F(\zeta)\equiv 0 \modu \pi^{p-1}$,
same as
$ (1-\zeta)^{p-1-\mu} F(\zeta)\equiv 0 \modu p$,
so
$ (1-\zeta)^{p-1-\mu} F(\zeta)=p \times f(\zeta)$ with $f(X)\in \Z[X]$,
hence
$(1-X)^{p-1-\mu} F(X)=p \times f(X)+\Phi(X)\times h(X)$, where $h(X)\in\Z[X]$ and
$\Phi(X)$ is the $p$-cyclotomic polynomial $\frac{X^p-1}{X-1}$.
By derivation of this relation, we obtain
\begin{displaymath}\label{ea190}
\begin{split}
&-(p-1-\mu)\times(1-X)^{p-1-\mu-1}\times F(X)+(1-X)^{p-1-\mu}
\times F_X^\prime(X)= \\
&p \times f_X^\prime(X) +\Phi_X^\prime(X)\times h(X)+\Phi(X)\times h_X^\prime(X).
\end{split}
\end{displaymath}
We have
$\Phi(\zeta)=0$ and $\Phi_X^\prime(\zeta)\equiv 0 \modu \pi^{p-2}$,
hence
\begin{displaymath}
- (p-1-\mu)(1-\zeta)^{p-1-\mu-1}\times F(\zeta)+(1-\zeta)^{p-1-\mu}F_X^\prime(\zeta)
\equiv 0 \modu \pi^{p-2},
\end{displaymath}
same  as
\begin{displaymath}
- (p-1-\mu)F(\zeta)+(1-\zeta)F_X^\prime(\zeta)\equiv 0 \modu \pi^{p-2-(p-1-\mu-1)},
\end{displaymath}
so
\begin{equation} \label{e14011}
- (p-1-\mu)F(\zeta)+(1-\zeta)F_X^\prime(\zeta)\equiv 0 \modu \pi^\mu,
\end{equation}
hence
\begin{displaymath}
(1-\zeta)F_X^\prime(\zeta)\equiv 0 \modu \pi^\mu
\end{displaymath}
 and finally
\begin{displaymath}
F_X^\prime(\zeta)\equiv 0 \modu \pi^{\mu-1}.
\end{displaymath}
From equation (\ref{e14011}) p.\pageref{e14011}, we deduce immediatly that if $\pi^\mu\| F(\zeta)$ then
$\pi^{\mu-1}\| F^\prime_X(\zeta)$, which achieves the proof.
\end{proof}
\end{lem}
%
\begin{lem}\label{p6}
Let $G(X)=\frac{F(X)}{H(X)}$ where $F(X), H(X)\in \Z[X]$.
If $F(\zeta) \equiv 0\modu \pi^\mu, \quad 0 <\mu \leq(p-1)$ and
$H(\zeta)\not\equiv 0 \modu \pi$
then $G_X^\prime(\zeta)\equiv 0 \modu \pi^{\mu-1}$.
\begin{proof}
From hypothesis, $F(\zeta)\equiv 0\modu \pi^\mu$. Then
$G_X^\prime(\zeta)=\frac{F_X^\prime(\zeta)}{H(\zeta)}-\frac{F(\zeta)H_X^\prime(\zeta)}{H(\zeta)^2}$, so
$G_X^\prime(\zeta)=\frac{F_X^\prime(\zeta) H(\zeta)-F(\zeta)H_X^\prime(\zeta)}{H(\zeta)^2}$,
 and finally
$G_X^\prime(\zeta)\equiv 0\modu \pi^{\mu-1}$.
\end{proof}
\end{lem}
%
\begin{lem}\label{l202201}
Let $F(X)=\sum c_i X^i\in\Q[X],\quad v_p(c_i)\geq 0$, of the indeterminate $X$. If $F(\zeta)\equiv 0\modu \pi^{p}$ and if $F(1)\equiv 0\modu p^2$ then
$F_X^{^\prime}(\zeta)\equiv 0\modu \pi^{p-1}$.
\begin{proof}
From hypothesis we get
\begin{displaymath}
\begin{split}
& F(\zeta)= p \times G(\zeta)\times (\zeta-1),\\
& G(\zeta)=\sum_{i=0}^{p-2} g_i \zeta^i,\quad v_p(g_i)\geq 0.
\end{split}
\end{displaymath}
Algebraically we get
\begin{displaymath}
\begin{split}
& F(X)=p\times G(X)\times (X-1)+h(X)\times \Phi_p(X),\\
& h(X)\in \Q[X],\quad h(X)=\sum h_i X^i,\quad v_p(h_i)\geq 0,\\
& \Phi_p(X)=X^{p-1}+X^{p-2}+\dots+X+1.
\end{split}
\end{displaymath}
We derive that
\begin{displaymath}
F(X)=p\times G(X)\times (X-1)+(h(X)-h(1))\times\Phi_p(X)+h(1)\times \Phi_p(X),
\end{displaymath}
also that
\begin{displaymath}
\begin{split}
& F(X)=p\times G(X)\times (X-1)+h_1(X)\times(X^p-1)+h(1)\times \Phi_p(X),\\
& h_1(X)\in\Q(X),\quad v_p(h_{1,i})\geq 0.
\end{split}
\end{displaymath}
But $F(1)\equiv 0\modu p^2$ implies that $h(1)\equiv 0\modu p$ because
$\Phi_p(1)=p$.
By derivation we get
\begin{displaymath}
\begin{split}
& F_X^\prime(X)= p\times G_X^\prime(X)\times (X-1)+p\times G(X)\\
& +h_1^\prime(X)\times (X^p-1)+h_1(X)\times p\times X^{p-1}
+h(1)\times\Phi_p^\prime(X),
\end{split}
\end{displaymath}
which leads to
$F_X^\prime(\zeta)\equiv 0\modu \pi^{p-1}$.
\end{proof}
\end{lem}
{\bf Remark:} This improvement $F_X^\prime(\zeta)\equiv \pi^{p-1}$ instead of
$F_X^\prime(\zeta)\equiv 0\modu \pi^{p-2}$ derived from lemma \ref{p6} p. \pageref{p6} shall be of great importance in the sequel of this monograph.
%
\subsection{Congruences $\modu p$ connected to sums
$\sum_{i=1}^f \frac{g_i l_i}{(x+\zeta^{g_i}y)}$.}\label{s02021a}
In this  subsection, we transform the congruences $\modu \pi^{p-1}$ related to products $\prod_{i=1}^f(x+\zeta^{g_i}y)^{l_i}$ in congruences $\modu \pi^{p-1}$ related to sums $\sum_{i=1}^f \frac{g_i l_i}{(x+\zeta^{g_i}y)}$.
\begin{lem}\label{p6a}
Let $P_{L,X=\zeta}^\prime(X)$ be  the value of the derivative  in the inderminate $X$  of
$P_L(X)$ for the value $X=\zeta$.
Let $P_{L,X=\zeta}^\prime(X^{-1})$ be  the value of the derivative  in the inderminate $X$  of
$P_L(X^{-1})$ for the value $X=\zeta$.
Then
\begin{displaymath}\label{e20a}
P_{L,X=\zeta}^\prime(X)-P_{L,X=\zeta}^\prime(X^{-1}) \equiv 0\modu \pi^{p-1}.
\end{displaymath}
\begin{proof}
From proposition \ref{p4} p.\pageref{p4}, we have

$P_L(\zeta)-\overline{P_L(\zeta)}\equiv P_L(\zeta)-P_L(\zeta^{-1})\equiv 0\modu \pi^{p+1},$

$P_L(\zeta)-P_L(\zeta^{-1})=p\times (1-\zeta)^2\times f(\zeta),$

$ f(\zeta)\in\Z[\zeta].$

Therefore, $P_L(X)-P_L(X^{-1})=p\times(1-X)^2\times f(X)+h(X)\times\Phi(X),
\quad f(X)\in\Z[X],\quad  h(X)=\frac{h_1(X)}{X^u}, \quad h_1(X)\in\Z[X],\quad u\in\N, \quad \Phi(X)\in\Z[X],\quad \Phi(X)=\frac{(X^p-1)}{(X-1)}$,
$p$-cyclotomic polynomial.
We have
$P_L(1)-P_L(1^{-1})=0\quad$, therefore we deduce the relation
\begin{displaymath}\label{ea200}
 P_L(X)-P_L(X^{-1})=p\times (1-X)^2\times f(X)+\frac{h_2(X)(1-X)}{X^u}\times \Phi(X),
\end{displaymath}
which leads to
\begin{equation}\label{e300113}
P_L(X)-P_L(X^{-1})=p\times (1-X)^2\times f(X)+\frac{h_2(X)}{X^u}\times (X^p-1).
\end{equation}
The proof is  similar to proposition \ref{p5} p.\pageref{p5}, noticing that, here,
$(X^p-1)_X^\prime=p \times X^{p-1}\equiv 0 \modu \pi^{p-1}$.
\end{proof}
\end{lem}
%
\begin{lem}\label{p7}
For $S_1(I_f)$ defined in relation (\ref{e18}) p.\pageref{e18}, we have
$S_1(I_f)\equiv 0\modu p$.
\begin{proof}
\begin{displaymath}\label{e20}
\prod_{i=1}^f (x+\zeta^{g_i} y)^{l_i} - \prod_{i=1}^f (x+\zeta^{-g_i}  y)^{l_i}
\equiv 0 \modu \pi^{p+1}.
\end{displaymath}
From previous notations,
\begin{displaymath}\label{e21}
\begin{split}
& (x+\zeta^{g_i} y)^{l_i}\equiv (x+y+y(\zeta^{g_i}-1))^{l_i}\\
&\equiv (x+y)^{l_i}+(x+y)^{l_i-1}y g_i l_i \lambda\modu \pi^2,\\
& (x+\zeta^{-g_i} y)^{l_i}\equiv (x+y+y(\zeta^{-g_i}-1))^{l_i}\\
&\equiv (x+y)^{l_i}+(x+y)^{l_i-1} y  (-g_i) l_i \lambda\modu \pi^2.
\end{split}
\end{displaymath}
Denote $S_0(I_f)\in\Z$  the sum $S_0(I_f) = \sum_{i=1}^f  l_i\modu p$.
Then, we deduce that
\begin{displaymath}\label{e22}
\begin{split}
& \prod_{i=1}^f (x+\zeta^{g_i} y)^{l_i} \equiv
(x+y)^{S_0(I_f)}+(x+y)^{S_0(I_f)-1}\times (\sum_{i=1}^f y g_i l_i\lambda) \modu \pi^2\\
& \prod_{i=1}^f (x+\zeta^{-g_i}  y)^{l_i} \equiv
(x+y)^{S_0(I_f)}+(x+y)^{S_0(I_f)-1}\times (-\sum_{i=1}^f y g_i l_i\lambda) \modu \pi^2\\
& \Rightarrow 2\sum_{i=1}^f y g_i l_i \lambda \equiv 0 \modu \pi^2
\Rightarrow S_1(I_f)\equiv 0\modu \pi,\\
& \Rightarrow S_1(I_f)\equiv 0 \modu p.
\end{split}
\end{displaymath}
\end{proof}
\end{lem}
%
\begin{lem}\label{p8}
\begin{displaymath}\label{e23}
\begin{split}
& \prod_{i=1}^f (x+\zeta^{g_i} y)^{l_i} \equiv \prod_{i=1}^f (x+\zeta^{-g_i}  y)^{l_i} \equiv \\
& \prod_{i=1}^f (y+\zeta^{g_i} x)^{l_i} \equiv \prod_{i=1}^f (y+\zeta^{-g_i}  x)^{l_i} \modu \pi^{p+1}.
\end{split}
\end{displaymath}
\begin{proof}
From proposition \ref{p7} p. \pageref{p7}, we have
$\sum_{i=1}^f g_i l_i\equiv 0\modu p$
hence $\zeta^{\sum_{i=1}^f g_i l_i} =1.$
Then we have
\begin{displaymath}\label{e24}
\begin{split}
& \prod_{i=1}^f (x+\zeta^{g_i} y)^{l_i}=
\zeta^{\sum_{i=1}^f g_i l_i}\prod_{i=1}^f (\zeta^{-g_i} x+y)^{l_i}\\
& \prod_{i=1}^f (x+\zeta^{-g_i}  y)^{l_i} =
\zeta^{\sum_{i=1}^f(-g_i l_i)}\prod_{i=1}^f (\zeta^{g_i} x+y)^{l_i}
\end{split}
\end{displaymath}
which leads to the result.
\end{proof}
\end{lem}
%
\begin{lem}\label{p9}
\begin{displaymath}\label{e25}
\sum_{i=1}^f \frac{g_i l_i\zeta^{g_i}}{(x+\zeta^{g_i}y)} +
\sum_{i=1}^f \frac{g_i l_i\zeta^{-g_i} }{(x+\zeta^{-g_i} y)} \equiv 0 \modu \pi^{p-1}.
\end{displaymath}
\begin{proof}
We compute the value of  the derivative in the indeterminate $X$ of $P_L(X)$ and $P_L(X^{-1})$ for $X=\zeta$.
We have
\begin{displaymath}\label{e26}
\begin{split}
& P_{L,X=\zeta}^\prime(X)=
P_{L}(\zeta)\times y\sum_{i=1}^f \frac{g_i l_i\zeta^{g_i-1}}{(x+\zeta^{g_i} y)},\\
& P_{L,X=\zeta}^\prime(X^{-1})=
P_L(\zeta^{-1})\times y\sum_{i=1}^f \frac{(-g_i l_i)\zeta^{-g_i-1}}{(x+\zeta^{-g_i} y)},
\end{split}
\end{displaymath}
and from relation (\ref{e4}) p. \pageref{e4}
\begin{equation}\label{e1109a}
\begin{split}
& P_L(\zeta)\equiv P_L(\zeta^{-1}) \modu \pi^{p+1}\\
& P_L(\zeta)\not\equiv 0 \modu \pi.
\end{split}
\end{equation}

Then from lemma \ref{p6a} p. \pageref{p6a}, we deduce that
$P_{L,X=\zeta}^\prime(X)-P_{L,X=\zeta}^\prime(X^{-1})\equiv 0 \modu \pi^{p-1}$.
Therefore, we have
\begin{displaymath}\label{ea260}
\begin{split}
& P_L(\zeta)\times y\sum_{i=1}^f \frac{g_i l_i\zeta^{g_i-1}}{(x+\zeta^{g_i}y)}
\equiv \\
& P_L(\zeta^{-1})\times y\sum_{i=1}^f \frac{(-g_i) l_i\zeta^{-g_i-1}}
{(x+\zeta^{-g_i} y)} \modu \pi^{p-1}.
\end{split}
\end{displaymath}
Then,  we have, from relation (\ref{e1109a}) p. \pageref{e1109a}
\begin{displaymath}\label{ea261}
\begin{split}
& P_L(\zeta)\times y\sum_{i=1}^f \frac{g_i l_i\zeta^{g_i-1}}{(x+\zeta^{g_i}y)} \equiv \\
& P_L(\zeta)\times y\sum_{i=1}^f \frac{(-g_i) l_i\zeta^{-g_i-1}}{(x+\zeta^{-g_i} y)}
\modu \pi^{p-1}.
\end{split}
\end{displaymath}
hence
\begin{displaymath}
\begin{split}
& \sum_{i=1}^f \frac{g_i l_i \zeta^{g_i-1}}{(x+\zeta^{g_i} y)}\equiv \\
& \sum_{i=1}^f \frac{(-g_i) l_i\zeta^{-g_i-1}}{(x+\zeta^{-g_i} y)} \modu \pi^{p-1}.\\
\end{split}
\end{displaymath}
which completes the proof.
\end{proof}
\end{lem}
%
\begin{lem}\label{p9a}
\begin{displaymath}\label{e30a}
y\sum_{i=1}^f \frac{g_i l_i\zeta^{g_i}}{(x+\zeta^{g_i} y)} -
x\sum_{i=1}^f \frac{g_i l_i\zeta^{g_i}}{(y+\zeta^{g_i} x)} \equiv 0 \modu \pi^{p-1}.
\end{displaymath}
\begin{proof}
Similar to lemma \ref{p9} p.\pageref{p9} proof, starting from lemma \ref{p8} p.\pageref{p8}.
\end{proof}
\end{lem}
%
\begin{thm}{ *** }\label{p9b}
With previous notations,
\begin{equation}\label{e35}
\begin{split}
&x\sum_{i=1}^f \frac{g_i l_i}{(x+\zeta^{g_i} y)} -
 y\sum_{i=1}^f \frac{g_i l_i}{(y+\zeta^{g_i} x)} \equiv 0 \modu \pi^{p-1}\\
& \sum_{i=1}^f \frac{g_i l_i}{(x+\zeta^{g_i} y)} +
  \sum_{i=1}^f \frac{g_i l_i}{(x+\zeta^{-g_i}  y)} \equiv 0 \modu \pi^{p-1}
\end{split}
\end{equation}
\begin{proof} $ $
\begin{itemize}
\item
From lemma \ref{p9a} p.\pageref{p9a}, we get
\begin{displaymath}\label{e36}
\sum_{i=1}^f \frac{g_i l_i\zeta^{g_i} y}{(x+\zeta^{g_i} y)} -
 \sum_{i=1}^f \frac{g_i l_i\zeta^{g_i} x}{(y+\zeta^{g_i} x)} \equiv 0 \modu \pi^{p-1},
\end{displaymath}
which leads to
\begin{displaymath}
 \sum_{i=1}^f g_i l_i(1 -\frac{x}{(x+\zeta^{g_i} y)}) -
  \sum_{i=1}^f g_i l_i(1 -\frac{y}{(y+\zeta^{g_i} x)}) \equiv 0 \modu \pi^{p-1}.
\end{displaymath}
But, from lemma \ref{p7} p.\pageref{p7},  $\sum_{i=1}^f g_i l_i\equiv 0 \modu p$.
\item
From proposition \ref{p9} p.\pageref{p9}, we have
\begin{displaymath}\label{e37}
\sum_{i=1}^f \frac{g_i l_i\zeta^{g_i} y}{(x+\zeta^{g_i} y)} +
 \sum_{i=1}^f \frac{g_i l_i\zeta^{-g_i} y}{(x+\zeta^{-g_i} y)} \equiv 0 \modu \pi^{p-1},
\end{displaymath}
and hence
\begin{displaymath}
 \sum_{i=1}^f g_i l_i(1 -\frac{x}{(x+\zeta^{g_i} y)}) +
  \sum_{i=1}^f g_i l_i(1 -\frac{x}{(x+\zeta^{-g_i} y)}) \equiv 0 \modu \pi^{p-1}.
\end{displaymath}
\end{itemize}
\end{proof}
\end{thm}
%
\subsection {Eichler's theorem}
We derive from theorem \ref{p9b} p.\pageref{p9b} the  proof of Eichler's theorem; see also for a proof, Washington \cite{was} theorem 6.23 p 107.
\begin{thm} { *** }\label{t011182}
Let $e_p^-$ be the power of $p$ dividing the relative class number $h^-$ of $\Q(\zeta)$.
If first case of FLT fails for $p$, then $e_p^-\geq [\sqrt{p}]-1$.
\begin{proof} $ $
\begin{itemize}
\item
Let $(x+\zeta y)\Z[\zeta]=\s^p$, where $\s$ is an integral ideal of $\Z[\zeta]$.
From $f<\frac{p-1}{2}$, we derive that
$\frac{\prod_{i=1}^f \s ^{l_i}}
{\prod_{i=1}^f \overline{\s}^{l_i}}$ is an ideal whose class belongs to the group $C_p^-=C_p/C_p^+$ of order  $p^{e_p^-}$.
\item
Observe that $e_p^-\geq [\sqrt{p}]-1 \Leftrightarrow e_p^- > \sqrt{p}-2$.
Suppose that $e_p^-+2<\sqrt{p}$ and search for a contradiction.
Let $f=e_p^-+2$. From relation (\ref{e35}) p.\pageref{e35}, we derive with $g_i=i,\quad i=1,\dots,f$,
\begin{displaymath}
\prod_{i=1}^f(x+\zeta^i y)(x+\zeta^{-i} y)\times
(\sum_{i=1}^f \frac{i l_i}{x+\zeta^i y}
+\frac{i l_i}{x+\zeta^{-i} y})\equiv 0 \modu p.
\end{displaymath}
\item
The higher power term in $\zeta$ is
\begin{displaymath}
l_1\times y^{f-1} \times x^f\times \zeta^{(1+2+\dots+f-1+f)-1}
= l_1\times  y^{f-1}\times  x^f \times \zeta^{f(f-1)/2-1}.
\end{displaymath}
In the same way, the smaller  power in $\zeta$ is
\begin{displaymath}
l_1 \times y^{f-1} \times x^f\times \zeta^{-((1+2+\dots+f-1+f)-1)}
= l_1 \times y^{f-1}\times  x^f\times  \zeta^{-(f(f-1)/2-1)}.
\end{displaymath}
\item
We have $e_p+2 < \sqrt{p} \Rightarrow f  < \sqrt{p}\Rightarrow f(f-1) < p
\Rightarrow f(f-1)-1 < p-1 \Rightarrow \frac{f(f-1)}{2} -1 < \frac{p-1}{2}$.
Therefore there is no term in $\zeta^{(p-1)/2}$.
\item
Thus we deduce that $l_1\equiv 0 \modu p$ and even that $l_1 =0$ because $0\leq l_1\leq p-1$; then the higher power in $\zeta \modu p$ is
\begin{displaymath}
2\times l_2\times y^{f-1}\times x^f\times\zeta^{(1+2+\dots+f-1+f)-2},
\end{displaymath}
so $l_2 = 0$ and so
on, $l_1 = l_2 = \dots= l_f =  0$, which contradicts the fact that, in $L(I_f)$,
the $l_i$ are not all simultaneously null.
Therefore, if first case of FLT fails for $p$  then
$e_p^-\geq [\sqrt{p}]-1$.
\end{itemize}
\end{proof}
\end{thm}
%
%
\clearpage
\section{FLT  first case : Some congruences on Mirimanoff polynomials  }\label{s4}
In this section, we generalize Mirimanoff's congruences. The next theorem is the first result towards this direction.
We have proved that if $x,y,z$ are solutions of the first case of FLT equation, whe should have
\begin{equation}\label{e23012}
 \sum_{i=1}^f \frac{g_i l_i}{(x+\zeta^{g_i} y)}
+\sum_{i=1}^f \frac{g_i l_i}{(x+\zeta^{-g_i} y)}\equiv 0 \modu \pi^{p-1}.
\end{equation}
The results of relation (\ref{e35}) p.\pageref{e35}, used in this article directly for the context of Fermat's equation can in fact be generally proved, in a  context {\bf strictly independant}  of Fermat's  Last Theorem  , with the only hypothesis:

\begin{equation}\label{e23013}
\begin{split}
& a,b\in\Z,\quad a\times b\times (a^2-b^2)\not\equiv 0 \modu p,\\
&\sum_{i=1}^f \frac{g_i l_i}{(a+\zeta^{g_i} b)}
+\sum_{i=1}^f \frac{g_i l_i}{(a+\zeta^{-g_i} b)}\equiv 0 \modu \pi^{p-1}.
\end{split}
\end{equation}
%
\subsection{The foundation Theorem}\label{s12113}
\begin{thm}{ *** }\label{p21a}
Let $I_f$ be the sets defined in the relation  (\ref{ea30bis}).
Let $\phi_{2m+1}(T),\quad m=1,\dots,\frac{p-3}{2}$, be the odd Mirimanoff polynomials of the indeterminate $T$.
Let $t\in\N,\quad 1<t\leq p-1,\quad t\equiv -\frac{y}{x}\modu p$.
Then, we have for all set $I_f$
\begin{displaymath}\label{e50g}
\begin{split}
&S_{1}(I_f)\equiv 0 \modu p,\\
&S_{2m+1}(I_f)\times \phi_{2m+1}(t)\equiv 0 \modu p,\quad m=1,\dots,\frac{p-3}{2}.
\end{split}
\end{displaymath}
\begin{proof}
From lemma \ref{p7} p.\pageref{p7}, we have $S_{1}(I_f)\equiv 0 \modu p$.
From proposition \ref{p9b} p.\pageref{p9b},
\begin{displaymath}\label{e51g}
 \sum_{i=1}^f \frac{g_i l_i}{(x+\zeta^{g_i} y)}
+\sum_{i=1}^f \frac{g_i l_i}{(x+\zeta^{-g_i} y)}\equiv 0 \modu \pi^{p-1}.
\end{displaymath}
We have the algebraic identities
\begin{displaymath}\label{e52g}
\begin{split}
&\frac{(x^p+y^p)}{(x+\zeta^{g_i} y)}=
(x^{p-1}-\zeta^{g_i} x^{p-2}y+\zeta^{2g_i} x^{p-3}y^2-\dots
-\zeta^{(p-2)g_i}x y^{p-2}+\zeta^{(p-1)g_i}y^{p-1}),\\
&\frac{(x^p+y^p)}{(x+\zeta^{-g_i}  y)}=
(x^{p-1}-\zeta^{-g_i}  x^{p-2}y+\zeta^{-2g_i} x^{p-3}y^2-\dots
-\zeta^{-(p-2)g_i}x y^{p-2}+\zeta^{-(p-1)g_i}y^{p-1}).
\end{split}
\end{displaymath}
We deduce of these relations that
\begin{equation}\label{e53g}
\begin{split}
&F(\zeta)=\\
&\sum_{i=1}^f g_i l_i
(x^{p-1}-\zeta^{g_i} x^{p-2}y+\zeta^{2g_i} x^{p-3}y^2-\dots
-\zeta^{(p-2)g_i}x y^{p-2}+\zeta^{(p-1)g_i}y^{p-1})\\
&+\sum_{i=1}^f g_i l_i
(x^{p-1}-\zeta^{-g_i}  x^{p-2}y+\zeta^{-2g_i} x^{p-3}y^2-\dots
-\zeta^{-(p-2)g_i}x y^{p-2}+\zeta^{-(p-1)g_i}y^{p-1})\\
&\equiv 0 \modu \pi^{p-1}.
\end{split}
\end{equation}
From relation (\ref{e53g}) p. \pageref{e53g}, we can write $F(\zeta)=\frac{f(\zeta)}{\zeta^u}$ whith $f(X)\in\Z[X]$ and $u\in\N$.
Then, let $F(X)\in\Q(X)$ be the rational function of the indeterminate $X$ where $\zeta$ is replaced by $X$ in the previous value of $F(\zeta)$, thus $F(X)=\frac{f(X)}{X^u}$.
Here, let us denote $b_X^\prime(X)$  the derivative of any function $b(X)$ for $X$. We define the functions $F_1(X),F_2(X),F_3(X),\dots$ by
\begin{equation}\label{e54g}
\begin{split}
&F_1(X)=X \times F_X^  \prime(X),\\
&F_2(X)=X \times (F_1)_X^\prime(X),\\
&F_3(X)=X \times (F_2)_X^\prime(X),\\
&\vdots
\end{split}
\end{equation}
From  proposition \ref{p6} p.\pageref{p6} we have :
\begin{equation}
\begin{split}
& F(\zeta)\equiv 0 \modu \pi^{p-1},\\
&F_1(\zeta)\equiv 0 \modu \pi^{p-2},\\
&F_2(\zeta)\equiv 0 \modu \pi^{p-3},\\
&F_3(\zeta)\equiv 0 \modu \pi^{p-4},\\
&\vdots\\
&F_{p-3}(\zeta)\equiv 0 \modu \pi^2.
\end{split}
\end{equation}
Then, we obtain
\begin{equation}\label{e56g}
\begin{split}
&F_1(\zeta)=\\
&\sum_{i=1}^f (g_i)^2 l_i
(-\zeta^{g_i} x^{p-2}y+2\zeta^{2 g_i} x^{p-3}y^2-\dots
-(p-2)\zeta^{(p-2) g_i}x y^{p-2}+(p-1)\zeta^{(p-1)g_i}y^{p-1})\\
&-\sum_{i=1}^f (g_i)^2 l_i
(-\zeta^{-g_i}  x^{p-2}y+2\zeta^{-2g_i} x^{p-3}y^2-\dots
-(p-2)\zeta^{-(p-2)g_i}x y^{p-2}+(p-1)\zeta^{-(p-1)g_i}y^{p-1}),
\end{split}
\end{equation}
thus, substituting $1$ to $\zeta$, we get
$F_1(1)=0.$
In the same way for $F_2(\zeta)$ we obtain,
\begin{equation}\label{e57g}
\begin{split}
&F_2(\zeta)=\\
& \sum_{i=1}^f (g_i)^3 l_i
(-\zeta^{g_i} x^{p-2}y+2^2\zeta^{2g_i} x^{p-3}y^2-\dots
-(p-2)^2\zeta^{(p-2)g_i}x y^{p-2}+(p-1)^2\zeta^{(p-1)g_i}y^{p-1})\\
&+\sum_{i=1}^f (g_i)^3 l_i
(-\zeta^{-g_i}  x^{p-2}y+2^2\zeta^{-2g_i} x^{p-3}y^2-\dots
-(p-2)^2\zeta^{-(p-2)g_i}x y^{p-2}+(p-1)^2\zeta^{-(p-1)g_i}y^{p-1}),
\end{split}
\end{equation}
so, by substitution of $1$ to $\zeta$,
\begin{displaymath}
2\times (\sum_{i=1}^f (g_i)^3 l_i
(- x^{p-2}y+2^2 x^{p-3}y^2-\dots
-(p-2)^2 x y^{p-2}+(p-1)^2 y^{p-1}))\equiv 0 \modu \pi,
\end{displaymath}
also
\begin{displaymath}
2\times(\sum_{i=1}^f g_i^3 l_i
(- \frac{y}{x}+2^2 \frac{y^2}{x^2}-\dots
-(p-2)^2 \frac{ y^{p-2}}{x^{p-2}}+(p-1)^2\frac{ y^{p-1}}{x^{p-1}}))\equiv 0 \modu p,
\end{displaymath}
also
\begin{displaymath}
2\times(\sum_{i=1}^f g_i^3 l_i)\times
((- \frac{y}{x}+2^2 \frac{y^2}{x^2}-\dots
-(p-2)^2 \frac{ y^{p-2}}{x^{p-2}}+(p-1)^2\frac{ y^{p-1}}{x^{p-1}}))\equiv 0 \modu p,
\end{displaymath}
hence
$F_2(1)=2\times S_3(I_f)\times \Phi_3(t)\equiv 0 \modu p$.
In the  same way
\begin{equation}\label{e57gb}
\begin{split}
&F_3(1)=0,\\
&F_4(1)=2\times S_5(I_f)\Phi_5(t)\equiv 0 \modu p,\\
&F_5(1)=0,\\
&\vdots\\
&F_{2m}(1)=2\times S_{2m+1}(I_f)\Phi_{2m+1}(t)\equiv 0 \modu p,\\
&\vdots\\
&F_{p-3}(1)=2\times S_{p-2}(I_f)\Phi_{p-2}(t)\equiv 0 \modu p,\\
&F_{p-2}(1)=0.
\end{split}
\end{equation}
which achieves the proof.
\end{proof}
\end{thm}
{\bf A definition :}
We define the polynomials $\phi_{2m+1}^*(T)$ by
\begin{equation}\label{e104051}
\begin{split}
& \phi_{1}^*(T)=1,\\
& \phi_{2m+1}^*(T)=\phi_{2m+1}(T),\quad m=1,\dots,\frac{p-3}{2}.
\end{split}
\end{equation}
Note that, with this definition, the previous theorem can be stated
\begin{displaymath}
\phi_{2m+1}^*(t)\times S_{2m+1}(I_f)\equiv 0\modu p,\quad m=0,\dots,\frac{p-3}{2}.
\end{displaymath}
The interest of this definition will appear more clearly in the sequel of this section.
%
The next lemma plays an important part in the sequel of this article.
\begin{lem}\label{propaa}
Let $S_{2m+1}(I_f),\quad m=0,\dots,\frac{p-3}{2}$, be the sums defined in the relation (\ref{e18}) p.\pageref{e18}.
Then, if $S_{1}(I_f)\equiv 0 \modu p$ and $S_{2m+1}(I_f)\times\Phi_{2m+1}(t)\equiv 0\modu p$ for $m=1,\dots,\frac{p-3}{2}$  then
(reciprocal result of  theorem \ref{p21a} p.\pageref{p21a}) :
\begin{displaymath}
\sum_{i=1}^f \frac{g_i l_i}{(x+\zeta ^{g_i} y)}+
\sum_{i=1}^f \frac{g_i l_i}{(x+\zeta^{-g_i}y)}
\equiv 0 \modu p.
\end{displaymath}
\begin{proof}
Let us  define $F(X)\in\Q(X)$   by:
\begin{equation}\label{e31101}
\begin{split}
&F(X)=\\
&\sum_{i=1}^f g_i l_i\times
(x^{p-1}-X^{g_i} x^{p-2}y+X^{2g_i} x^{p-3}y^2-\dots
-X^{(p-2)g_i}x y^{p-2}+X^{(p-1)g_i}y^{p-1})\\
&+\sum_{i=1}^f g_i l_i\times
(x^{p-1}-X^{-g_i}  x^{p-2}y+X^{-2g_i} x^{p-3}y^2-\dots
-X^{-(p-2)g_i}x y^{p-2}+X^{-(p-1)g_i}y^{p-1}).
\end{split}
\end{equation}
We have
\begin{displaymath}
F(\zeta)=(x^p+y^p)\times (\sum_{i=1}^f \frac{g_i l_i}{(x+\zeta ^{g_i} y)}+
\sum_{i=1}^f \frac{g_i l_i}{(x+\zeta^{-g_i}y)})
\end{displaymath}
It is enough to prove that $F(\zeta)\equiv 0 \modu p$.
The function $F(X)$ is of the form $F(X)=\frac{a(X)}{X^u}$ where $a(X)\in\Z[X]$ and $u\in\N$.
With  $\pi$-adic developpment, we have
\begin{displaymath}
\begin{split}
&F(\zeta)=F(1)+(\zeta-1)F_X^{\prime}(1)+\frac{(\zeta-1)^2}{2!}F_X^{(2)}(1)+\dots\\
&+\frac{(\zeta-1)^{p-2}}{(p-2)!}F_X^{(p-2)}(1)
+\sum_{i=p-1}^{+\infty} \frac{(\zeta-1)^i}{i!} F_X(1),
\end{split}
\end{displaymath}
where we verify directly that, for $i\geq p-1$,
we  obtain $\frac{(\zeta-1)^i}{i!}F^{(i)_X(1)}\equiv 0\modu \pi^i,\quad i\geq p-1$:
\begin{enumerate}
\item
For $i=p-1$, it is clear.
\item
For $i>p-1$:
\begin{itemize}
\item
If $j\geq i$:
let  the term $c_jX^j$ of $F(X),\quad c_j\in\Z$. We get
$\frac{(X^j)_X^{(i)}}{i!}=\frac{j(j-1)\dots(j-i+1)}{i!}X^{j-i}
=C_j^i \times X^{j-i}$.
Thus, we get $C_j^i\equiv 0\modu p^\alpha,\quad \alpha\geq 0$.
\item
If $i>j\geq 0$ then $(X^j)^{(i)}_X=0$.
\item
If $j<0$ then let $j_1=-j$. We get
\begin{displaymath}
\frac{(X^j)^{(i)}_X}{i!}=\pm \frac{j_1(j_1+1)\dots (j_1+i-1)}{i!}
=\pm C_{j_1+i-1}^i\in\Z.
\end{displaymath}
We  obtain also $C_{j_1+i-1}^i\equiv 0\modu p^\alpha,\quad \alpha\geq 0$.
\end{itemize}
which leads to the result.
\end{enumerate}
So we obtain
\begin{equation}\label{e104291}
\begin{split}
&F(\zeta)\equiv F(1)+(\zeta-1)F_X^{(1)}(1)
+\frac{(\zeta-1)^2}{2!}F_X^{(2)}(1)+\dots\\
&+\frac{(\zeta-1)^{p-2}}{(p-2)!}F_X^{(p-2)}(1)\modu \pi^{p-1}.
\end{split}
\end{equation}
We have $F(1)\equiv 0\modu p.$
Let us denote, with $F_X^\prime(X)$  the derivative of $F(X)$ for $X$,
\begin{displaymath}\label{e}
\begin{split}
&F_1(X)=X \times F_X^  \prime(X),\\
&F_2(X)=X \times (F_1)_X^\prime(X),\\
&F_3(X)=X \times (F_2)_X^\prime(X),\\
&\vdots
\end{split}
\end{displaymath}
We have
\begin{displaymath}
\begin{split}
& F_2(X)=X^2 F_X^{(2)}(X)+X F_X^\prime(X)\\
& F_3(X)=X^3 F_X^{(3)}(X)+3X^2 F_X^{(2)}(X)+X F_X^\prime(X)\\
& \vdots
\end{split}
\end{displaymath}
Then,
\begin{equation}\label{e22021}
\begin{split}
& F_X^\prime(1)=F_1(1),\\
& F_X^{(2)}(1)=F_2(1)-F_1(1),\\
& F_X^{(3)}(1)=F_3(1)-3F_2(1)+2F_1(1),\\
& \vdots
\end{split}
\end{equation}
The same computation as in relations (\ref{e53g}) p.\pageref{e53g},
(\ref{e56g}) p.\pageref{e56g}, (\ref{e57g}) p.\pageref{e57g} and
(\ref{e57gb}) p. \pageref{e57gb} of  proposition \ref{p21a} p. \pageref{p21a}, leads by derivation of $F(X)$, to
\begin{equation}\label{lab02}
\begin{split}
& F(1)\equiv 2 \times S_1(I_f)\times\phi_1(t) +2\sum_{i=1}^f g_i l_i=2 \times S_1(I_f)\times (\phi_1(t)+1)\equiv 0\modu p,\\
& F_1(1)=0,\\
& F_2(1)\equiv 2\times  S_3(I_f)\times \phi_3(t)\equiv 0\modu p,\\
& F_3(1)=0,\\
& F_4(1)\equiv 2\times  S_5(I_f)\times \phi_5(t)\equiv 0\modu p,\\
& \vdots\\
& F_{p-3}(1)\equiv 2\times  S_{p-2}(I_f)\times \Phi_{p-2}(t)\equiv 0\modu p\\
& F_{p-2}(1)=0.
\end{split}
\end{equation}
From hypothesis, relations (\ref{e22021}) p.\pageref{e22021} and
(\ref{lab02}) p.\pageref{lab02}, we have
$F(1)\equiv  F_X^\prime(1)\equiv F_X^{(2)}(1)\equiv \dots \equiv F_X^{(p-2)}(1)
\equiv 0\modu p$ and then
$F(\zeta)\equiv 0 \modu p$.
\end{proof}
\end{lem}
%
%
%
\subsection{Definitions connected with sets $I_f$}\label{edefr}
The definitions of this subsection are important for the sequel of this note.
Let $h=p^{e_p}\times h_2,\quad h_2\not\equiv 0 \modu p$, be the class number of $\Q(\zeta)$. Let
$f=min(e_p+2,\frac{p-1}{2})$.
Consider the set $I_f=\{ g_i \quad | \quad i=1,\dots,f\}$  defined in theorem \ref{t2009} p.\pageref{t2009}.
There exists  a {\bf minimum  rank} $r(I_f)\in\N,\quad r(I_f)\leq f$, or more shortly $r$ for $r(I_f)$ if context allows it,  which verifies simultaneously  the relations :
\begin{equation}\label{eqdefr}
\begin{split}
& r(I_f)\leq f= \mbox{ min }(e_p+2, \frac{p-1}{2}),\\
& J_r =\{ g_i\quad| \quad i=1,\dots,r\},\quad J_r\subset
I_f=\{g_i\quad |\quad i=1,\dots,f\}, \\
& J_r^\prime=\{p-g_i \quad|\quad i=1,\dots,r\}\not=J_r, \\
& N(J_r)=\{n_1,\dots,n_r\},\quad n_1,\dots,n_r \mbox{ not all } 0,\quad n_r=1,\\
& \sum_{i=1}^r \frac{g_{i} n_i}{(x+\zeta^{g_{i}} y)}
+\sum_{i=1}^r \frac{g_{i} n_i}{(x+\zeta^{-g_{i}} y)}\equiv 0 \modu p,
\end{split}
\end{equation}
these relations being impossible for any strict subset $J_{r^\prime}\subset J_r,
\quad r^\prime < r$.
Observe that
\begin{itemize}
\item
We have $r(I_f)>1$ because
$\frac{g_1n_1}{(x+\zeta^{g_1}y)}+\frac{g_1 n_1}{(x+\zeta^{-g_1} y)}\not\equiv 0\modu p$.
\item
For this value of $r=r(I_f)$, recall the notation
\begin{displaymath}
S_{2m+1}(J_r)\equiv \sum_{i=1}^r g_i^{2m+1}n_i\modu p.
\end{displaymath}
\end{itemize}
%
\subsection{Definitions connected with sets $U_f$.}\label{s31101}
Let $u$ be a primitive root $\modu p$; recall that, for $0\leq i\leq p-1$, we note indifferently
\begin{displaymath}
u_i\equiv u^i\modu p, \quad 1\leq u_i\leq p-1.
\end{displaymath}
\begin{lem}\label{pa210}
Let $u\in\N,\quad 1 < u < p-1,$ be a  primitive root $\modu p$.
Then,
\begin{displaymath}
\begin{split}
& U_f=\{u, u_2 ,\dots, u_{f} \} \not =\\
& U_f^\prime=\{p-u, p-u_2 ,\dots p-u_{f} \}.
\end{split}
\end{displaymath}
and the set $U_f$ verifies  the conditions of definition of $I_f$ sets in subsection \ref{edefr} p. \pageref{edefr}.
\begin{proof}

We shall prove that the set $U_f$ verifies the hypothesis assumed for sets $I_f$:
$u$ is a primitive root $\modu p$, so  we have $u^{(p-1)/2}+1\equiv 0 \modu p$.
Let $d\in \N$ with $d \quad|\quad \frac{(p-1)}{2},\quad d\not=\frac{(p-1)}{2}$.
We have $u^d\not\equiv \pm 1 \modu p$ : if not, we would have

$u^{2d}\equiv 1 \modu p$, which would contradict the fact that $u$ is a primitive root $\modu p$.
\begin{itemize}
\item
For $1\leq k<k^\prime \leq f$, we have $u_k\not\equiv u_{k^\prime}\modu p$ : if not, we should have
$u^{k-k^\prime}\equiv 1\modu p, \quad 1 <k^\prime-k <p-1$, a contradiction.
\item Show that $U_f\not=U_f^\prime $ : Suppose that $U_f=U_f^\prime$ and search for a contradiction: for all $k$, there should exist $k^\prime$ such that
$u_k\equiv p-u_{k^\prime}\modu p$,
so $u^{k^\prime-k}\equiv -1\modu p$,
so $k^\prime-k\equiv 0 \modu \frac{p-1}{2}$,
hence $k^\prime=k+\frac{p-1}{2}$,
that contradicts $1\leq k<k^\prime \leq f\leq \frac{p-1}{2}$.
\end{itemize}
Therefore, the set $U_f=\{u,u_2,\dots,u_{f-1} \}$  verifies the relations
assumed for the sets $I_f$.
\end{proof}
\end{lem}
This definition of $U_f$ set is important  for  the sequel of this note.
%
 We have shown that the  definition for  the set $I_f$ given in relation(\ref{eqdefr}) p.\pageref{eqdefr} can be applied, {\it mutatis mutandis}, to the set $U_f$:

There exists  a minimum {\bf rank }$r\in\N,\quad r\leq f,\quad$ with
\begin{equation}\label{e31102}
\begin{split}
& r \leq f =  \mbox{ min }(e_p+2, \frac{p-1}{2})\\
& u \in \N,\quad \mbox { primitive root }\modu p\\
& J_r =\{ u_{i} \quad| \quad i=1,\dots,r\}, \\
& J_r^\prime=\{p-u_{i} \quad|\quad i=1,\dots,r\}\not=J_r \\
& N_r=\{n_1,\dots,n_r\},\quad n_1\not\equiv 0\modu p,\quad n_r=1\\
& \sum_{i=1}^r \frac{u_{i} n_i}{(x+\zeta^{u_{i}} y)}
+\sum_{i=1}^r \frac{u_{i} n_i}{(x+\zeta^{-u_{i}} y)}\equiv 0 \modu p,
\end{split}
\end{equation}
these relations being impossible for any $r^\prime< r$.
For this value $r$, let us denote $S_{2m+1}(J_r)$ or more brievely $S_{2m+1}(r)$ defined by
\begin{displaymath}
S_{2m+1}(J_r)\equiv \sum_{i=1}^r u_{i(2m+1)} n_i\modu p.
\end{displaymath}
The value $r$ is called the  $p$-rank of FLT congruences $\modu p$.
{\bf Remarks :}
\begin{itemize}
\item We can always suppose that $n_1\not\equiv 0\modu p$ and that $n_r=1$, if not $r$ would not be the minimal rank.
\item
It is clear that set  $U_r$ is a  $I_f$ set.
\item
These definitions play a {\bf crucial part} in the sequel of this chapter.
\item
Observe that the two next congruences are equivalent, we shall use indifferently one or other form, corresponding to a conjugation $\sigma$ which does not affect congruences
$\modu \pi^{p-1}$ :
\begin{displaymath}
\begin{split}
& \sum_{i=1}^r \frac{u_{i} n_i}{(x+\zeta^{u_{i}} y)}
+\sum_{i=1}^r \frac{u_{i} n_i}{(x+\zeta^{-u_{i}} y)}\equiv 0,\\
& \sum_{i=1}^r \frac{u_{i-1} n_i}{(x+\zeta^{u_{i-1}} y)}
+\sum_{i=1}^r \frac{u_{i-1} n_i}{(x+\zeta^{-u_{i-1}} y)}\equiv 0.
\end{split}
\end{displaymath}
\end{itemize}
%
\subsection{Connection between sum and product congruences $\modu \pi^{p-1}$ of sets $U_f$}\label{s19022a}
The next results connect product congruences and sum congruences of minimal rank.
We shall show in the two next lemmas,  that $f=r$ is the smallest value such that there exists simultaneously $n_1,\dots,n_f$ not all zero with
\begin{displaymath}
\begin{split}
&\sum_{i=1}^f \frac{u_i n_i}{(x+\zeta^{u_i}y)}
+\sum_{i=1}^f\frac{u_i n_i}{(x+\zeta^{u_i}y)}\equiv 0 \modu \pi^{p-1},\\
&\prod_{i=1}^f (x+\zeta^{u_i}y)^{n_i}-\prod_{i=1}^f(x+\zeta^{-u_i}y)^{n_i}
\equiv 0 \modu \pi^{p-1}.
\end{split}
\end{displaymath}

%
The next lemma is a purely technical result.
\begin{lem}\label{l09112}
For every set $U_f$, the congruence
\begin{displaymath}
\sum_{i=1}^{f} \frac{u_i n_i}{(x+\zeta^{u_i}y)}
+\sum_{i=1}^f \frac{u_i n_i}{(x+\zeta^{-u_i} y)}
\equiv 0 \modu p
\end{displaymath}
implies the congruence
\begin{displaymath}
\sum_{i=1}^{f} \frac{u_i n_i\zeta^{u_i}}{(x+\zeta^{u_i}y)}
+\sum_{i=1}^f \frac{u_i n_i\zeta^{-u_i}}{(x+\zeta^{-u_i} y)}
\equiv 0 \modu p.
\end{displaymath}
\begin{proof}
Note at first that, from hypothesis, we get $\sum_{i=1}^f u_i n_i\equiv 0\modu \pi$, hence $\sum_{i=1}^f u_i n_i\equiv 0\modu p$.
\begin{displaymath}
\sum_{i=1}^f \frac{u_i n_i}{(x+\zeta^{u_i}y)}
+\sum_{i=1}^f \frac{u_i n_i}{(x+\zeta^{-u_i} y)}
\equiv 0 \modu p
\end{displaymath}
leads to
\begin{displaymath}
\begin{split}
&\sum_{i=1}^f \frac{u_i n_ix}{(x+\zeta^{u_i}y)}
+\sum_{i=1}^f \frac{u_i n_ix}{(x+\zeta^{-u_i}y)}
\equiv 0 \modu p\\
&2\sum_{i=1}^f u_i n_i
-(\sum_{i=1}^f \frac{y u_i n_i\zeta^{u_i}}{(x+\zeta^{u_i}y)}
+\sum_{i=1}^f \frac{y u_i n_i\zeta^{-u_i}}{(x+\zeta^{-u_i}y)})
\equiv 0 \modu p,
\end{split}
\end{displaymath}
which achieves the proof.
\end{proof}
\end{lem}
%
The relation (\ref{e31102}) p.\pageref{e31102}  gives a minimal $r$ for the congruences with the sums $\sum_{i=1}^r\frac{u_i n_i}{(x+\zeta^{u_i}y)}$. The next lemma shows that this minimal $r$ is also for the product congruences with $\prod_{i=1}^{r} (x+\zeta^{u_i}y)$.
%
\begin{lem}\label{l28053a}
Let $r$ and $n_1,\dots,n_{r}\in\N$  defined in the relation (\ref{e31102}) p.\pageref{e31102} of section (\ref{s31101}) p.\pageref{s31101}. Then

$\prod_{i=1}^{r} (x+\zeta^{u_i}y)^{n_i}
-\prod_{i=1}^{r} (x+\zeta^{-u_i}y)^{n_i}\equiv 0\modu p.$
Moreover, there does not exist $t<r$ and $l_1,\dots,l_t\in \N$ such that
$\prod_{i=1}^{t} (x+\zeta^{u_i}y)^{l_i}-\prod_{i=1}^{t} (x+\zeta^{u_i}y)^{l_i}\equiv 0\modu p$.
\begin{proof} $ $
\begin{itemize}
\item
From relation (\ref{e31102}) p.\pageref{e31102}, we have
$\sum_{i=1}^{r} \frac{u_i n_i}{(x+\zeta^{u_i}y)}
+\sum_{i=1}^{r} \frac{u_i n_i}{(x+\zeta^{-u_i}y)}\equiv 0 \modu \pi^{p-1}$.
Note that $\prod_{i=1}^r (x+\zeta^{u_i}y)^{n_i}-\prod_{i=1}^r (x+\zeta^{-u_i}y)^{n_i}
\equiv 0 \modu \pi$.

Suppose that for $\nu \in \N,\quad  1\leq \nu \leq p-2,$ we have $\pi^\nu
\parallel (\prod_{i=1}^{r}(x+\zeta^{u_i} y)^{n_i}-\prod_{i=1}^{r}(x+\zeta^{-u_i}y)^{n_i})$ and search for a contradiction.
Then, by derivation, (with the substitution, from lemma \ref{p6} p.\pageref{p6}, of $\zeta$  by the indeterminate $X$,  we should deduce that
\begin{displaymath}
\pi^{\nu-1}\parallel
 (\prod_{i=1}^{r}(x+\zeta^{u_i}y)^{n_i}\times
\sum_{i=1}^{r} \frac{u_i n_i\zeta^{u_i}}{(x+\zeta^{u_i}y)}-
\prod_{i=1}^{r}(x+\zeta^{-u_i}y)^{n_i}\times
\sum_{i=1}^{r} \frac{-u_i n_i\zeta^{-u_i}}{(x+\zeta^{-u_i}y)}),
\end{displaymath}
which contradicts the two simultaneous relations
\begin{displaymath}
\begin{split}
&\pi^\nu\|(\prod_{i=1}^{r}(x+\zeta^{u_i}y)^{n_i}-\prod_{i=1}^{r} (x+\zeta^{-u_i}y)^{n_i})\\
&\sum_{i=1}^{r} \frac{u_i n_i \zeta^{u_i}}{(x+\zeta^{u_i}y)}
-\sum_{i=1}^{r} \frac{-u_i n_i\zeta^{-u_i}}{(x+\zeta^{-u_i}y)}\equiv 0 \modu \pi^{p-1},
\end{split}
\end{displaymath}
the second one obtained from lemma \ref{l09112} p.\pageref{l09112}.
Therefore, we have $\nu\geq p-1$.
\item Suppose that for $t<r$, there would exist $l_1,\dots,l_t$ with
$\prod_{i=1}^t (x+\zeta^{u_i}y)^{l_i}-\prod_{i=1}^t (x+\zeta^{-u_i}y)^{l_i}\equiv 0 \modu p$.
Then, we would get
$\sum_{i=1}^{t} \frac{u_i l_i}{(x+\zeta^{u_i}y)}
+\sum_{i=1}^{t}\frac{u_i l_i}{(x+\zeta^{-u_i}y)}\equiv 0 \modu p$,  which would contradict the minimality of $r$.
\end{itemize}
\end{proof}
\end{lem}
To summarize, at the end of this section,we thus assert that, for the minimal $r$ of relation (\ref{e31102}) p.\pageref{e31102}, simultaneously:
\begin{displaymath}
\begin{split}
&\mbox{There exists} \quad n_1,\dots,n_r\\
&\sum_{i=1}^r\frac{u_i n_i}{(x+\zeta^{u_i}y)}
+\sum_{i=1}^r\frac{u_i n_i}{(x+\zeta^{-u_i}y)}\equiv 0\modu \pi^{p-1}\\
&\prod_{i=1}^{r} (x+\zeta^{u^i}y)^{n_i}
-\prod_{i=1}^{r} (x+\zeta^{-u^i}y)^{n_i}\equiv 0\modu \pi^{p-1}.
\end{split}
\end{displaymath}
%
\subsection{A canonical form of congruences $\modu p$ for FLT first case-I}
The Mirimanoff  congruences assert that

$\phi_{2m+1}(t) \times B_{p-1-2m}\equiv 0 \modu p, \quad m=1,\dots \frac{p-3}{2}$,
where $B_{p-1-2m}$ are even Bernoulli numbers.
But they do not give any information on the possibility to simultaneously have
$\phi_{2m+1}(t)\equiv 0 \modu p$ and $B_{p-1-2m}\equiv 0 \modu p$.
In the congruences
$\phi_{2m+1}(t)\times S_{2m+1}(J_r)\equiv 0 \modu p,\quad m=1,\dots,\frac{p-3}{2}$,
that we obtain, we shall show that {\bf else} $\phi_{2m+1}(t)\equiv 0 \modu p$, {\bf else} $S_{2m+1}(J_r)\equiv 0 \modu p$.

Recall that, for the indeterminate $T$, we define $\phi_{2m+1}^*(T)$ by $\phi_{1}^*(T)=1$ for $m=0$ and $\phi_{2m+1}^*(T)=\phi_{2m+1}(T)$ for $m=1,\dots,\frac{p-3}{2}$.
%
\begin{thm}{ *** }\label{pa212}
Let $t\equiv -\frac{y}{x}\modu p$. Let $r$ defined in the relation (\ref{e31102}) p.\pageref{e31102}.
If the First Case of  Fermat's Last Theorem fails for $p$,  then there exist exactly
$\frac{(p-1)}{2}-r+1$ different values
$m,\quad m\in\N,\quad 1\leq m \leq \frac{p-3}{2}$, such that we have the Mirimanoff polynomial congruences:
\begin{displaymath}\label{ea570}
\phi_{2m+1}(t)\equiv 0 \modu p.
\end{displaymath}
Moreover, $r$ does not depend on the primitive root $u \modu p$ taken to define $U_f$.
\begin{proof}
We have $J_r=\{u,u_2,\dots,u_{r} \}$ where $u$ is a primitive root
$\modu p$.
\begin{itemize}
\item
Show at first that there are at least $\frac{p-1}{2}-r+1$ different values
$m,\quad 0\leq m \leq \frac{p-3}{2}$, such that $\phi_{2m+1}^*(t)\equiv 0 \modu p$, where the notation $\phi_{2m+1}^*(T)$ has been previously defined. Suppose that this assertion is false and search for a contradiction :
then there exists at most $\frac{p-1}{2}-r$ different values such that
$\phi_{2m+1}^*(t)\equiv 0 \modu p,\quad 0\leq m\leq\frac{p-3}{2}$.
Then there exists at least $r$ values
$m_1=0,m_2\not=0,\dots,m_r\not=0$ such that
\begin{displaymath}\label{ea571}
\begin{split}
&\phi_{2m_1+1}^*(t)=1\not\equiv 0 \modu p,\\
&\phi_{2m_2+1}^*(t)\not\equiv 0 \modu p,\\
&\vdots\\
&\phi_{2m_r+1}^*(t)\not\equiv 0 \modu p.
\end{split}
\end{displaymath}
From theorem \ref{p21a} p.\pageref{p21a} and definitions (\ref{e31102}) p.\pageref{e31102}, {\it mutatis mutandis},
\begin{displaymath}
\begin{split}
& f\rightarrow r, \\
& I_f\rightarrow J_r,\\
& L(I_f)\rightarrow N(J_r),
\end{split}
\end{displaymath}
then we get
\begin{displaymath}\label{ea572}
\begin{split}
&S_{2m_1+1}(J_r)\equiv 0 \modu p,\\
&S_{2m_2+1}(J_r)\equiv 0 \modu p,\\
&\vdots\\
&S_{2m_r+1}(J_r)\equiv 0\modu p,
\end{split}
\end{displaymath}
hence
\begin{displaymath}
\begin{split}
&n_1 + n_2 u^{2m_1+1}+\dots+n_r u^{(r-1)(2m_1+1)}\equiv 0 \modu p,\\
&n_1 + n_2 u^{2m_2+1}+\dots+n_r u^{(r-1)(2m_2+1)}\equiv 0 \modu p,\\
&\vdots\\
&n_1 + n_2 u^{2m_r+1}+\dots+n_r u^{(r-1)(2m_r+1)}\equiv 0 \modu p.
\end{split}
\end{displaymath}
The determinant of this system of $r$ congruences in the $r$ indeterminates
$n_1, n_2,\dots, n_r$,
\begin{equation}\label{ea575}
\begin{array}{|llll|}
1 & u^{2m_1+1} & \dots & u^{(r-1)(2m_1+1)} \\
1 & u^{2m_2+1} & \dots & u^{(r-1)(2m_2+1)} \\
\vdots&\vdots &&\vdots \\
1 & u^{2m_r+1} & \dots & u^{(r-1)(2m_r+1)} \\
\end{array}
\end{equation}
is a Vandermonde determinant with factors $(u^{2m_k+1}-u^{2m_l+1})\not\equiv 0 \modu p$ because $u$ is a primitive root $\modu p$.
Therefore $n_1= n_2= \dots= n_r = 0$, which contradicts hypothesis on $N(J_r)$.
\item
Show secondly that there are at most $\frac{p-1}{2}-r+1$ different values such that

$\phi_{2m+1}^*(t)\equiv 0 \modu p,\quad 0\leq m\leq \frac{p-3}{2}$ : Suppose that this assertion is false and search for a contradiction :
Then there exist at least $\frac{p-1}{2}-r+2$ different values
$m, \quad 0 \leq m\leq\frac{p-3}{2}$, such that $\phi_{2m+1}^*(t)\equiv 0 \modu p$.
Then there exists  at most $\nu\leq r-2$ values
$m_1=0,m_2\not=0,\dots,m_{\nu}\not=0$, such that
\begin{displaymath}\label{ea571}
\begin{split}
&\phi_{2m_1+1}^*(t)=1\not\equiv 0 \modu p\\
&\vdots\\
&\phi_{2m_{\nu}+1}^*(t)\not\equiv 0 \modu p.
\end{split}
\end{displaymath}
It is possible to find $l_1,l_2,\dots, l_\nu, l_{\nu+1}=1\in\N$, not all zero such that
\begin{displaymath}
\begin{split}
&l_1+l_2 u^{2m_1+1}+\dots+l_\nu u^{(\nu-1)(2m_1+1)}
+l_{\nu+1} u^{\nu(2m_1+1)}\equiv 0\modu p,\\
&l_1+l_2 u^{2m_2+1}+\dots+l_\nu u^{(\nu-1)(2m_2+1)}
+l_{\nu+1} u^{\nu(2m_2+1)}\equiv 0\modu p,\\
&\vdots\\
&l_1 + l_2 u^{2m_{\nu}+1}+\dots+l_\nu u^{(\nu-1)(2m_{\nu}+1}
+l_{\nu+1} u^{\nu(2m_{\nu}+1)}\equiv 0 \modu p.
\end{split}
\end{displaymath}
because the determinant
\begin{equation}\label{ea575bis}
\begin{array}{|llll|}
1 & u^{2m_1+1} & \dots & u^{(\nu-1)(2m_1+1)} \\
1 & u^{2m_2+1} & \dots & u^{(\nu-1)(2m_2+1)} \\
\vdots&\vdots &&\vdots \\
1 & u^{2m_{\nu}+1} & \dots & u^{(\nu-1)(2m_{\nu}+1)} \\
\end{array}
\end{equation}
is a Vandermonde determinant with factors
$(u^{2m_k+1}-u^{2m_l+1})\not\equiv 0 \modu p$,
because $u$ is a primitive root $\modu p$.
Therefore we have,
\begin{displaymath}
\begin{split}
& J_{\nu+1}=\{u,u^2,\dots,u^{\nu+1}\}\\
& L(J_{\nu+1})=\{l_1,l_2,\dots,l_{\nu+1}\}\\
& S_{2m_1+1}(J_{\nu+1})\equiv 0\modu p\\
& S_{2m_2+1}(J_{\nu+1})\equiv 0\modu p\\
&\vdots\\
& S_{2m_\nu+1}(J_{\nu+1})\equiv 0 \modu p
\end{split}
\end{displaymath}
Therefore, for $m\in M = \{m_1=0,m_2,\dots,m_\nu\}$,  we have $S_{2m+1}(J_{\nu+1})\equiv 0 \modu p$ and for $m\in\{0,1,\dots,\frac{p-3}{2}\}-M$, we have
$\phi_{2m+1}^*(t)\equiv 0 \modu p$ and hence
\begin{displaymath}\label{ea575t}
S_{2m+1}(J_{\nu+1})\times \phi_{2m+1}^*(t)\equiv 0\modu p,
\quad m=0,1,\dots,\frac{p-3}{2}.
\end{displaymath}
From lemma \ref{propaa} p.\pageref{propaa}, this leads to the relation
\begin{displaymath}
 \sum_{i=1}^{\nu+1} \frac{u_i n_i}{(x+\zeta^{u_i} y)}
+\sum_{i=1}^{\nu+1} \frac{u_i n_i}{(x+\zeta^{-u_i} y)}\equiv 0 \modu p,
\end{displaymath}
which contradicts the minimal value of $r=r(U_f)>\nu+1$.
\item From first part of the proof, $r$ depends only of the number of Mirimanoff polynomials $\phi_{2m+1}(T),\quad 0\leq m\leq \frac{p-3}{2}$, with $\phi_{2m+1}(t)\equiv 0 \modu p$ and does not depend on the primitive root $u \modu p$ taken to define $U_f$.
\item From $\phi_1^*(t)=1$ and $\phi_{2m+1}^*(t)=\phi_{2m+1}(t)$ for $m=1,\dots,\frac{p-3}{2}$, we deduce that the number of $m, \quad 1\leq m \leq \frac{p-3}{2}$ with
$\phi_{2m+1}(t)\equiv 0 \modu p$ is the number of $m, \quad 0\leq m\leq \frac{p-3}{2}$ with $\phi_{2m+1}^*(t)\equiv 0 \modu p$, which completes the proof.

\end{itemize}
\end{proof}
\end{thm}
%
\begin{cor}{ *** }\label{c02101}
Let $h=p^{e_p}\times h_2,\quad h_2\not\equiv 0 \modu p$, be the class number of $\Q(\zeta)$.
Suppose that first case of FLT fails for $p$.
Let $t,\quad t\in\N,\quad t\equiv -\frac{x}{y}\modu p$.
If $e_p <\frac{p-1}{2}$, then there are at least $\frac{p-3}{2}-e_p$ Mirimanoff polynomials of odd indice $\phi_{2m+1}(T),\quad 1\leq m \leq \frac{p-3}{2}$, with
\begin{displaymath}
\phi_{2m+1}(t)\equiv 0 \modu p.
\end{displaymath}
\begin{proof}
From relation (\ref{eqdefr}) p.\pageref{eqdefr}, we have $r\leq e_p+2$. The result is then an immediate consequence of  theorem \ref{pa212} p.\pageref{pa212}, where the number of $\phi_{2m+1}(t)\equiv 0$ is $\frac{p-1}{2}-r+1$.
\end{proof}
\end{cor}
%

\subsection{A comparison with Mirimanoff, Herbrand and Ribet/Kolyvagin}\label{s25101}
$ $
\begin{itemize}
\item
Recall and fix some notations:
\begin{itemize}
\item
Recall that the class number $h$ of $\Q(\zeta)$ is $h=p^{e_p}\times h_2,\quad h_2\not\equiv 0 \modu p$.
\item
Let $r_p$ be the $p$-rank of the class group of  $\Q(\zeta)$. Let $r_p^+$ be the $p$-rank of the class group of
$\Q(\zeta+\zeta^{-1})$. Let $r_p^-=r_p-r_p^+$ be the relative $p$-rank.
\item
Let, as previously, $t\equiv -\frac{x}{y}\modu p$.
\item
Let $G=Gal(\Q(\zeta_p)/\Q)\simeq (\Z/p\Z)^*$.
\item
Let  $A$ be the $p-$Sylow subgroup of the ideal class group of $\Q(\zeta)$. So $A$ is a $\Z_p[G]$-module for $G$ acting on $A$.
Let $A=\oplus_{i=0}^{p-2} A_i$ be the decomposition of $A$ as $\Z_p[G]-$module.
(see for instance \cite {was})  6.3 p 100).
\item
Let $E_A=\{A_{2m+1} \quad|\quad 1\leq m\leq \frac{p-3}{2}, \quad A_{2m+1}\not=0\}$.
\item
Let $i_A=Card(E_A)$.
\item
Let $B_{p-2m-1}, \quad 1\leq m\leq \frac{p-3}{2}$, be even Bernoulli numbers.
\item
Let $E_B=\{ B_{p-2m-1} \quad|\quad 1\leq m\leq \frac{p-3}{2},
\quad B_{p-2m-1}\equiv 0 \modu p\}$.
\item
We note $i_p=Card(E_B)$ the index of irregularity of $p$.
\item
Let $\phi_{2m+1}(T), \quad  1\leq m\leq \frac{p-3}{2}$, be the odd Mirimanoff polynomials.
\item
Let $E_\phi(t)=\{ \phi_{2m+1}(T)\quad |\quad 1\leq m\leq \frac{p-3}{2}, \quad \phi_{2m+1}(t)\equiv 0 \modu p\}$.
\item
We note $\phi_c(t)=Card(E_{\phi(t)})$.
\end{itemize}
\item
We compare our approach to known results:
\begin {itemize}
\item
Mirimanoff congruences :
Let $t\equiv -\frac{x}{y}\modu p$.
We have $\phi_{2m+1}(t)\times B_{p-2m-1}\equiv 0 \modu p,
\quad 1\leq m \leq \frac{p-3}{2}$:
see for instance \cite{rib} (1B) p 145, which implies $\phi_c(t)\geq \frac{p-3}{2}-i_p$.
\item
Trivially $r_p\leq e_p$.
\item Herbrand theorem : $i_A\leq i_p$  (see for instance \cite{was} thm 6.17 p 101). Note that Herbrand theorem does not allow to conclude on $\phi_c(t)$, because if $B_{p-2m-1}\equiv 0 \modu p$, the Mirimanoff congruence does not allow to determine $\phi_{2m+1}(t)\modu p$.
\item
Ribet/Kolyvagin converse of Herbrand theorem : $i_A=i_p$ and  also $i_p\leq r_p$ (see for instance \cite{was} thm 6.18 p 102). Then these results allow to obtain $i_p\leq r_p\leq e_p$.
\item
Then, from Mirimanoff's congruences and Ribet theorem simultaneously, we get
\begin{displaymath}
\phi_c(t)\geq \frac{p-3}{2}-i_p\geq \frac{p-3}{2}-r_p\geq \frac{p-3}{2}-e_p.
\end{displaymath}
which is the result that we have found in corollary \ref{c02101} p.\pageref{c02101} by an {\bf elementary method}.
\item
This result can also be obtained from corollary 10.15 p 198 of \cite{was} and a Mazur-Wiles  theorem $|\varepsilon_i|=p-$part of $B_{1,\omega^{-i}}$, see \cite{was} p 300 and p 348.
\item
It is important to note that our proof uses only strictly elementary properties of Dedekind algebraic number field $\Q(\zeta)$; Proofs of Herbrand and  Ribet/Kolyvagin converse of Herbrand theorem and Mazur-Wiles theorem uses
$p-$adic theory, Hilbert class field theory, Modular forms, Euler systems, ... totally beyond material used in this paper.
\item
Note that our  result gives  $\phi_c(t)\geq \frac{p-3}{2}-e_p$. It is even possible to prove the stronger result $\phi_c(t)\geq \frac{p-3}{2}-r_p$ obtained from Mirimanoff/Herbrand/Ribet/Kolyvagin and even $\phi_c(t)\geq \frac{p-3}{2}-r_p^-$, always in our {\bf strictly elementary approach},  (stronger  because it is possible that $r_p^-<r_p<e_p$), see corollary \ref{c09113} p.\pageref{c09113}.
\item
It is also possible to prove the stronger result $\phi_c(t)\geq \frac{p-3}{2}-r_p^-$ see for instance theorem \ref{t011182} p. \pageref{t011182}, which is better if $h^+\equiv 0\modu p$.
\end{itemize}
\end{itemize}
%
\subsection{A canonical form of congruences $\modu p$ for FLT first case-II}
\label{s108251}
In this subsection, we give some other structure theorems on congruences $\modu \pi^{p-1}$ of Fermat's equation.

\begin{thm} { *** }\label{t12092}
If the first case of FLT fails for $p$,  then there exists exactly $\mu=r-1$ different values
$m,\quad m\in\N,\quad 0\leq m\leq \frac{p-3}{2}$, such that
\begin{displaymath}
S_{2m+1}(J_r)\equiv 0\modu p.
\end{displaymath}
\begin{proof}
Consider always the set $U_f=\{u,\dots,u_{f}\}$ with $u$ primitive root $\modu p$ and the minimal set $J_r =\{u,\dots,u_{r}\}$. Let $\mu$ be the number of $m,\quad m\in\N,\quad 0\leq m \leq \frac{p-3}{2}$, such that $S_{2m+1}(J_r)\equiv 0 \modu p$.
\begin{itemize}
\item
If $\mu<r-1$ then there are more than $\frac{p-1}{2}-r+1$ Mirimanoff polynomials with $\phi_{2m+1}^*(t)\equiv 0 \modu p,\quad 0\leq m \leq \frac{p-3}{2}$, which contradicts previous theorem  \ref{pa212} p.\pageref{pa212}.
\item
If $\mu>r-1$ then, from the first $r$ congruences
\begin{displaymath}
\begin{split}
& S_{2m_1+1}(J_r)\equiv 0 \modu p,\\
& S_{2m_2+1}(J_r)\equiv 0 \modu p,\\
&\vdots\\
& S_{2m_r+1}(J_r)\equiv 0 \modu p,
\end{split}
\end{displaymath}
we get
\begin{displaymath}
\begin{split}
&n_1 + n_2 u^{2m_1+1}+\dots+n_r u^{(r-1)(2m_1+1)}\equiv 0 \modu p,\\
&n_1 + n_2 u^{2m_2+1}+\dots+n_r u^{(r-1)(2m_2+1)}\equiv 0 \modu p,\\
&\vdots\\
&n_1 + n_2 u^{2m_r+1}+\dots+n_r u^{(r-1)(2m_r+1)}\equiv 0 \modu p.
\end{split}
\end{displaymath}

The determinant of this system of $r$ congruences in the $r$ indeterminates
$n_1, n_2,\dots, n_r$ is
\begin{displaymath}\label{ea575}
\begin{array}{|llll|}
1 & u^{2m_1+1} & \dots & u^{(r-1)(2m_1+1)} \\
1 & u^{2m_2+1} & \dots & u^{(r-1)(2m_2+1)} \\
\vdots&\vdots &&\vdots \\
1 & u^{2m_r+1} & \dots & u^{(r-1)(2m_r+1)} \\
\end{array}
\end{displaymath}
a Vandermonde determinant with factors $(u^{2m_k+1}-u^{2m_l+1})\not\equiv 0 \modu p$ because $u$ is a primitive root $\modu p$.
Therefore $n_1\equiv n_2\equiv \dots\equiv n_r \equiv 0 \modu p$, which contradicts hypothesis and achieves the proof.
\end{itemize}
\end{proof}
\end{thm}
%
We use in the sequel that value $\mu=r-1$.
Therefore, with  $n_r=1$ assumed without loss of generality, we are in the situation where we have a system of $r-1$ linear congruences in the $r-1$
indeterminates $n_1,\dots,n_{r-1}$, not all zero, with no null second member,
\begin{equation}\label{e14101}
\begin{split}
&n_1 + n_2 u^{2m_1+1}+\dots+n_{r-1} u^{(r-2)(2m_1+1)}\equiv
- u^{(r-1)(2m_1+1)} \modu p\\
&n_1 + n_2 u^{2m_2+1}+\dots+n_{r-1} u^{(r-2)(2m_2+1)}\equiv
- u^{(r-1)(2m_2+1)} \modu p\\
&\vdots\\
&n_1 + n_2 u^{2m_{r-1}+1}+\dots+n_{r-1} u^{(r-2)(2m_{r-1}+1)}\equiv
- u^{(r-1)(2m_{r-1}+1)} \modu p.
\end{split}
\end{equation}
with determinant
\begin{equation}\label{e14102}
\Delta=
\begin{array}{|llll|}
1 & u^{2m_1+1} & \dots & u^{(r-2)(2m_1+1)} \\
1 & u^{2m_2+1} & \dots & u^{(r-2)(2m_2+1)} \\
\vdots&\vdots &&\vdots \\
1 & u^{2m_{r-1}+1} & \dots & u^{(r-2)(2m_{r-1}+1)} \\
\end{array}
\not\equiv 0 \modu p.
\end{equation}
This result shall be used in section \ref{s14101} p.\pageref{s14101}.
%
Recall that the definitions used in the next corollary are in subsection \ref{s31101} p.\pageref{s31101}.
%
\begin{cor}{ *** }\label{cor1209}
If the first case of  Fermat's Last Theorem fails for $p$, then there exists $r\in \N$ depending only on $p$ and $t$ such that, for any primitive root $u\modu p$,
$S_1(J_r)\equiv 0 \modu p$ and
for $m=1,\dots,\frac{p-3}{2}$:
\begin{itemize}
\item
$\phi_{2m+1}(t)\times S_{2m+1}(J_r) \equiv 0 \modu p$,
\item
If $\phi_{2m+1}(t)\equiv 0\modu p$ then $S_{2m+1}(J_r)\not\equiv 0 \modu p$,
\item
If $S_{2m+1}(J_r)\equiv 0\modu p$ then $\phi_{2m+1}(t)\not\equiv 0 \modu p$.
\end{itemize}
\begin{proof}
Immediate consequence of the two previous theorems and of relation (\ref{eqdefr}) p.\pageref{eqdefr}  where $r(U_f)$ is defined.
\end{proof}
\end{cor}
%
In the next theorem, we shall show that it is possible to compute explicitly $S_{2m+1}(J_r)$
and we shall give an explicit  solution $\{l_1,\dots,l_r\}$
of the congruence $\modu p$
\begin{displaymath}
\sum_{i=1}^r\frac{u_i l_i}{x+\zeta^{u_i} y}+\sum_{i=1}^r \frac{u_i l_i}{x+\zeta^{-u_i} l_i}\equiv 0 \modu p.
\end{displaymath}
%
\begin{thm}{ *** } { {\bf An explicit computation of $S_{2m+1}(J_r)$.}}\label{t23011}

Let the set

$I_\phi=\{2m_i+1 \quad |\quad i=1,\dots,r-1,\quad \phi^*_{2m_i+1}(t)\not\equiv 0 \modu p,\quad 0\leq m_i\leq \frac{p-3}{2}\}$.
\begin{enumerate}
\item
Then an  explicit formula for $S_{2n+1}(J_r)$ is :
\begin{equation}\label{e23014}
S_{2n+1}(J_r)=u_{2n+1}\times\prod_{i=1}^{r-1} (u_{2n+1}-u_{2m_i+1}).
\end{equation}
\item
We have the equivalence
\begin{displaymath}
\begin{split}
& \sum_{i=1}^r\frac{u_i l_i}{x+\zeta^{u_i}y}
+\sum_{i=1}^r\frac{u_i l_i}{x+\zeta^{-u_i} y}\equiv 0 \modu \pi^{p-1} \\
& \Longleftrightarrow\\
& S_{2n+1}(J_r)=u_{2n+1}\times\prod_{i=1}^{r-1} (u_{2n+1}-u_{2m_i+1}) \\
& \Longleftrightarrow\\
& l_1 \equiv (-1)^{r-1}\times \prod_{i=1}^{r-1} u_{2m_i+1} \modu p,\\
&\vdots \\
& l_{r-1}\equiv (-1)\times \sum_{i=1}^{r-1} u_{2m_i+1} \modu p, \\
& l_r=1.
\end{split}
\end{displaymath}
\end{enumerate}
\begin{proof} $ $
\begin{enumerate}
\item
Let $P_{2n+1}=\prod_{i=1}^{r-1} (u_{2n+1}-u_{2m_i+1})$.
With a direct computation, we get
\begin{displaymath}
\begin{split}
& P_{2n+1}=\sum_{i=1}^r l_i\times u_{(i-1)(2n+1)},\\
& l_1\equiv (-1)^{r-1}\times\prod_{i=1}^{r-1} u_{2m_i+1} \modu p,\\
& \vdots\\
& l_{r-1}\equiv  (-1)\times \sum_{i=1}^{r-1} u_{2m_i+1}\modu p,\\
& l_{r}=1.
\end{split}
\end{displaymath}
We have $S_{2n+1}(J_r)=\sum_{i=1}^{r} n_i\times u_{i(2n+1)}$, where, without loss of generality, we can assume $n_{r}=1$.
Then, suppose that $S_{2n+1}(J_r)\not=u_{2n+1}\times P_{2n+1}$ and search for a contradiction:
Let us consider the expression  $S_{2n+1}^\prime =S_{2n+1}(J_r)-u_{2n+1}\times P_{2n+1}$.
We have seen, theorem  \ref{p21a} p.\pageref{p21a}, that
\begin{displaymath}
\phi_{2n+1}^*(t)\times S_{2n+1}(J_r)\equiv 0 \modu p,\quad n=0,\dots,\frac{p-3}{2}.
\end{displaymath}
The definition of $P_{2n+1}$ implies that
\begin{displaymath}
\phi^*_{2n+1}(t)\times u_{2n+1}\times P_{2n+1}\equiv 0\modu p,\quad n=0,\dots,\frac{p-3}{2}.
\end{displaymath}
Therefore
\begin{displaymath}
\phi^*_{2n+1}(t)\times S_{2n+1}^\prime\equiv 0\modu p,\quad n=0,\dots,\frac{p-3}{2}.
\end{displaymath}
and so
\begin{displaymath}
\phi^*_{2n+1}(t)\times \sum_{i=1}^{r-1} (n_i-l_i)\times u_{i(2n+1)}\equiv 0\modu p,
\quad n=0,\dots,\frac{p-3}{2}.
\end{displaymath}
From lemma \ref{propaa} p.\pageref{propaa} we get
\begin{displaymath}
\sum_{i=1}^{r-1} \frac{u_i(n_i-l_i)}{(x+\zeta^{u_i} y)}
+\sum_{i=1}^{r-1}\frac{u_i(n_i-l_i)}{(x+\zeta^{-u_i} y)}
\equiv 0 \modu p,
\end{displaymath}
which contradicts minimality of $r$.
\item
Immediate consequence of 1)
\end{enumerate}
\end{proof}
\end{thm}
%
\subsection{A comparison with Mirimanoff}\label{s31103}
Let $t\equiv -\frac{x}{y}\modu p$.
We have for $m=1,2,\dots,\frac{p-3}{2}$ the two congruences $\modu p$ :
\begin{itemize}
\item
$\phi_{2m+1}(t)\times B_{p-1-2m}\equiv 0 \modu p$, see \cite{rib} p145.
\item
$\phi_{2m+1}(t)\times S_{2m+1}(J_r)\equiv 0 \modu $, corollary \ref{cor1209} p.\pageref{cor1209}, where
$S_{2m+1}(V_r)=\sum_{i=1}^r u_{i (2m+1)}n_i\modu p$ is, for $u$ primitive root $\modu p$ given,  explicitly computable.
\end{itemize}
Note that, in our result, it is not possible to have simultaneously $\phi_{2m+1}(t)\equiv 0\modu p $ and
$S_{2m+1}(J_r)\equiv 0 \modu p$. In opposite, we have not seen in the literature if it is possible or no to have simultaneously   $\phi_{2m+1}(t)\equiv 0 \modu p$ and $B_{p-1-2m}\equiv 0 \modu p$. The section \ref{s14101} p.\pageref{s14101} of this article gives an example of importance of the property: the two congruences
$\phi_{2m+1}(t)\equiv 0\modu p $ and
$S_{2m+1}(J_r)\equiv 0 \modu p$ are not simultaneously possible.

Related to this subject, we have only found the beautiful results of F. Thaine, relating even Mirimanoff polynomials and even Bernoulli numbers:
\begin{itemize}
\item
Let $h^+$ be the class number of $\Q(\zeta+\zeta^{-1})$. If $h^+\not\equiv 0 \modu p$ (in particular case of Vandiver conjecture) and if $B_{p-1-2m}\equiv 0 \modu p,\quad 1\leq m \leq \frac{p-3}{2}$ and if the first case of Fermat's Last Theorem fails, then  $\phi_{p-1-2m}(t)\equiv 0 \modu p$.
see (\cite{th1}) prop 3 p 141.
This results allows a symmetric formulation : if $h^+\not\equiv 0 \modu p$ then
\begin{itemize}
\item
if $B_{p-1-2m}\equiv 0 \modu p$ then $\phi_{p-1-2m}(t)\equiv 0 \modu p$;
\item
if $B_{p-1-2m}\not\equiv 0 \modu p$ then $\phi_{2m+1}(t)\equiv 0 \modu p$.
\end{itemize}
Observe that, from Mirimanoff, we know already that $\phi_{2m+1}(t)\times\phi_{p-1-2m}(t)\equiv 0\modu p$, but Thaine's result is more precise.
\item
If $\prod_{i=1}^{p-1} (1-\zeta^i)^{i^{2m}}$ is not a $p$th power in $\Z[\zeta]$, if $B_{p-1-2m}\equiv 0 \modu p$ and if the first case of Fermat's Last Theorem fails,
then $\phi_{p-1-2m}\equiv 0 \modu p$, see (\cite{th2}) cor. p 299.
\end{itemize}

%
\subsection{A connection between Mirimanoff congruences and $p-$rank of the class group  of
$\Q(\zeta)/\Q$}
Let $t\in\N,\quad  t\equiv -\frac{x}{y}\modu p $.
Let $\phi_c(t)$ be the number of odd Mirimanoff polynomials
$\phi_{2m+1}(T),\quad m=1,\dots,\frac{p-3}{2}$, with $\phi_{2m+1}(t)\equiv 0 \modu p$.
Let $r_p$ be the $p-$rank of the class group of
$\Q(\zeta)$. We shall prove that $\phi_c(t)\geq \frac{p-3}{2} -r_p$.
%
\begin{thm} { *** }\label{t09111}
Let $r_p$ be the $p-$rank of the class group $C_p$ of $\Q(\zeta)$. Let $r$ be the $p$-rank of Fermat congruences  defined in the section (\ref{s31101}). Then $r\leq r_p+2$.
\begin{proof}
Suppose that $r> r_p+2$ and search for a contradiction :
Recall that $(x+\zeta y)\Z[\zeta]=\s^p$ where $\s$ is an integral ideal. Let $r_1$ be the $p$-rank of the $p$-subgroup $C_p$ of the class group of $\Q(\zeta)/\Q$ generated by action of $Gal(\Q(\zeta)/\Q)$ on the class group $<Cl(\sigma(\s)>$. We have $r_1\leq r_p$.
From elementary properties of the $p$-class group, it is possible to find $l_1,l_2,\dots,l_{r_1+1}$, such that $\prod_{i=1}^{r_1+1}(x+\zeta^{u_i}y)^{l_i}=\zeta^u\times \eta\times \gamma^p$ with $u\in\N$ and $\gamma\in \Z[\zeta]$. It is then possible to find $n_i,\quad i=1,2,\dots,r_1+2$, such that
\begin{displaymath}
\begin{split}
&\prod_{i=1}^{r_1+2} (x+\zeta^{u_i}y)^{n_i}=\eta_1\times \gamma_1^p,
\quad \eta_1\in\Z[\zeta+\zeta^{-1}]^*,\quad \gamma_1\in\Z[\zeta],\\
&\sum_{i=1}^{r_1+2} \frac{u_i l_i\zeta^{u_i}}{(x+\zeta^{u_i}y)}
+\sum_{i=1}^{r_1+2}\frac {u_i l_i\zeta^{-u_i}}{(x+\zeta^{-u_i}y)}\equiv 0 \modu p.
\end{split}
\end{displaymath}
Therefore, from the  minimality of $r$ in the definition \ref{s31101} p.\pageref{s31101}, we get
\begin{equation}\label{e012091}
r_1+2\geq r,
\end{equation}
and hence
\begin{equation}\label{e011182}
r_p+2\geq r_1+2\geq r.
\end{equation}
\end{proof}
\end{thm}
%
\begin{cor} { *** }\label{c09113}
Let $r_p$ be the $p-$rank of the class group of $\Q(\zeta)$.
Suppose that FLT first case fails for $p$.
Let $t\in \N,\quad t\equiv  -\frac{x}{y}\modu p$.
Let $\phi_{2m+1}(T),\quad m=1,\dots,\frac{p-3}{2}$, be the odd Mirimanoff polynomials.
Let $\phi_c(t)$ be the number of $\phi_{2m+1}(t)\equiv 0 \modu p,\quad m=1,\dots,\frac{p-3}{2}$.
Then $\phi_c(t)\geq \frac{p-3}{2}-r_p$.
\begin{proof}
From theorem \ref{pa212} p.\pageref{pa212}, we have $\phi_c(t)=\frac{p-1}{2}-r+1$. Then we apply theorem
\ref{t09111} p.\pageref{t09111}.
\end{proof}
\end{cor}
%
\clearpage
%
\section{FLT first case: Mirimanoff in intermediate fields $K$ with $\Q\subset K\subset\Q(\zeta)$.
Part I.}\label{s108252}
\label{s104284}
In this section we shall obtain some relations on Mirimanoff polynomials with    a systematic  study of FLT in  intermediate fields between $\Q$ and $\Q(\zeta)$.
\subsection{Some definitions}\label{s104091}
Let us recall or give some new notations used in this section.
\begin{itemize}
\item
$u$ is a primitive root $\modu p$; we note $u_i$ for $u^i\modu p$.
\item
$\sigma:\Q(\zeta)\rightarrow \Q(\zeta)$ is a $\Q$-isomorphism defined by $\sigma(\zeta)=\zeta^u$, where $<\sigma>=Gal(\Q(\zeta)/\Q)$.
\item
Let $g,\quad g\in \N,\quad g>1,\quad \frac{p-1}{2}\equiv 0 \modu g$.
\item
We have defined the odd Mirimanoff's polynomials by
\begin{displaymath}
\phi_{2m+1}(T)=T+2^{2m}\times T^2+\dots+(p-1)^{2m}\times T^{p-1},
\quad m=0,\dots,\frac{p-3}{2}.
\end{displaymath}
\item
We recall the definition of  $\phi_{2m+1}^*(T)$ in relation (\ref{e104051}) p.\pageref{e104051} by
\begin{displaymath}
\phi^*_{2m+1}(T)=\phi_{2m+1}(T),\quad m=1,\dots,\frac{p-3}{2},\quad
\phi_1^*(T)=1.
\end{displaymath}
\item
Let $M_g=\frac{p-1}{2 g}$. Therefore, from hypothesis on $g$,  $M_g\in\N$.
\item
Let us suppose, in this section,  that
\begin{equation}\label{e112156}
M_g=\frac{p-1}{2 g}\geq e_p+2.
\end{equation}
where $e_p$ is defined, from the class number $h$ of $\Q(\zeta)$,  by the relation $h=p^{e_p}\times h_2,\quad h_2\not\equiv 0\modu p$.
\item
Let $K$ be the field,
$\Q\subset K\subset\Q(\zeta),\quad [\Q(\zeta):K]=\frac{p-1}{g},\quad [K:\Q]=g$.
\item
Let us consider the $K$-isomorphism
$\sigma^g : \Q(\zeta) \rightarrow \Q(\zeta),\quad \zeta\rightarrow\zeta^{u_g}$.
\item
Let us consider the set of ideals $\sigma^{g i}(\s),\quad i=0,\dots, \frac{p-1}{g}-1$.
\item
Recall that, for the Fermat equation $x^p+y^p+z^p=0$, we set $t\equiv-\frac{y}{x}\modu p$.
\end{itemize}

%
\subsection{ Results on Mirimanoff in intermediate fields $K,\quad \Q\subset K\subset \Q(\zeta)$}
\label{s209282}
The proof of  theorem \ref{t2009} p.\pageref{t2009}, that we have explained, when
$\frac{p-1}{2}\geq f \geq e_p+2$,
for the $\Q$-isomorphism $\sigma :\Q(\zeta)\rightarrow \Q(\zeta)$ and the set of ideals $\sigma^i(\s)\quad i=0,\dots,p-2$, can be strictly applied, {\it mutatis mutandis}, when $\frac{p-1}{2 g}\geq f_g \geq e_p+2$, to the $K$-isomorphism
$\sigma^g :\Q(\zeta)\longrightarrow \Q(\zeta)$ and to the set of ideals $\sigma^{gi}(\s),\quad i=0,\dots,\frac{p-1}{g}-1$. Then the theorem \ref{t2009} p.\pageref{t2009} gives in that case:
%
\begin{thm}\label{t10101}
Let $h=p^{e_p} \times h_2,\quad h_2\not\equiv 0 \modu p$, be the class number of $\Q(\zeta)/\Q$.
Let $g\in \N,\quad g\geq 1,\quad \frac{p-1}{2}\equiv 0 \modu g$ with $\frac{p-1}{2 g}\geq e_p+2$.
Let  $f_g\in\N, \quad \frac{p-1}{2 g}\geq f_g\geq e_p+2$.
Then, there exists at least one set of $f_g$ natural numbers $L(I_{f_g})$ not all simultaneously null
\begin{equation}\label{e4bis}
L(I_{f_g})=\{\ l_i \ | \quad l_i \in\N,\quad 1\leq i\leq f_g,\quad 0\leq l_i\leq p-1\ \}
\end{equation}
such that simultaneously,
\begin{displaymath}\label{e5bis}
\begin{split}
&\prod_{i=1}^{f_g} \frac{(x+\zeta^{ u_{g i}}y)^{l_i}}{(x+\zeta^{- u_{g i}}y)^{l_i}}=\beta^p,
\quad\beta\in\Q(\zeta),\quad\beta\not=1,\\
&\prod_{i=1}^{f_g} (x+\zeta^{ u_{g i}} y)^{l_i}
- \prod_{i=1}^{f_g}(x+\zeta^{-u_{g i}}  y)^{l_i} \equiv 0 \modu \pi^{p+1}.
\end{split}
\end{displaymath}
\begin{proof}
The same proof as proposition \ref{p4} p.\pageref{p4} and theorem \ref{t2009} p.\pageref{t2009} can be applied here; we observe only that
$U_{f_g}=\{u^{gi}\quad | \quad i=0,\dots, f_g-1\}$ verifies the conditions assumed for $I_{f}$ sets:
\begin{displaymath}
\begin{split}
& u^{g i_1}\not\equiv u^{g i_2}\modu p , \quad i_1\not = i_2,\\
& \{u_{gi}\quad|\quad i=1,\dots,f_g-1\}
\not=\{p-u_{gi}\quad|\quad i=1,\dots,f_g-1\}.
\end{split}
\end{displaymath}
because $u$ is a primitive root $\modu p$ and
$g i_1 < g\times\frac{p-1}{2 g}=\frac{p-1}{2}$ and
$g i_2 < g\times\frac{p-1}{2 g}=\frac{p-1}{2}$.
\end{proof}
\end{thm}
{\bf Remark:} Observe that if $g=1$ it is the theorem \ref{t2009} p.\pageref{t2009}.
%
\subsubsection{Definitions connected with sets $U_{f_g}$.}

This definition is important  for  the sequel of this section. This definition is the translation for $\sigma^g$ of the definition given in relation (\ref{e31102}) p.\pageref{e31102}  for $\sigma$.
There exists  a minimum Fermat congruence  {\bf rank }
$r_g\in\N$ verifying
\begin{equation}\label{e104012}
\begin{split}
& r_g \leq   e_p+2\leq  \frac{p-1}{2g}),\\
& u \in \N,\quad \mbox { primitive root }\modu p,\\
& J_{r_g} =\{ u_{g i} \quad| \quad i=1,\dots,r_g\}, \\
& N(J_{r_g})=\{n_1,\dots,n_{r_g}\},\quad n_1\not\equiv 0\modu p ,\quad
n_{r_g}=1,\\
& \sum_{i=1}^{r_g} \frac{u_{g i} n_i}{(x+\zeta^{u_{g i}} y)}
+\sum_{i=1}^{r_g} \frac{u_{g i} n_i}{(x+\zeta^{-u_{g i}} y)}\equiv 0 \modu p,
\end{split}
\end{equation}
these relations being impossible for  $r_g^\prime <r_g$.

For this value $r_g$, let us denote
\begin{equation}\label{e112141}
S_{g(2m+1)}(J_{r_g})\equiv \sum_{i=1}^{r_g} u_{g i (2m+1)} \times n_i\modu p.
\end{equation}
%
\subsubsection{Instanciation of $\sigma$-theorem to $\sigma^g$-theorems}
\label{s108253}
The instanciation of {\it foundation} theorem (\ref{p21a}) becomes in that case:

%
\begin{thm}\label{t104012}
\begin{displaymath}
\begin{split}
& S_g(J_{r_g})\equiv 0\modu p,\\
&  S_{g(2m+1)}(J_{r_g})\times \phi_{2m+1}^*(t)\equiv 0 \modu p,
\quad m=0,\dots,\frac{p-3}{2}.
\end{split}
\end{displaymath}
\begin{proof}$ $
\begin{enumerate}
\item
Similar to lemma \ref{p7} p.\pageref{p7} with $\sigma^g$ in place of $\sigma$ and $u_g$ in place of $u$.
\item
Similar to the  proof of theorem \ref{p21a} p.\pageref{p21a} with $\sigma^g$ in place of $\sigma$ and $u_g$ in place of $u$.
\end{enumerate}
\end{proof}
\end{thm}
%
%
We show that theorems \ref{pa212} p.\pageref{pa212}  can be {\it partially} applied here, with {\it mutatis mutandis}, $f \rightarrow f_g$, $r\rightarrow r_g$.
%

\begin{thm}{ *** }\label{t104013}
Let $t\equiv -\frac{y}{x}\modu p$.
Let $r_g$ defined in the relation (\ref{e104012}) p \pageref{e104012}.
If the First Case of  Fermat's Last Theorem fails for $p$,  then
\begin{enumerate}
\item
there exists {\bf exactly}
$\frac{p-1}{2g}-r_g+1$ different values
$m,\quad m\in\N,\quad 1\leq m \leq \frac{p-1}{2g}-1$,
such that we have for $m$ the $g$  Mirimanoff polynomial congruences:
\begin{displaymath}\label{ea570}
\phi_{2m+1+\alpha(p-1)/g}^{*}(t)\equiv 0 \modu p,\quad \alpha=0,\dots,g-1.
\end{displaymath}
\item
there exists {\bf exactly}
$r_g-1$ different values
$m,\quad m\in\N,\quad 1\leq m \leq \frac{p-1}{2g}-1$,
such that we have for $m$ {\bf at least} one $\alpha\in\N, \quad 0\leq \alpha\leq g-1$ with the  Mirimanoff polynomial congruence:
\begin{displaymath}\label{ea570}
\phi_{2m+1+\alpha(p-1)/g}^{*}(t)\not\equiv 0 \modu p.
\end{displaymath}
\end{enumerate}
\begin{proof}$ $
We duplicate, {\it mutatis mutandis}, the proof of  theorem \ref{pa212} p.\pageref{pa212}.
\begin{enumerate}
\item
We have $J_{r_g}=\{u_g,u_{2 g},\dots,u_{g r_g} \}$ where $u$ is a primitive root
$\modu p$.
\begin{itemize}
\item
Observe  that, from relation (\ref{e104012}) p.\pageref{e104012}, we have
$r_g\leq f_g\leq \frac{p-1}{2g}=M_g$. Let $\beta=2\alpha\times M_g$.
\item
Show  that there are at least $M_g-r_g+1$ different values
$m,\quad 0\leq m \leq M_g-1$, such that $S_{g(\beta+2m+1)}(J_{r_g})\not\equiv 0 \modu p$. Suppose that this assertion is false and search for a contradiction :
then there exists at most $M_g-r_g$ different values $m,\quad 0\leq m\leq M_g-1$ such that
$S_{g(\beta+2m+1)}(J_{r_g})\not\equiv 0 \modu p$.
Then there exists at least $r_g$ different values
$m_i,\quad 0\leq m_i\leq M_g-1,\quad i=1,\dots,r_g$, such that
\begin{displaymath}
\begin{split}
&S_{g(\beta+2m_1+1)}(J_{r_g})\equiv 0 \modu p,\\
&S_{g(\beta+2m_2+1)}(J_{r_g})\equiv 0 \modu p,\\
&\vdots\\
&S_{g(\beta+2m_{r_g}+1)}(J_{r_g})\equiv 0\modu p,
\end{split}
\end{displaymath}
hence
\begin{displaymath}
\begin{split}
&n_1 + n_2 u^{g(\beta+2m_1+1)}
+\dots+n_{r_g} u^{g(r_g-1)(\beta+2m_1+1)}\equiv 0 \modu p,\\
&n_1 + n_2 u^{g(\beta+2m_2+1)}
+\dots+n_{r_g} u^{g(r_g-1)(\beta+2m_2+1)}\equiv 0 \modu p,\\
&\vdots\\
&n_1 + n_2 u^{g(\beta+2m_{r_g}+1)}
+\dots+n_{r_g} u^{g(r_g-1)(\beta+2m_{r_g}+1)}\equiv 0 \modu p.
\end{split}
\end{displaymath}
The determinant of this system of $r_g$ congruences in the $r_g$ indeterminates
$n_1, n_2,\dots, n_{r_g}$, after simplification by a power of
$u\not\equiv 0\modu p$ not intervening in $\modu p$ congruences,
\begin{equation}\label{e104013}
\begin{array}{|llll|}
1 & u^{g(2m_1+1)} & \dots & u^{g(r_g-1)(2m_1+1)} \\
1 & u^{g(2m_2+1)} & \dots & u^{g(r_g-1)(2m_2+1)} \\
\vdots&\vdots &&\vdots \\
1 & u^{g(2m_{r_g}+1)} & \dots & u^{g(r_g-1)(2m_{r_g}+1)} \\
\end{array}
\end{equation}
is a Vandermonde determinant with factors
$(u^{g(2m_k+1)}-u^{g(2m_l+1)})\not\equiv 0 \modu p$ because $u$ is a primitive root $\modu p$ and because
$g(2m_i+1)<p-1,\quad i=1,\dots,r_g$: in fact $g(2m_i+1)\leq g(2(M_g-1)+1)
=g(2(\frac{p-1}{2g}-1)+1)=p-1-g$.
Therefore $n_1= n_2= \dots= n_{r_g} = 0$, which contradicts hypothesis on $N(J_{r_g})$.
\item
Then, for $m$ varying through  a set $M\subset \{1,\dots,\frac{p-1}{2g}-1\},\quad Card(M)=
\frac{p-1}{2g}-r_g+1$ and for $\alpha=0,\dots,g-1$, we have
\begin{displaymath}
S_{g(2m+1+\alpha(p-1)/g)}(J_{r_g})
\equiv S_{g(2m+1)}(J_{r_g})\not\equiv 0 \modu p,
\end{displaymath}
and from foundation theorem
\begin{displaymath}
\phi_{2m+1+\alpha(p-1)/g}^*(t)\equiv 0\modu p,\quad \alpha=0,\dots,g-1.
\end{displaymath}
\end{itemize}
\item
\begin{itemize}
Show secondly that there are at most $\frac{p-1}{2g}-r_g+1$ different values $m,\quad 0\leq m\leq M_g-1$, such that
\begin{displaymath}
\phi_{2m+1+\alpha(p-1)/g}^*(t)\equiv 0 \modu p,\quad 0\leq m\leq M_g-1,\quad \alpha=0,\dots,g-1.
\end{displaymath}
Suppose that this assertion is false and search for a contradiction :
Then there exist at least $\frac{p-1}{2g}-r_g+2$ different values $m\in\N$ such that
\begin{displaymath}
\phi_{2m+1+\alpha(p-1)/g}^*(t)\equiv 0 \modu p,\quad 0\leq m\leq M_g-1,\quad \alpha=0,\dots,g-1.
\end{displaymath}
Then there exists  at most $\nu\leq r_g-2$ values
$m_1=0,m_2\not=0,\dots,m_{\nu}\not=0$ with at least one
$\alpha_i\in\N,\quad 0\leq\alpha_i\leq g-1,\quad i=1,\dots,\nu$,  such that
\begin{displaymath}\label{ea571}
\begin{split}
&\phi_{2m_1+1+\alpha_1(p-1)/g}^*(t)\not\equiv 0 \modu p\\
&\vdots\\
&\phi_{2m_{\nu}+1+\alpha_\nu(p-1)/g}^*(t)\not\equiv 0 \modu p.
\end{split}
\end{displaymath}
It is possible to find $l_1,l_2,\dots, l_\nu, l_{\nu+1}=1\in\N$, not all zero such that
\begin{displaymath}
\begin{split}
&l_1+l_2 u^{g(2m_1+1)}+\dots+l_\nu u^{(\nu-1)g(2m_1+1)}
+l_{\nu+1} u^{\nu g(2m_1+1)}\equiv 0\modu p,\\
&l_1+l_2 u^{g(2m_2+1)}+\dots+l_\nu u^{(\nu-1)g(2m_2+1)}
+l_{\nu+1} u^{\nu g(2m_2+1)}\equiv 0\modu p,\\
&\vdots\\
&l_1 + l_2 u^{g(2m_{\nu}+1)}+\dots+l_\nu u^{(\nu-1)g(2m_{\nu}+1)}
+l_{\nu+1} u^{\nu g(2m_{\nu}+1)}\equiv 0 \modu p.
\end{split}
\end{displaymath}
because the determinant
\begin{equation}\label{ea575bis}
\begin{array}{|llll|}
1 & u^{g(2m_1+1)} & \dots & u^{(\nu-1)g(2m_1+1)} \\
1 & u^{g(2m_2+1)} & \dots & u^{(\nu-1)g(2m_2+1)} \\
\vdots&\vdots &&\vdots \\
1 & u^{g(2m_{\nu}+1)} & \dots & u^{(\nu-1)g(2m_{\nu}+1)} \\
\end{array}
\end{equation}
is a Vandermonde determinant with factors
$(u^{g(2m_k+1)}-u^{g(2m_l+1)})\not\equiv 0 \modu p$,
because $u$ is a primitive root $\modu p$ and $g(2m_i+1) <p-1,\quad i=1,\dots,\nu$.
Therefore we have,
\begin{displaymath}
\begin{split}
& J_{\nu+1}=\{u^g,u^{2g},\dots,u^{(\nu+1)g}\}\\
& L(J_{\nu+1})=\{l_1,l_2,\dots,l_{\nu+1}\}\\
& S_{g(2m_1+1)}(J_{\nu+1})\equiv 0\modu p\\
& S_{g(2m_2+1)}(J_{\nu+1})\equiv 0\modu p\\
&\vdots\\
& S_{g(2m_\nu+1)}(J_{\nu+1})\equiv 0 \modu p
\end{split}
\end{displaymath}
Therefore, for $m\in M = \{m_1=0,m_2,\dots,m_\nu\}$,  we have $S_{g(2m+1)}(J_{\nu+1})\equiv 0 \modu p$ and for $m\in\{0,1,\dots,\frac{p-1}{2g}-1-M \}$, we have
\begin{displaymath}
\phi_{2m+1+\alpha(p-1)/g}^*(t)\equiv 0 \modu p,\quad \alpha=0,\dots,g-1,
\end{displaymath}
and hence
\begin{displaymath}
\begin{split}
& S_{g(2m+1)}(J_{\nu+1})\times \phi_{2m+1+\alpha(p-1)/g}^*(t)\equiv 0\modu p,\\
& m=0,1,\dots,\frac{p-1}{2g}-1,\quad \alpha=0,\dots,g-1,\\
\end{split}
\end{displaymath}
which leads to
\begin{displaymath}
S_{g(2m+1)}(J_{\nu+1})\times \phi_{2m+1}^*(t)\equiv 0\modu p,
\quad m=0,1,\dots,\frac{p-3}{2}.
\end{displaymath}
From lemma \ref{propaa} p.\pageref{propaa}, this leads to the relation
\begin{displaymath}
 \sum_{i=1}^{\nu+1} \frac{u_{gi} n_i}{(x+\zeta^{u_{gi}} y)}
+\sum_{i=1}^{\nu+1} \frac{u_{gi} n_i}{(x+\zeta^{-u_{gi}} y)}\equiv 0 \modu p,
\end{displaymath}
which contradicts the minimal value of $r_g>\nu+1$.
\item From first part of the proof, $r_g$ depends only of the number of Mirimanoff polynomials $\phi_{2m+1}(T),\quad 0\leq m\leq \frac{p-3}{2}$, with $\phi_{2m+1}(t)\equiv 0 \modu p$ and does not depend on the primitive root $u \modu p$.
\item From $\phi_1^*(t)=1$ and $\phi_{2m+1}^*(t)=\phi_{2m+1}(t)$ for $m=1,\dots,\frac{p-3}{2}$, we deduce that the number of $m, \quad 1\leq m \leq \frac{p-3}{2}$ with
$\phi_{2m+1}(t)\equiv 0 \modu p$ is the number of $m, \quad 0\leq m\leq \frac{p-3}{2}$ with $\phi_{2m+1}^*(t)\equiv 0 \modu p$, which completes the proof.
\end{itemize}
\end{enumerate}
\end{proof}
\end{thm}
%
%
\clearpage
%
\section {FLT First case: Mirimanoff in intermediate fields $K$ with
 $\Q\subset K\subset\Q(\zeta)$. Part II}\label{s14101}
This section deals, with another point of view, of Mirimanoff polynomials in intermediate fields $K,\quad \Q\subset K\subset \Q(\zeta)$.
\begin{itemize}
\item
The first subsection deals of Gauss periods for FLT.
\item
The second subsection gives an example of application: if $p$ is prime, $p\equiv 1 \modu 3$,  does not divide the class number of $K/\Q,\quad [K:\Q]=\frac{p-1}{3}$, then first case of FLT holds for $p$.
\item
The third subsection is a study of FLT first case in quadratic subfield of $\Q(\zeta)$.
\end{itemize}
%

\subsection  { First case of FLT and Gauss periods}\label{s108254}
In this subsection, we apply Gauss periods of cyclotomic equations to generalize foundation theorem to intermediate fields $\Q\subset K\subset \Q(\zeta)$, see Ribenboim \cite{rib},
p 115, to the Fermat's equation $x^p+y^p+z^p=0, \quad x y z\not\equiv 0 \modu p$ assumed as hypothesis.
%
We have shown the foundation theorem \ref{s12113} \pageref{s12113}
\begin{displaymath}
\phi_{2m+1}(t)\times\sum_{i=0}^{r-1} u_{(2m+1)i} \times l_i\equiv 0 \modu p,\quad m=1,\dots,\frac{p-3}{2}.
\end{displaymath}
where we recall that $t\equiv -\frac{x}{y}\modu p$ and $\phi_{2m+1}(T), \quad
m=1,\dots,\frac{p-3}{2}$ are the odd Mirimanoff polynomials.
In this section, we shall generalize foundation theorem to intermediate fields $K,\quad
\Q\subset K \subset \Q(\zeta)$, with $[K:\Q]=g,\quad p-1=g\times f$, where $f$ is odd.
%
\begin{thm}{ ***  Translation of Mirimanoff congruences in intermediate fields}\label{t10071a}

Let $K$ be an intermediate field, $\Q\subset K\subset \Q(\zeta)$, with $[K:\Q]=g,
\quad p-1=g\times f$, where $f$ is odd. Let $r_K$ be the $p$-rank of the class group of $K/\Q$. There exists a minimal $r_A\in\N,\quad r_A\leq r_K,\quad r_A<g$ such that
the odd Mirimanoff polynomials verify the system of  congruences $\modu p$
\begin{equation}\label{e105261}
\phi_{f(2n+1)}(t)\times \{\sum_{j=0}^{r_A} u_{j f^2 (2n+1)}\times  l_j\}\equiv 0\modu p,
\quad 2n+1=1,3,\dots, g-1,
\end{equation}
where $l_i\in\N,\quad 0\leq l_i\leq p-1,\quad l_0\not=0,\quad l_{r_A}\not=0$.
\begin{proof} $ $
\begin{itemize}
\item
Let us consider the algebraic integer $A=\prod_{i=0}^{f-1} (x+\zeta^{u_{g i}} y)\in\Z[\zeta]$.
We have $A\in K$.
\item
Let $O_K= K\cap \Z[\zeta]$  be the ring of integers of $K$.
We have $A O_K =U^p$, where $U$ is an ideal of $O_K$.
Let $\rho_A$ be the $p$-rank of the action of $G=Gal(\Q(\zeta)/\Q)$ on the  ideal class group $<Cl(U \Z[\zeta])>$ subgroup of the $p$-group $C_p$. Clearly $\rho_A\leq r_K$ and $\rho_A<g$ because, for all ideals $\mathbf b\subset K$, then ${\it N}_{K/\Q}(\mathbf b)$ is principal. It is possible to find $l_j,\quad 0\leq l_j\leq p-1,\quad l_0=1,\quad \quad l_{\rho_A}\not\equiv 0\modu p,\quad j=0,\dots, \rho_A$, such that
$\prod_{j=0}^{\rho_A} \sigma^{((p-1)/g))\times j}(U)^{l_j}$ is a principal ideal of $O_K$ and so, observing that $f=\frac{p-1}{g}$,
\begin{displaymath}
\begin{split}
& Q=\prod_{j=0}^{\rho_A} \sigma^{f  j} (A)^{l_j}=\zeta^v \times\eta\times \gamma^p,\\
& \gamma\in\Z[\zeta], \quad \eta\in \Z[\zeta+\zeta^{-1}]^*,\quad v\in \N,
\quad 0\leq v\leq p-1.
\end{split}
\end{displaymath}
\item
Then, we obtain
\begin{displaymath}
Q=\prod_{j=0}^{\rho_A} \prod_{i=0}^{f-1} (x+\zeta^{u_{ g  i+ f  j }} y)^{l_j} =\zeta^v\times \eta\times\gamma^p.
\end{displaymath}
Then it is possible to find a minimal natural integer $\modu \pi^{p+1}$ Fermat congruence rank $r_A\leq \rho_A$ such that
\begin{displaymath}
\begin{split}
&\prod_{j=0}^{r_A} \prod_{i=0}^{f-1} (x+\zeta^{u_{ g  i+ f  j }} y)^{l_j} - \prod_{j=0}^{r_A} \prod_{i=0}^{f-1} (x+\zeta^{-u_{ g  i+ f  j }} y)^{l_j}
\equiv 0\modu \pi^{p+1},\\
&l_i\in\N,\quad 0\leq l_i\leq p-1,\quad l_0\not=0,\quad l_{r_A}\not=0.
\end{split}
\end{displaymath}
Then, we apply foundation theorem \ref{s12113} p.\pageref{s12113}, mutatis mutandis  in this situation, to get
for $m=1,\dots, \frac{p-3}{2}$,
\begin{displaymath}
\phi_{2m+1}(t)\times (\sum_{j=0}^{r_A}
\sum_{i=0}^{f-1} u_{(2m+1) (g i + f  j)} \times l_j)\equiv 0 \modu p,
\end{displaymath}
so
\begin{equation}\label{e105264}
\phi_{2m+1}^*(t)\times (\sum_{i=0}^{f-1} u_{(2m+1) g i})
\times (\sum_{j=0}^{r_A} u_{(2m+1) f j} \times l_j)\equiv 0 \modu p,
\quad m=0,\dots,\frac{p-3}{2}.
\end{equation}
which is the translation of foundation theorem \ref{pa212} p.\pageref{pa212} in intermediate field $K$.
\item
Let us evaluate $\sum_{i=0}^{f-1} u_{(2m+1) g i}$:
\begin{itemize}
\item
If $(2m+1) \times g\equiv 0 \modu (p-1)$, then $\sum_{i=0}^{f-1} u_{(2m+1) g i} \equiv f \not\equiv 0 \modu p$. In that case $2m+1= f(2n+1), \quad 1\leq 2n+1 <g$.
\item
If $(2m+1) \times g\not \equiv 0 \modu (p-1)$, then $\sum_{i=0}^{f-1} u_{(2m+1) g i}\equiv \frac{u^{(2m+1) g f}-1}{u^{(2m+1) g}-1}
\equiv \frac{u^{(2m+1)(p-1)}-1}{u^{(2m+1) g}-1}\equiv 0 \modu p$.
\end{itemize}
\item
Then we obtain $\phi_{f(2n+1)}(t)\times (\sum_{j=0}^{r_A } u_{j (2n+1) f^2 } \times l_j)\equiv 0 \modu p,\quad 1\leq 2n+1< g$.
\end{itemize}
\end{proof}
\end{thm}
%
{\bf Example : }
Suppose that $g=6$, that $p-1\not\equiv 0 \modu 3^2,\quad p-1\not\equiv 0 \modu 2^2$  and that $p\| h(K/\Q)$, so $r_A=1$. Then this  generalization of foundation theorem (\ref{t10071a}) gives, with $f=\frac{p-1}{6}$,
\begin{displaymath}
\begin{split}
&\phi_{(p-1)/6}(t)\times (l_0+ u_{((p-1)/6)^2} \times l_1)\equiv 0 \modu p,\\
&\phi_{3(p-1)/6}(t)\times (l_0+ u_{3((p-1)/6)^2}\times l_1)\equiv 0 \modu p,\\
&\phi_{5(p-1)/6}(t)\times (l_0+ u_{5((p-1)/6)^2} \times l_1)\equiv 0 \modu p.
\end{split}
\end{displaymath}
Then, this implies that at least two of the values $\phi_{(p-1)/6}(t), \phi_{3(p-1)/6}(t),
\phi_{5(p-1)/6}(t)$ are null $ \modu p$ because, from hypothesis on $p-1$, we get $2\times(\frac{p-1}{6})^2\not\equiv 0 \modu p-1$, so $u_{3(p-1)/6)^2}-u_{(p-1)/6)^2}\not\equiv 0\modu p$  and also $4\times(\frac{p-1}{6})^2\not\equiv 0 \modu p-1$, so $u_{5(p-1)/6)^2}-u_{(p-1)/6)^2}\not\equiv 0\modu p$.

{\bf Remark :} Let a field $K\subset \Q(\zeta)$. Let $r_K$ be the $p$-rank of the class group $K/\Q$; let $r_p$ be the $p$-rank of the class group $\Q(\zeta)/\Q$; the interest of this theorem is  that $r_K\leq r_p$ and that $r_K \leq [K:\Q]$ : so we have now to consider system  of linear congruences of rank smaller than $r_K$, very  more easily solvable that a system of $r_p$ conguences in $\Q(\zeta)/\Q$. We have translated Mirimanoff congruences in intermediate fields $K/\Q$.
%
\begin{cor} { *** } \label{c010161}
Let $p-1=f\times g$ with  $f$ odd. Let $K$ be the intermediate field $\Q\subset K \subset \Q(\zeta)$ and $[K:\Q]=g$. Suppose that $p$ does not divide the class number of $K/\Q$. Then we have  the Mirimanoff's polynomials congruences
\begin{displaymath}
\phi_{ f (2 n+1)}(t)\equiv 0 \modu p,\quad 2n+1=1,3,5,\dots, g-1.
\end{displaymath}
\begin{proof}
Apply theorem \ref{t10071a} p.\pageref{t10071a} with here $r_A=0$.
\end{proof}
\end{cor}
%
\begin{thm}{ *** }\label{t09071a}
Let $f,g \in\N$, with $f\times g=p-1$, $Gcd(f,g)=d$ and $f$ odd.
Let $r_A$ be the rank defined in theorem \ref{t10071a} p.\pageref{t10071a}.
Then \begin{equation}\label{e105262}
\begin{split}
& Card( \{\phi_{f(2n+1)}(t)\not\equiv 0\modu p \quad
|\quad n=0,\dots, \frac{g-2 d}{2 d}\})=r_A,\\
& Card( \{\phi_{f(2n+1)}(t)\equiv 0\modu p \quad
|\quad n=0,\dots, \frac{g-2d}{2d}\})=\frac{g-2 d}{2 d}+1-r_A.
\end{split}
\end{equation}
Moreover $r_A\leq r_K$ where $r_K$ is the $p$-rank of the  class group
of $K,\quad [K:\Q]=g$.
\begin{proof}
Let us denote
\begin{displaymath}
S_{f(2n+1)}\equiv\sum_{j=0}^{r_A} u_{jf^2(2n+1)}\times l_j\modu p,\quad
n=0,\dots,\frac{g-2d}{2d},
\end{displaymath}
where $r_A$ is the minimal natural integer which verifies the relation (\ref{e105264}) p.\pageref{e105264}.
The foundation theorem \ref{s12113} p.\pageref{s12113} is the relation (\ref{e105264}) p.\pageref{e105264} which implies for $2m+1=f(2n+1),\quad n=0,\dots,\frac{g-2d}{2d}$:
\begin{displaymath}
\phi_{f(2n+1)}(t)\times S_{f(2n+1)}\equiv 0\modu p.
\end{displaymath}
Observe that
$u_{f^2(2n+1)}\not\equiv u_{f^2(2n^\prime+1)},\quad n=0,\dots,\frac{g-2d}{2d}\modu p$,
because $f^2(2n+1)\not\equiv 0\modu f^2(2n^\prime+1)\modu p-1
\Leftrightarrow f(2n+1)\not\equiv f(2n^\prime+1)\modu g
\Leftrightarrow 2n+1\not\equiv 2n^\prime+1\modu \frac{g}{d}$, which is not possible from hypothesis
$0\leq n,n^\prime\leq \frac{g-2d}{2d}$.
It is not possible that
$S_{f(2n_i+1)}=\sum_{j=0}^{r_A} u_{jf^2(2n_i+1)}\not\equiv 0\modu p, \quad i=1,\dots,r_A$,
because $\phi_{f(2n_i+1)}(t)\not\equiv 0\modu p,\quad i=1,\dots,r_A$, and because the Vandermonde determinant $\Delta=|u_{jf^2(2n_i+1)}|_{1\leq i,j\leq r_A}$
 verifies  $\Delta\not\equiv 0\modu p$:
indeed the foundation theorem \ref{s12113} p.\pageref{s12113}, its reciprocal \ref{propaa} p\pageref{propaa}, the theorems \ref{pa212} p.\pageref{pa212}, \ref{t12092} p.\pageref{t12092} and the corollary \ref{cor1209} p.\pageref{cor1209} can be applied here
starting from relation (\ref{e105264}) p.\pageref{e105264}, which can be applied for
$m=1,\dots,\frac{p-3}{2}$.
As in corollary \ref{cor1209} p.\pageref{cor1209}, we have $\phi_{f(2n+1)}(t)\equiv 0\modu p$ else
$S_{f(2n+1)}\equiv 0\modu p$, the two congruences being simultaneously impossible,
because $r_A$ is the minimal natural integer verifying relation (\ref{e105264}).
This leads to relation (\ref{e105262}) p.\pageref{e105262}.
\end{proof}
\end{thm}
%
\begin{cor}{ *** }\label{c209261}
Let $f,g\in\N$ with $f\times g=p-1,\quad Gcd(f,g)=d > 1$ and $f$ odd.
If there exists $n\in\N,\quad n\leq \frac{g-2d}{2d}$ with $\phi_{f(2n+1)}(t)\equiv 0\modu p$
then
\begin{equation}\label{e209261}
\phi_{f(2n+1+j g/d)}(t)\equiv 0\modu p, \quad j=0,\dots, d-1.
\end{equation}
\begin{proof}
Observe the periodicity
\begin{displaymath}
S_{f(2n+1+jg/d)}\equiv \sum_{k=0}^{r_A} u_{k f^2(2n+1+j g/d)}\equiv \sum_{k=0}^{r_A} u_{k f^2(2n+1)}
\equiv S_{f(2n+1)} \modu p
\end{displaymath}
and apply theorem \ref{t09071a} p. \pageref{t09071a}.
\end{proof}
\end{cor}
%
\begin{cor}\label{c209271}
Let $f,g\in\N$ with $f\times g=p-1,\quad Gcd(f,g)=d>1$ and $f$ odd.
Then

{\bf else}
\begin{displaymath}
B_{p-f(2n+1)}\equiv 0\modu p,\quad n=0,\dots,\frac{g-2d}{2d},
\end{displaymath}
where
$B_{p-f(2n+1)}$ are even Bernoulli Numbers.

{\bf else} there exists at least one $n\in\N,\quad 0\leq n\leq \frac{g-2d}{2}$ such that
\begin{displaymath}
\phi_{f(2n+1+jg/d)}(t)\equiv 0 \modu p,\quad j=0,\dots,d-1.
\end{displaymath}
\begin{proof}
If $\phi_{f(2n+1)}(t)\not\equiv 0\modu p$ for $n=0,\dots,\frac{g-2d}{2d}$ then, from Mirimanoff congruences we derive that $B_{p-f(2n+1)}\equiv 0\modu p, \quad n=0,\dots,\frac{g-2d}{2d}$. If there exists one
$n,\quad 0\leq n\leq\frac{g-2d}{2d}$ with $\phi_{f(2n+1)}(t)\equiv 0\modu p$ apply corollary
\ref{c209261} p. \pageref{c209261}.
\end{proof}
\end{cor}

%


%

%

\subsection{ FLT results in intermediate fields $\Q\subset K\subset\Q(\zeta)$.}
\label{s108255}
Let, as previously, denote $u$  a primitive root $\modu p$.
Let $\sigma :\Q(\zeta)\rightarrow \Q(\zeta)$ be the $\Q$-isomorphism defined by
$\sigma : \zeta \rightarrow \zeta^u$. Let $g\in\N,\quad g>2,\quad g$ dividing $p-1$. Let $K$ be the intermediate field, $\Q\subset K\subset \Q(\zeta),\quad [K:\Q]=g$. We suppose, in this subsection that $p$ does not divides $h(K/\Q)$, class number of the extension $K/\Q$. Let $\s$ be the integral ideal of $\Z[\zeta]$ with $(x+\zeta y)\Z[\zeta]=\s^p$. Let $C=\prod_{i=0}^{(p-1)/g-1} \sigma^{g i} (\s)$.
We have $C=c\Z[\zeta]$, where $ c$ is an integral ideal of $K/\Q$.
Then,
\begin{displaymath}
\prod_{i=0}^{(p-1)g-1}(x+\zeta^{u_{gi}} y)\Z[\zeta]=
\prod_{i=0}^{(p-1)/g-1}\sigma^{g i}( \s^p)= c^p\Z[\zeta].
\end{displaymath}
From hypothesis, $p$ does not divide $h(K/\Q)$, and so we get
\begin{equation}\label{e12101}
\begin{split}
&\prod_{i=0}^{(p-1)/g-1} (x+\zeta^{u_{g i}} y)= \eta\zeta^v\gamma^p,\\
& \gamma\in\Z[\zeta],\quad \eta\in\Z[\zeta+\zeta^{-1}]^*,\quad v\in\Z,\quad  0\leq v\leq p-1.
\end{split}
\end{equation}
We shall apply this results when $3$ divides $p-1$ and $g=\frac{p-1}{3}$.
%
\begin{thm}\label{t12101}
Suppose that $3$ divides $p-1$.
Let $K$ be the field $\Q\subset K\subset\Q(\zeta)$ with $[K:\Q]=\frac{p-1}{3}$.
If $p$ does not divide $h(K/\Q)$, class number of the extension $K/\Q$, then the first case of FLT holds for $p$.
\begin{proof}
Let $g=\frac{p-1}{3}$.
From equation (\ref{e12101}) p.\pageref{e12101}, we get, similarly to proof of theorem \ref{p21a} p.\pageref{p21a},
for $m\geq 1$,
\begin{displaymath}
\phi_{2m+1}(t)\times\sum_{i=0}^{(p-1)/g-1} u_{gi}^{2m+1}\equiv 0 \modu p.
\end{displaymath}
Then, for $m=1$, we obtain
\begin{displaymath}
\phi_3(t)\times\sum_{i=0}^2 u^{3 g i}\equiv 0 \modu p.
\end{displaymath}
But $3 g= p-1$ and so $u^{3 g i}\equiv 1 \modu p$.
Then, we should have $3\times \phi_3(t)\equiv 0 \modu p$, which is not possible, because
$\phi_3(t)\not \equiv 0 \modu p$ : in fact $\phi_3(t)\equiv 0 \Longrightarrow P_3(-t)\equiv 0\modu p$, where $P_3(-t)=-t(1+t)$ is a Kummer Polynomial, see for instance \cite{rib}, 4D p 125; but $(1+t)=\frac{y-x}{y}\not\equiv 0 \modu p$ because we assume in this article that $p>5$ and therefore that $y-x\not\equiv 0 \modu p$.
\end{proof}
\end{thm}
%
{\bf Remark:} In comparison, when $g=\frac{p-1}{2}$, Vandiver's conjecture asserts that $h(K/\Q)=h(\Q(\zeta+\zeta^{-1})/\Q) \not\equiv 0 \modu p$ and there was several failed tentatives of proof of first case of FLT assuming Vandiver's conjecture.
%
\subsection{  Fermat equation in  quadratic subfield of $\Q(\zeta)/\Q$.}
Let us define $\delta_p\in\{-1,1\}$ by $\delta_p=1$ if $p\equiv 1\modu 4$ and
$\delta_p=-1$ if $p\equiv 3\modu 4$. Then the quadratic field $\Q(\sqrt{ \delta_p p})\subset\Q(\zeta)$.
%
\begin{lem}\label{l104281}
For $1 \leq j <\frac{p-1}{2}$, we have
$P=\prod_{j=0}^{(p-3)/2} \sigma_2^{j}(\s)= Q\Z[\zeta]$, where $Q$ is a principal integral ideal
of the quadratic field  $\Q(\sqrt{\delta_p p})$.
\begin{proof}
Note that $\sigma_2^{(p-1)/2} :\zeta\rightarrow \zeta^{u_2^{(p-1)/2}}$ is the identity map.
Note that $\sigma_2\circ\sigma_2^j=\sigma_2^{j+1}$.
Compute $S=\sigma_2(\prod_{j=0}^{(p-3)/2} \sigma_2^{j}(\s))$ where we recall that
$(x+\zeta y)\Z[\zeta]=\s^p$.
We have
\begin{equation}
\begin{split}
&S=\prod_{j=0}^{(p-3)/2} \sigma_2^{j+1}(\s)\\
&\Longrightarrow\\
&S=(\prod_{j=1}^{(p-3)/2} \sigma_2^{j}(\s))\times\sigma_2^{(p-1)/2}(\s)\\
&\Longrightarrow\\
&S=(\prod_{j=1}^{(p-3)/2} \sigma_2^{j}(\s))\times\s,\\
&\Longrightarrow\\
&S=P.
\end{split}
\end{equation}
Therefore, from $\sigma_2(P)=P$, there exists an integral ideal $Q$ of $\Q(\sqrt{\delta_p p})$ such that $P=Q\Z[\xi]$. The order of the class of the ideal
$P$ of $\Z[\xi]$ is  $1$ or $p$.
Then, the order of the ideal class  of $Q$ in $\Q(\sqrt{\delta_p p})$ is $1$ or $p$. From classical upper bound of the class number of the quadratic number fields, the ideal $Q$ is principal in $\Q(\sqrt{\delta_p p})$.
\end{proof}
\end{lem}
%
\begin{lem}\label{l104282}
Let  $p\equiv 3\modu 4$.
Let $A$ be the ring of integers of $\Q(\sqrt{-p})$.
Then the Fermat equation in $A$ implies
\begin{equation}\label{e23015}
\prod_{i=0}^{(p-3)/2}(x+\zeta^{u_{2i}} y)= \gamma^p,\quad \gamma\in A.
\end{equation}
\begin{proof}
We have $\prod_{i=0}^{(p-3)/2}(x+\zeta^{u_{2i}}y)\in A$.
From lemma (\ref{l104281}) we get

\begin{displaymath}
\prod_{i=0}^{(p-3)/2} (x+\zeta^{u_{2i}}y)=\varepsilon\times\gamma^p,\quad
\varepsilon\in \Z[\zeta]^*,\quad \gamma\in A.
\end{displaymath}
Therefore $\varepsilon\in A^*$.
We have  $\Q(\sqrt{-p})\subset \Q(\zeta)$.
Then the result is an immediate consequence of the structure of units group of imaginary quadratic extension of $\Q(\sqrt{-p})$ with $I\not\in A$.
\end{proof}

\end{lem}
%
\begin{lem}\label{l1040283}
Let  $p\equiv 1\modu 4$.
Then
\begin{equation}\label{e23016}
\prod_{i=0}^{(p-3)/2}(x+\zeta^{u_{2i}} y)= \varepsilon\times \gamma^p, \quad \gamma\in A,\quad\varepsilon \in A^*.
\end{equation}
\begin{proof}
Here $A^*$ is generated by a fundamental unit $\frac{\varepsilon_0+\sqrt{p}\varepsilon_1}{2}$.
\end{proof}
\end{lem}
%
%
\clearpage
\section{FLT first case: Elementary congruences on $t\equiv -\frac{x}{y}\modu p$.}\label{s18121}
Let $t\equiv -\frac{x}{y}\modu p$.
The definitions of section (\ref{s3s1}) with $f=min((e_p+2),\frac{p-1}{2})$, the sets $I_f=\{g_i\quad |\quad i=1,\dots,f\}$ and with the sets $L(I_f)=\{l_i\quad |\quad i=1,\dots,f\}$ used in proposition (\ref{p4}) are available in this section.
\begin{itemize}
\item
The two first subsections  gives  a system of polynomial congruences of $t$: we conjecture that this system of congruences is strong and contains a proof of first case of FLT. We prove also that order $a$ of $t \modu p$, the smaller $a\in\N,\quad a>1,\quad t^a\equiv 1\modu p$, verifies $a> e_p+2$. The two subsections gives two method of proof of these congruences.
\item
In the third subsection, we prove, from symmetries on Mirimanoff congruences,
that $\phi_{p-1}(t^2)\equiv 0\modu p,$ and we derive  a proof of this particular case of Fermat-Wiles theorem : if
$p \| (x-y)(y-z)(z-x),\quad x,y,z\in\N-\{0\}$, then $x^p+y^p+z^p\not=0$.
\item
In the fourth subsection, we  obtain a set of $p-1$ explicit congruences
$P_m(t)\equiv 0\modu p,\quad P_m(T)\in{\bf F}_p[T],\quad m=1,\dots,p-1$ and also as a consequence $\phi_{p-2}(t^2)\equiv 0\modu p$.
\end{itemize}
\subsection { A First elementary approach}
%
\begin{thm}{ *** }\label{p06031}
Let $t\in\N,\quad t\equiv -\frac{y}{x} \modu p$.
Let {\bf any}  set $I_f=\{g_i\ |\ i=1,\dots,f$, and $L(I_f)$.
Let us, for $i=1, \dots , f,\quad k=1, \dots , p-1$, define $m_{ik}\in\N$  by the congruence
\begin{equation}
g_i\times m_{ik}\equiv k \modu p,\quad 1\leq m_{ik}\leq p-1.
\end{equation}
Then, we have  the $p-2$ congruences:
\begin{equation}\label{e208072}
\sum_{i=1}^f (l_i g_i) \times
\{(t^{m_{ik}}+t^{p-m_{ik}})-(t^{m_{i 1}}+t^{p-m_{i 1}})\}\equiv 0 \modu p,
\quad k=2,\dots p-1.
\end{equation}
Reciprocally, if

\begin{displaymath}
\sum_{i=1}^f (l_i g_i) \times
\{(t^{m_{ik}}+t^{p-m_{ik}})-(t^{m_{i 1}}+t^{p-m_{i_1}})\}\equiv 0 \modu p,
\quad k=2,\dots p-1.
\end{displaymath}
then
\begin{displaymath}
\sum_{i=1}^f \frac{g_i l_i}{(x+\zeta^{g_i}y)}+\sum_{i=1}^f \frac{g_i l_i}{(x+\zeta^{g_i} y)}
\equiv 0\modu p
\end{displaymath}
\begin{proof}
From proof of theorem \ref{p21a} p.\pageref{p21a} the relation
\begin{displaymath}
\sum_{i=1}^f \frac{g_i l_i}{(x+\zeta^{g_i}y)}+\sum_{i=1}^f \frac{g_i l_i}{(x+\zeta^{g_i} y)}
\equiv 0\modu p
\end{displaymath}
is equivalent to the relation (\ref{e53g}) p.\pageref{e53g}
\begin{equation}\label{e208061}
\begin{split}
&F(\zeta)=\\
&\sum_{i=1}^f g_i l_i\times
(x^{p-1}-\zeta^{g_i} x^{p-2}y+\zeta^{2g_i} x^{p-3}y^2-\dots
-\zeta^{(p-2)g_i}x y^{p-2}+\zeta^{(p-1)g_i}y^{p-1})\\
&+\sum_{i=1}^f g_i l_i\times
(x^{p-1}-\zeta^{-g_i}  x^{p-2}y+\zeta^{-2g_i} x^{p-3}y^2-\dots
-\zeta^{-(p-2)g_i}x y^{p-2}+\zeta^{-(p-1)g_i}y^{p-1})\\
&\equiv 0 \modu \pi^{p-1}.
\end{split}
\end{equation}
Compute the coefficient  $\alpha_k$ of $\zeta^k$ in expression of $F(\zeta)$; there exists one $\alpha_k\in\Z$
such that  we have, for all $k,\quad 1\leq k\leq p-1$ :
\begin{equation}\label{ea583}
\begin{split}
&\sum_{i=1, g_i m_{ik}\equiv k \modu p}^{i=f} (l_i g_i) \times((-1)^{m_{ik}} (\frac{y}{x})^{m_{ik}}) \\
&+\sum_{i=1,(-g_i) m_{ik}^\prime\equiv k \modu p}^{i=f} (l_i g_i) \times((-1)^{m_{ik}^\prime} (\frac{y}{x})^{m_{i k}^\prime})\equiv \alpha_k \modu p,\\
&\Longleftrightarrow\\
&\sum_{i=1,g_i m_{ik}\equiv k \modu p}^{i=f} (l_i g_i) \times t^{m_{ik}}
+\sum_{i=1,(-g_i) m_{ik}^\prime\equiv k \modu p}^{i=f} (l_i g_i) \times t^{m_{ik}^\prime}
\equiv \alpha_k \modu p,\\
&\Longleftrightarrow\\
&\sum_{i=1,g_i m_{ik}\equiv k \modu p}^{i=f} (l_i g_i) \times t^{m_{ik}}
+\sum_{i=1,g_i m_{ik}\equiv k \modu p}^{i=f} (l_i g_i) \times t^{p-m_{ik}}
 \equiv \alpha_k \modu p,\\
&\Longleftrightarrow\\
&S_k=\sum_{i=1,g_i m_{ik}\equiv k \modu p}^{i=f} (l_i g_i) \times (t^{m_{ik}}+t^{p-m_{ik}})
\equiv \alpha_k \modu p,\quad k=1,\dots,p-1.
\end{split}
\end{equation}
From relation \ref{e208061} p. \pageref{e208061}, we derive that $S_k-S_1\equiv 0\modu p$.
then we obtain the result, computing $S_k-S_1$.

\end{proof}
\end{thm}
%
%
\begin{thm}{ *** }\label{t08031} $ $

Let $f$ be an integer with $p-2\geq f\geq min(e_p+2,\frac{p-1}{2})$.

Let $E_f$ be  {\bf any } set

$E_f=\{g_i\in\N\quad | \quad i=1,\dots,f,\quad 1\leq g_i\leq p-1 \}$.

Let $K_f$ be {\bf any} set

$K_f=\{k_j\in\N\quad |\quad j=1,\dots,f,\quad 1\leq k_j\leq p-1 \}$.

Let us note, for $i=1,\dots f,\quad j=1,\dots,f$,
\begin{displaymath}
\begin{split}
& m_{i j}\equiv g_i^{-1}\times k_j\modu p,\quad  1\leq m_{i j}\leq p-1,\\
&\theta_{i j}= (t^{m_{i j}}+t^{p-m_{i j}})-(t^{m_{i 1}}+t^{p-m_{i 1}}).
\end{split}
\end{displaymath}

If the first case of FLT would fail for $p$ then

the determinant
\begin{displaymath}
\Delta_{f}(t)=
\begin{array}{|lllll|}
\theta_{1 1} & \dots & \theta_{i 1} & \dots & \theta_{f 1}\\
\vdots & & \vdots & & \vdots \\
\theta_{1 j} & \dots & \theta_{i j} & \dots & \theta_{f j}\\
\vdots & & \vdots & & \vdots \\
\theta_{1 f}& \dots & \theta_{i f} & \dots & \theta_{f f}\\
\end{array}
\equiv 0 \modu p.
\end{displaymath}
\begin{proof} $ $
\begin{itemize}
\item
Note, at first, that if $g_{i_1}=g_{i_2}$ for some $1\leq i_1\leq i_2\leq p-1$, then
$m_{i_1 j}=m_{i_2 j}$ for $j=1,\dots,f$ and $\Delta_{f}(t)=0$ because it has two identical columns.
In the same way,
if $k_{j_1}=k_{j_2}$ for some $1\leq j_1\leq j_2\leq p-1$, then
$m_{i j_1}=m_{i j_2}$ for $i=1,\dots,f$ and $\Delta_{f}(t)=0$ because it has two identical rows.
\item
Suppose then that
\begin{displaymath}
\begin{split}
& E_f=\{g_i\in\N,\quad i=1,\dots,f,\quad 1\leq g_i\leq p-1,
\quad i\not= i^\prime\Longrightarrow g_i\not=g_{i^\prime}\}. \\
& K_f=\{k_j\in\N,\quad j=1,\dots,f,\quad 1\leq k_j\leq p-1,
\quad j\not=j^\prime\Longrightarrow k_j\not=k_{j^\prime}\}.
\end{split}
\end{displaymath}
It is then an immediate consequence of the system of $f$ first  congruences of the previous proposition \ref{p06031} p.\pageref{p06031} on the $f$ indeterminates $X_i= l_i g_i\modu p,\quad i=1,\dots,f$ with $I_f=E_f$.

\end{itemize}
\end{proof}
\end{thm}
%
\subsection{ A second elementary approach }
\begin{thm}{ *** }\label{t20021a}
Let $t\in\N,\quad t\equiv -\frac{y}{x} \modu p$.
Let {\bf any}  set $I_f$ and $L(I_f)$.
Let us, for $i=1, \dots , f,\quad k=1, \dots , p-1$, define $m_{ik}\in\N$  by the congruence
\begin{equation}
g_i\times m_{ik}\equiv k \modu p,\quad 1\leq m_{ik}\leq p-1.
\end{equation}
Then, we have  the $p-1$ congruences:
\begin{equation}
\sum_{i=1}^f (l_i g_i) \times
(t^{m_{ik}}+t^{p-m_{ik}})\equiv 0 \modu p,
\quad k=1,\dots p-1.
\end{equation}
Reciprocally, if
\begin{displaymath}
\sum_{i=1}^f (l_i g_i) \times
(t^{m_{ik}}+t^{p-m_{ik}})\equiv 0 \modu p,
\quad k=1,\dots p-1.
\end{displaymath}
then
\begin{displaymath}
\sum_{i=1}^f \frac{g_i l_i}{(x+\zeta^{g_i}y)}+\sum_{i=1}^f \frac{g_i l_i}{(x+\zeta^{g_i} y)}
\equiv 0\modu p
\end{displaymath}
\begin{proof}
Let us consider, for the indeterminate $X$, the rational function  $F(X)\in \Q(X)$ defined in relation (\ref{e53g}) p.\pageref{e53g} and (\ref{e31101}) p.\pageref{e31101}.
\begin{displaymath}
F(X)=\sum_{i=1}^f g_i l_i\times
(\sum_{m=0}^{p-1} (-1)^m x^{p-1-m}y^m X^{m g_i}
+\sum_{m=0}^{p-1} (-1)^m x^{p-1-m} y^m X^{-m g_i})
\end{displaymath}
From theorem (\ref{t02021a}) p.\pageref{t02021a} we get
$F(1)\equiv  F^\prime_X(1)\equiv   F^{(2)}_X(1)\equiv\dots\equiv  F^{(p-1)}\equiv 0 \modu p$.
Considering congruences $\modu p$ we can use the lighter notation

\begin{displaymath}
F(X)=\sum_{i=1}^f g_i l_i\times
(\sum_{m=0}^{p-1} (-1)^m (-t)^m X^{m g_i}
+\sum_{m=0}^{p-1} (-1)^m (-t)^m X^{-m g_i})
\end{displaymath}
and so
\begin{displaymath}
F(X)=\sum_{i=1}^f g_i l_i\times
(\sum_{m=0}^{p-1} t^m X^{m g_i}
+\sum_{m=0}^{p-1} t^m X^{-m g_i})
\end{displaymath}
We compute $F^{(p-1)}_X(X)$.
We have
\begin{displaymath}
(X^{m g_i})^{(p-1)}_X = m g_i\times (m g_i-1)\times\dots (m g_i- (p-2))\times X^{m g_i -(p-1)}
= a_i\times  X^{m g_i -(p-1)}.
\end{displaymath}
If $a_i\not\equiv 0 \modu p$ then $m g_i-(p-1)\equiv 0 \modu p$  and, from Wilson theorem, $a_i\equiv -1\modu p$.
Then,
\begin{displaymath}
(X^{m g_i})^{(p-1)}_{X=1}= - g_i l_i,\quad m g_i\equiv p-1\modu p.
\end{displaymath}
In the same way, we find that
\begin{displaymath}
(X^{-m g_i})^{(p-1)}_{X=1}=- g_i l_i,\quad m g_i\equiv 1\modu p.
\end{displaymath}
and finally, from $F_X^{(p-1)}(1)\equiv 0\modu p$,  that
\begin{displaymath}
\sum_{i=1, \quad m_i g_i\equiv 1 \modu p}^f  g_i l_i \times (t^{m_i }+t^{(p-m_i)})\equiv 0\modu p.
\end{displaymath}
The same computation for $F^{(p-2)}(X)$ leads to the congruence
\begin{displaymath}
\begin{split}
&\sum_{i=1,\quad m_{i,2}\times g_i \equiv 2 \modu p}^f
g_i l_i \times (t^{m_{i,2}}+t^{p-m_{i,2}})\\
&-(\sum_{i=1,\quad m_{i,1}\times  g_i\equiv 1 \modu p}^f
g_i l_i \times (t^{m_{i,1}}+t^{p-m_{i,1}}))\times \zeta
\equiv 0 \modu p,
\end{split}
\end{displaymath}
so
\begin{displaymath}
\sum_{i=1, \quad m_{i,2}\times g_i \equiv 2 \modu p }^f
g_i l_i \times (t^{m_{i,2}}+t^{p-m_{i,2}})
\equiv 0 \modu p,
\end{displaymath}
and so on.
The proof of reciprocal result is similar to theorem \ref{p06031} p.\pageref{p06031}.
\end{proof}
\end{thm}
%
\begin{thm}{ *** }\label{t21021a} $ $

Let $f$ be an integer with $f\geq min(e_p+2,\frac{p-1}{2})$.

Let $E_f$ be  {\bf any } set

$E_f=\{g_i\in\N,\quad i=1,\dots,f,\quad 1\leq g_i\leq p-1 \}$.

Let $K_f$ be {\bf any} set

$K_f=\{k_j\in\N,\quad j=1,\dots,f,\quad 1\leq k_j\leq p-1 \}$.

Let us define for $i=1,\dots f,\quad j=1,\dots,f$
\begin{displaymath}
\begin{split}
& m_{i j}\equiv g_i^{-1}\times k_j\modu p,\quad  1\leq m_{i j}\leq p-1,\\
& g_i\in E_f,\quad k_j\in K_f,\\
&\theta_{i j}= (t^{m_{i j}}+t^{p-m_{i j}}\modu p).
\end{split}
\end{displaymath}
If the first case of FLT would fail for $p$ then
the determinant
\begin{displaymath}
\Delta_{f}(t)=
\begin{array}{|lllll|}
\theta_{1 1} & \dots & \theta_{i 1} & \dots & \theta_{f 1}\\
\vdots & & \vdots & & \vdots \\
\theta_{1 j} & \dots & \theta_{i j} & \dots & \theta_{f j}\\
\vdots & & \vdots & & \vdots \\
\theta_{1 f}& \dots & \theta_{i f} & \dots & \theta_{f f}\\
\end{array}
\equiv 0 \modu p.
\end{displaymath}
\begin{proof} $ $
Similar to theorem \ref{t08031} p.\pageref{t08031} starting of theorem \ref{t20021a} p.\pageref{t20021a}.
\end{proof}
\end{thm}
%

\begin{cor}{ *** }\label{c010051} $ $

Let $t\equiv -\frac{x}{y}\modu p$.

Let $f$ be an integer with $f\geq min(e_p+2,\frac{p-1}{2})$.

Let $u$ be a primitive root $\modu p$.

Let $l\in\N,\quad 1\leq l\leq p-1$.

Let $m\in\N,\quad 1\leq l\leq p-1$.

If the first case of FLT would fail for $p$ then

the determinant $\Delta_{f}(t)=$
\begin{displaymath}
\begin{array}{|lllll|}
t^{u_{m -l}}+t^{p-u_{m-l}} & \dots &  t^{u_{m -l\times i}}+t^{p- u_{m-l\times i}}
& \dots & t^{u_{m -l \times f}}+t^{p- u_{m-l\times  f}}\\
\vdots & & \vdots & & \vdots \\
t^{u_{m\times j -l}}+t^{p-u_{m\times j-l}} & \dots &  t^{u_{m\times j -l\times i}}
+t^{p- u_{m\times j -l\times i}}
& \dots & t^{u_{m\times j -l \times f}}+t^{p- u_{m\times j -l\times  f}}\\
\vdots & & \vdots & & \vdots \\
t^{u_{m\times f -l}}+t^{p-u_{m\times f-l}} & \dots &  t^{u_{m\times f -l\times i}}
+t^{p- u_{m\times f -l\times i}}
& \dots & t^{u_{m\times f -l \times f}}+t^{p- u_{m\times f -l\times  f}}\\
\end{array}
\end{displaymath}
verifies the congruence $\Delta_{f}(t)\equiv 0 \modu p$.
\begin{proof} $ $
We apply theorem \ref{t21021a} p.\pageref{t21021a}
with the sets $E_f=\{ u_{l \times i}\quad |\quad i=1,\dots,f\}$ and with the set
$K_f=\{ u_{m \times j}\quad |\quad j=1,\dots,f \}$.
\end{proof}
\end{cor}
%
We give now two examples of application of theorem \ref{t21021a} p.\pageref{t21021a}.
%
\begin{thm}\label{t20022a}
Suppose that $p\equiv 3 \modu 4$.
Let as usual $u$ be a primitive root $\modu p$.
Then we get the $p-1$ congruences
\begin{equation}\label{e20021a}
\begin{split}
&\sum_{i=1}^{(p-1)/2} u_{2i} \times(t^{m_{i,k}}+t^{p-m_{i,k}})\equiv 0 \modu p,\\
&k=1,\dots,p-1,\quad u_{2i}\times m_{i,k}\equiv k \modu p, \quad 1\leq m_{i,k}\leq p-1.
\end{split}
\end{equation}
which gives  for  $k=1$,
\begin{equation}\label{e009121}
\sum_{i=1}^{(p-1)/2} u_{2i}\times (t^{u_{-2i}}+t^{p-u_{-2i}})\equiv 0 \modu p,
\end{equation}
\begin{proof}
Apply  theorem  \ref{t20021a} p.\pageref{t20021a}.
\end{proof}
\end{thm}
%
\begin{thm}\label{t20024a}{ *** On order of $t \modu p$}

Let $a\in\N$ be the order of $t \modu p$ (the smaller integer different of $0$ such that
$t^a\equiv 1\modu p$). Then, if first case of FLT fails for $p$, $a>\frac{p-1}{e_p+2}$.

\begin{proof}
Suppose that $a\leq \frac{p-1}{e_p+2}$ and search for a contradiction:
It is always possible to choose $f= e_p+2$ and $g_i,\quad i=1,\dots, f-1$ such that
\begin{displaymath}
\begin{split}
& g_i\equiv i^{-1} a^{-1}\modu p,\quad i=1,\dots,f-1,\\
& k_1=1.
\end{split}
\end{displaymath}
Then, we have, from hypothesis on $a$ applied to  theorem \ref{t20021a} p.\pageref{t20021a},
\begin{displaymath}
\begin{split}
& m_{i 1}= i\times a,\quad i=1,\dots,f-1, \\
& 0 < m_{i 1} < p, \\
& t^{m_{i 1}}=t^{i\times a}\equiv 1 \modu p.
\end{split}
\end{displaymath}

It is also possible to choose
$g_f$ and  $g_f\times m_{f 1}\equiv 1 \modu p$, so   such that $t^{m_{f 1}}+t^{p-m_{f 1}}\not \equiv 1+t\modu p$;
then, from theorem \ref{t20021a} p.\pageref{t20021a}, we get
\begin{displaymath}
\sum_{i=1}^f l_i g_i \times(1+t)+ l_f g_f \times(t^{m_{f 1}}+t^{p-m_{f 1}})\equiv 0\modu p.
\end{displaymath}
But, from lemma \ref{p7} p.\pageref{p7},  we have $\sum_{i=1}^f l_i g_i\equiv 0 \modu p$ which would lead to $t^{m_{f 1}}+t^{p-m_{f 1}} \equiv 1+t\modu p$, contradiction.
\end{proof}
\end{thm}
%
\subsection{A partial proof of FLT first case with symmetric sums of Mirimanoff polynomials}\label{s110071}
Recall some results seen previously
\begin{itemize}
\item
Let $x, y, z\in\Z-\{0\}$ verifying Fermat's equation $x^p+y^p+z^p=0$.
Let, without loss of generality, $t=-\frac{x}{y}\modu p,\quad t\in\N,\quad 1< t \leq p-1$.
\item
Let $T$ be an indeterminate. Let $\phi_m(T),\quad m=1,\dots, p-1$, be the
Mirimanoff polynomials:
\begin{equation}\label{e110071}
\phi_m(T)=\sum_{m=1}^{p-1} i^{m-1}\times T^i.
\end{equation}
\item
The Mirimanoff's congruences used here are
\begin{displaymath}
\phi_m(t)\times \phi_{p-m}(t)\equiv 0\modu p,\quad m=1,\dots, p-1.
\end{displaymath}
\item
Let us define $S(T)=\sum_{m=1}^{p-1} \phi_m(T)\times\phi_{p-m}(T)$.
\end{itemize}
%

\begin{lem}\label{l110071}
\begin{displaymath}
S(T)\equiv -\phi_{p-1}(T^2)\modu p.
\end{displaymath}
\begin{proof}
From previous definitions, we get
\begin{displaymath}
S(T)=\sum_{m=1}^{p-1}(\sum_{i=1}^{p-1} i^{m-1}\times T^i)\times
(\sum_{j=1}^{p-1} j^{p-m-1}\times T^j).
\end{displaymath}
With a reordering of finite sums, we get
\begin{displaymath}
S(T)=\sum_{m=1}^{p-1}\sum_{i=1}^{p-1}\sum_{j=1}^{p-1} i^{m-1}\times
j^{p-m-1}\times T^{i+j}.
\end{displaymath}
We collect the terms with $i=j$ to obtain
\begin{displaymath}
 S(T)=\sum_{m=1}^{p-1}\sum_{i=1}^{p-1}\sum_{j=1,\quad j\not=i}^{p-1} i^{m-1}\times
j^{p-m-1}\times T^{i+j} +\sum_{m=1}^{p-1}\sum_{i=1}^{p-1} i^{m-1}\times
i^{p-m-1}\times T^{2i}.
\end{displaymath}
But
\begin{displaymath}
\begin{split}
& \sum_{m=1}^{p-1}\sum_{i=1}^{p-1}\sum_{j=1,\quad j\not=i}^{p-1} i^{m-1}\times
j^{p-m-1}\times T^{i+j}=\\
& \sum_{i=1}^{p-1}\sum_{j=1,\quad j\not= i}^{p-1}\sum_{m=1}^{p-1} i^{m-1}\times
j^{p-m-1}\times T^{i+j}.
\end{split}
\end{displaymath}
But, from $i\not=j$,
\begin{displaymath}
\begin{split}
& \sum_{m=1}^{p-1} i^{m-1}\times
j^{p-m-1}\times T^{i+j}= j^{p-2}\times (\sum_{m=1}^{p-1} i^{m-1}\times
j^{-(m-1)})\times T^{i+j}=\\
&\sum_{m=1}^{p-1} i^{m-1}\times
j^{p-m-1}\times T^{i+j}= j^{p-2}\times (\frac{i^{p-1}\times j^{-(p-1)}-1}
{i\times j^{-1}-1})\times T^{i+j} \equiv 0\modu p.
\end{split}
\end{displaymath}
Therefore, we get
\begin{displaymath}
S(T)\equiv \sum_{m=1}^{p-1}\sum_{i=1}^{p-1} i^{m-1}\times
i^{p-m-1}\times T^{2i} \modu p,
\end{displaymath}
also
\begin{displaymath}
S(T)\equiv \sum_{m=1}^{p-1}\sum_{i=1}^{p-1} i^{p-2} \times T^{2i} \modu p,
\end{displaymath}
and finally
\begin{displaymath}
S(T)\equiv (p-1)\times \phi_{p-1}(T^2),
\end{displaymath}
which achieves the proof.
\end{proof}
\end{lem}
%
%
\begin{lem}\label{l110071}
Let $S_2(T)=\sum_{n=1}^{(p-1)/2}\phi_{2n}(T)\times\phi_{p-2n}(T)$. Then
\begin{displaymath}
S_2(T)\equiv \frac{p-1}{2}\times\phi_{p-1}(T^2)\modu p.
\end{displaymath}
\begin{proof}
From previous definitions, we get
\begin{displaymath}
S_2(T)=\sum_{n=1}^{(p-1)/2}(\sum_{i=1}^{p-1} i^{2n-1}\times T^i)\times
(\sum_{j=1}^{p-1} j^{p-2n-1}\times T^j).
\end{displaymath}
With a reordering of finite sums, we get
\begin{displaymath}
S_2(T)=\sum_{n=1}^{(p-1)/2}\sum_{i=1}^{p-1}\sum_{j=1}^{p-1} i^{2n-1}\times
j^{p-2n-1}\times T^{i+j}.
\end{displaymath}
We collect the terms with $i^2=j^2$ to obtain
\begin{displaymath}
\begin{split}
& S_2(T)=\sum_{n=1}^{(p-1)/2}\sum_{i=1}^{p-1}\sum_{j=1,\quad j^2\not\equiv i^2 (p)}^{p-1} i^{2n-1}\times
j^{p-2n-1}\times T^{i+j}\\
& +\sum_{n=1}^{(p-1)/2}\sum_{i=1}^{p-1} i^{2n-1}\times i^{p-2n-1}\times T^{2i}\\
& +\sum_{n=1}^{(p-1)/2}\sum_{i=1}^{p-1} i^{2n-1}\times (-i)^{p-2n-1}\times T^p.
\end{split}
\end{displaymath}
But
\begin{displaymath}
\begin{split}
& \sum_{n=1}^{(p-1)/2}\sum_{i=1}^{p-1}\sum_{j=1,\quad j^2\not=i^2}^{p-1} i^{2n-1}\times
j^{p-2n-1}\times T^{i+j}=\\
& \sum_{i=1}^{p-1}\sum_{j=1,\quad j^2\not= i^2}^{p-1}\sum_{n=1}^{(p-1)/2} i^{2n-1}\times
j^{p-2n-1}\times T^{i+j}.
\end{split}
\end{displaymath}
But, from $i^2\not\equiv j^2\modu p$,
\begin{displaymath}
\sum_{n=1}^{(p-1)/2} i^{2n-1}\times
j^{p-2n-1}\times T^{i+j}
= j^{p-2}\times (\sum_{n=1}^{(p-1)/2} i^{2n-1}\times j^{-(2n-1)})\times T^{i+j}.
\end{displaymath}
We have, from $i^2\not\equiv j^2\modu p$,
\begin{displaymath}
\begin{split}
&\sum_{n=1}^{(p-1)/2} i^{2n-1}j^{-(2n-1)}+\sum_{n=1}^{(p-1)/2} i^{2n}j^{-2n}\equiv 0\modu p,\\
&\sum_{n=1}^{(p-1)/2} i^{2n}j^{-2n}\equiv 0\modu p,
\end{split}
\end{displaymath}
Therefore $\sum_{n=1}^{(p-1)/2}i^{2n-1}j^{-(2n-1)}\equiv 0\modu p,$ so for $i^2\not\equiv 0\modu p$,
\begin{displaymath}
\sum_{n=1}^{(p-1)/2} i^{2n-1}\times
j^{p-2n-1}\times T^{i+j}\equiv 0\modu p.
\end{displaymath}
Therefore, we get
\begin{displaymath}
S_2(T)\equiv \sum_{n=1}^{(p-1)/2}\sum_{i=1}^{p-1} i^{p-2}\times T^{2i}
+\sum_{n=1}^{(p-1)/2}\sum_{i=1}^{p-1} i^{p-2} T^p \modu p,
\end{displaymath}
also
\begin{displaymath}
S_2(T)\equiv \sum_{n=1}^{(p-1)/2}\sum_{i=1}^{p-1} i^{p-2} \times T^{2i} \modu p,
\end{displaymath}
and finally
\begin{displaymath}
S_2(T)\equiv \frac{p-1}{2}\times \phi_{p-1}(T^2),
\end{displaymath}
which achieves the proof.
\end{proof}
\end{lem}
%
\begin{lem}\label{l110072}
Let $x^p+y^p+z^p=0, \quad x\times y\times z\not\equiv 0\modu p$, be the Fermat equation, first case.
\begin{enumerate}
\item
$\phi_{p-1}(-t)\equiv\phi_{p-1}(t)\equiv\phi_{p-1}(t^2)\equiv 0\modu p.$
\item
If $x\equiv y\modu p$ then $x\equiv y\modu p^2$.
\item
If $x\not\equiv y\modu p$ then $(x-y)^{p-1}\equiv 1\modu p^2$.
\end{enumerate}
\begin{proof}
From lemma \ref{l110071} p. \pageref{l110071} and from Mirimanoff congruences
\ref{e110071} p. \pageref{e110071}, we derive that $\phi_{p-1}(t^2)\equiv 0\modu p$.
From Ribenboim \cite{rib} relation (1.23) p 144,  we get
\begin{displaymath}
\phi_{p-1}(t^2)\equiv \frac{1-t^{2p}-(1-t^2)^p}{p}\equiv 0\modu p,
\end{displaymath}
so
\begin{displaymath}
1-t^{2p}-(1-t^2)^p\equiv 0\modu p^2,
\end{displaymath}
and so
\begin{displaymath}
(1-t^p)\times (1+t^p)-(1-t)^p\times (1+t)^p\equiv 0\modu p^2.
\end{displaymath}
From Ribenboim \cite{rib} p 145 we get $\phi_{p-1}(t)\equiv 0\modu p$, so
\begin{displaymath}
\frac{1-t^p-(1-t)^p}{p}\equiv 0\modu p,
\end{displaymath}
so
\begin{displaymath}
1-t^p-(1-t)^p\equiv 0\modu p^2.
\end{displaymath}
From these relations, we get
\begin{displaymath}
1+t^p-(1+t)^p\equiv 0\modu p^2,
\end{displaymath}
because $1-t=\frac{x+y}{y}\not\equiv 0\modu p$.
Thus
\begin{displaymath}
\frac{(1+t^p)-(1+t)^p}{p}\equiv\phi_{p-1}(-t)\equiv 0\modu p.
\end{displaymath}
Also
\begin{displaymath}
1-(\frac{x}{y})^p-(1-\frac{x}{y})^p\equiv 0\modu p^2,
\end{displaymath}
and so
\begin{displaymath}
y^p-x^p\equiv (y-x)^p\modu p^2.
\end{displaymath}
But, from $x^p\equiv x\modu  p^2$ and $y^p\equiv y\modu p^2$, classical Fermat result, we get
\begin{displaymath}
y-x\equiv (y-x)^p\modu p^2,
\end{displaymath}
also
\begin{displaymath}
(y-x)\times( 1- (y-x)^{p-1})\equiv 0\modu p^2.
\end{displaymath}
Therefore $y-x\equiv 0\modu p$ implies that $y-x\equiv 0\modu p^2$ and $y-x\not\equiv 0\modu p$ implies that $(y-x)^{p-1}-1\equiv 0\modu p^2$, which achieves the proof.
\end{proof}
\end{lem}
%
{\bf Remark:} observe that the assertion $y-x\not\equiv 0\modu p \Rightarrow (y-x)^{p-1}-1\equiv 0\modu p^2$ can be directly derived of second Futwangler theorem,
see for instance Ribenboim, \cite{rib} (3C) p. 169.
\begin{thm} {*** }{A particular case of Fermat-Wiles}\label{t110071}

Let $x,y,z\in \Z-\{0\}$, pairwise coprime.
Then $p\|(x-y)\times (y-z)\times (z-x)$ implies that
$x^p+y^p+z^p\not=0$.
\begin{proof}$ $
\begin{itemize}
\item
Suppose that second case holds and search for a contradiction:
If $p\|(x-y)$ then $x\times y\not\equiv 0\modu p$, so $p|z$, and also $p|2x^p\Rightarrow p|x$, which is not possible because $gcd(x,z)=1$.
\item
If $ x\times y\times z\not\equiv 0\modu p$, then the theorem is an immediate consequence of lemma \ref{l110072} p. \pageref{l110072}.
\end{itemize}
\end{proof}
\end{thm}
%

%
\subsection{ Some elementary congruences derived of Mirimanoff's congruences- II}
Recall that in exponential congruences, we denote $a^{b (p)}$ and not
$a^{b \modu p}$.

In this subsection, starting from the Mirimanoff's congruences
$\phi_m(t)\times\phi_{p-m}(t)\equiv 0\modu p,\quad m=1,\dots p-1$, we prove the set of elementary explicit congruences $\modu p$:
%
\begin{thm} { *** }\label{t110181}

$t=-\frac{x}{y}\modu p$ verifies the strong set of $p-1$ congruences
\begin{displaymath}
\sum_{j=1}^{p-1} j^{p-2}\times t^{k j (p)+j}\equiv 0\modu p,\quad k=1,\dots,p-1.
\end{displaymath}
\begin{proof}
We have, from Mirimanoff's congruences, see Ribenboim \cite{rib} (1B) p 145,
\begin{displaymath}
\phi_m(t)\times\phi_{p-m}(t)\equiv 0\modu p,\quad m=1,\dots,p-1,
\end{displaymath}
also, expliciting Mirimanoff's polynomials,
\begin{displaymath}
(\sum_{i=1}^{p-1} i^{m-1} t^i)\times (\sum_{j=1}^{p-1} j^{p-m-1} t^j)
\equiv 0\modu p,\quad m=1,\dots,p-1.
\end{displaymath}
This relation leads to
\begin{displaymath}
\sum_{i=1}^{p-1}\sum_{j=1}^{p-1} (i j^{-1})^{m-1}\times j^{-1}\times t^{i+j},
\quad m=1,\dots,p-1.
\end{displaymath}
Let $i j^{-1}\equiv k\modu p,\quad 1\leq k\leq p-1$, so with $i\equiv k j\modu p$.
We get
\begin{displaymath}
\sum_{k=1}^{p-1} k^{m-1} \{ \sum_{j=1}^{p-1} j^{-1}\times t^{k j (p)+j} \}\equiv 0\modu p.
\end{displaymath}
Let $X_k\in\N,\quad X_k\equiv \sum_{k=1}^{p-1} j^{-1}\times t^{ k j (p)+j} \modu p$.
We obtain the set of $p-1$ linear congruences of the $p-1$ indeterminates
$X_k,\quad k=1,\dots,p-1$,
\begin{displaymath}
\sum_{k=1}^{p-1} k^{m-1}\times X_k\equiv 0\modu p,\quad m=1,\dots,p-1.
\end{displaymath}
This system has only the trivial solution $X_k=0,\quad k=1,\dots,p-1$, because
the Vandermonde determinant
\begin{displaymath}
| k^{m-1}|_{k=1,\dots, p-1,\quad m=1,\dots,p-1}\not\equiv 0\modu p,
\end{displaymath}
which achieves the proof.
\end{proof}
\end{thm}
%

\begin{cor}\label{c110181}
$\phi_{p-2}(t^2)\equiv 0\modu p.$
\begin{proof}
Immediate consequence of theorem \ref{t110181} p. \pageref{t110181}.
\end{proof}
\end{cor}
%

{\bf Remark:}
In Numerical MAPLE investigations, we have  noticed that, for the indeterminate
$T$,
\begin{displaymath}
\begin{split}
&\sum_{j=1}^{p-1}j^{m-1}\times T^{kj\  (p)}\equiv c_{k,m}\times\phi_m(T) \modu p,
\quad m=1,\dots,p-2,\quad k=1,\dots,p-1,\\
& c_{k,m}\in{\bf F}_p[m,k].
\end{split}
\end{displaymath}
For instance
\begin{displaymath}
\sum_{j=1}^{p-1}j^{p-2}\times T^{kj\  (p)}\equiv (p-k)\times\phi_{p-1}(T) \modu p,
\quad k=1,\dots,p-1.
\end{displaymath}

%

\subsection { Numerical MAPLE applications and a conjecture implying first case of FLT}\label{s010032}
The theorem \ref{t21021a} p.\pageref{t210021a} gives  $N_f= (C_{p-1}^{f})^2$ polynomial exponential congruences $\Delta_{f}(t)\equiv 0 \modu p$ where the  polynomials  $\Delta_{f}(T)\in \Z[T]$ are explicitly  computable   on $I_f=\{g_1,\dots,g_f\}$.
In our opinion, these congruences are a strong tool for a {\bf full} proof of first case of FLT (try to combine some of them to give a contradiction). We {\bf conjecture} that these  congruences are never all simultaneously possible.
Towards this direction, let for the prime $p$, the polynomials $\Delta_l(T) \modu p,\quad l=1,\dots,N_f$. We give a verbatim MAPLE list of  a program that factorize
$Gcd(\Delta_l(T) \modu p,\quad l=1,\dots, N_f)$ in the finite field ${\bf F}_p[T]$. We have in fact in this example computed the Gcd  of a subset of this set of $N_f$ polynomials for $f=[\frac{p}{3}]$, (observe that for a lot of $p$-cyclotomic fields $\Q(\zeta)$, and even probably all, we have $e_p<[\frac{p}{3}]$, so it is possible to take $f=[\frac{p}{3}]$.
We have defined $\Delta_l(T) \modu p$ by :
\begin{displaymath}
\begin{split}
& f= [\frac{p}{3}], \\
& l_{max}= [\frac{p}{2}]-f, \\
& g_i =i,\quad i=1,\dots,f, \\
& k_j =j+l,\quad j=1,\dots,f, \\
& m_{ij}\equiv g_i^{-1}\times k_j \modu p.
\end{split}
\end{displaymath}
The program computes $Gcd(\Delta_l(T) \modu p,\quad l=1,\dots,l_{max})$.
We have found for $16 < p < 70$, that $Gcd(\Delta_l(T) \modu p = (T+1)\times (T+p-1)$.
We have not pursued because our CPU computer is slow to compute determinants of rank $>70$.
Nevertheless, the numerical evidence that we observe in these computations lead us to :
%


{\bf Conjecture}
\begin{displaymath}
Gcd(\Delta_l(T)) \modu p,\quad l=1,\dots,N_f)= (T+1)\times (T+p-1).
\end{displaymath}
\begin{itemize}
\item
Observe that in all examples computed for $16< p < 70$, it was enough to take the gcd of {\bf two} polynomials $\Delta_1(T)\mod p$ and $\Delta_2(T) \modu p$, which reinforces strongly our congruence.
\item
Observe that this conjecture and theorem (\ref{t21021a})  implies First case of FLT, because $(t+1)\times (t+p-1)\not\equiv 0 \modu p$; it is the same proof than in subsection \ref{s010031} p.\pageref{s010031}.
\end{itemize}

\begin{verbatim}
# This MAPLE program  gives a numerical "obviousness"   that  # theorem 9.4 implies FLT
> #-----------------------
> #  the m_ij g_i k_j theta_ij are those of Theorem 9.4
> # Delta_k(T) mod p, k=1,...,N_f
> restart:
> with(linalg);
> p_inf:=7:
> p_sup:=50:
> p:=p_inf:
> while p < p_sup do
>   p:=nextprime(p):
>   f:=trunc(p/3):
>   lm:=trunc(p/2):
>   print (`p =`,p,` f = `,f,` lm = `,lm):
>   print(`----------`):


>   T_T:=expand((T+1)*(T+p-1)) mod p:
>   r:=array(1..lm):
>   k_j:=array(1..f):
>   g_i:=array(1..f):
>   Delta:=array(1..f,1..f):
>   for i from 1 to f do
>     for j from 1 to f do
>       delta[i,j]:=0:
>     od:
>   od:
>   #----------------------------------
>   for l from 1 to lm do
>     ind_r_gcd:=0:
>     for i from 1 to f do
>       g_i[i]:=i&^(p-2) mod p:
>     od:
>     for j from 1 to f do
>       k_j[j]:= j+l mod p:
>     od:
>     for i from 1 to f do
>       for j from 1 to f do
>         tem_i:= g_i[i]&^(p-2) mod p:
>         m_ij:= k_j[j]*tem_i mod p:
>         theta_ij:=T^(m_ij)+T^(p-m_ij):
>         Delta[i,j]:=theta_ij:
>       od:
>     od:
>     #-------------------------
>     #print(Delta):
>     v_Delta:=Det(Delta) mod p:
>     v_Delta:=expand(v_Delta) mod p:
>     #print(`p =`,p,` l =`,l,` v_Delta =`,v_Delta):
>     #-------------------------
>     if v_Delta = 0 then
>       print(`p =`,p,` l =`,l,` v_Delta =`,v_Delta):
>     fi:
>     #----------------------------
>     if v_Delta <> 0 then
>       b_T:=(T^(p-1)-1) mod p :
>       r[l]:=Rem(v_Delta,b_T,T) mod p:
>       r[l]:=expand(r[l]) mod p:
>       r[l]:=Gcd(r[l],b_T) mod p:
>       #--------------------------
>       if l = 1 then
>         r_gcd:= r[1] mod p:
>         print(`p =`,p,` l = `,l,` r[1] =`,Factors(r[l]) mod p):
>         print(`p =`,p,` l = `,l,` r_gcd =`,Factors(r_gcd) mod p):
>       fi:
>       #--------------------------
>       if l > 1 then
>         r_gcd:=Gcd(r_gcd,r[l]) mod p:
>         print(`p =`,p,` l = `,l,` r[l] =`,Factors(r[l]) mod p):
>         print(`p =`,p,` l = `,l,` r_gcd =`,Factors(r_gcd) mod p):
>         if r_gcd = T_T then
>           print(` ***** criterion ok **** `,\
>           `p =`,p,` l = `,l,` r_gcd = `,\
>           Factors(r_gcd) mod p):
>           break:
>         fi:
>       fi:
>     fi:
>     #---------------------------
>   od:
>   if r_gcd <> T_T then
>     print(` ***** criterion not ok **** `,\
>     `p =`,p, ` r_gcd = `,\
>     Factors(r_gcd) mod p):
>   fi:
> print(`#############################`):
> od:
\end{verbatim}
%


%

\clearpage
\section {FLT first case: Representations of $Gal(\Q(\zeta)/\Q)$ approach}\label{s20113}
We give as applications :
\begin{enumerate}
\item
A connection between Mirimanoff's polynomials and representations
\item
A connection between relative class group $C_p^-$ structure and congruences on Mirimanoff's polynomials
\item
A relation between Fermat's congruence $p$-rank $r$ and $p$-rank $r_1$ of the action of $G$ on the ideal $\s,\quad \s^p=(x+\zeta y)\Z[\zeta]$.
\item
An application of Eichler's theorem in intermediate fields
$\Q\subset K\subset \Q(\zeta)$.
\end{enumerate}
The notations and results of previous sections are used in this subsection.
\begin{itemize}
\item
Recall that the ideal $\s$ is such that $\s^p=(x+\zeta y)\Z[\zeta]$.
\item
In proposition \ref{p4} p.\pageref{p4} , we proved that there exists
$L_f=\{l_i\quad |\quad i=1,\dots,f\}$,
 such that
the ideal $\prod_{i=0}^{f-1} \sigma^i(\s)^{l_i}$ is principal.
\item
We shall translate this result , {\it mutatis mutandis,} to the principal ideal
$\prod_{i=0}^{r_1} \sigma^{i} (\s)^{(-1)^{r_1-i}\times S_{r_1-i}(1)}$

found in relation (\ref{e04112}) p.\pageref{e04112}.
\item
With the language of representations, the polynomial
$P_{r_1}(U)=\prod_{i=1}^{r_1}(U-\mu_i)$ on the indeterminate $U$ is  the minimal polynomial such that $P_{r_1}(\sigma)$ annihilates the class of ideal $\s$, also
denoted $\s^{P_{r_1}(\sigma)}\simeq\Z[\zeta]$.
\item
The rank $r_1$ has two components $r_1^-$ and $r_1^+$ corresponding to ranks of subgroups of the $p$-class groups $C_p^-$ and $C_p^+$ of $\Q(\zeta)$ with $r_1=r_1^-+r_1^+$.
\item
Recall that, for $r$ defined in relation  (\ref{e31102}) p.\pageref{e31102}, we have
\begin{displaymath}
\begin{split}
& r-1 = Card(\{\phi_{2m+1}^*(t)\not \equiv 0 \modu p\quad
|\quad 0\leq m\leq \frac{p-3}{2}\}),\\
& r-2 = Card(\{\phi_{2m+1}(t)\not \equiv 0 \modu p\quad
|\quad 1\leq m\leq \frac{p-3}{2}\}).
\end{split}
\end{displaymath}
\item
$I_\phi=\{2m+1\quad |\quad 0\leq m\leq \frac{p-3}{2},
\quad \phi_{2m+1}^*(t)\not\equiv 0 \modu p\}$
with $Card(I_\phi)=r-1$.
\item
$M_{r_1}= \{ \mu_i\quad |\quad i=1,\dots, r_1\}$,
with the meaning of $\mu_i$ of theorem (\ref{t31101a})
\item
With the language of representations,
Let $P_{r}(U)$ be the minimal polynomial of the indeterminate $U$ such that
\begin{displaymath}
(\frac{x+\zeta y}{x+\zeta^{-1}y})^{P_{r}(\sigma)}\equiv 1\modu \pi^{p+1},
\quad P_{r}(U)=\prod_{j=1}^{r-1}(U-\mu_j).
\end{displaymath}
equivalent to
\begin{displaymath}
\prod_{i=0}^{r-1}(\frac{x+\zeta^{u_i} y}{x+\zeta^{-u_i} y})^{n_i}
\equiv 1\modu \pi^{p+1},
\quad n_0\not\equiv 0\modu p,\quad n_{r-1}\not\equiv 0\modu p.
\end{displaymath}
\end{itemize}
%
\subsection{Connection between  Mirimanoff's polynomials and representations} \label{s011191}
\begin{thm} \label{t03111}
Suppose that first case of FLT fails for $p$.
Then, for each $2 m_j+1\in I_\phi-1,\quad j=1,\dots,r-2$, there exists $i_j, \quad 1\leq i_j\leq r_1$, with
\begin{equation}\label{e05121}
u_{2m_j+1} = \mu_{i_j},
\end{equation}
where $\mu_i,\quad i=1,\dots,r_1,$ are those of representation theorem (\ref{t31101a}).
\begin{proof} $ $
\begin{itemize}
\item
From relation (\ref{e04112}) p.\pageref{e04112} with $d=1$, the ideal
$\prod_{i=0}^{r_1} \sigma^{i}(\s)^{(-1)^{r_1-i} \times S_{r_1-i}(1)}$ is a principal ideal; in relation (\ref{e5}) p.\pageref{e5}, we considered  the principal ideal
$\prod_{i=0}^{f-1} \sigma^i(\s)^{l_i}$;
 the relation  (\ref{e5}) p. \pageref{e5} is replaced
by the relation (\ref{e04112}) p.\pageref{e04112},  where $\phi_{2m+1}(t)\not\equiv 0 \modu p$ exactly for the {\bf same} values of $2m+1$ depending only on $p$ and $t$.
\item
Then
\begin{equation}\label{e104011}
\begin{split}
&\prod_{i=0}^{r_1} (x+\zeta^{u_i}y)^{(-1)^{r_1-i} \times S_{r_1-i}(1)}= \zeta^v\times\eta\times\gamma^p,\\
& v\in \N,\quad 0\leq v\leq p-1,\quad \eta\in\Z[\zeta+\zeta^{-1}]^*,\quad \gamma\in\Z[\zeta].
\end{split}
\end{equation}
where, by comparison, in relation (\ref{e5}) p.\pageref{e5}, we had
\begin{displaymath}
\prod_{i=0}^{f-1} (x+\zeta^{u_i} y)^{l_i}=\eta_2\times\gamma^p,\quad
\eta_2\in\Z[\zeta+\zeta^{-1}]^*,\quad \gamma\in\Z[\zeta],
\end{displaymath}
(with here $\zeta^{v_2}=1$).
\item
Then, from foundation theorem \ref{p21a} p.\pageref{p21a} ,
we obtain the polynomial congruences:
\begin{displaymath}
\phi_{2m+1}(t)\times\sum_{i=0}^{r_1} (-1)^{r_1-i} \times S_{r_1-i}(1) \times u^{(2m+1)\times i} \equiv 0\modu p,\quad
m=1,\dots,\frac{p-3}{2},
\end{displaymath}
and so
\begin{equation}\label{e20121}
\sum_{i=0}^{r_1} (-1)^{r_1-i} \times S_{r_1-i}(1) \times u^{(2m_j+1)\times i} \equiv 0\modu p,\quad
j=1,\dots,r-2.
\end{equation}
\item
The fact that $v$ can be here different of zero had induced to disgard the case $m=0$ (or $2m+1=1$), corresponding to $\phi_{2m+1}^*(t)=1$, in theorem \ref{p21a} p.\pageref{p21a}, so to disgard the congruence corresponding to $m=0$ involving the sum
\begin{displaymath}
\sum_{i=0}^{r_1} (-1)^{r_1-i} \times S_{r_1-i}(1) \times u^i \modu p.
\end{displaymath}
\item
Then the  relation (\ref{e20121}) p.\pageref{e20121}  leads algebraically to
\begin{displaymath}
\prod_{i=1}^{r_1} (u^{2m_j+1 }-\mu_i)\equiv 0 \modu p,
\end{displaymath}
so there exists $i_j,\quad 1\leq i_j\leq r_1$ such that
$u_{2m_j+1 }-\mu_{i_j}\equiv 0 \modu p$, which achieves the proof.
\end{itemize}
\end{proof}
\end{thm}
%
Another point of view for a similar result, the representations approach : let us consider the Fermat's minimal multiplicative congruence
\begin{displaymath}
\prod_{i=0}^{r-1} (\frac{x+\zeta^{u_i}y}{x+\zeta^{-u_i}y})^{n_i}\equiv 1\modu\pi^{p+1}.
\end{displaymath}
This can be written, in representations notations,
with $P_{r}(\sigma)=\prod_{j=1}^{r-1}(\sigma-\mu_j)\in{\bf F}_p[G]$,
\begin{equation}\label{e103051}
(\frac{x+\zeta y}{x+\zeta^{-1}y})^{P_r(\sigma)}\equiv 1\modu \pi^{p+1}.
\end{equation}
%
\begin{lem} \label{t03111}
Let us note $\rho=r_1^-$.
The polynomials $P_{\rho}(U)$ and $P_r(U)$ of the indeterminate $U$ verify the relation
\begin{displaymath}
P_r(U) \ | \ P_{\rho}(U)\times (U-u).
\end{displaymath}
\begin{proof} $ $
The ideal $(\frac{\s}{\overline{\s}})^{P_{\rho}(\sigma)}$ is principal. In terms of congruences, it induces that
\begin{equation}\label{e103052}
(\frac{x+\zeta y}{x+\zeta^{-1} y})^{P_{\rho}(\sigma)}\equiv \zeta^{v_1}\modu\pi^{p+1}, \quad v_1\in\N.
\end{equation}
Observe that $U-u$ does not divides $P_{\rho}(U)$ from lemma \ref{l103021} p.\pageref{l103021} derived of Stickelberger theorem; this implies that $v_1\not\equiv 0\modu p$ in  this relation.
We derive that
\begin{equation}\label{e109141}
(\frac{x+\zeta y}{x+\zeta^{-1} y})^{P_{\rho}(\sigma)\times (\sigma-u)}
\equiv 1\modu\pi^{p+1}.
\end{equation}
From the definition of $r$ and from the relations (\ref{e103051}) p.\pageref{e103051} and (\ref{e109141}) p.\pageref{e109141}, we get $r-1\leq \rho+1$.
Observe at first that  definition of $P_r(U)$ implies that $(U-u)\ | \ P_r(U)$. Let us make the euclidean division of $P_{\rho}(U)\times (U-u)$
by $P_r(U)$ in ${\bf F}_p[U]$:
\begin{displaymath}
P_{\rho}(U)\times (U-u)=P_r(U)\times B(U)+R(U),\quad B(U), R(U)\in{\bf F}_p[U],
\quad deg(R(U)) <deg(P_r(U)).
\end{displaymath}
The remainder $R(U)$ cannot be different of zero : if not, we should have
\begin{displaymath}
(\frac{x+\zeta y}{x+\zeta^{-1} y})^{R(\sigma)}\equiv 1\modu \pi^{p+1},
\end{displaymath}
which contradicts the minimality of $P_r(U)$.
Therefore $P_r(U)$ divides $P_{\rho}(U)\times (U-u)
=(\prod_{i=1}^{\rho} (U-\mu_i))\times (U-u)$.
\end{proof}
\end{lem}
%

\subsection {On rank $r$ of Fermat congruences $\modu p$ }
Let us note $G$ the Galois group of the field $\Q(\zeta)$. Recall that the ideal
$\s$ is defined by $\s^p=(x+\zeta y)\Z[\zeta]$. We show that, when $h^+\not\equiv 0\modu p$, there is an explicit relation between:
\begin{enumerate}
\item
the $p$-rank $r$ of FLT congruences $\modu \pi^{p+1}$, previously defined in subsection \ref{s31101} p.\pageref{s31101},
\item
the $p$-rank $r_1$ of the ${\bf F}_p[G]$ module resulting of  $G$-action  on the subgroup
$<Cl(\s)>$ of $C_p$, defined in subsection \ref{s20111} p.\pageref{s20111} of this chapter.
\end{enumerate}
This subsection appears to be a {\bf strong} gateway between the  $\modu \pi^{p+1}$ congruences approach of chapters \ref{s3} p.\pageref{s3}, \ref{s4} p.\pageref{s4}, and representation approach of this chapter \ref{s09121} p.\pageref{s09121}.
%
\begin{thm} { *** }\label{t011241}
Let $r$ bet the FLT congruences $p$-rank defined in  relation (\ref{e31102}) p.\pageref{e31102}.
Let $r_1$ be the $p$-rank of the ${\bf F}_p[G]$-module resulting of $G$-action on the subgroup
$<Cl(\s)>$ of $C_p$, defined in  relation (\ref{e011251}) p.\pageref{e011251}.
Assume that $p$ does not divide $h^+$.
Then the two rank verify the relation
\begin{equation}\label{e011241}
r_1= r-2.
\end{equation}
\begin{proof}
We have shown in relation (\ref{e31102}) p.\pageref{e31102} that
\begin{displaymath}
\sum_{i=1}^r \frac{ u_i l_i}{x+\zeta^{u_i} y}
+\sum_{i=1}^r \frac{u_i l_i}{x+\zeta^{-u_i} y} \equiv 0 \modu \pi^{p-1}.
\end{displaymath}
We have shown, from lemma \ref{l28053a} p.\pageref{l28053a} that this sum kind of  congruence induces a product kind of  congruence
\begin{displaymath}
\prod_{i=1}^r (x+\zeta^{u_i} y)^{l_i}
-\prod_{i=1}^r (x+\zeta^{-u_i} y)^{l_i}\equiv 0 \modu \pi^{p-1}.
\end{displaymath}
This congruence implies the polynomial identity, for the indeterminate $X$,
\begin{displaymath}
\begin{split}
&\prod_{i=1}^r (x+ X^{u_i} y)^{l_i}
-\prod_{i=1}^r (x+X^{p-u_i} y)^{l_i} = p\times f(X)+(X^p-1) \times g(X),\\
& g(X)\in\Z[X],
\end{split}
\end{displaymath}
which annihilates for $X=1$, and so gives
\begin{displaymath}
\begin{split}
& \prod_{i=1}^r (x+ X^{u_i} y)^{l_i}
-\prod_{i=1}^r (x+X^{p-u_i} y)^{l_i} = p\times (X-1)\times  f_1(X)
+(X^p-1) \times g(X),\\
& f_1(X)\in\Z[X],
\end{split}
\end{displaymath}
and also
\begin{displaymath}
\prod_{i=1}^r (x+ \zeta^{u_i} y)^{l_i}
-\prod_{i=1}^r (x+\zeta^{-u_i} y)^{l_i} = p\times (\zeta-1)\times  f_1(\zeta).
\end{displaymath}
Let us denote
\begin{displaymath}
\varpi = \prod_{i=1}^r\frac{ (x+ \zeta^{u_i} y)^{l_i}}
{ (x+\zeta^{-u_i} y)^{l_i}}
\end{displaymath}
From these relations, we get immediately that
\begin{displaymath}
\begin{split}
&\varpi \in\Q(\zeta),\\
&\varpi \equiv 1 \modu \pi^p, \\
&\overline{\varpi}=\varpi^{-1}.
\end{split}
\end{displaymath}
From Washington, \cite{was}, lemma 9.1 p 169, we derive
that the extension $\Q(\zeta,\varpi^{1/p})$ is unramified at $\pi$ and, from \cite{was} exercice 9.1 p 182,  that it is unramified at the others primes. Then, from \cite{was} lemma 9.2 p 170, we derive that
$\varpi = \beta^p, \quad \beta\in\Q(\zeta).$
Therefore
\begin{displaymath}
\varpi = \prod_{i=1}^r\frac{ (x+ \zeta^{u_i} y)^{l_i}}
{ (x+\zeta^{-u_i} y)^{l_i}} = \beta^p,\quad \beta\in\Q(\zeta)
\end{displaymath}
and so
\begin{displaymath}
\prod_{i=1}^r\frac{ \sigma^{u_i}(\s)^{l_i}}
{ \sigma^{-u_i}(\s)^{l_i}} = \beta\Z[\zeta],\quad  \beta\in\Q[\zeta]
\end{displaymath}
also,   we derive of this relation
\begin{displaymath}
\frac{\prod_{i=1}^r \sigma^{u_i}(\s)^{2 l_i}}
{(\prod_{i=1}^r \sigma^{u_i}(\s)^{l_i})\times
(\prod_{i=1}^r \sigma^{-u_i}(\s)^{l_i})} = \beta\Z[\zeta],
\end{displaymath}
and finally, because $h^+\not\equiv 0 \modu p$,
\begin{displaymath}
\prod_{i=1}^r \sigma^{u_i}(\s)^{2 l_i}
 = \gamma\Z[\zeta],\quad \gamma \in\Z(\zeta).
\end{displaymath}
or
\begin{displaymath}
\prod_{i=1}^r \sigma^{u_i}(\s)^{ l_i}
 = \gamma_1\Z[\zeta],\quad \gamma_1 \in\Z(\zeta).
\end{displaymath}
Therefore, from definition of $r_1$, we deduce that $r\geq r_1+1$.
From relation (\ref{e012091}) p.\pageref{e012091}, we get $r_1\geq r-2$, and finally $r-1\geq r_1\geq r-2$.
Then we take the determination  $r_1=r-2$ because here $r_1=r_1^-$ and $r_1^-$ is the smaller natural integer such that
\begin{displaymath}
(\frac{x+\zeta y}{x+\zeta^{-1} y})^{P_{r_1^-}(\sigma)}\equiv \zeta^v\modu \pi^{p+1},
\quad v\in \N,
\end{displaymath}
and we have seen that $v\not\equiv 0\modu p$: the occurence of this extraneous $\zeta^v$ factor leads to chose the determination $r_1=r-2$ and not
$r_1=r-1$.
\end{proof}
\end{thm}
{\bf Remark: be carefull,} the $p$-rank $r_1^-$ is  the dimension  of the
${\bf F}_p[G]$ module generated by $G$-action on $< Cl(\s)>$. Thus, the smaller not trivial relation involves $r_1^-+1$ elements.  In an other part, what we call {\it minimal Fermat's congruence rank}  $r$ is the smallest  number of elements for a  $\modu \pi^{p-1}$ congruence, and not the greatest $r$ without not trivial congruences $\modu \pi^p$. Therefore a residual difference of $1$  between $r_1^-$ and $r-1$ is explained by these {\it terminology conventions} taken in this monograph.
%

%
\subsection { On Eichler's Theorem in intermediate fields }\label{s011192}
\begin{thm} \label{t1011181}
Let $d\in\N, \quad p-1\equiv 0 \modu d, \quad d < [\sqrt{p}]-1$. If first case of FLT would fail for $p$ then $r_d > \frac{[\sqrt{p}]-1}{d}$.
\begin{proof}
From  relation (\ref{e011181}) p.\pageref{e011181} we get
$r_d\times d\geq r_1$.
From Eichler's theorem \ref{t011182} p.\pageref{t011182}, we get $r_1 >[\sqrt{p}]-1$.
Therefore $r_d\times d\geq [\sqrt{p}]-1$ which leads to the result.
\end{proof}
\end{thm}
%
\clearpage
%
%
\clearpage
\section{FLT first case: Representations and  congruences $\phi_{m}(t)\equiv 0 \modu p$.}\label{s103112}
We have studied up to now the congruences $\phi_{2m+1}(t)\equiv 0\modu p,\quad
0\leq m\leq \frac{p-3}{2}$. In this section, we show that it is possible to give some important generalizations to congruences
$\phi_m(t)\equiv 0\modu p,\quad m=1,\dots, p-1$.
%
\subsection{General properties of congruences $\phi_m(t)\equiv 0\modu p$.}
\label{s104281}
\begin{itemize}
\item
Recall that $u$ is a primitive root $\modu p$.
Our approach, in proposition \ref{p4} p.\pageref{p4} and  theorem \ref{t2009} p.\pageref{t2009}, was to show that, for
$ f\in \N$, with
$\frac{p-1}{2}\geq f\geq min(e_p+2,\frac{p-1}{2})$, there exists
$L_{I_{f}}=\{l_i\quad |\quad l_i\in \N,\quad i=1,\dots,f \}$,
such that
\begin{equation}\label{e202026}
\begin{split}
&\prod_{i=1}^{f} (x+\zeta^{u_{i}}y)^{l_i}=\eta\times\gamma^p, \quad
\eta\in\Z[\zeta+\zeta^{-1}]^*,\quad \gamma\in\Z[\zeta],\\
& l_i\in\N,\quad 0\leq i\leq p-1,\quad l_0\not=0,\quad l_f\not=0,\\
&\prod_{i=1}^f (\frac{x+\zeta^{u_i} y}{x+\zeta^{-u_i}y})^{l_i}\equiv 1\modu \pi^{p+1},\\
& \prod_{i=1}^{f} (x+\zeta^{u_{i}}y)^{l_i} -
\prod_{i=1}^{f} (x+\zeta^{-u_{i}}y)^{l_i}\not=0.
\end{split}
\end{equation}
\item
Classically, Mirimanoff congruences assert that
\begin{displaymath}
\phi_{2m+1}(t)\times B_{p-1-2m}\equiv 0 \modu p.
\end{displaymath}
Observe that these congruences are only for {\bf odd} Mirimanoff Polynomials $\phi_{2m+1}(t)$ because the odd Bernoulli Number are trivially divisible by $p$. We shall show that, more generally,  we can obtain no trivial results also for  Mirimanoff Polynomials
$\phi_{m}(T),\quad m=1,\dots, p-1$.
\item
We can show that there exists
$f^\prime\in \N,\quad f< f^\prime\leq f +\frac{p-3}{2}$,
such that
\begin{equation}\label{e202025}
\begin{split}
& \prod_{i=1}^{f^\prime} (x+\zeta^{u_{ i}}y)^{n_i}\equiv c^p \modu \pi^{p+1},
\quad c\in\Z,\quad c\not\equiv 0\modu p,\\
& n_i\in\N,\quad 0\leq n_i\leq p-1,\quad n_0\not=0,\quad n_{f^\prime}\not=0,\\
& \prod_{i=1}^{f^\prime} (x+\zeta^{u_{i}}y)^{n_i} -
\prod_{i=1}^{f^\prime} (x+\zeta^{-u_{i}}y)^{n_i}\not = 0.
\end{split}
\end{equation}
To see these congruences, we apply  Dirichlet's boxes principle to unit $\eta\in\Z[\zeta+\zeta^{-1}]^*$
\begin{displaymath}
\eta=\pm\prod_{i=1}^{(p-3)/2}\eta_i^{\nu_i},\quad \nu_i\in\Z,
\end{displaymath}
where $\eta_i,\quad i=1,\dots,\frac{p-3}{2}$, are a system of fundamental units of $\Z[\zeta+\zeta^{-1}]^*$, noticing that
$\eta_i^{\nu_i}\equiv \eta_i^{\nu_i+\alpha_i p}\modu \pi^{p+1},\quad\alpha_i\in\Z$.
\item
Therefore, we get
\begin{displaymath}
\sum_{i=1}^{f^\prime} \frac{u_{i}\times n_i}{(x+\zeta^{u_{i}} y)}\equiv 0 \modu \pi^{p-1}
\end{displaymath}
\item
Then there exists a {\bf minimal} Fermat congruence rank
$r^\prime\leq \frac{p-3}{2}+e_p+2$, such that
\begin{equation}\label{e04031a}
\sum_{i=1}^{r^\prime} \frac{u_{i} n^\prime_i}{(x+\zeta^{u_{i}}y)}\equiv 0 \modu \pi^{p-1}.
\end{equation}
\item  We get now, similarly to foundation theorem \ref{p21a} p.\pageref{p21a},
\begin{equation}\label{e23011}
\begin{split}
& \phi_m^*(t)\times S_m(J_{r^\prime})\equiv 0 \modu p,\quad m=1,\dots,p-1,\\
& S_m(J_{r^\prime})\equiv \sum_{i=1}^{r^\prime} u_{m i}\times  n_i^\prime\modu p.
\end{split}
\end{equation}
\item
The correspondance, {\it mutatis mutandis}, between previous theory and this subsection is summarized in notations:
\begin{displaymath}
\begin{split}
& f\longleftrightarrow f^\prime, \\
& r\longleftrightarrow r^\prime, \\
& \phi^*_{2m+1}(t),\quad  m=0,\dots \frac{p-3}{2} \longleftrightarrow \phi^*_m(t),
\quad m=1,\dots,p-1,\\
& I_{\phi^*}=\{  2m+1\quad |\quad 0\leq m\leq \frac{p-3}{2},
\quad \phi^*_{2m+1}(t)\not\equiv 0 \modu p\}\\
& \longleftrightarrow \\
&I^\prime_{\phi^*}=\{  m\quad |\quad 1\leq m\leq p-1,\quad \phi^*_m(t)\not\equiv 0 \modu p \}
\end{split}
\end{displaymath}
\end{itemize}
%
\subsection{The theorems on $\phi_m(t)\modu p,\quad,m=1,\dots,p-1$}\label{s104282}
We give now,   for  Mirimanoff polynomials $\phi_{m}(t),\quad m=1,\dots, p-1$, the  theorems corresponding to theorem \ref{pa212} p.\pageref{pa212}, theorem \ref{t12092} p.\pageref{t12092}, corollary \ref{cor1209} p.\pageref{cor1209} and theorem \ref{t23011} p. \pageref{t23011} for odd polynomials $\phi_{2m+1}(t),\quad m=0,\dots,\frac{p-3}{2}$.
%
\begin{thm}{ *** }\label{t104285}
Let $t\equiv -\frac{y}{x}\modu p$. Let $r^\prime$ be the minimal rank of Fermat congruences defined in the relation (\ref{e04031a}) p.\pageref{e04031a}.
If FLT first case  fails for $p$,  then there exists exactly
$p-r^\prime$ different values $m,\quad m\in\N,\quad 1\leq m \leq p-1$, such that we have the Mirimanoff polynomial congruences:
\begin{displaymath}\label{ea570}
\phi_{m}(t)\equiv 0 \modu p.
\end{displaymath}
Moreover $r^\prime\leq\frac{p+1}{2}$.
\begin{proof}
See theorem \ref{pa212} p.\pageref{pa212}. From Mirimanoff, see Ribenboim \cite{rib} (1B) p 145, we get $p-r^\prime\geq \frac{p-1}{2}$, which leads to the result.
\end{proof}
\end{thm}
%
\begin{thm}{ *** }\label{t104201}
If FLT  first case  fails for $p$,  then there exists exactly $\mu^\prime=r^\prime-1$ different values
$m,\quad m\in\N,\quad 1\leq m\leq p-1$, such that
\begin{displaymath}
S_m(J_{r^\prime})\equiv 0\modu p.
\end{displaymath}
\end{thm}
%
\begin{cor}{ *** }\label{c104201}
If FLT  first case  fails for $p$, then there exists
$r^\prime\in \N$, depending only on $p$ and $t$,  such that for any primitive root $u\modu p$,
$S_1(J_r)\equiv 0 \modu p$ and
for $m=2,\dots,p-1$:
\begin{itemize}
\item
$\phi_m(t)\times S_m(J_{r^\prime}) \equiv 0 \modu p$,
\item
If $\phi_m(t)\equiv 0\modu p$ then $S_m(J_{r^\prime})\not\equiv 0 \modu p$,
\item
If $S_m(J_{r^\prime})\equiv 0\modu p$ then $\phi_m(t)\not\equiv 0 \modu p$.
\end{itemize}
\end{cor}
%

\begin{thm}{ *** } { {\bf An explicit computation of $S_{m}(J_r)$.}}\label{t108313}

Let the set $I^\prime_\phi$ defined by

$I^\prime_\phi=\{m_i \quad |\quad i=1,\dots,r^\prime-1,\quad \phi^*_{m_i}(t)\not\equiv 0 \modu p,\quad 1\leq m_i\leq p-1 \}$.
\begin{enumerate}
\item
Then an  explicit formula for $S_{n}(J_r)$ is :
\begin{equation}\label{e23014}
S_{n}(J_r)=u_{n}\times\prod_{i=1}^{r-1} (u_{n}-u_{m_i}).
\end{equation}
\item
We have the equivalence
\begin{displaymath}
\begin{split}
& \sum_{i=1}^{r^\prime}\frac{u_i l_i}{x+\zeta^{u_i}y} \equiv 0 \modu \pi^{p-1} \\
& \Longleftrightarrow\\
& S_{n}(J_r)=u_{n}\times\prod_{i=1}^{r^\prime-1} (u_{n}-u_{m_i}) \\
& \Longleftrightarrow\\
& l_1 \equiv (-1)^{r^\prime-1}\times \prod_{i=1}^{r-1} u_{m_i} \modu p,\\
&\vdots \\
& l_{r^\prime-1}\equiv (-1)\times \sum_{i=1}^{r^\prime-1} u_{m_i} \modu p, \\
& l_{r^\prime}=1.
\end{split}
\end{displaymath}
\end{enumerate}
\end{thm}
%
%
\subsection{The representation approach for Fermat-Wiles}\label{s106021}

In this subsection, we take,  for the study of odd and even  Mirimanoff polynomials
$\phi_m(T),\quad m=1,\dots,p-1$, the approach taken in section \ref{s09121} p.\pageref{s09121} for the study of odd Mirimanoff polynomials $\phi_{2m+1}(T),\quad m=1,\dots,\frac{p-3}{2}$.
We shall give a description of Fermat's equation obtained from action of Galois group $G=Gal(\Q(\zeta)/\Q)$ on the class group and also  on the unit group of $\Q(\zeta)$.
%
\subsubsection{A summary on the different  $p$-groups of the Fermat's equation}\label{s106041}
Recall or fix some notations used in this subsection.
\begin{itemize}
\item
Let $r_p$ be the $p$-rank of the class group of the extension  $\Q(\zeta)/\Q$. We have $r_p=r_p^++r_p^-$,  where  $r_p^+$ is the $p$-rank of the class group of
the field $\Q(\zeta+\zeta^{-1})$, and $r_p^-$ is the relative class $p$-group of $\Q(\zeta)$.
\item
The ideal $\s$ of Fermat's equation is defined by $\s^p=(x+\zeta y)\Z[\zeta]$.
\item
The action of $G=Gal(\Q(\zeta)/\Q)$ on ideal class group  $<Cl(\s)>$ yields
\begin{equation}\label{e106021}
\begin{split}
& \prod_{i=0}^{r_1}(x+\zeta^{u_i}y)^{l_i}=\zeta^v\times\eta\times\gamma^p,\\
& v\in\N,\quad \eta\in\Z[\zeta+\zeta^{-1}]^*,\quad\gamma\in\Q(\zeta)
\quad l_1\not\equiv 0\modu p,\quad l_{r_1}\not\equiv 0\modu p.
\end{split}
\end{equation}
Observe that in first case of FLT then $v\not\equiv 0\modu p$.
The corresponding polynomial $P_{r_1}(V)=\prod_{i=1}^{r_1}(V-\mu_i),\quad \mu_i\in{\bf F}_p$, is the minimal polynomial annihilating the ideal $\s$,
so with exponential notation $\s^{P_{r_1}(\sigma)}\simeq \Z[\zeta]$.
\item
In the same way, the $p$-rank $r_1$ can be written $r_1=r_1^++r_1^-$, where $r_1^+$ is the $p$-rank of the action of the Galois group of $Gal(\Q(\zeta)/\Q)$ on the
ideal class group $<Cl(\s\times\overline{\s})>$, and $r_1^-$, relative $p$-rank, corresponds to the action of $Gal(\Q(\zeta)/\Q)$ on the ideal class group  $<Cl(\s\overline{\s}^{-1})>$.
\item
The $p$-rank $r_1^-$ is the minimal natural number such that
\begin{displaymath}
\prod_{i=0}^{r_1^-}(\frac{x+\zeta^{u_i} y}{x+\zeta^{-u_i} y})^{l_i}=\zeta^v\times\gamma^p,
\quad v\in\N,\quad v\not\equiv 0\modu p,\quad \gamma\in\Q(\zeta).
\end{displaymath}
or in exponential notation
\begin{displaymath}
(\frac{x+\zeta y}{x+\zeta^{-1}y})^{P_{r_1^-}(\sigma)}\simeq \Z[\zeta],
\end{displaymath}
\item
The $p$-rank $r_1^+$ is the minimal natural number such that
\begin{equation}\label{e112182}
\prod_{i=1}^{r_1^+} (x+\zeta^{u_i}y)^{l_i}\times (x+\zeta^{-u_i}y)^{l_i}
=\eta\times\gamma^p,\quad \eta\in\Z[\zeta+\zeta^{-1}]^*,
\quad \gamma\in\Z[\zeta+\zeta^{-1}],
\end{equation}
or in exponential notation,
\begin{displaymath}
((x+\zeta y)\times (x+\zeta^{-1}y))^{P_{r_1^+}(\sigma)}\simeq \Z[\zeta+\zeta^{-1}].
\end{displaymath}
\item
The minimal Fermat's congruence rank $r$ defined in relation (\ref{e31102}) p.\pageref{e31102} is the minimal natural integer such that
\begin{displaymath}
\prod_{i=0}^{r-1} (\frac{x+\zeta^{u_i}y}{x+\zeta^{-u_i}y})^{n_i}\equiv 1\modu\pi^{p+1},
\end{displaymath}
where $n_0\not\equiv 0\modu p,\quad n_{r-1}\not\equiv 0\modu p$,
or in exponential notation

\begin{displaymath}
(\frac{x+\zeta y}{x+\zeta^{-1}y})^{P_r(\sigma)}\equiv 1\modu\pi^{p+1}.
\end{displaymath}
Observe that we have prefered in this article the not too cumbersome notation $r$ instead of $r^-$ because this notation is frequently used.
\item
The minimal Fermat's congruence rank $r^+$ is the minimal rational integer such that
\begin{equation}\label{e112183}
\prod_{i=0}^{r^+-1} ((x+\zeta^{u_i} y)\times (x+\zeta^{-u_i}y))^{n_i^+}\equiv c^p\modu\pi^{p+1},\quad c\in\Z.
\end{equation}
or in representations notation
\begin{displaymath}
((x+\zeta^u y)\times (x+\zeta^{-u}y))^{P_{r^+}(\sigma)}\equiv c^p\modu\pi^{p+1},\quad
 c\in\Z.
\end{displaymath}
\item
The minimal Fermat's congruence rank $r^\prime$ is the minimal natural integer such that
\begin{displaymath}
\begin{split}
&\prod_{i=0}^{r^\prime-1} (x+\zeta^{u_i} y)^{n_i^\prime}
-c^{(p\times\sum_{i=1}^{r^\prime} n_i^\prime)}\equiv 0\modu\pi^{p-1},
\quad x+y=c^p,\quad c\in\Z,\\
& \prod_{i=0}^{r^\prime-1} (x+\zeta^{u_i} y)^{n_i^\prime} -
\prod_{i=0}^{r^\prime-1} (x+\zeta^{-u_i} y)^{n_i^\prime}\not =  0.
\end{split}
\end{displaymath}
This relation can also be written
\begin{equation}\label{e201251}
 (x+\zeta^{u_i} y)^{P_{r^\prime}(\sigma)}
-c^{(p\times\sum_{i=0}^{r^\prime-1} n_i^\prime)}\equiv 0\modu\pi^{p+1},
\quad x+y=c^p,
\end{equation}
where $P_{r^\prime}(\sigma)$ is the minimal corresponding polynomial.
\item
Let $r_s\in \N$, ($\s$ like first letter of sum)  be the minimal natural integer such that
\begin{displaymath}
\prod_{i=0}^ {r_s-1} (x+\zeta^{u_i} y)^{n_i}=\gamma^p,\quad\gamma\in\Z[\zeta],
\end{displaymath}
which can be written, for the minimal polynomial $P_{r_s}(\sigma)=
\prod_{i=1}^{r_s}(\sigma -\mu_i)$.
\item
For the unit $\eta$ verifying the relation (\ref{e106021}) p.\pageref{e106021}, there exists a {\bf minimal} unit rank $r_\eta$, defined in relation \ref{e203051} p. \pageref{e203051} in the section dealing of unit group $\Z[\zeta+\zeta^{-1}]^*$,
such that we have congruence in form
\begin{equation}\label{e203052}
\eta^{\nu_0}\times\sigma(\eta)^{\nu_1}\times\dots\times
\sigma^{r_\eta-1}(\eta)^{\nu_{r_\eta-1}}
\times\sigma^{r_\eta}(\eta)^{\nu_{r_\eta}}=\varepsilon^p,
\quad\varepsilon\in\Z[\zeta+\zeta^{-1}]^*,
\end{equation}
or, for the minimal polynomial $P_{r_\eta}(V)=\prod_{i=1}^{r_\eta}(V-\mu_i)$, the relation  $\eta^{P_{r_\eta}(\sigma)}=\varepsilon^p$.
\end{itemize}
%
\subsubsection{Some results on Representations approach}\label{s106042}
In this subsection we shall  prove some properties of minimal polynomials.
%
\begin{lem}\label{l106032}
Let $\mathbf b$ be an ideal of $\Z[\zeta]$, not principal and with $\mathbf b^p$ principal. Suppose that for $\mu\in{\bf F}_p^*$, the ideal $\mathbf b^{\sigma-\mu}$ is a principal ideal. Then
\begin{displaymath}
\begin{split}
&\mathbf b\times\overline{\mathbf b} \mbox{ principal}\Leftrightarrow
\mu^{(p-1)/2}=-1,\\
&\frac{\mathbf b}{\overline{\mathbf b}} \mbox{ principal}\Leftrightarrow
\mu^{(p-1)/2}=1.
\end{split}
\end{displaymath}
corresponding, for the minimal polynomial $P_{r_1^-}(V)$ and $P_{r_1^+}(V)$, to
\begin{displaymath}
\begin{split}
& P_{r_1^-}(V)=\prod_{i=1}^{r_1^-}(V-\mu_i),\quad \mu_i=u_{2 m_i+1},\quad m_i\in\N,
\quad 1\leq m_i\leq\frac{p-3}{2},\\
&P_{r_1^+}(V)=\prod_{i=1}^{r_1^+} (V-\mu_i),\quad \mu_i=u_{2 m_i},\quad m_i\in\N,
\quad 1\leq m_i\leq\frac{p-1}{2}.
\end{split}
\end{displaymath}
\begin{proof} $ $
\begin{itemize}
\item
If $\mathbf b\times \overline{\mathbf b}$ is principal then
$\mathbf b^{\sigma^{(p-1)/2}+1}$ is principal. From hypothesis
$\mathbf b^{\sigma-\mu}$ is principal,
so $\mathbf b^{\sigma^{(p-1)/2}-\mu^{(p-1)/2}}$ is principal. $\mu^{(p-1)/2}=1$ should imply that $\mathbf b^{\sigma^{(p-1)/2}-1}$ is principal, thus,
 with $\mathbf b^{\sigma^{(p-1)/2}+1}$ principal, it  should imply that
$\mathbf b$ is principal, contradiction.
\item
If $\mu^{(p-1)/2}=-1$ then from $\mathbf b^{\sigma-\mu}$ principal, we get
$\mathbf b^{\sigma^{(p-1)/2}-\mu^{(p-1)/2}}$ principal, and so
$\mathbf b^{\sigma^{(p-1)/2}+1}$ principal or
$\mathbf b\times \overline{\mathbf b}$ principal.
\item
Similar proof if $\frac{\mathbf b}{\overline{\mathbf b}}$ is principal.
\end{itemize}
\end{proof}
\end{lem}
%
\begin{lem}\label{l112181}
$r^+\leq r_p^++r_\eta.$
\begin{proof}
From relation (\ref{e112182}) p. \pageref{e112182} of definition of the rank $r_1^+$, we get
\begin{displaymath}
\prod_{i=1}^{r_1^+} (x+\zeta^{u_i}y)^{n_i}(x+\zeta^{-u_i}y)^{n_i}
=\eta\times\gamma_1^p,\quad \eta\in\Z[\zeta+\zeta^{-1}]^*,
\quad \gamma_1\in\Z[\zeta+\zeta^{-1}].
\end{displaymath}
From definition of the unit ranks $r_\eta$, we get
\begin{displaymath}
\prod_{j=0}^{r_\eta} \sigma^j(\eta)^{n_j^\prime}\equiv c^p \modu\pi^p,\quad
c\in\Z,\quad c\not\equiv 0\modu p.
\end{displaymath}
Gathering these two relations we get
\begin{displaymath}
\prod_{i=1}^{r_1^++r_\eta} (x+\zeta^{u_i}y)^{m_i}(x+\zeta^{-u_i}y)^{m_i}
\equiv c_2^p\modu \pi^p,\quad c_2\in\Z,\quad c_2\not\equiv 0\modu p.
\end{displaymath}
The minimal congruence rank $r^+$ defined in relation (\ref{e112183}) p. \pageref{e112183} is the minimal integer such that
\begin{displaymath}
\prod_{i=1}^{r^+} (x+\zeta^{u_i}y)^{l_i}(x+\zeta^{-u_i}y)^{l_i}
\equiv c^p\modu\pi^{p},\quad c\in\Z.
\end{displaymath}
From these definitions, we get $r^+\leq r_1^++r_\eta$.
Clearly $r_1^+\leq r_p^+$  which achieves the proof.
\end{proof}
\end{lem}
%
\subsubsection{A summary on the  ranks and minimal polynomials}
Here we summarize some relations  obtained between  the different ranks,  the different minimal polynomials.
\begin{itemize}
\item
The minimal corresponding polynomials of the indeterminate $V$, all dividing
the polynomial $V^{p-1}-1=\prod_{i=1}^{p-1}(V-i)$, can be written in form
\begin{equation}\label{e108312}
\begin{split}
& P_{r_1}(V)=\prod_{i=1}^{r_1} (V-\mu_i),\quad \mu_i\in {\bf F}_p,\quad
i_1\not=i_2\Rightarrow \mu_{i_1}\not=\mu_{i_2},\\
& P_{r_1^-}(V)=\prod_{i=1}^{r_1^-} (V-\mu_i),\quad \mu_i\in {\bf F}_p,\quad
i_1\not=i_2\Rightarrow \mu_{i_1}\not=\mu_{i_2},\quad \mu_i=u_{2m_i+1},\\
& P_{r_1^+}(V)=\prod_{i=1}^{r_1^+} (V-\mu_i),\quad \mu_i\in {\bf F}_p,\quad
i_1\not=i_2\Rightarrow \mu_{i_1}\not=\mu_{i_2},\quad \mu_i=u_{2m_i},\\
& P_{r}(V)=\prod_{i=1}^{r-1} (V-\mu_i),\quad \mu_i\in {\bf F}_p,\quad
i_1\not=i_2\Rightarrow \mu_{i_1}\not=\mu_{i_2},\quad \mu_i=u_{2m_i+1},\\
& P_{r^+}(V)=\prod_{i=1}^{r^+} (V-\mu_i),\quad \mu_i\in {\bf F}_p,\quad
i_1\not=i_2\Rightarrow \mu_{i_1}\not=\mu_{i_2},\quad \mu_i=u_{2m_i},\\
& P_{r^\prime}(V)=\prod_{i=1}^{r^\prime-1} (V-\mu_i),\quad \mu_i\in {\bf F}_p,\quad
i_1\not=i_2\Rightarrow \mu_{i_1}\not=\mu_{i_2},\\
& P_{r_s}(V)=\prod_{i=1}^{r_s} (V-\mu_i),\quad \mu_i\in {\bf F}_p,\quad
i_1\not=i_2\Rightarrow \mu_{i_1}\not=\mu_{i_2},\\
& P_{r_\eta}(V)=\prod_{i=1}^{r_\eta} (V-\mu_i),\quad \mu_i\in {\bf F}_p,\quad
i_1\not=i_2\Rightarrow \mu_{i_1}\not=\mu_{i_2},\quad \mu_i=u_{2m_i}.
\end{split}
\end{equation}
\item
The minimal polynomials verify the following divisibility properties :
\begin{equation}\label{e108313}
\begin{split}
& P_{r_1^-}(V)\times P_{r_1^+}(V)=P_{r_1}(V),\\
& P_r(V)\ |\    P_{r_1^-}(V)\times (V-u),\\
& P_r(V)\ |\ P_{r^\prime}(V),\\
& P_{r^\prime}(V)\ |\ P_{r_1}(V)\times (V-u)\times P_{r_\eta}(V)=P_{r_s}(V).
\end{split}
\end{equation}
\item
The different ranks verify the inequalities:
\begin{equation}\label{e108314}
\begin{split}
& r_1=r_1^-+r_1^+,\\
& r_1^-\leq r_p^-,\\
& r_1^+\leq r_p^+,\\
& r-1\leq r_1^-+1,\\
& r^+\leq r_1^++r_\eta,\\
& r^\prime-1 \leq r_s = r_1+r_\eta+1.\\
\end{split}
\end{equation}
\item
In this relation (\ref{e108314}), the difference $r_1^--(r-1)$ corresponds to the singular primary $C_i$ among the $r_1^-$ terms  $C_i$
of the formula
\begin{displaymath}
\frac{x+\zeta y}{x+\zeta^{-1}y}=\zeta^{2v}\times \prod_{i=1}^{r_1^-} C_i^{\nu_i}\times\gamma^p,
\quad \gamma\in\Q(\zeta),\quad v\in{\bf F}_p^*,
\end{displaymath}
where $C_i$ are defined in relation (\ref{e108191}) p. \pageref{e108191}.
In fact these $C_i\equiv 1\modu\pi^p$ can be neglected in  $\modu\pi^p$ congruences.
\item
In this relation (\ref{e108314}), the difference $r_s^--(r^\prime-1)$ corresponds to the singular primary $\eta_i$ among the $r_\eta$ units   $\eta_i$
of the relation
\begin{displaymath}
x+\zeta y=\zeta^{v}\times \prod_{i=1}^{r_1} C_i^{\nu_i}
\times\prod_{j=1}^{r_\eta} \eta_j^{\delta_j}\times \gamma^p,
\quad \gamma\in\Q(\zeta),\quad \nu_i, \delta_j, v\in{\bf F}_p^*,
\end{displaymath}
where units $\eta_i$ are those of the relation \ref{e201274} p. \pageref{e201274}.
In fact these units $\eta_j\equiv c^p\modu\pi^p,\quad c\in{\bf F}_p^*$ can be neglected in  $\modu\pi^p$ congruences.
\end{itemize}
%
%
\subsubsection{On singular integers and  minimal polynomial $P_{r_s}(V)$.}
\label{s201134}
Recall that the polynomial $P_{r_s}(V)$ of the indeterminate $V$ is the minimal polynomial such that
\begin{displaymath}
(x+\zeta y)^{P_{r_s}(\sigma)}=\gamma^p,\quad \gamma\in\Z[\zeta].
\end{displaymath}
Let us consider the method of singular integers explained in Ribenboim, \cite{rib} paragraph 4. p 170 : an element $b\in\Z[\zeta]$ is said to be a singular integer if there exists an ideal $\mathbf b$ such that $b\Z[\zeta]=\mathbf b^p$. In that sense,
in the Fermat's equation with $z\not\equiv 0\modu p$ then $x+\zeta y$ with
$(x+\zeta y)\Z[\zeta]= \s^p$, is a singular integer.
Let $V$ be an indeterminate.
Let the polynomial $H(V)\in {\bf F_p}[V]$ defined by
\begin{displaymath}
H(V)=\sum_{i=0}^{p-2}\frac{u\times u_{-i}-u_{-i+1}}{p}\times V^i.
\end{displaymath}
Let us define also the $p-2$ polynomials
\begin{displaymath}
I_a(V)=\sum_{i=0}^{p-2}([\frac{(a+1)\times u_i}{p}]-[\frac{a\times u_i}{p}])
\times V^{p-2-i},\quad a=1,\dots, p-2.
\end{displaymath}
We get the  lemma on the polynomial $P_{r_s}(V)$.
%
\begin{lem}\label{l11006}
\begin{displaymath}
\begin{split}
& P_{r_s}(V)\ |\ (V-u)\times H(V).\\
& P_{r_s}(V)\ |\ (V-u)\times I_a(V),\quad a=1,\dots,p-2.
\end{split}
\end{displaymath}
\begin{proof}
$P_{r_s}(V)$ is the minimal polynomial such that
\begin{displaymath}
(x+\zeta y)^{P_{r_s}(\sigma)}=\gamma_1^p,\quad \gamma_1\in\Z[\zeta].
\end{displaymath}
From Inkeri \cite{ink}, cited in Ribenboim \cite{rib}, (4A)-2 p 171, we get
\begin{displaymath}
(x+\zeta y)^{H(\sigma)}=\zeta^{v_2}\times\gamma_2^p,\quad \gamma_2\in\Z[\zeta].
\end{displaymath}
We derive
\begin{displaymath}
(x+\zeta y)^{(\sigma-u)\times H(\sigma)}=(\gamma^\prime_2)^p,
\quad \gamma^\prime_2\in\Z[\zeta],
\end{displaymath}
which implies that the minimal polynomial $P_{r_s}(V)$ divides the polynomial
$(V-u)\times H(V)$.
The end of the proof is similar starting of the result of a result of Fueter \cite{fue} cited in Ribenboim \cite{rib} (4C)-2 p 172.
\end{proof}
\end{lem}
{\bf Remarks:}
\begin{itemize}
\item
Observe that the coefficients $\delta_i(a)=[\frac{(a+1)\times u_i}{p}]
-[\frac{a\times u_i}{p}]$ of the polynomial $I_a(V)$ verify
\begin{displaymath}
\delta_i(a)\in \{0,1\},\quad i=0,\dots, p-2,\quad a=1,\dots,p-2.
\end{displaymath}
\item
Observe, as an example, that for $a=1$ then
\begin{displaymath}
I_1(V)=\sum_{i=0,\quad u_i> p/2}^{p-2} V^{p-2-i}.
\end{displaymath}
\item
A numerical MAPLE observation has allowed to observe:
Let $d(V)$ be the polynomial in ${\bf F}_p[V]$ of the indeterminate $V$ defined by
\begin{displaymath}
d(V)=Gcd(I_a(V)\ |\ a=1,\dots,p-2).
\end{displaymath}
Let the two polynomials $P^-(V)$ and $P^+(V)$ of the indeterminate $V$ defined by
\begin{displaymath}
\begin{split}
& P^-(V)=Gcd(V^{(p-1)/2}+1, d(V)), \\
& P^+(V)=Gcd(V^{(p-1)/2}-1,d(V)).
\end{split}
\end{displaymath}
We have observed that $P^+(V)=1+V+\dots+V^{p-3}{2}$,
for $h^+\not\equiv 0 \modu p$,  the relations
\begin{equation}\label{e111043}
\begin{split}
& \mbox{degree}_V(P^-(V))= \mbox{ index of irregularity of\ } p,\\
& \mbox{degree}_V(P^+(V))= \frac{p-3}{2}.
\end{split}
\end{equation}
Therefore this should lead to  $d(V)=P_{r_1^-}(V)$,  general observation independant of Fermat which could also be a strong tool for first case of FLT.
\item
Observe that these results given in the case of Fermat's equation for the singular integer $x+\zeta y,\quad (x+\zeta y)\Z[\zeta]=\s^p$,  can be formulated for the minimal polynomial of any singular integer $b$.
\end{itemize}
%
%
\subsubsection{Congruences on Mirimanoff polynomials}\label{s108256}
We give here an  application to Mirimanoff polynomials in Fermat context.
Let $p$ be a prime, $p>5$.
%
\begin{thm}\label{t108311}
Let $p>5$ be a prime. If FLT first case would fail for $p$ there should be
{\bf at most} $r_\eta+r_p^+$ Mirimanoff's polynomials $\phi_{2m}(T)$ with
$\phi_{2m}(t)\not\equiv 0\modu p, \quad 1\leq m \leq \frac{p-1}{2}$.
\begin{proof} $ $
\begin{itemize}
\item
Recall that $(x+\zeta y)^{P_{r_1}(\sigma)}=\zeta^v\times \eta\times\gamma^p,
\quad \eta\in\Z[\zeta+\zeta^{-1}]^*,\quad \gamma\in \Q(\zeta)$, and that
$E(\eta) =E(\eta_1)\times \dots\times E(\eta_{r_\eta}),
\quad E(\sigma(\eta_i))=E(\eta_i^{\mu_i}),\quad i=1,\dots,r_\eta.$
\item
Let the minimal polynomial $P_{r^\prime}(V)$ defined in relation (\ref{e108312}) p.\pageref{e108312}. We have seen in relation (\ref{e108313}) p.\pageref{e108313} that
\begin{equation}\label{e201252}
P_{r^\prime}(V)\ |\ P_{r_1^-}(V)\times P_{r_1^+}(V)\times (V-u)\times P_{r_\eta}(V).
\end{equation}
Let $\mu=u_{2m}$ with $(V-\mu)\ |\ P_{r^\prime}(V)$.
We have $(V-\mu)\not| P_{r_1^-}(V)$ because $P_{r_1^-}(V)=\prod_{i=1}^{r_1^-}(V-\mu_i),\quad \mu_i=u_{2m_i+1}$.
We have also $V-\mu\not| V-u$.
Therefore $V-\mu\ |\ P_{r_1^+}(V)\times P_{r_\eta}(V)$ and even
\begin{displaymath}
V-\mu\ |\ lcm(P_{r_1^+}(V), P_{r_\eta}(V)).
\end{displaymath}
\end{itemize}
\end{proof}
\end{thm}
%
%

%
%
%
\clearpage
\section{ FLT first case:  Hilbert's class field of $\Q(\zeta)$ approach}\label{s204261}
This section contains :
\begin{enumerate}
\item
a general  application of properties of $p$-Hilbert's class field of $\Q(\zeta)$ to FLT,
\item
some properties of  unit $\varepsilon$ of Fermat's equation in $p$-Hilbert class field $H$ of $\Q(\zeta)$.
\end{enumerate}
%
\subsection{ Application of Hilbert's class field to FLT}
We apply these explicit computations of $p$-elementary Hilbert's class field to FLT context.
A first application when $h^+\not\equiv 0\modu p$ is derived directly from
theorem \ref{t201061} p. \pageref{t201061}.
%
\begin{thm}\label{t201062}
Suppose that $h^+\not\equiv 0\modu p$.
If $\prod_{i=1}^r(\frac{x+\zeta^{u_i}y}{x+\zeta^{-u_i}y})^{l_i}\equiv 1\modu \pi^p$ then
$\prod_{i=1}^r(\frac{x+\zeta^{u_i}y}{x+\zeta^{-u_i}y})^{l_i}=\gamma^p,\quad
\gamma\in\Q(\zeta)$.
\end{thm}
%
Let $H$ be the $p$-elementary Hilbert field of $\Q(\zeta)$. Let $B_H$ be the ring of integers of $H$. Let $B_H^*$ be the group of units of $B_H$. Let $Gal(H/\Q(zeta)$ be the Galois group of the extension of
$H/\Q(\zeta)$.
With a second approach, from Principal Ideal Theorem, the Fermat's equation with $z\not\equiv 0 \modu p$ becomes  in $B_H$:
\begin{displaymath}
x+\zeta y = \varepsilon\times \gamma^p, \quad \gamma\in B_H,\quad
\varepsilon\in B_H^*,
\end{displaymath}
because the ideal $\s$ with $\s^p=(x+\zeta y)\Z[\zeta]$  becomes principal in $B_H$.
%
\begin{thm} { *** }\label{t012161}
Let $x,y,z$ be a solution of Fermat's equation with $z\not\equiv 0 \modu p$.
The Fermat's equation becomes in the ring of integers $B_H$ of the $p$-elementary Hilbert class field $H$
\begin{equation}\label{e012174}
x+\zeta y=\varepsilon\times \gamma^p,\quad \gamma\in B_H,\quad
\varepsilon\in B_H^*.
\end{equation}
with
\begin{equation}\label{e012161}
\frac{\theta(\varepsilon)}{\varepsilon}
=\varepsilon_1^p,\quad \varepsilon_1\in B_H^*,\quad
\forall\ \theta\in Gal(H/\Q(\zeta)).
\end{equation}
\begin{proof}
Let us consider the relation seen before
\begin{displaymath}
x+\zeta y=\varepsilon\times \gamma^p,\quad \gamma\in B_H,\quad
\varepsilon\in B_H^*.
\end{displaymath}
We have
\begin{displaymath}
\theta(x+\zeta y)=x+\zeta y
\end{displaymath}
and so
\begin{equation}\label{e012162}
\theta(\varepsilon)\times
\theta( \gamma^p)
=\varepsilon\times\gamma^p.
\end{equation}
We have $\frac{\theta(\varepsilon)}{\varepsilon}\in B_H^*$ and
$\frac{\theta(\varepsilon)}{\varepsilon}=\frac{\gamma^p}{\theta(\gamma^p)}\in H^p$;
thus $\frac{\theta(\varepsilon)}{\varepsilon}\in B_H^*\cap H^p$ and finally
$\frac{\theta(\varepsilon)}{\varepsilon}=\varepsilon_1^p, \quad
\varepsilon_1\in B_H^*$,
which achieves the proof.
\end{proof}
\end{thm}
{\bf Question:} Does there exists a way to caracterize units verifying strong relation (\ref{e012161}) p.\pageref{e012161}?
%
%

%
\subsection {Some congruences  on unit $\varepsilon$ of Fermat's equation in
$p$-Hilbert class field  $H$.}
This subsection contains a study of the units of the Fermat's equation in the $p$-Hilbert class field:
\begin{itemize}
\item
with a representations of $Gal(\Q(\zeta)/\Q)$ approach,
\item
with an application of the Dirichlet's theorem on units.
\end{itemize}
%
\subsubsection{ A representation approach}
\begin{itemize}
\item
Let $r_p$ be the $p$-rank of the class group $\Q(\zeta)$.
\item
Let $H$ be the $p$-elementary Hilbert class field of $\Q(\zeta)$,
\begin{displaymath}
H=\Q(\zeta,\omega_1,\dots,\omega_{r_p}).
\end{displaymath}
Let $B_H$ be the ring of integers of $H$.
\item
Let us consider the Fermat's equation in $B_H$,
\begin{displaymath}
x+\zeta y =\varepsilon\times\gamma^p,\quad \varepsilon\in B_H^*,\quad
\gamma\in B_H.
\end{displaymath}
\item
The extension $H/\Q$ is Galois.
Let $\Theta$ be the $\Q$-isomorphisms of $H/\Q$. Let $\theta$ be the
$\Q(\zeta)$-isomorphisms of $H/\Q(\zeta)$. The $\Q$-isomorphism $\sigma$ of $H$, extending the $\Q$-isomorphism $\sigma :\zeta\rightarrow \zeta^u$ of $\Q(\zeta)$, has been defined in relations (\ref{e102091}) p.\pageref{e102091} and (\ref{e102092}) p.\pageref{e102092}.
\begin{displaymath}
\Theta=\theta\circ\sigma^k,\quad\forall\theta\in Gal(H/\Q(\zeta),\quad
k=0,\dots,p-2.
\end{displaymath}
\item
From relation (\ref{e012161}) p.\pageref{e012161}, we have
\begin{displaymath}
\theta\circ\sigma^k(\frac{\varepsilon}{\overline{\varepsilon}})
\equiv\sigma^k(\frac{\varepsilon}{\overline{\varepsilon}})\modu (B_H^*)^p,
\quad \forall\theta\in Gal(H/\Q(\zeta)),\quad k=0,\dots,p-2,
\end{displaymath}
or
\begin{displaymath}
\theta\circ\sigma^k(\frac{\varepsilon}{\overline{\varepsilon}})=
\frac{\varepsilon}{\overline{\varepsilon}}\times \varepsilon_2^p,\quad\varepsilon_2\in B_H^*.
\end{displaymath}
\item
Let $r_1\leq r_p$ be the $p$-rank of the action of $Gal(\Q(\zeta)/\Q$ on the $p$-class group   $<Cl(\s)>,\quad \s^p=(x+\zeta y)\Z[\zeta]$.
\item
We have proved in theorem \ref{t2009} p.\pageref{t2009} that
there exists $l_i,\quad i=0,\dots,r_1$, such that
\begin{displaymath}
\prod_{i=0}^{r_1} (\frac{x+\zeta^{u_i} y}{x+\zeta^{-u_i} y})^{l_i}=
\gamma_\zeta^p,\quad \gamma_\zeta\in \Q(\zeta).
\end{displaymath}
Therefore, from the Fermat's equation  $x+\zeta y=\varepsilon\times\gamma^p$ in $H$,
we get
\begin{displaymath}
\prod_{i=0}^{r_1} \sigma^i (\frac{\varepsilon}{\overline{\varepsilon}})^{l_i}
\times \prod_{i=0}^{r_1}(\frac{\gamma}{\overline{\gamma}})^{p l_i}
=\gamma_\zeta^p,\quad\gamma_\zeta\in \Q(\zeta).
\end{displaymath}
\end{itemize}
Thus, we have proved
%
\begin{lem}\label{l102081}
Let the Fermat's equation in $B_H$
\begin{displaymath}
x+\zeta y=\varepsilon\times\gamma^p,\quad \varepsilon\in B_H^*,
\quad \gamma\in B_H.
\end{displaymath}
There exists a minimal $r_\varepsilon\in\N,\quad 1\leq r_\varepsilon\leq r_1$, such that
\begin{displaymath}
\prod_{i=0}^{r_\varepsilon} \sigma^i(\frac{\varepsilon}{\overline{\varepsilon}})^{m_i}=
\varepsilon_H^p,\quad \varepsilon_H\in B_H^*,
\quad m_i\in \N,\quad m_0\not=0,\quad m_{r_\varepsilon}\not=0.
\end{displaymath}
\end{lem}
%

It is possible to formulate  this result in terms of representations.
\begin{lem}\label{l106081}
Let the Fermat's equation in $B_H$
\begin{displaymath}
x+\zeta y=\varepsilon\times\gamma^p,\quad \varepsilon\in B_H^*,
\quad \gamma\in B_H.
\end{displaymath}
There exists a minimal rank $r_\varepsilon\leq r_1\leq r_p$ and a minimal polynomial
$P_{r_\varepsilon}(\sigma)\in {\bf F}_p[G]$ such that
\begin{displaymath}
(\frac{\varepsilon}{\overline{\varepsilon}})^{P_{r_\varepsilon}(\sigma)}
=\varepsilon_H^p,\quad \varepsilon_H\in B_H^p.
\end{displaymath}
For the indeterminate $V$, $P_{r_\varepsilon}(V)|P_{r_1}(V)$.
\begin{proof}
We have shown that $(\frac{\varepsilon}{\overline{\varepsilon}})^{P_{r_1}(\sigma)}\in (B_H^*)^p$ in lemma \ref{l102081} p.\pageref{l102081}. We apply the representation of $G$ like in section \ref{s09121} p.\pageref{s09121} to get the result.
\end{proof}
\end{lem}
%
\subsubsection{An application of Dirichlet's units theorem in $H/\Q$}

We shall give a criterion on system of fundamental units of Hilbet class field $H$ of $\Q(\zeta)$ verified when  FLT holds for $p$.
\begin{itemize}
\item
Let $\rho=\frac{(p-1)\times p^{r_p}}{2}-1$ be the rank of the group of units $B_H^*$.
\item
Let $\eta_i,\quad i=1,\dots,\rho$, be a system of fundamental units of $B_H^*$.
\item
For each $\theta\in Gal(H/\Q(\zeta))$ we have a system of relations
\begin{displaymath}
\begin{split}
& \theta(\eta_1)=\eta_1^{\mu_{1,1}}\times\dots\times\eta_{\rho}^{\mu_{1,\rho}},\\
& \vdots \\
& \theta(\eta_i)=\eta_1^{\mu_{i,1}}\times\dots\times\eta_{\rho}^{\mu_{i,\rho}},\\
& \vdots \\
& \theta(\eta_\rho)=\eta_1^{\mu_{\rho,1}}\times\dots\times
\eta_{\rho}^{\mu_{\rho,\rho}}.
\end{split}
\end{displaymath}
\item
Let $\varepsilon$ be the unit of the Fermat's equation:
\begin{displaymath}
\varepsilon=\eta_1^{\nu_1}\times\dots\times\eta_\rho^{\nu_{\rho}}.
\end{displaymath}
We have seen that $\frac{\theta(\varepsilon)}{\varepsilon}=\varepsilon_1^p,\quad \varepsilon_1\in B_H^*$. From these relations, we get the linear system of $\rho$ congruences $\modu p$ on the $\rho$ indeterminates $\nu_j,\quad j=1,\dots,\rho$:
\begin{displaymath}
\begin{split}
& \mu_{1,1}\times\nu_1+\dots+\mu_{j,1}\times\nu_j+\dots+
\mu_{\rho,1}\times\nu_{\rho}\equiv \nu_1 \modu p,\\
&\vdots \\
& \mu_{1,i}\times\nu_1+\dots+\mu_{j,i}\times\nu_j+\dots+
\mu_{\rho,i}\times\nu_{\rho}\equiv \nu_i \modu p,\\
&\vdots \\
& \mu_{1,\rho}\times\nu_1+\dots+\mu_{j,\rho}\times\nu_j+\dots+
\mu_{\rho,\rho}\times\nu_{\rho}\equiv \nu_\rho \modu p.
\end{split}
\end{displaymath}
\item
This implies that the determinant $M(p,\theta)$ depending only on $p$ and $\theta\in Gal(H/\Q(\zeta))$
\begin{equation}\label{e102172}
M(p,\theta)=
\begin{array}{|lllll|}
\mu_{1,1}-1 & \dots & \mu_{j,1}& \dots & \mu_{\rho,1} \\
\vdots &&\vdots && \\
\mu_{1,j} & \dots & \mu_{j,j}-1& \dots & \mu_{\rho,j} \\
\vdots &&\vdots && \\
\mu_{1,\rho} & \dots & \mu_{j,\rho}& \dots & \mu_{\rho,\rho}-1 \\
\end{array}
\end{equation}
\end{itemize}
verifies the congruence $M(p,\theta)\equiv 0\modu p$. We have obtained the result :
%
\begin{thm}\label{t102171}
If there exists a sytem $\eta_i,\quad i=1,\dots,\rho$, of fundamental units of
$B_H^*$ and one $\theta\in Gal(H/\Q(\zeta))$, such that the determinant $M(p,\theta)\not\equiv 0\modu  p$, then FLT first case holds for $p$.
\end{thm}
%
{\bf Remark :} Observe that, in this theorem,  we give a criterion for first case of FLT depending only on $p$ and  structure of the units group of the $p$-elementary Hilbert class field of $\Q(\zeta)$. It does not depends on $x,y,z$ verifying Fermat's equation.
%
%

%
%

\clearpage
\begin{section} {FLT first case: An improvement of  Kummer's result on Jacobi resolvents}\label{seckum}

In this chapter, we improve some classical Kummer's results on FLT resting on Jacobi resolvents.
It contains:
\begin{enumerate}
\item
some general definitions, subsection  independant of FLT,
\item
some properties of Kummer's method of annihilation of ideals classes of the relative $p$-class group,
subsection independant of FLT,
\item
some applications in FLT context.
\end{enumerate}
\subsection{Some definitions}\label{s22031a}
These results are in relation with Kummer results \cite{kum} mentionned in Ribenboim, \cite{rib} p119.
These definitions faithfully follows (\cite{rib}) page 118-119.
\begin{itemize}
\item
Let $a,b\in\R$; recall that $[a,b]=\{ \alpha\quad | \alpha\in \R,\quad a\leq \alpha\leq b\}$ and that
$[a,b[=\{ \alpha\quad | \alpha\in \R,\quad a\leq \alpha<  b\}$.
\item
In this section,  when $r\equiv a \modu p$, then $r$ verifies $0\leq r\leq p-1$.
\item
Let $u\in\N$ be a primitive root $\modu p$.
Let $\sigma :\Q(\zeta)\rightarrow \Q(\zeta)$ be the $\Q$-isomorphism defined by $\sigma(\zeta)=\zeta^u$.
\item
Let $u_i\equiv u^i \modu p,\quad 1\leq u_i\leq p-1$.
Note that, if $i<0$, then this is to be understood as $u_i\times u^{-i}\equiv 1 \modu p$.
Note that $u_{ij}\equiv u^{ij}\equiv (u^i)^j\equiv u_i^j\equiv (u_i)_j\modu p$.
Note that $u_i\equiv u_{i \mbox{ mod } p-1}$.
\item
Let $d\in\N,\quad 1\leq d \leq p-2$. Let $\nu=\frac{p-1}{2}$.
Let $s=ind_u(d)\in\N,\quad 0\leq s \leq p-2$ be defined by $d\equiv u^s \modu p$.
There exists always one and only one $s$ because $u$ is a primitive root $\modu p$ and therefore $ind_u(d)$ is well defined.
\item
Let $I(u,d)$, and if there  is no ambiguity $I(d)$,  be the set defined by
\begin{equation} \label{e24031a}
I(d)=\{i \quad|\quad 1\leq i \leq p-1, \quad u_{\nu-i}+u_{\nu-i+ind_u(d)}>p\}.
\end{equation}
\item
Let $q\in\N,\quad q\equiv 1\modu p$, be a prime.
Let $\mathbf q\subset\Z[\zeta]$ be a prime ideal above $q$. The relation (2.5) p119 of \cite{rib} says that
$\prod_{i\in I(u,d)} \sigma^i(\mathbf q)$
is a principal integral ideal of $\Z[\zeta]$.
\item
We shall prove here also two properties of $I(d)$ used in the sequel of this article:
\begin{itemize}
\item
Let $I^\prime(d)=\{i^\prime\quad|\quad 1\leq i^\prime \leq p-1,\quad i^\prime\not\in I(d)\}$. From
\cite{rib} p 121, (line after relation 4.5),
$i^\prime\in I^\prime(d)\Longleftrightarrow i^\prime\equiv i+\frac{p-1}{2} \modu p-1$
where $i\in I(d)$.  Observe that we write $i^\prime\equiv i+\frac{p-1}{2}\modu (p-1)$ and not $\modu p$ that we think to be a print error of Ribenboim book, we shall detail under why:
\item
$Card(I(d))=Card(I^\prime(d))=\frac{p-1}{2}$.
\end{itemize}
Let $i\in I(d)$; Let $i^\prime=i\pm \frac{p-1}{2},\quad 1\leq i^\prime \leq p-1$.
Let us consider the quantities
\begin{equation}\label{e204101}
\begin{split}
& u_{\nu-i}=r_1,\quad 1\leq r_1\leq p-1,\\
& u_{\nu-i+ind_u(d)}=r_2,\quad 1\leq r_2\leq p-1,\\
& u_{\nu-i}+u_{\nu-i+ind_u(d)}=r_1+r_2,\\
& u_{\nu-i^\prime}= p-r_1,\\
& u_{\nu-i^\prime+ind_u(d)}=p-r_2,\\
& u_{\nu-i^\prime}+u_{\nu-i^\prime+ind_u(d)}=p-r_1+p-r_2.
\end{split}
\end{equation}
Show that
\begin{equation}\label{e009161}
i\in I(d)\Leftrightarrow i^\prime\equiv i+\frac{p-1}{2}\modu p-1\in I^\prime(d).
\end{equation}
\begin{itemize}
\item
If $i\in I(d)$ then $r_1+r_2>p$ so $p-r_1+p-r_2<p$ and $i^\prime\not\in I(d)$.
\item
If $i^\prime\not\in I(d)$ then $p-r_1+p-r_2\leq p$ so $r_1+r_2\geq p$.
\begin{itemize}
\item
if $r_1+r_2=p$ then $u_{\nu-i}+(d\times u_{\nu-i}\modu p)=p$ so
$u_{\nu-i}+d\times u_{\nu-i}\equiv 0 \modu p$, thus $d+1\equiv 0 \modu p$ which contradicts hypothesis assumed  $1\leq d\leq p-2$.
\item
if $r_1+r_2>p$ then $i\in I(d)$ and we are done.
\end{itemize}
\end{itemize}
Note that it implies that $Card(I(d))=Card(I^\prime(d))=\frac{p-1}{2}$.
\item
We complete these notations by the definition of the set
\begin{equation}\label{eqkum01}
K(d)=\{k \quad|\quad k=p-1-i,\quad i\in I(d)\}.
\end{equation}
\end{itemize}


%
\subsection{Elementary properties of sets $I(d)$}\label{s20031a}
This subsection deals of cyclotomic fields independently of FLT.
\begin{lem} \label{l08041a}
For $d=1,\dots,p-2$, we have
\begin{equation}\label{e08041a}
I(d)= I(p-1-d).
\end{equation}
\begin{proof}
We have
\begin{equation}\label{e009153}
I(d)=\{i \quad | \quad u_{\nu-i}+d\times u_{\nu-i}\modu p>p \},
\end{equation}
and
\begin{equation}\label{e009154}
I(p-1-d)=\{i \quad | \quad u_{\nu-i}+((p-1-d)\times u_{\nu-i}\modu p)>p \}.
\end{equation}
Let $i\in I(p-1-d)$.
Let us consider $r_1$ and $r_2$ defined in relation (\ref{e204101}) p. \pageref{e204101}. If $(r_1+r_2)\leq p$ then
\begin{displaymath}
(p-1-d)\times u_{\nu-i}\modu p = p-r_1-r_2,
\end{displaymath}
and then, from equation (\ref{e009154}) p.\pageref{e009154}, we obtain $r_1+p-r_1-r_2>p$ contradiction,
thus $r_1+r_2>p$ which implies that $i\in I(d)$.
\end{proof}
\end{lem}
%
\begin{lem}  \label{l21031a}
$I(1)=\{ i\quad |\quad u_i\equiv -j^{-1}\modu p, \quad |\quad \frac{p}{2} < j <p\}$.
\begin{proof}
Let us  evaluate $I(1)$.
We have $u^{ind_u(1)}\equiv 1\mod p$ and so
\begin{displaymath}
I(1)=\{i\quad |\quad 0\leq i\leq p-2,\quad (u_{\nu-i}+u_{\nu-i})>p\}.
\end{displaymath}
Therefore,
we have
\begin{displaymath}
K(1)=\{k \quad |\quad 1\leq k\leq p-1,\quad (u_{\nu+k}+u_{\nu+k})>p\}.
\end{displaymath}
so $p>j=u_{\nu+k}=u_{\nu-i}>\frac{p}{2}$.
But  $u_i\equiv -(u_{\nu-i})^{-1}\equiv -j^{-1}\modu p$, because $u^\nu\equiv -1\modu p$,
and finally
\begin{displaymath}
I(1)=\{i\quad |\quad u_i\equiv -j^{-1}\modu p,\quad \frac{p}{2}< j <p\}.
\end{displaymath}
\end{proof}
\end{lem}
{\bf Remark :}Observe that we have $Card(I(1))=\frac{p-1}{2}$.

%
\begin{lem} \label{l21032a}
$I(2)=\{ i\quad |\quad u_i\equiv  j^{-1}\modu p, \quad j\in [\frac{p}{3}, \frac{p}{2}]
\cup [ \frac{2p}{3}, p-1]\}$.
\begin{proof}
From $u^{ind_u(2)}\equiv 2\mod p$,
we have
\begin{displaymath}
I(2)=\{i\quad |\quad 0\leq i\leq p-2,\quad (u_{\nu-i}+(2u_{\nu-i}\modu p))>p\}.
\end{displaymath}
Therefore,
we have
\begin{displaymath}
K(2)=\{k\quad |\quad 1\leq k\leq p-1,\quad (u_{\nu+k}+(2u_{\nu+k}\modu p))>p\}.
\end{displaymath}
If $u_{\nu+k}<\frac{p}{2}$ then $(2u_{\nu+k}\modu p)=2u_{\nu+k}$
leads to
$u_{\nu+k}>\frac{p}{3}$, hence $\frac{p}{3} < j=u_{\nu+k} <\frac{p}{2}$.
If $u_{\nu+k}>\frac{p}{2}$ then $(2u_{\nu+k}\modu p)=2u_{\nu+k}-p$
leads to
$p>u_{\nu+k}>\frac{2p}{3}$ then to $\frac{2p}{3} <j=u_{\nu+k} <p$.

But $u_i\equiv -(u_{\nu-i})^{-1}\equiv -j^{-1}\modu p$, because $u^\nu\equiv -1\modu p$
and finally
\begin{displaymath}
I(2)=\{ i \quad |\quad u_i\equiv -j^{-1}\modu p, \quad  j \in[\frac{p}{3}, \frac{p}{2}]
\cup [\frac{2p}{3}, p-1] \}.
\end{displaymath}
\end{proof}
\end{lem}
%
\begin{lem}\label{l21033a}
\begin{displaymath}
I(k)=\{ i \quad | \quad u_i\equiv -j^{-1}\modu p,\quad  1\leq j\leq p-1, \quad  j\in \cup_{l=1}^k
[ \frac{l \times p}{k+1}, \frac{l\times p}{k}[\}.
\end{displaymath}
\begin{proof}
generalization of lemmas \ref{l21031a} p.\pageref{l21031a} and \ref{l21032a}
p.\pageref{l21032a}.
\end{proof}
\end{lem}
%
\begin{lem}\label{l21034a}
Let $J_k=\cup_{i=1}^k I(k)$.
Then
\begin{displaymath}
J_k=\{ i\quad |\quad u_i\equiv -j^{-1} \modu p,\quad  \frac{p}{k+1} <j \leq p-1 \}.
\end{displaymath}
\begin{proof}
Apply lemmas \ref{l21031a} p.\pageref{l21031a}, \ref{l21032a} p.\pageref{l21032a} and \ref{l21033a} p.\pageref{l21033a}.
More precisely,
we have
\begin{displaymath}
I(1)=\{i \quad |\quad u_i\equiv -j^{-1}\modu p, \quad j\in[\frac{p}{2},p-1]\}.
\end{displaymath}
and
\begin{displaymath}
I(2)=\{i \quad |\quad u_i\equiv -j^{-1}\modu p, \quad j\in[\frac{p}{3},\frac{p}{2}]\cup [\frac{2p}{3}, p-1]\}.
\end{displaymath}
Then, we get
\begin{displaymath}
I(1)\cup I(2)=
\{ i\quad |\quad u_i\equiv -\j^{-1}\modu p \quad |\quad j\in[\frac{p}{3},p-1]\}
\end{displaymath}
and so on.
\end{proof}
\end{lem}
%
\begin{lem}\label{l02041a}
Let $J_k=\cap_{i=1}^k I(k)$.
Then
\begin{displaymath}
J_k=\{ i \quad |\quad u_i\equiv -j^{-1}\modu p,\quad  \frac{k p}{(k+1)}<j\leq p-1  \}.
\end{displaymath}
\begin{proof}
We have
\begin{displaymath}
I(1)=\{ i \quad |\quad u_i\equiv -j^{-1}\modu p, \quad j\in[\frac{p}{2},p-1]\}.
\end{displaymath}
We have
\begin{displaymath}
I(2)=\{ i \quad |\quad u_i\equiv -j^{-1}\modu p, \quad j\in[\frac{p}{3},\frac{p}{2}]\cup [\frac{2p}{3}, p-1]\}.
\end{displaymath}
Then, we get
\begin{displaymath}
I(1)\cap I(2)=\{ i \quad |\quad u_i\equiv -\j^{-1}\modu p, \quad j\in[\frac{2p}{3},p-1]\}
\end{displaymath}
and so on.
\end{proof}
\end{lem}
%
\subsection{Some Kummer's principal ideals}
This subsection deals of cyclotomic fields independently of FLT.
Let us denote $I(d_1,\dots,d_k)=\cup_{i=1}^k I(d_i)$.
Let us denote $J(d_1,\dots,d_k)=\cap_{i=1}^k I(d_i)$.
Let $\mathbf q$ be a prime ideal of $\Z[\zeta]$ of inertial degree $1$.
Kummer showed that the ideal $\prod_{i\in I(d)}\sigma^i(\mathbf q),\quad d=1,\dots,p-2$, is a principal ideal of $\Z[\zeta]$.
Let $C_p$ be the class group of $\Q(\zeta)/\Q$. Let $C_p^+$ be the class group of $\Q(\zeta+\zeta^{-1})/\Q$. Let us consider the group $C_p^- = C_p/C_p^+$. In this section, we give some more general results in the group $C_p^- = C_p/C_p^+$.
Let us note $\mathbf q \simeq \Z[\zeta]$ when the ideal $\mathbf q$ is principal in $C_p^-$.
%
\begin{thm}  \label{t26031a}
Let $d_1,d_2 \in \N$ with  $1\leq d_1<d_2\leq p-2$; then
\begin{displaymath}
\prod_{i\in I(d_1)\cap I(d_2)} \sigma^i(\mathbf q)\simeq\Z[\zeta].
\end{displaymath}
\begin{proof}
The ideal  $P=\prod_{i\in I(d_1)} \sigma^i(\mathbf q)\times
\prod_{i\in I(d_2)}\sigma^i (\mathbf q)$ is principal.
Let $J(d_1,d_2)=I(d_1)\cap I(d_2)$.
Let $M(d_1,d_2)= (I(d_1)\cup I(d_2)) - J(d_1,d_2)$.
We get
\begin{displaymath}
P=\prod_{i\in J(d_1,d_2)} (\sigma^i(\mathbf q))^2 \times \prod_{i\in M(d_1,d_2)} \sigma^i(\s).
\end{displaymath}
If $i\in M(d_1,d_2)$ then:
\begin{itemize}
\item
else $i\in I(d_1)$ and $i\not\in I(d_2)$; in that case, from relation (\ref{e009161}) p.\pageref{e009161}, there exists $\varepsilon_i\in
\{-1,1\}$ with $1\leq i+\varepsilon_i\frac{p-1}{2}\leq p-1$ and $i+\varepsilon_i\frac{p-1}{2}\in I(d_2),\quad i+\varepsilon_i\frac{p-1}{2}\not\in I(d_1)$, therefore
$\{i,i+\varepsilon_i\frac{p-1}{2}\}\subset M(d_1,d_2)$.
\item
else $i\in I(d_2)$ and $i\not\in I(d_1)$; in that case there exists
$\varepsilon_i\in \{-1,1\}$ with $1\leq i+\varepsilon_i\frac{p-1}{2}\leq p-1$ and $i+\varepsilon_i\frac{p-1}{2}\in I(d_1),\quad i+\varepsilon_i\frac{p-1}{2}\not\in I(d_2)$, therefore
$\{i,i+\varepsilon_i\frac{p-1}{2}\}\subset M(d_1,d_2)$.
\end{itemize}
Therefore
\begin{displaymath}
P=\prod_{i\in J(d_1,d_2)} (\sigma^i(\mathbf q))^2 \times
\prod_{\{i,i+\varepsilon_i (p-1)/2\}\subset K(d_1,d_2)} \sigma^i(\mathbf q)\times\overline{\sigma^i(\mathbf q)}
\end{displaymath}
because $\sigma^{i+\varepsilon_i (p-1)/2} (\mathbf q )=\overline{\sigma^i(\mathbf q)}$.
Then we deduce that $P\cong \prod_{i\in J(d_1,d_2)} (\sigma^i(\mathbf q))^2$ in the class group $C_p^-$
and finally that $\prod_{i\in J(d_1,d_2)} \sigma^i(\mathbf q)^2$ is principal in $C_p^-$, which leads to the result.
\end{proof}
\end{thm}
%
\begin{lem}\label{l01041a} { ****}
Let $k\in \N,\quad 1 \leq k \leq p-2$.
Let $d_l\in\N,\quad l=1,\dots,k,\quad 1\leq d_1<d_2<\dots < d_k\leq p-2$.
Then
\begin{equation}\label{e01041a}
\prod_{i\in I(d_1,\dots, d_k)} \sigma^i(\mathbf q)
\simeq \prod_{i\in J(d_1,\dots, d_k)} \sigma^i(\mathbf q).
\end{equation}
\begin{proof}
Let any $i\in I(d_1,\dots,d_k)$. There exists $j_i\in \N,\quad  1\leq j_i\leq k,$ such that, with a certain reordering of $d_1,\dots,d_k$, and without loss of generality,
\begin{displaymath}
\begin{split}
& i\in I(d_1),\dots  i\in I(d_{j_i}),\\
& i\not\in I(d_{j_i+1}),\dots i\not\in I(d_k).
\end{split}
\end{displaymath}
\begin{itemize}
\item
Suppose at first that $j_i<k$ :
then there exists one and only one $\varepsilon_i\in\{-1,1\}$ with $1\leq i+\varepsilon_i\frac{p-1}{2}\leq p-1$ and with $i+\varepsilon_i\frac{p-1}{2}\in I(d_{j_i+1})$.
Therefore, $\{i,i+\varepsilon_i\frac{p-1}{2}\}\subset I(d_1,\dots,d_k)$.
Therefore, the term $\sigma^i(\mathbf q)\times\sigma^{i+\varepsilon_i (p-1)/2}(\mathbf q)
=\sigma^i(\mathbf q)\times \overline{\sigma^i(\mathbf q)}$ which is principal in $C_p^-$ can be omitted in the computation of class  of $\prod_{i\in I(d_1,\dots,d_k)} \sigma^i(\mathbf q)$ in $C_p^-$.
\item
Suppose then that $j_i=k$: it means that $i\in I(d_1),\dots, i\in I(d_k)$ and so
$i\in J(d_1, \dots,d_k)$, which achieves the proof.
\end{itemize}
\end{proof}
\end{lem}
%


\subsection {On  $p$-rank of the matrix of Kummer's system}\label{s27031a}
Let us consider $P_d=\prod_{i\in I_{u,d}} \sigma^i(\s)$. With the language of representations of section \ref{s09121} p.\pageref{s09121}, it means that, if $X\cong Cl(\s)$, then we have the system of relations
\begin{displaymath}
P_i(X)=\sum_{i\in I(d)} X^i=0,\quad d=1,\dots,p-2.
\end{displaymath}
This system can be reduced to the system
\begin{displaymath}
\begin{split}
& P_i(X)=\sum_{i=1}^{(p-1)/2} \varepsilon_{i,d} X^i=0,\quad d=1,\dots,p-2,\\
& \varepsilon_{i,d} =1 \Leftrightarrow i\in I(d),\\
& \varepsilon_{i,d}=-1\Leftrightarrow i\pm (p-1)/2\in I(d).
\end{split}
\end{displaymath}
A direct computation show that $\sum_{i=1}^{(p-1)/2} \varepsilon_{i,d} X^i$ is identically similar to $\sum_{i=1}^{(p-1)/2} \varepsilon_{i,p-1-d} X^i$. Observe also that
$\sigma^{i+(p-1)/2}(\s)=\overline{\sigma^i(\s)}$ and that $\sigma^i(\s)\times\overline{\sigma^i(\s)}$ is principal in $C_p^-$.
Therefore the Kummer's $p-2$ relations can {\it canonically} be reduced to a system of
$\frac{p-1}{2}$ relations :
\begin{equation}\label{e27031a}
\begin{split}
& P_i(X)=\sum_{i=1}^{(p-1)/2} \varepsilon_{i,d} X^i =0,\quad d=1,\dots,\frac{p-1}{2}, \\
& \varepsilon_{i,d}=1
\Leftrightarrow (u_{(p-1)/2-i} +(d\times u_{(p-1)/2-i} \modu p))>p, \\
& \varepsilon_{i,d}=-1
\Leftrightarrow (u_{-i} +(d\times u_{-i} \modu p))>p.
\end{split}
\end{equation}
The square matrix $M_p= (\varepsilon_{i,d})_{1\leq d,i \leq (p-1)/2}$  is then well defined : we shall study this matrix carefully.
First we give as an example from a MAPLE program the value of the matrix and of its determinant for $p=13$.
\begin{displaymath}
M_p\equiv
\begin{array}{|rrrrrr|}
-1& -1& 1& -1& -1& -1\\
1& -1& -1& -1& -1& -1\\
-1& -1& 1& 1& -1& -1\\
1& 1& 1& -1& -1& -1\\
-1& -1& -1& -1& -1& -1\\
1&-1& 1& 1& 1& -1\\
\end{array}
\modu p
\end{displaymath}
From a MAPLE program, we have computed the $p$-rank of $M_p$ for $2<p<200$.
The results are given in the following table
\begin{displaymath}
\begin{array}{ |r|r|r|r|r|r|r|r|r|r|}
p&k_p&p&k_p&p&k_p&p&k_p&p&k_p \\
 3&1& 5&2& 7&3& 11&5& 13&6\\
17&8&19&9& 23&11& 29&14& 31&15\\
{\bf 37}&{\bf 17}&41&20&43&21& 47&23& 53&26\\
{\bf 59}&{\bf 28}& 61&30&{\bf 67}&{\bf 32}& 71&35&73&36\\
 79&39& 83&41&89&44&97&48&{\bf 101}&{\bf 49}\\
{\bf 103}&{\bf 50}& 107&53&109&54&113&56& 127&63\\
 {\bf 131}&{\bf 64}&137&68& 139&69&{\bf 149}&{\bf 73}&151&75\\
{\bf 157}&{\bf 76}&163&81&167&83&173&86&179&89\\
181&90&191&95&193&96& 197&98&199&99\\
\end{array}
\end{displaymath}
All the $p$-rank $k_p<(p-1)/2$ bold printed correspond to {\bf irregular} primes. In this table, the irregular prime  all have rank $\frac{p-1}{2}-1$ except the prime $p=157$ with $k_p=\frac{p-1}{2}-2$.Observe that these results are in accordance with table of relative class numbers in Washington \cite{was}, p 412.

We can now give, as an application to FLT, the following theorem:
%
\begin{thm} {A criterion on $M_p$ for First case of FLT}\label{t28031a}

If the first case of FLT would fail for $p$, then the $p$-rank $k_p$ of the  matrix $M_p$ should verify
\begin{equation}\label{e28031a}
k_p<\frac{p-1}{2}-\sqrt{p}-1.
\end{equation}
\begin{proof}
Suppose that $k_p>\frac{p-1}{2}-\sqrt{p}-1$ and search for a contradiction:
Let us consider the $\frac{p-1}{2}$ relations, with $X^i\cong Cl(\sigma^i(\s)$,
\begin{displaymath}
\sum_{i\in I(d)} \varepsilon_{i,d} X^i=0,\quad d=1,\dots, p-2.
\end{displaymath}
It would be possible to find, by elimination, a polynomial $\mathbf P_d(X)$ of degree $d<\sqrt{p}-1$ with $\mathbf P_d(X)=0$; from theorem \ref{t31101a} p.\pageref{t31101a} we should have $d\geq r_1$, degree of the representation $\rho_1 :  G\rightarrow {\bf GL}({\bf F}_p.V_1)$ and also degree of the minimal associate polynomial $P_{r_1}(X)$, because we should have $\mathbf P_d(X)=P_{r_1}(X)\times Q(X),\quad Q(X)\in\bf F_p[X]$. But, from Eichler's, we deduce that $r_1>\sqrt{p}-1$.
\end{proof}
\end{thm}
%
\subsection{Generalizations of Kummer congruences - I}\label{s23091}
This subsection is in the Fermat's equation context.
\begin{itemize}
\item
Let $q\in\N,\quad q\equiv 1\modu p$, be a prime.

Let $\mathbf q\subset\Z[\zeta]$ be a prime ideal above $q$. Recall that the relation (2.5) p119 of \cite{rib} says that
$\prod_{i\in I(d)} \sigma^i(\mathbf q)$
is a principal integral ideal of $\Z[\zeta]$.
\item
Using this relation for all the prime ideals, all different of $\pi=(1-\zeta)\Z[\zeta]$, dividing $(x+\zeta y)$, We can generalize for the ideal $\s,\quad \s^p=(x+\zeta y)\Z[\zeta]$, (First and second case of Fermat equation) to,
\begin{equation}\label{eqkum02}
P_d=\prod_{i\in I(d)} \sigma_i(\s)
\end{equation}
is a principal ideal.
All the results of the sections \ref{s3} p.\pageref{s3}, \ref{s4} p.\pageref{s4}  start from the existence of the principal ideal $\prod_{i=1}^f \sigma_i(\s)^{l_i}$. The sketch of proof of this section of the article is exactly the same with, here, the principal ideal of the relation
(\ref{eqkum02}), $P_d=\prod_{i\in I(d)} \sigma_i(\s)$.
Using this equation and results of previous sections, we shall give several generalizations   of Kummer congruences $\modu p$.
\end{itemize}
%
\begin{prop}\label{thmkum01}
Let $u\in\N$ be a primitive root $\modu p$.
Let $d\in\N,\quad 0\leq d\leq p-2$.
Then there exists $v\in\N, \quad 0\leq v\leq p-1$ such that :
\begin{equation}\label{eqkum03}
-2v+\sum_{i\in I(d)}\frac{u_i}{(x+\zeta^{u_i}y)}
+\sum_{i\in I(d)}\frac{u_i}{(x+\zeta^{-u_i} y)}\equiv 0 \modu \pi^{p-1}.
\end{equation}
\begin{proof}
From relation (\ref{eqkum02}) p.\pageref{eqkum02} we deduce, similarly to proposition \ref{p4} p.\pageref{p4},  that
\begin{displaymath}
 \prod_{i\in I(d)} (x+\zeta^{u_i} y) =\zeta^m\eta\gamma^p,
\quad \eta\in\Z[\zeta+\zeta^{-1}]^*,\quad \gamma\in\Z[\zeta],
\end{displaymath}
then, by conjugation,
\begin{displaymath}
 \prod_{i\in I(d)} (x+\zeta^{-u_i} y) =\zeta^{-m}\eta\overline{\gamma^p},
\quad \eta\in\Z[\zeta+\zeta^{-1}]^*,\quad \overline{\gamma}\in\Z[\zeta]\,
\end{displaymath}
hence
\begin{displaymath}
\zeta^{-m}\prod_{i\in I(d)} (x+\zeta^{u_i} y)
\equiv \zeta^{m}\prod_{i\in I(u,d)} (x+\zeta^{-u_i} y)
\modu \pi^{p+1}.
\end{displaymath}
Then, similarly to proposition \ref{p9} p.\pageref{p9}, we deduce that:
\begin{displaymath}
\begin{split}
&(-m) \zeta^{-m-1}\prod_{i\in I(d)} (x+\zeta^{u_i} y)
+\zeta^{-m}\prod_{i\in I(d)} (x+\zeta^{u_i} y)\sum_{i\in I(u,d)}\frac{u_i \zeta^{u_i-1}y}
{(x+\zeta^{u_i} y)} \equiv\\
& (m)\zeta^{m-1}\prod_{i\in I(d)} (x+\zeta^{-u_i} y)
+\zeta^{m}\prod_{i\in I(d)} (x+\zeta^{-u_i}y)\sum_{i\in I(d)}
\frac{(-u_i) \zeta^{-u_i-1} y}{(x+\zeta^{-u_i} y)} \modu \pi^{p-1},
\end{split}
\end{displaymath}
then, multiplying by $\zeta$,
\begin{displaymath}
\begin{split}
&(-m)\zeta^{-m}\prod_{i\in I(d)} (x+\zeta^{u_i} y)
+\zeta^{-m}\prod_{i\in I(d)} (x+\zeta^{u_i} y)\sum_{i\in I(d)}
\frac{u_i \zeta^{u_i}y} {(x+\zeta^{u_i}y)} \equiv\\
&(m)\zeta^{m}\prod_{i\in I(d)} (x+\zeta^{-u_i} y)
+\zeta^{m}\prod_{i\in I(d)} (x+\zeta^{-u_i}y)\sum_{i\in I(d)}
\frac{(-u_i) \zeta^{-u_i} y}{(x+\zeta^{-u_i} y)} \modu \pi^{p-1},
\end{split}
\end{displaymath}
so
\begin{displaymath}
(-m)+\sum_{i\in I(d)}\frac{u_i \zeta^{u_i} y}{(x+\zeta^{u_i}y)}\equiv
 (m)+\sum_{i\in I(d)}\frac{(-u_i) \zeta^{-u_i} y}{(x+\zeta^{-u_i} y)}\modu \pi^{p-1},
\end{displaymath}
also
\begin{displaymath}
(-2m)+\sum_{i\in I(d)}\frac{u_i \zeta^{u_i} y}{(x+\zeta^{u_i}y)}
+\sum_{i\in I(d)}\frac{u_i \zeta^{-u_i} y}{(x+\zeta^{-u_i} y)}\equiv 0\modu \pi^{p-1},
\end{displaymath}
hence
\begin{displaymath}
\begin{split}
& -2m+2( \sum_{i\in I(d)} u_i)\frac{y}{(x+y)}\\
&-x\times(\sum_{i\in I(d)} \frac{u_i}{(x+\zeta^{u_i} y)}
+\sum_{i\in I(d)}\frac{u_i}{(x+\zeta^{-u_i} y)})\equiv 0\modu \pi^{p-1}.
\end{split}
\end{displaymath}
\end{proof}
\end{prop}
%
\begin{thm} { *** }\label{t31101}
Let $\phi_{2m+1}(T),\quad m=1,\dots,\frac{p-3}{2}$ be the odd Mirimanoff polynomials.
Let $t\equiv -\frac{y}{x}\modu p$.
Let $r$ be the value defined in theorem (\ref{pa212}).
Let $C_M=\{m_1,m_2,\dots,m_{r-1}\}\subset \{1,2,\dots,\frac{p-3}{2}\}$ be the set such that $m\in C_M \Longleftrightarrow \phi_{2m+1}(t)\not\equiv 0 \modu p$, see corollary \ref{cor1209} p.\pageref{cor1209}.
Let $u\in\N$ be any primitive root $\modu p$. Let $d\in\N,\quad 1\leq d\leq p-2$.
Let $I(u,d)=\{i,\quad 1\leq i\leq p-1,\quad (u_{\nu-i}+u_{\nu-i+ind_u(d)})>p\}$.

Then for all $m\in C_M$, we have the congruence :
\begin{displaymath}
\sum_{i\in I(d)} u_i^{2m+1}\equiv 0 \modu p.
\end{displaymath}
\begin{proof}
From proposition \ref{thmkum01} p.\pageref{thmkum01} and similarly to theorem \ref{p21a} p.\pageref{p21a}, we have for
$1\leq m\leq \frac{p-3}{2}$ the congruence
$\phi_{2m+1}(t)\times \sum_{i\in I(d)}u_i^{2m+1}\equiv 0 \modu p$.
\end{proof}
\end{thm}
%
\begin{cor} \label{c21091}

Let $\mu\in\{3,5,7\}$.
Let $d\in\N,\quad 1\leq d\leq p-2$.

Let $I(u,d)=\{i,\quad 1\leq i\leq p-1,\quad (u_{\nu-i}+u_{\nu-i+ind_u(d)})>p\}$.
Then, for each primitive root $u \modu p$, we have
\begin{displaymath}
\sum_{i\in I(d)} u_i^\mu\equiv 0 \modu p.
\end{displaymath}
\begin{proof}
From \cite{rib}, relation (1.12) p143, we get
\begin{displaymath}
\phi_{2m+1}(t)\equiv -(1-t)^{p-1-2m} \times P_{2m+1}(-t)\modu p.
\end{displaymath}
where $P_{2m+1}(T)$ is the Kummer polynomial corresponding to the Mirimanoff polynomial $\phi_{2m+1}(T)$.
From \cite{rib} p125 and p 140, we have
\begin{displaymath}
\begin{split}
&P_3(-t)=-t(1+t),\\
&P_5(-t)=-t(1+t)(1+10t+t^2),\\
&P_7(-t)=-t(1+t)(1+56t+246t^2+56t^3+t^4).
\end{split}
\end{displaymath}

Let $G=\{-\frac{x}{y},-\frac{y}{x},-\frac{x}{z},-\frac{z}{x},-\frac{x}{z},
-\frac{y}{z},\frac{z}{y}\}\modu p$.
It is possible to prove by direct computation that:
\begin{itemize}
\item
$P_3(-t)\equiv 0 \modu p$ is not simultaneously possible for all $t\in G$,
\item
$P_5(-t)\equiv 0 \modu p$ is not simultaneously possible for all $t\in G$,
\item
$P_7(-t)\equiv 0 \modu p$ is not simultaneously possible for all $t\in G$.
\end{itemize}
(see \cite{rib} p 140, (1A)).
Then , this corollary results immediatly from theorem \ref{t31101} p.\pageref{t31101}.
\end{proof}
\end{cor}
%
{\bf Remark:} for $\mu>3$, the proof use two or three of the Barlow-Abel relations
\begin{displaymath}
\frac{x^p+y^p}{x+y}=t_1^p,\quad
 \frac{y^p+z^p}{x+y}=r_1^p,\quad
 \frac{z^p+x^p}{x+y}=s_1^p
\end{displaymath}
and then we cannot apply directly this result to Generalized Fermat-Wiles equation, see section \ref{s08081} p.\pageref{s08081}.
%

\begin{cor} \label{c22091}
Let $\mu\in\{3,5,7\}$.

If $p$ does not divide the sum
\begin{displaymath}
S=\sum_{j\in [\frac{p}{3},\frac{p}{2}]\cup [\frac{2p}{3},p-1]} \frac{1}{j^\mu},
\end{displaymath}
then the First case of Fermat's Last Theorem holds for the exponent $p$.
\begin{proof} $ $
\begin{itemize}
\item
Let $u^{ind_u(2)}\equiv 2\mod p$.
We have
\begin{displaymath}
I(d)=\{i,\quad 0\leq i\leq p-2,\quad (u_{\nu-i}+(2u_{\nu-i}\modu p))>p\}.
\end{displaymath}
Therefore,
We have
\begin{displaymath}
K(d)=\{k,\quad 1\leq k\leq p-1,\quad (u_{\nu+k}+(2u_{\nu+k}\modu p))>p\}.
\end{displaymath}
\item
If $u_{\nu+k}<\frac{p}{2}$ then $(2u_{\nu+k}\modu p)=2u_{\nu+k}$
leads to

$u_{\nu+k}>\frac{p}{3}$, hence $\frac{p}{3} < j=u_{\nu+k} <\frac{p}{2}$.
\item
If $u_{\nu+k}>\frac{p}{2}$ then $(2u_{\nu+k}\modu p)=2u_{\nu+k}-p$
leads to
$u_{\nu+k}>\frac{2p}{3}$ then to $\frac{2p}{3} <j=u_{\nu+k} \leq p-1$.
\item
But $u_i^\mu\equiv u^{i\mu} \equiv
(u^{-1})^{-i\mu}\equiv (u^{-i})^{-\mu}$

so
$u_i^\mu\equiv -(u_{\nu-i})^{-\mu}\equiv -j^{-\mu}\modu p$, because $u^\nu\equiv -1\modu p$.
\item
the result is then an immediate consequence of the previous proposition \ref{c21091} p.\pageref{c21091}.
\end{itemize}
\end{proof}
\end{cor}
%
\begin{cor} \label{c09081}$ $
Let $\mu\in\{3,5,7\}$.
If $p$ does not divide the sum
\begin{displaymath}
S=\sum_{j\in[\frac{p}{4},\frac{p}{3}]\cup [\frac{2p}{4},\frac{2p}{3}]\cup
[\frac{3p}{4}, p-1]}\frac{1}{j^\mu},
\end{displaymath}
then the First case of Fermat's Last Theorem holds for the exponent
$p$.
\begin{proof} $ $
\begin{itemize}
\item
Let  $d=3$. We have $I(d)=\{i,\quad 0\leq i\leq p-2,
\quad (u_{\nu-i}+(3u_{\nu-i}\modu p))>p\}$.
\item
If $u_{\nu-i}<\frac{p}{3}$ then $(3u_{\nu-i}\modu p)=3u_{\nu-i}$
so
$u_{\nu-i}>\frac{p}{4}$,
hence
$\frac{p}{4} < j=u_{\nu-i} <\frac{p}{3}$.
\item
If $u_{\nu-i}<\frac{2p}{3}$ and $u_{\nu-i}>\frac{p}{3}$
then $(3u_{\nu-i}\modu p)=3u_{\nu-i}-p$
so
$u_{\nu-i}>\frac{2p}{4}$,
hence $\frac{2p}{4} <j=u_{\nu-i} < \frac{2p}{3}$.
\item
If $u_{\nu-i}>\frac{2p}{3}$ then
$3u_{\nu-i} \modu p = 3u_{\nu-i}-2p$
so
$u_{\nu-i}>\frac{3p}{4}$,
hence
$\frac{3p}{4}<j=u_{\nu-i}\leq p-1$.
\item
But $u_i^\mu\equiv u^{i\mu} \equiv
(u^{-1})^{-i\mu}\equiv (u^{-i})^{-\mu}$
so
$u_i^\mu\equiv -(u_{\nu-i})^{-\mu}\equiv -j^{-\mu}\modu p$ because $u^\nu\equiv -1\modu p$.
\item
the result is then an immediate consequence of the  proposition \ref{c21091} p.\pageref{c21091}.
\end{itemize}
\end{proof}
\end{cor}
%

\begin{cor}  \label{c03104}
Let $\mu\in\{3,5,7\}$.
If $p$ does not divide the sum $\sum_{j\in[1,\frac{p}{4}]} \frac{1}{j^\mu}$ then the first case of Fermat's Last Theorem holds for the exponent $p$.
\begin{proof}
If the First case of Fermat fails, then, from the two previous propositions by addition of congruences,
\begin{displaymath}
\sum_{j>p/3}^{j<p/2} j^{-\mu}
+\sum_{j>2p/3}^{j= p-1} j^{-\mu}
+\sum_{j>p/4}^{j<p/3} j^{-\mu}+\sum_{j>p/2}^{j<2p/3} j^{-\mu}
+\sum_{j>3p/4}^{j= p-1}j^{-\mu}\equiv 0 \modu p,
\end{displaymath}
so, by permutation of the sums,
\begin{displaymath}
\sum_{j>p/4}^{j<p/3} j^{-\mu}+\sum_{j>p/3}^{j<p/2} j^{-\mu}+\sum_{j>p/2}^{j<2p/3} j^{-\mu}
+\sum_{j>2p/3}^{j= p-1} j^{-\mu}
+\sum_{j>3p/4}^{j= p-1}j^{-\mu}\equiv 0 \modu p,
\end{displaymath}
so grouping sums,
\begin{displaymath}
\sum_{j>p/4}^{j=p-1} j^{-\mu}+\sum_{j>3p/4}^{j= p-1} j^{-\mu}\equiv 0 \modu p,
\end{displaymath}
also, noticing that $\sum_{j=1}^{p-1}j^{-\mu}\equiv 0 \modu p$,
\begin{displaymath}
\sum_{j>3p/4}^{j=p-1} j^{-\mu} -\sum_{j=1}^{j<p/4} j^{-\mu}\equiv 0 \modu p,
\end{displaymath}
and
\begin{displaymath}
\sum_{j=1}^{j<p/4} j^{-\mu}\equiv 0\modu p.
\end{displaymath}
which completes the proof.
\end{proof}
\end{cor}
%
\begin{cor} { *** $3(p-2)$ congruences $\modu p.$ }\label{c2009}

Let $\mu\in\{3,5,7\}$.
Let $d\in\N, \quad 1\leq d \leq p-2$.
If $p$ does not divide the sum
\begin{displaymath}
S= \sum_{j\in\cup_{l=1}^d [\frac{lp}{d+1},\frac{lp}{d}[} \frac{1}{j^\mu}
\end{displaymath}
then the first case of Fermat's Last Theorem holds for the exponent $p$.
\begin{proof}
Same proof  as \ref{c22091} p.\pageref{c22091}.
\end{proof}
\end{cor}
%

{\bf Remarks :}
\begin{itemize}
\item
These results generalize  the Cauchy congruence

$\sum_{j=1}^{(p-1)/2} \frac{1}{j^3}\equiv 0 \modu p$ (see for instance \cite{cau} or \cite{rib} (4A) p 120)
and the other summation criteria of the bibliography (see the survey in section\ref{s06101} p.\pageref{s06101}).
Corollary (\ref{c2009}) p.\pageref{c2009} generalizes  Cauchy congruences depending only on $p$ to $3(p-2)$ explicit congruences depending only on $p$ for $k=1,\dots,p-2$.
\item
With the meaning of $r$ and $C_M$ given in theorem \ref{t31101} p.\pageref{t31101}, we can generalize these congruences involving $\sum \frac{1}{j^{\mu}}$ for $\mu\in\{3,5,7\}$ to all sums of the form
$\sum_j \frac{1}{j^{2m+1}}$ for $m\in C_M$.

In this direction, we obtain the result :

Let $m\in C_M =\{m\quad |\quad \phi_{2m+1}(t)\not\equiv 0 \modu p\}$.
Let $d\in\N, \quad 1\leq d \leq p-2$. If $p$ does not divide the sum
$\quad \sum_{l=1}^d \sum_{j>(lp)/(d+1)}^{j<(lp)/d} \frac{1}{j^{2m+1}}$ then the first case of Fermat's Last Theorem holds for the exponent $p$. The previous results being a particular case of this theorem.
\item
Note that our proof does not use the Bernoulli numbers $B_k$. It is also possible to obtain a proof connected to Bernoulli Numbers : From proposition (5B) p 108 in \cite{rib}, derived of a result of Vandiver \cite{va3} resting on class field theory, we have for $d\in\N,\quad d=2,\dots,p-2$, the congruences:
\begin{displaymath}
(1-d^{-2m})B_{p-2m-1}\equiv
(p-2m-1)d^{-2m-1}\sum_{l=1}^{d-1}\sum_{1\leq j<\frac{lp}{d}} \frac{1}{j^{2m+1}}\modu p.
\end{displaymath}
Then, we get in the same way, with $d+1$ in place of $d$, the second relations
\begin{displaymath}
 (1-(d+1)^{-2m})B_{p-2m-1}\equiv
(p-2m-1)(d+1)^{-2m-1}\sum_{l=1}^{d}\sum_{1\leq j<\frac{lp}{(d+1)}} \frac{1}{j^{2m+1}}\modu p.
\end{displaymath}
From Ribenboim \cite{rib}, (4E) p 125 and  (1A) p 140, we have $B_{p-3}\equiv B_{p-5}\equiv B_{p-7}\modu p$. Then,  by difference of the two previous congruences for $m=1,m=2,m=3$, we get the results of the corollary \ref{c2009} p.\pageref{c2009}.
\end{itemize}
From remark and corollary \ref{c09113} p.\pageref{c09113} , we state the theorem :
%
\begin{thm}{***} \label{t22021}

Let $r_p$ be the $p-$rank of the class group of $\Q(\zeta)$.
Let $C_M=\{m\quad |\quad 1\leq m\leq \frac{p-3}{2},\quad\phi_{2m+1}(t)\not\equiv 0 \modu p\}$.
If the first case of FLT fails for $p$ then
\begin{itemize}
\item
$3\leq Card(C_M)\leq r_p$ and $\{1,2,3\}\subset C_M$.
\item $Card(C_M)\times (p-2)$ congruences :

For every $m\in C_M$ and every $d\in\N,\quad 1\leq d\leq p-2$,
\begin{displaymath}
\sum_{j\in\cup_{l=1}^d [\frac{lp}{d+1},\frac{lp}{d}[} \frac{1}{j^{2m+1}}\equiv 0 \modu p
\end{displaymath}
\end{itemize}
\end{thm}
%
\begin{cor} { $[Log(p)^{1/3}](p-2)$ congruences $\modu p$ }\label{c15031}

Let $n_0=(45!)^{88}$.
Let $m\in\N, \quad m=1,\dots,[Log(p)^{1/3}]$.
Let $d\in\N, \quad 1\leq d \leq p-2$.
If $p$ does not divide the sum
\begin{displaymath}
S=\sum_{j\in\cup_{l=1}^d [\frac{lp}{d+1},\frac{lp}{d}[} \frac{1}{j^{2m+1}}
\end{displaymath}
then the first case of Fermat's Last Theorem holds for the exponent $p$.
\begin{proof}
Is an immediate consequence of \cite{kra}, of the Vandiver congruences \cite{va3} and of a previous remark between our results and Vandiver \cite{va3}.
\end{proof}
\end{cor}
%
\subsection{ Bernoulli Numbers and
sums $\sum_{ j\in\cup_{l=1}^d[\frac{l p}{d+1},\frac{ l p}{d}[} \frac{1}{j^{2m+1}}$.}
This subsection deals with a connection, in a general context independant of FLT, between even  Bernoulli numbers $B_{p-1-2m}$ and sum congruences
$ \sum_{j\in\cup_{l=1}^d [\frac{l p}{d+1},\frac{l p}{d}]}\frac{1}{j^{2m+1}}$ for $d=1,\dots,p-2$.

Let $S(d-1)=\sum_{l=1}^{d-1}\sum_{1\leq j<lp/d} \frac{1}{j^{2m+1}}\modu p$.
From Vandiver in Ribenboim \cite{rib}, (5B) p 108, we get for $m=1,\dots,\frac{p-3}{2}$
\begin{displaymath}
\begin{split}
& S(d-1)\equiv \frac{d\times (d^{2m}-1)\times B_{p-1-2m}}{p-1-2m}\modu p,\\
& S(d)\equiv \frac{(d+1)\times ((d+1)^{2m}-1)\times B_{p-1-2m}}{p-1-2m}\modu p.
\end{split}
\end{displaymath}
Then we have
\begin{displaymath}
S(d)\equiv \sum_{l=1}^{d-1}\sum_{1\leq j<l p/(d+1)}\frac{1}{j^{2m+1}}
+\sum_{1\leq j< d p/(d+1)} \frac{1}{j^{2m+1}} \modu p.
\end{displaymath}
This leads to
\begin{displaymath}
S(d)\equiv \sum_{l=1}^{d-1} (\sum_{1\leq j < l p/d}\frac{1}{j^{2m+1}}
-\sum_{j >l p/(d+1)}^{j< l p/d} \frac{1}{j^{2m+1}})
+\sum_{1\leq j < d p/(d+1)} \frac{1}{j^{2m+1}}\modu p,
\end{displaymath}
and so
\begin{displaymath}
S(d)\equiv  S(d-1)
-\sum_{l=1}^{d-1}\sum_{j >l p/(d+1)}^{j< l p/d} \frac{1}{j^{2m+1}}
+\sum_{1\leq j < d p/(d+1)} \frac{1}{j^{2m+1}}\modu p,
\end{displaymath}
also
\begin{displaymath}
\begin{split}
&S(d)\equiv  S(d-1)
-\sum_{l=1}^{d}\sum_{j >l p/(d+1)}^{j< l p/d} \frac{1}{j^{2m+1}}\\
&+\sum_{d p/(d+1) <j< p}\frac{1}{j^{2m+1}}
+\sum_{1\leq j < d p/(d+1)} \frac{1}{j^{2m+1}}\modu p,
\end{split}
\end{displaymath}
also
\begin{displaymath}
\begin{split}
&S(d)\equiv  S(d-1)
-\sum_{l=1}^{d}\sum_{j >l p/(d+1)}^{j< l p/d} \frac{1}{j^{2m+1}}\\
&+\sum_{1\leq j< p}\frac{1}{j^{2m+1}} \modu p,
\end{split}
\end{displaymath}
and finally
\begin{displaymath}
S(d)\equiv S(d-1)-\sum_{l=1}^{d}\sum_ {l p/(d+1)<j}^{j< l p /d}\frac{1}{j^{2m+1}}\modu p.
\end{displaymath}
Then with previous relations obtained with Bernoulli numbers we get
%
\begin{thm} \label{t12031a}
For $d\in\N,\quad 1\leq d\leq p-2$ and for $m=1,\dots,\frac{p-3}{2}$ , we have  for the even Bernoulli Numbers
\begin{equation}\label{e12031a}
\sum_{j\in\cup_{l=1}^d[\frac{ l p}{d+1},\frac{ l p}{d}[}\frac{1}{j^{2m+1}}\equiv
\frac{d(d^{2m}-1)-(d+1)((d+1)^{2m}-1)}{p-1-2m}\times B_{p-1-2m}\modu p
\end{equation}
\end{thm}
%
{\bf Remarks:}
\begin{itemize}
\item
Observe that these general $p-2$ congruences without any connection with FLT  are a  tool to test $B_{p-1-2m}\equiv 0 \modu p$.

\item
We verified numerically this congruence with a Maple Program, which is an indirect insurance that our previous results are correct.
\item
Observe that else $\frac{1}{j^{2m+1}}$ else $\frac{1}{(p-j)^{2m+1}}$ intervenes in the sum in the first member of the relation (\ref{e12031a}) p.\pageref{e12031a} which corresponds to
$i+\frac{p-1}{2}\not\in I(u,d)$: this shows that the class of the ideal
$\sigma^i(\s)\sigma^{i+(p-1)/2}(s)=\sigma^i(\s)\overline{\sigma^i(s)}$ of
$\Q(\zeta+\zeta^{-1})$ does not intervenes in this computation of $B_{p-1-2m}$ and that
$B_{p-1-2m}\equiv 0 \modu p$ caracterizes $h^-$ and finally that $p$ irregular implies that
$p | h^-$.
\end{itemize}
%
These congruences can be explicitly written in the form
\begin{equation}\label{e15031a}
\begin{split}
&\gamma_{1,2}\frac{1}{2^{2m+1}}+\dots+\gamma_{1,j}\frac{1}{j^{2m+1}}+\dots+
\gamma_{1,p-1}\frac{1}{(p-1)^{2m+1}}\equiv 0 \modu p,\\
&\vdots\\
&\gamma_{d,2}\frac{1}{2^{2m+1}}+\dots+\gamma_{d,j}\frac{1}{j^{2m+1}}+\dots+
\gamma_{d,p-1}\frac{1}{(p-1)^{2m+1}}\equiv 0 \modu p,\\
&\vdots\\
&\gamma_{p-2,2}\frac{1}{2^{2m+1}}+\dots+\gamma_{p-2,j}\frac{1}{j^{2m+1}}+\dots+
\gamma_{p-2,p-1}\frac{1}{(p-1)^{2m+1}}\equiv 0 \modu p,\\
&\gamma_{d,j}\in \{0,1\},\quad d=1,\dots,p-2,\quad j=2,\dots,p-1,\\
&\gamma_{d,j}=1 \ {\mbox if}\  \exists l\in\N,\quad  \frac{l p}{(d+1)}< j < \frac{ l p}{d},\\
&\gamma_{d,j}=0 \ {\mbox if}\  \not\exists l\in\N,\quad  \frac{l p}{(d+1)}< j < \frac{ l p}{d}.
\end{split}
\end{equation}
We have the property
\begin{displaymath}
Card\{ \gamma_{d,j}\quad |\quad j=2,\dots,p-1,\quad \gamma_{d,j}\not=0 \}=\frac{p-1}{2},
\end{displaymath}
because it is not possible to have simultaneously
\begin{displaymath}
\frac{ l_1 p}{(d+1)}<j<\frac{l_1 p}{d}, \quad \frac{l_2 p}{(d+1)}< p-j<\frac{l_2 p}{d};
\end{displaymath}
If not we should have $\frac{(l_1+l_2) p}{d+1}<p<\frac{(l_1+l_2)p}{d}$ and so
$\frac{l_1+l_2}{d+1}<1<\frac{l_1+l_2}{d}$ so $l_1+l_2\leq d$ and $(d+1)\leq l_1+l_2$ contradiction.
Then observing that $\frac{1}{j^{2m+1}}\equiv j^{p-1-(2m+1)}\equiv -(p-j)^{p-1-(2m+1)}
\equiv -\frac{1}{(p-j)^{2m+1}}\modu p$, from linear system (\ref{e15031a}) p.\pageref{e15031a}
we get the reduced system
\begin{equation}\label{e15032a}
\begin{split}
&\delta_{1,1}\frac{1}{1^{2m+1}}+\dots+\delta_{1,j}\frac{1}{j^{2m+1}}+\dots+
\delta_{1,(p-1)/2}\frac{1}{((p-1)/2)^{2m+1}}\equiv 0 \modu p,\\
&\vdots\\
&\delta_{d,1}\frac{1}{1^{2m+1}}+\dots+\delta_{d,j}\frac{1}{j^{2m+1}}+\dots+
\delta_{d,(p-1)/2}\frac{1}{((p-1)/2)^{2m+1}}\equiv 0 \modu p,\\
&\vdots\\
&\delta_{(p-1)/2,1}\frac{1}{2^{2m+1}}+\dots+\delta_{(p-1)/2,j}\frac{1}{j^{2m+1}}+\dots+
\delta_{(p-1)/2,(p-1)/2}\frac{1}{((p-1)/2^{2m+1}}\equiv 0 \modu p,\\
&\delta_{d,j}\in \{-1,1\},\quad d=1,\dots,(p-1)/2,\quad j=1,\dots,(p-1)/2,\\
&\delta_{d,j}=1 \ {\mbox if}\  \exists l\in\N,\quad  \frac{l p}{(d+1)}< j < \frac{ l p}{d},\\
&\delta_{d,j}=-1 \ {\mbox if}\  \not\exists l\in\N,\quad
\frac{l p}{(d+1)}< j < \frac{ l p}{d}.
\end{split}
\end{equation}
If the prime $p$ is irregular then the determinant $\Delta= |\delta_{d,j}|,
\quad d=1,\dots,\frac{p-1}{2},\quad j=1,\dots,\frac{p-1}{2}$ verifies
$\Delta\equiv 0 \modu p$.
Then let $i_p$ be the index of irregularity of the prime $p$, (the number
of Bernoulli numbers $B_{p-1-2m}\equiv 0 \modu p,\quad m=1,\dots,\frac{p-3}{2})$
then we get the theorem
%
\begin{thm}\label{t15031a}
Let $\Delta_p$ be the determinant
\begin{equation}\label{e16031a}
\begin{split}
&\Delta_p=|\delta_{d,j}|,\quad d=1,\dots,\frac{(p-1}{2},\quad j=1,\dots,
\frac{p-1}{2},\\
&\delta_{d,j}\in \{-1,1\},\quad d=1,\dots,(p-1)/2,\quad j=1,\dots,(p-1)/2,\\
&\delta_{d,j}=1 \ {\mbox if}\  \exists l\in\N,\quad  \frac{l p}{(d+1)}< j
 < \frac{ l p}{d},\\
&\delta_{d,j}=-1 \ {\mbox if}\  \not\exists l\in\N,\quad
\frac{l p}{(d+1)}< j < \frac{ l p}{d}.
\end{split}
\end{equation}
If the prime $p$ has index of irregularity $i_p\geq 1$ then all the minor determinant $M$ of $\Delta_p$ of rank $r_M\geq \frac{p+1}{2}-i_p$ are divisible by $p$.
\begin{displaymath}
\Delta
\end{displaymath}
\begin{proof}
To each Bernoulli Number $B_{p-1-2m}\equiv 0\modu p$ corresponds a system of congruences
(\ref{e15032a}) p.\pageref{e15032a} with $\frac{1}{i^{2m+1}},\quad i=1,\dots,\frac{p-1}{2}$; two such systems are not equivalent.
\end{proof}
\end{thm}
%
As an example we have, for $p=17$, irregular prime,   from an easy MAPLE computation:
\begin{displaymath}
\Delta_{17} =
\begin{array}{|rrrrrrrr|}
-1   & -1  &  -1  &  -1   & -1    &-1   & -1  &  -1\\
-1   & -1  &  -1  &  -1   & -1    & 1   &  1  &   1\\
-1   & -1  &  -1  &  -1   &  1    &-1   & -1  &  -1\\
-1   & -1  &  -1  &   1   & -1    &-1   &  1  &   1\\
-1   & -1  &   1  &  -1   & -1    & 1   & -1  &  -1\\
-1   & -1  &  -1  &  -1   &  1    &-1   & -1  &   1\\
-1   & -1  &  -1  &  -1   & -1    &-1   &  1  &  -1\\
-1   &  1  &  -1  &   1   & -1    & 1   & -1  &   1\\
\end{array}
\equiv 9 \modu 17
\end{displaymath}
The same computation for $p=37$, first irregular prime, show that $\Delta_{37}\equiv 0 \modu 37$.
%
\subsection{Generalizations of Kummer congruences - II} \label{s23092}

We have seen that $P_1=\prod_{i\in I(d_1)} (x+\zeta^{u_i} y)\Z[\zeta]$ is a principal ideal. For another value $d_2\not= d_1$, we have also $P_2=\prod_{i\in I(d_2)} (x+\zeta^{u_i} y)\Z[\zeta]$ is a principal ideal. Here, we use simultaneously these two properties. From \cite{rib} p121 we have $Card(I(d_1))=Card(I(d_2))=\frac{p-1}{2}$.
%
\begin{prop}\label{t20092}
Let $u\in\N$ be a primitive root $\modu p$.
Let $d_1,d_2\in\N,\quad 0\leq d_1<d_2\leq p-2$.

Let $I(d_1)=\{i,\quad 1\leq i\leq p-1,\quad (u_{\nu-i}+u_{\nu-i+ind_u(d_1)})>p\}$.

Let $I(d_2)=\{i,\quad 1\leq i\leq p-1,\quad (u_{\nu-i}+u_{\nu-i+ind_u(d_2)})>p\}$.

Then, there exists $m\in\N, \quad 0\leq m\leq p-1$ such that :

\begin{equation}\label{eqkum03}
-2m+\sum_{i\in I(d_1)\cap I(d_2)}\frac{u_i}{(x+\zeta^{u_i}y)}
+\sum_{i\in I(d_1)\cap I(d_2)}\frac{u_i}{(x+\zeta^{-u_i} y)}\equiv 0 \modu \pi^{p-1}.
\end{equation}
\begin{proof}
Let $P_1=\prod_{i\in I(d_1)} (x+\zeta^{u_i}y)$.
Let $P_2=\prod_{i\in I(d_2)} (x+\zeta^{u_i}y)$.
We have
$P_1\times P_2=\zeta^{m_1+m_2}\eta_1\eta_2\gamma_1^p\gamma_2^p$
where
$m_1,m_2\in \N,\quad 0\leq m_1\leq p-1,\quad 0\leq m_2\leq p-1$,
where $\eta_1,\eta_2\in\Z[\zeta+\zeta^{-1}]^*$,
where $\gamma_1,\gamma_2\in\Z[\zeta]$.
Then we obtain

\begin{displaymath}
\begin{split}
&P_1P_2=\zeta^m\eta\gamma^p,\\
& m\in\N,\quad 0\leq m \leq p-1,\quad \eta\in\Z[\zeta+\zeta^{-1}]^*,\quad \gamma\in\Z[\zeta].
\end{split}
\end{displaymath}
Let $J(d_1,d_2)=(I(d_1)\cup I(d_2))-(I(d_1)\cap I(d_2))$.
We have
\begin{displaymath}
P_1P_2=\prod_{i\in I(d_1)\cap I(d_2)} (x+\zeta^{u_i}y)^2
\prod_{i\in J(d_1,d_2)} (x+\zeta^{u_i} y)(x+\zeta^{-u_i} y)
\end{displaymath}
because, from \cite{rib} p 121, if
$i\in I(d_1)$ and $i\not\in I(d_2)$ then $i+\frac{p-1}{2}\modu p\in I(d_2)$
and $u_{i+(p-1)/2}\equiv -u_i \modu p$.
Then,
\begin{displaymath}
\begin{split}
&\zeta^{-m}\prod_{i\in I(d_1)\cap I(d_2)} (x+\zeta^{u_i}y)^2
\prod_{i\in J(d_1,d_2)} (x+\zeta^{u_i} y)(x+\zeta^{-u_i} y)\equiv\\
&\zeta^{m}\prod_{i\in I(d_1)\cap I(u,d_2)} (x+\zeta^{-u_i}y)^2
\prod_{i\in J(d_1,d_2)} (x+\zeta^{u_i} y)(x+\zeta^{-u_i} y)\modu \pi^{p+1},
\end{split}
\end{displaymath}
so
\begin{displaymath}
\zeta^{-m}\prod_{i\in I(d_1)\cap I(d_2)} (x+\zeta^{u_i}y)^2\equiv
\zeta^{m}\prod_{i\in I(d_1)\cap I(d_2)} (x+\zeta^{-u_i}y)^2\modu \pi^{p+1},
\end{displaymath}
also
\begin{displaymath}
\zeta^{-2m}\prod_{i\in I(d_1)\cap I(d_2)} (x+\zeta^{u_i}y)^2\equiv
\prod_{i\in I(d_1)\cap I(d_2)} (x+\zeta^{-u_i}y)^2\modu \pi^{p+1},
\end{displaymath}

\begin{displaymath}
\zeta^{-2m(p+1)/2}\prod_{i\in I(d_1)\cap I(d_2)} (x+\zeta^{u_i}y)^{2(p+1)/2}\equiv
\prod_{i\in I(d_1)\cap I(d_2)} (x+\zeta^{-u_i}y)^{2(p+1)/2}\modu \pi^{p+1},
\end{displaymath}
hence
\begin{displaymath}
\zeta^{-m}\prod_{i\in I(d_1)\cap I(d_2)} (x+\zeta^{u_i}y)\equiv
\prod_{i\in I(d_1)\cap I(d_2)} (x+\zeta^{-u_i}y)\modu \pi^{p+1},
\end{displaymath}
because $(x+\zeta^{u_i}y)^p\equiv (x+\zeta^{-u_i}y)^p \modu \pi^{p+1}$.
Then, the end of the proof is the same than the proof of the theorem \ref{thmkum01} p.\pageref{thmkum01}.
\end{proof}
\end{prop}
%
\begin{cor}{  $\frac{3(p-2)(p-3)}{2}$ congruences $\modu p.$}\label{c05101}

Let $\mu\in \{3,5,7\}$.
Let $d_1,d_2\in\N,\quad 1\leq d_1<d_2\leq p-2$.

Let $I(d_1)=\{i,\quad 1\leq i\leq p-1,\quad (u_{\nu-i}+u_{\nu-i+ind_u(d_1)})>p\}$.

Let $I(d_2)=\{i,\quad 1\leq i\leq p-1,\quad (u_{\nu-i}+u_{\nu-i+ind_u(d_2)})>p\}$.

Then, for each primitive root $u \modu p$, we have
$\sum_{i\in I(d_1)\cap I(d_2)} u_i^\mu\equiv 0 \modu p$.
\begin{proof}
Proof similar to corollary \ref{c21091} p.\pageref{c21091}.
\end{proof}
\end{cor}
%
\begin{thm} { $\frac{3(p-2)(p-3)}{2}$ congruences $\modu p.$}\label{t05101}

Let $\mu\in\{3,5,7\}$.
Let $d_1,d_2\in\N, \quad 1\leq d_1<d_2 \leq p-2$.
If $p$ does not divide the sum
\begin{displaymath}
\begin{split}
& S= \sum_{l=1}^{d_1}\sum_{m=1}^{d_2}
\sum_{j\in [\frac{lp}{d_1+1}, \frac{lp}{d_1}[\cap [\frac{mp}{d_2+1},\frac{mp}{d_2}[}
\frac{1}{j^\mu},\\
& S=\sum_{j\in\cup_{l=1}^{d_1}\cup_{m=1}^{d_2} ([\frac{lp}{d_1+1}, \frac{lp}{d_1}[\cap [\frac{mp}{d_2+1},\frac{mp}{d_2}[\ )}
\frac{1}{j^\mu}.\\
\end{split}
\end{displaymath}
then the first case of Fermat's Last Theorem holds for the exponent $p$.
\begin{proof}$ $
\begin{itemize}
\item
Observe at first that the two forms of $S$ are equivalent because
$l_1\not=l_2\Rightarrow [\frac{l_1 p}{d_1+1},\frac{l_1 p}{d_1}]
\cap[\frac{l_2 p}{d_1+1},\frac{l_2 p}{d_1}]=\emptyset$; therefore
\begin{displaymath}
([\frac{l_1 p}{d_1+1},\frac{l_1 p}{d_1}]\cup[\frac{m_1 p}{d_2+1},\frac{m_1 p}{d_2}])\cap
([\frac{l_2 p}{d_1+1},\frac{l_2 p}{d_1}]\cup[\frac{m_2 p}{d_2+1},\frac{m_2 p}{d_2}])\not=\emptyset
\end{displaymath}
implies that $l_1=l_2$ and $m_1=m_2$.
\item
From corollary \ref{c05101} p.\pageref{c05101}, similarly to corollary \ref{c21091} p.\pageref{c21091}, we have

\begin{displaymath}
\sum_{i\in I(d_1)\cap I(d_2)} u_i^\mu\equiv 0 \modu p.
\end{displaymath}
Then the result is similar to corollary \ref{c2009} p.\pageref{c2009}.
\end{itemize}
\end{proof}
\end{thm}
%


\begin{cor} { $3(p-2)$ congruences $\modu p.$ }\label{c09101}

Let $\mu\in\{3,5,7\}$.
Let $d\in\N,\quad 1\leq d \leq p-2$.
Then, we have

$\sum_{j>p/2}^{p-1} \sum_{m=1}^{d}\sum_{j\in[\frac{mp}{(d+1)},\frac{mp}{d}[} \frac{1}{j^\mu}
\equiv 0 \modu p$.
\begin{proof}
Immediate consequence of theorem \ref{t05101} p.\pageref{t05101} with $d_1=1$ and $d_2=d$.
\end{proof}
\end{cor}
%

%
\subsection{ Generalization to exponential congruences  $f(p,t)\equiv 0\modu p$.}
All the results of this section can be generalized {\it mutatis mutandis} to exponential congruences  on $p$ and $t\equiv -\frac{y}{x}\modu p$: we give, as {\bf an} example the generalization of Cauchy criterion $\sum_{j=(p+1)/2}^{p-1} \frac{1}{j^3}\equiv 0 \modu p$.
\begin{thm}{ $p-1$ congruences, generalization of Cauchy criterion to exponential congruences.}\label{t21023a}

For $k=1,\dots, p-1$
\begin{displaymath}
\sum_{j=(p+1)/2}^{p-1} \frac{ 1}{j^3} \times
(t^{(k\times j^3 \modu p)}+t^{(p-k\times j^3\modu p)})\equiv 0 \modu p.
\end{displaymath}
\begin{proof}
Apply theorem \ref{t20021a} p.\pageref{t20021a} to the set
\begin{displaymath}
I_f=I(d)=\{u_i=\frac{1}{j^3}\modu p \quad |\quad j=\frac{p-1}{2},\dots,p-1 \}.
\end{displaymath}
\end{proof}
\end{thm}
\end{section}
%
%
\clearpage
\section{FLT first case: A comparative survey with bibliography}\label{s19102}
\label{s06101}
This table  locate our results in the context of a survey of some  significant necessary conditions on $x,y, p$ or on $p$ for the first case of Fermat's Last Theorem seen in the literature.
Here, we recall some notations used in this table :
\begin{itemize}
\item
Let $a,b\in \R$; recall the classic notations
\begin{itemize}
\item
$[a,b]=\{\alpha\quad |\quad \alpha\in\R,\quad a\leq \alpha\leq b\}$,
\item
$[a,b[=\{\alpha\quad |\quad \alpha\in\R,\quad a\leq \alpha<  b\}$.
\end{itemize}
\item
$p\in\N, \quad p>5$ prime.
\item
$p-$cyclotomic number field $\Q(\zeta)$.
\item
$\Z[\zeta]$, ring of integers of $\Q(\zeta)$.
\item
$u$ is a primitive root $\modu p$. For $i\in\N,\quad u_i\equiv u^i \modu p$ with
$1\leq u_i\leq p-1$; for $i< 0$, it is to be understood as: $u_i$ is defined by $u_i\times u_{-i}\equiv 1 \modu p$ and $1\leq u_i \leq p-1$.
\item
$h=p^{e_p}\times h_2,\quad e_p\in \N,\quad h_2\not\equiv 0 \modu p$, class number of $\Q(\zeta)$,
where $e_p<p-2$, from theorem (\ref{thmhb}).
\item
$h=h^*\times h^+$, where $h^*$ is the first factor of the class number.
\item
$r_p$, $p-$rank of the class group of the cyclotomic number field $\Q(\zeta)$.
\item
$x^p+y^p+z^p=0,\quad $
$x,y,z\in \Z-\{0\}$  coprime by pair: Fermat equation.
\item
$xyz\not\equiv 0 \modu p$ : First case.
\item
$t\in\N,\quad t\equiv -\frac{x}{y}\modu p$.
\item
$\phi_{2m+1}(T)=1^{2m}T+2^{2m}T^2+3^{2m}T^3+\dots+(p-1)^{2m}T^{p-1}$,
odd Mirimanoff polynomial of the indeterminate $T$, for $m=1,\dots,\frac{p-3}{2}$.
\item
$P_{2m+1}(T)$ : explicitly computable Kummer's  polynomials : $P_1(T)=T,\quad P_3(T)=T(1-T),\quad P_5(T)=T(1-T)(1-10T+T^2),\dots$.
\item
$B_{p-1-2m},\quad m\in\N$, \quad even Bernoulli number for $m=1,\dots,\frac{p-3}{2}$.

\item
$i_p=Card\{1\leq m\leq\frac{p-3}{2}\quad |\quad B_{p-1-2m}\equiv 0 \modu p\}$, index of irregularity of $p$.

\item
{\bf gfw} : This mention joined to the formulation of a criterion on the following table means that the generalization of Fermat equation
\begin{displaymath}
\begin{split}
& x^p+y^p+c\times z^p=0,\\
& c\in\Z-\{0\},\\
& x\times y \times z\times c\not\equiv 0 \modu p,\\
& x-y\not\equiv 0 \modu p^2,
\end{split}
\end{displaymath}

where the prime factors $c_i$ of $c$ verify $c_i\not\equiv 1 \modu p$, as explained in section \ref{s08081} p.\pageref{s08081}, can be applied to this criterion.
\end{itemize}
%
$ $
\clearpage

\begin{tabular}{|p{3.in}|l|l|l|}\hline
\multicolumn{4}{|c|}{\bf (1/8) If the first case of Fermat's Last Theorem would fail for $p$ then:}\\ \hline
\multicolumn{1}{|l|}{\it Condition}&{\it Author}&{\it date}&{\it Reference}\\ \hline\hline
{\bf gfw} $\sum_{j=1}^{(p-1)/2} \frac{1}{j^3}\equiv 0 \modu p. $&Cauchy&1847&\cite{rib} (4A) p 120\\ \hline
{\bf gfw}For Kummer Polynomials explicitly computable

$P_1(T)=T$,

$P_3(T)=T(1-T)$,

$P_5(T)=T(1-T)(1-10T+T^2),\dots$

and for $m=1,\dots,\frac{p-3}{2}$ then

$P_{2m+1}(-t)B_{p-1-2m}\equiv 0 \modu p.$

&Kummer&1857&\cite{rib} (4D) p 125\\
\hline
Mirimanoff congruences :

$ \phi_{p-1}(t)\equiv 0 \modu p$,

$\phi_{p-2}(t)\phi_{2}(t)\equiv 0 \modu p$,

$\phi_{(p+1)/2}\phi_{(p-1)/2}(t)\equiv 0 \modu p.$

&Mirimanoff&1909&\cite{rib}, (1B) p 145\\

\hline
$ 2^{p-1}-1\equiv 0 \modu p^2.$
&Wieferich& 1909&\cite{rib}, (3A) p151\\
\hline
{\bf gfw}For $m=1,\dots,\frac{p-3}{2}$ then

$\phi_{2m+1}(t)B_{p-1-2m}\equiv 0 \modu p.$&Mirimanoff&1909&\cite{rib} p 145 \\
\hline

$ 5^{p-1}-1 \equiv 0 \modu p^2 $

and $\sum_{j=1}^{[p/5]} \frac{1}{j}\equiv 0 \modu p$&Vandiver&1914& \cite{van}\\

\hline
$\sum_{j=1}^{[p/3]}\frac{1}{j^2}\equiv 0 \modu p.$&Vandiver&1925&\cite{rib} (5B) p 157 \\ \hline
For $m=1,\dots,31$ then

$m^{p-1}\equiv 1 \modu p^2$&Morishima&1931&\cite{rib} (6E) p 160 \\
\hline

$\sum_{j=1}^{[p/6]}\frac{1}{j^2}\equiv 0 \modu p.$&Schwindt&1933&\cite{rib} (5C) p 157 \\ \hline
Let $n_0=(45!)^{88}.$

If $p>n_0$ then for the even Bernoulli Number

$B_{p-1-2m}$ for $m=1,\dots, [Log(p)^{1/3}]$

verify $B_{p-1-2m}\equiv 0 \modu p.$
&Krasner&1934&\cite{kra}, \cite{rib} (2A) p 149\\
\hline
For $d=2,\dots,p-1$

$(1-d^{-2m})B_{p-1-2m}\equiv$

$(p-1-2m)d^{-2m-1}\sum_{l=1}^d\sum_{1\leq j<\frac{lp}{d}} \frac{1}{j^{2m+1}}\modu p$. &Vandiver &1937&\cite{va3}, \cite{rib} (5B) p 108.
\\
\hline
for $n=2,3,4,6$ then

$\sum_{j=1}^{[p/n]} \frac{1}{j}\equiv 0 \modu p. $&Yamada&1941&\cite{rib}, (5A) p 156 \\ \hline

\end{tabular}

\clearpage


\begin{tabular}{|p{3.in}|l|l|l|}\hline

\multicolumn{4}{|c|}{\bf (2/8) If the first case of  Fermat's Last Theorem would fail for $p$ then:}\\ \hline
\multicolumn{1}{|l|}{\it Condition}&{\it Author}&{\it date}&{\it Reference}\\ \hline\hline

{\bf gfw}$\quad p^{[\sqrt{p}]-1}$ divides the first factor $h*$&Eichler&1965&\cite{rib}, (7A) p185 \\ \hline

If $s\in\Z, s\not\equiv 0 \modu p$ denote $\overline{s}$ the unique integer such that
$s\equiv\overline{s} \modu p$ and $1\leq \overline{s}\leq p-1$.
Let $s^\prime$ be any integer such that $s^\prime s\equiv 1 \modu p$. For each $n=1,2,\dots,p-2$, let $E_n$ be the set of all integers $s,\quad 1\leq s \leq p-2$ such that $\overline{(n+1)^\prime n s}<s$.
The Le Lidec polynomials are
$\Lambda_n(T)=\sum_{s\in E_n} \overline{s^\prime}T^{p-s}$.

The Le Lidec congruences are :

$\phi_{p-1}(t)\equiv 0 \modu p$

$\Lambda_2(t)\equiv 0 \modu p$

$\vdots$

$\Lambda_{p-2}(t)\equiv 0 \modu p$.
&Le Lidec&1967&\cite{rib},(1D) p148\\
\hline
$i(p)>[\sqrt{p}]-1$&  Skula&1977&\cite{sk1}\\ \hline
For $\mu=1,3,5,7,9,11,13$ then

$\sum_{j=1}^{[p/3]} \frac{1}{j^\mu}\equiv 0 \modu p,$

$\sum_{j=1}^{[p/6]} \frac{1}{j^\mu}\equiv 0 \modu p.$
&Ribenboim&1979& \cite{rib},(5E) p 159\\
\hline
$q^{p-1}-1\equiv 0 \modu p^2$ (Granville-Monagan)

for all prime $q\leq 89$
&
Granville
&1988& \cite{gra}\\
\hline

For $k,N\in\N$,

$N\in \{2,3,4,5,6,7,8,9,10\}\cup \{12\}$

$0\leq k \leq N-1$ then

$S(k,N)=
\sum_{j=[kp/N]+1}^{[(k+1)p/N]} \frac{1}{j}\equiv 0 \modu p$&  Skula& 1992
& \cite{sk2} ,\cite{sk5}\\
\hline
If the Kummer/Mirimanoff system of congruences has

non trivial solutions, then

$i(p)\geq [(\frac{p}{2})^{1/3}]$ &  Skula& 1994& \cite{sk6},\cite{sk7}\\
\hline
$q^{p-1}-1\equiv 0 \modu p^2$

for all prime $q\leq 113$
&
Susuki
&1994& \cite{sus}\\
\hline

There exists  $L\in\N$ such that for all $p>L$, then :

for $k,N\in\N, \quad 2\leq  N\leq 94,\quad 0\leq k \leq N-1$

$S(k,N)=\sum_{j=[kp/N]+1}^{[(k+1)p/N]} \frac{1}{j}\equiv 0 \modu p$& Cir\`anek& 1994&\cite{cir}\\
\hline
\end{tabular}
.


.
\begin{tabular}{|p{3.in}|l|l|l|}\hline
\multicolumn{4}{|c|}{\bf (3/8) If the first case of  Fermat's Last Theorem would fail for $p$ then:}\\ \hline

\multicolumn{1}{|l|}{\it Condition}&{\it Author}&{\it date}&{\it Reference}\\ \hline\hline

For $k,N\in\N, \quad 2\leq  N\leq 46,\quad 0\leq k \leq N-1$

$S(k,N)=\sum_{j=[kp/N]+1}^{[(k+1)p/N]} \frac{1}{j}\equiv 0 \modu p$& Dilcher Skula& 1995
& \cite{dil} \\
\hline
{\bf gfw} For $f\in \N,\quad e_p+2\leq f\leq p-1$, for each set

$I_f=\{g_1,g_2,\dots,g_f \quad| \quad i\not=i^\prime\Rightarrow g_i\not=g_i^\prime \}$,

there exists a set

$L(I_f)=\{l_i, \quad l_i\in\N,\quad i=1,\dots,f \}$

such that $\sum_{i=1}^f g_i \times l_i\equiv 0 \modu p$ and

$\phi_{2m+1}(t)\sum_{i=1}^f g_{i(2m+1)}\times l_i\equiv 0 \modu p$

for $m=1,\dots,\frac{p-3}{2}$.&Qu\^eme&1998&pro \ref{p4} p. \pageref{p4}\\
\hline
{\bf gfw} There exists $r\in \N,\quad r\leq min(e_p+2,\frac{p-1}{2})$, depending only on $p$ and $t$,   such that, for each primitive root $u \modu p$, there exists a set

$N_r(u)=\{n_1,n_2,\dots,n_r\}$

verifying $\sum_{i=1}^r u_i \times n_i\equiv 0 \modu p$

and for all $m=1,2,\dots,\frac{p-3}{2}$
\begin{itemize}
\item
$\phi_{2m+1}(t)\times(\sum_{i=1}^r u_{(2m+1)i}\times  n_i)\equiv 0 \modu p$,
\item

$p$ does not divide simultaneously $\phi_{2m+1}(t)$ and $\sum_{i=1}^r u_{(2m+1)i}\times n_i$.
\end{itemize}
&Qu\^eme&1998&cor \ref{cor1209} p. \pageref{cor1209}\\
\hline

{\bf gfw} $r_p$  $p-$rank of class group of $\Q(\zeta)$.

There are at least $\frac{p-3}{2}-r_p$ different $m\in\N$,

$1\leq m\leq \frac{p-3}{2}$ with $\phi_{2m+1}(t)\equiv 0\modu p.$
&Qu\^eme&1998&cor \ref{c09113} p. \pageref{c09113} \\
\hline
\end{tabular}
\clearpage



\begin{tabular}{|p{3.in}|l|l|l|}\hline

\multicolumn{4}{|c|}{\bf (4/8) If the first case of  Fermat's Last Theorem would fail for $p$ then:}\\ \hline
\multicolumn{1}{|l|}{\it Condition}&{\it Author}&{\it date}&{\it Reference}\\ \hline\hline
{\bf gfw}
Let

$I_\phi=\{2m_i+1 \quad |\quad i=1,\dots,r-1,$

$\phi^*_{2m_i+1}(t)\not\equiv 0 \modu p,\quad 0\leq m_i\leq \frac{p-3}{2}\}$.

We have the equivalence

$\sum_{i=1}^r\frac{u_i l_i}{x+\zeta^{u_i}y}
+\sum_{i=1}^r\frac{u_i l_i}{x+\zeta^{-u_i} y}\equiv 0 \modu \pi^{p-1}$

$ \Longleftrightarrow $

$ l_1 \equiv (-1)^{r-1}\times \prod_{i=1}^{r-1} u_{2m_i+1} \modu p,$

$\vdots $

$ l_{r-1}\equiv (-1)\times \sum_{i=1}^{r-1} u_{2m_i+1} \modu p,$

$ l_r=1$.
& Qu\^eme & 2000 & thm \ref{t23011} p. \pageref{t23011} \\
\hline

Let $K$ be an intermediate field, with $\Q\subset K\subset \Q(\zeta),\quad [K:\Q]=g, \quad
p-1 =g\times f$, where $f$ is odd and $gcd(f,g)=1$.
Let $r_K$ be the $p$-rank of the class group of $K/\Q$. Then, if $r_K<\frac{g-1}{2}$, there are at least
$\frac{g-1}{2}-r_K$ Mirimanoff polynomials congruences
\begin{displaymath}
\phi_{f (2n+1)}(t)\equiv 0 \modu p, \quad 1\leq 2n+1 < g.
\end{displaymath}
& Qu\^eme & 2000 & thm \ref{t09071a} p. \pageref{t09071a} \\
\hline
\end{tabular}
\clearpage


\begin{tabular}{|p{3.in}|l|l|l|}\hline
\multicolumn{4}{|c|}{\bf (5/8) If the first case of  Fermat's Last Theorem would fail for $p$ then:}\\ \hline

\multicolumn{1}{|l|}{\it Condition}&{\it Author}&{\it date}&{\it Reference}\\ \hline\hline
Let $p-1=f\times g$ with $f$ odd. Let $K$ be the intermediate field $\Q\subset K \subset \Q(\zeta)$ and $[K:\Q]=g$. Suppose that $p$ does not divide the class number of $K/\Q$. Then we have  the Mirimanoff's polynomials congruences

$\phi_{ f (2 n+1)}(t)\equiv 0 \modu p,$

$2n+1=1,3,5,\dots, g-1.$
& Qu\^eme & 2000 & cor \ref{c010161} p. \pageref{c010161} \\
\hline
{\bf gfw}
Let $f\in \N,  \quad p-1\equiv 0 \modu f,\quad f\not\equiv 0 \modu 2$.
Let $g=\frac{p-1}{f}$.
Let $K$ be the field $\Q\subset K\subset \Q(\zeta),
\quad [K:\Q]=g$.
If  $p$ does not divide the class number of the field extension $K/\Q$ then

$\sum_{i=0, \quad i^{2 f}\equiv 1 \ mod \ p}^{p-1} (i^{f-1}\times t^i)\equiv 0 \modu p.$
& Qu\^eme & 2000 & cor \ref{c010162} p. \pageref{c010162} \\
\hline

{\bf gfw} Suppose that $p-1\equiv 0 \modu 3.$
Let $K$ be the field
$\Q\subset K\subset\Q(\zeta),\quad [K:\Q]=\frac{p-1}{3}.$

then $p$ should divide the class number of $K/\Q$.
&Qu\^eme&1999& thm \ref{t12101} p. \pageref{t12101} \\
\hline
{\bf gfw}
Let $t\equiv -\frac{y}{x}\modu p$.
Let $r_g$ be the minimal rank defined in the relation (\ref{e104012}) p \pageref{e104012}.
If FLT  First Case fails for $p$,  then

1.) there exists {\bf exactly}
$\frac{p-1}{2g}-r_g+1$ different values
$m,\quad m\in\N,\quad 1\leq m \leq \frac{p-1}{2g}-1$,
such that we have for $m$ the $g$  Mirimanoff polynomial congruences:

$\phi_{2m+1+\alpha(p-1)/g}^{*}(t)\equiv 0 \modu p,$

for $\alpha=0,\dots,g-1.$

2.) there exists {\bf exactly}
$r_g-1$ different values
$m,\quad m\in\N,\quad 1\leq m \leq \frac{p-1}{2g}-1$,
such that we have for $m$ {\bf at least} one $\alpha\in\N, \quad 0\leq \alpha\leq g-1$ with the  Mirimanoff polynomial congruence:

$\phi_{2m+1+\alpha(p-1)/g}^{*}(t)\not\equiv 0 \modu p.$
&Qu\^eme&2001& thm \ref{t104013} p. \pageref{t104013} \\
\hline
\end{tabular}
\clearpage


\begin{tabular}{|p{3.in}|l|l|l|}\hline
\multicolumn{4}{|c|}{\bf (6/8) If the first case of  Fermat's Last Theorem would fail for $p$ then:}\\ \hline
\multicolumn{1}{|l|}{\it Condition}&{\it Author}&{\it date}&{\it Reference}\\ \hline\hline
{\bf gfw}
Let $u$ be a primitive root $\modu p$.

Then $p-1$ congruences:

$\begin{array}{l}
\displaystyle \sum_{i=1}^{(p-1)/2} u_i \times(t^{m_{i,k}}+t^{p-m_{i,k}})
    \equiv 0 \modu p, \\
k=1,\dots,p-1,\quad u_i\times m_{i,k}\equiv k \modu p, \\
1\leq m_{i,k}\leq p-1. \end{array}$
&Qu\^eme& 2000 & thm \ref{t20022a} p.\pageref{t20022a}\\
\hline
{\bf gfw}

Let $f$ be an integer with $f\geq min(e_p+2,\frac{p-1}{2})$.

Let $E_f$ be  {\bf any } set $E_f=$

$\{g_i\in\N,\quad i=1,\dots,f,\quad 1\leq g_i\leq p-1 \}$.

Let $K_f$ be {\bf any} set $K_f=$

$\{k_j\in\N,\quad j=1,\dots,f,\quad 1\leq k_j\leq p-1 \}$.

Let us define for $i=1,\dots f,\quad j=1,\dots,f$

$ m_{i j}\equiv g_i^{-1}\times k_j\modu p,\quad  1\leq m_{i j}\leq p-1,$

$ g_i\in E_f,\quad k_j\in K_f,$

$\theta_{i j}= (t^{m_{i j}}+t^{p-m_{i j}})\modu p.$

If the first case of FLT would fail for $p$ then

the $f\times f$ determinant
\begin{displaymath}
\Delta_{f}(t)= | \theta_{i j}|_{1\leq i\leq p-1,\quad 1\leq j\leq p-1}
\equiv 0 \modu p.
\end{displaymath}
&Qu\^eme&1998&thm \ref{t21021a} p. \pageref{t21021a}\\

\hline

\end{tabular}


\begin{tabular}{|p{3.in}|l|l|l|}\hline
\multicolumn{4}{|c|}{\bf (7/8) If the first case of  Fermat's Last Theorem would fail for $p$ then:}\\ \hline
\multicolumn{1}{|l|}{\it Condition}&{\it Author}&{\it date}&{\it Reference}\\ \hline\hline

{\bf gfw}

Let $f$ be an integer with $f\geq min(e_p+2,\frac{p-1}{2})$.

Let $u$ be a primitive root $\modu p$.

Let $l\in\N,\quad 1\leq l\leq p-1$.

Let $m\in\N,\quad 1\leq l\leq p-1$.

the determinant $\Delta_{f}(t)= $

$ | t^{u_{m\times j -l\times i}}+
t^{p- u_{m\times j -l\times i}}|_{1\leq i\leq f,\quad 1\leq j\leq f}$

$\equiv 0 \modu p$
&Qu\^eme&2000& cor \ref{c010051} p. \pageref{c010051}\\
\hline

{\bf gfw}

Let $a\in\N$ be the order of $t \modu p$. Then $a>\frac{p-1}{e_p+2}$.
&Qu\^eme & 2000 & thm \ref{t20024a} p.\pageref{t20024a} \\
\hline
\hline
({\bf gfw} if $\mu=3)$. Let $\mu\in\{3,5,7\}$. Then

$\sum_{j=1}^{[p/4]} \frac{1}{j^\mu}\equiv 0 \modu p.$
&Qu\^eme&1998& cor \ref{c03104} p.\pageref{c03104} \\ \hline
({\bf gfw} if $\mu=3)$. $3(p-2)$ congruences mod $p$: For all $\mu\in\{3,5,7\}$.

For all $d\in\N,\quad 1\leq d \leq p-2$ then

$\sum_{j\in\cup_{l=1}^d[\frac{lp}{d+1},\frac{lp}{d}[} \frac{1}{j^\mu}\equiv 0 \modu p.$
&Qu\^eme&1998&cor \ref{c2009} p. \pageref{c2009} \\
\hline
{\bf gfw} For every $1\leq m\leq \frac{p-3}{2}$ with $\phi_{2m+1}(t)\not\equiv 0 \modu p$,

for every $d\in\N,\quad 1\leq d \leq p-2$ then

$\sum_{j\in\cup[\frac{lp}{d+1},\frac{lp}{d}[} \frac{1}{j^{2m+1}}\equiv 0 \modu p.$
&Qu\^eme&1998&cor \ref{t22021} p.\pageref{t22021} \\
\hline
{\bf gfw} Let $p>(45!)^{88}$.

$[Log(p)^{1/3}](p-2)$ congruences mod $p$:

For all $m=1,\dots,[Log(p)^{1/3}]$.

For all $d\in\N,\quad 1\leq d \leq p-2$ then

$\sum_{j\in\cup_{l=1}^d [\frac{lp}{d+1},\frac{lp}{d}]} \frac{1}{j^{2m+1}}\equiv 0 \modu p.$
&Qu\^eme&1998&cor \ref{c15031} p. \pageref{c15031} \\
\hline
({\bf gfw} if $\mu=3)$. $3(p-2)(p-3)$ congruences $\modu p$: let $\mu\in\{3,5,7\}$.

For all $d_1,d_2\in\N, \quad 1\leq d_1<d_2 \leq p-2$,

$\quad \sum_{l=1}^{d_1}\sum_{m=1}^{d_2}
\sum_{j\in [\frac{lp}{d_1+1}, \frac{lp}{d_1}]\cap [\frac{mp}{d_2+1},\frac{mp}{d_2}[}
\frac{1}{j^\mu}$

$\equiv 0 \modu p,$

 &Qu\^eme&1998&thm \ref{t05101} p. \pageref{t05101} \\
\hline
\end{tabular}
\clearpage


\begin{tabular}{|p{3.in}|l|l|l|}\hline
\multicolumn{4}{|c|}{\bf (8/8) If the first case of  Fermat's Last Theorem would fail for $p$ then:}\\ \hline
\multicolumn{1}{|l|}{\it Condition}&{\it Author}&{\it date}&{\it Reference}\\ \hline\hline

({\bf gfw} if $\mu=3)$. $3(p-2)$ congruences $\modu p$: let $\mu\in\{3,5,7\}$.

For all $d\in\N,\quad 1\leq d \leq p-2$,

$\sum_{j>p/2}^{p-1} \sum_{m=1}^{d}\sum_{j\in[\frac{mp}{(d+1)},\frac{mp}{d}[} \frac{1}{j^\mu}
\equiv 0 \modu p$,

&Qu\^eme&1998& cor \ref{c09101} p. \pageref{c09101}\\
\hline
{\bf gfw}
$p-1$ congruences, generalization of Cauchy criterion to exponential congruences.

For $k=1,\dots, p-1$

$\begin{array}{l} \displaystyle
\sum_{j=(p+1)/2}^{p-1} \frac{ 1}{j^3} \times
(t^{(k\times j^3 \modu p)}+t^{(p-k\times j^3\modu p)}) \\
 \equiv 0 \modu p.
\end{array}$
&Qu\^eme & 2000 & thm \ref{t21023a} p. \pageref{t21023a}\\
\hline
{\bf A particular case of Fermat-Wiles:}

Let $x,y,z\in \Z-\{0\}$, pairwise coprime.
Then $p\|(x-y)\times (y-z)\times (z-x)$ implies that

$x^p+y^p+z^p\not=0$.
&Qu\^eme & 2001 & thm \ref{t110071} p. \pageref{t110071}\\
\hline

\end{tabular}
\clearpage
%
\section{FLT first case: A generalization   to  diophantine equations
$\frac{x^p+y^p}{x+y}=t_1^p$ and $x^p+y^p+c\times z^p=0$}\label{s08081}
\subsection{Introduction}\label{s210162}
Let $p>5$ be a prime as assumed in all this paper. Let $x,y,z\in\Z-\{0\}$ mutually coprime. Let us consider the two generalizations of the Fermat's Equation  first case
\begin{equation}\label{e05037a}
\frac{x^p+y^p}{x+y}=t_1^p, \quad x,y,t_1\in\Z-\{0\},\quad xy(x^2-y^2)\not\equiv 0 \modu p,
\end{equation}
and
\begin{equation}\label{e05038a}
 x^p+y^p+c\times z^p=0,\quad
 x,y,z,c\in\Z-\{0\},\quad x y (x^2-y^2)\not \equiv 0 \modu p.
\end{equation}
Recall that $p>5$ is assumed in this paper.
\begin{itemize}
\item
When $x\times y \times z\not \equiv 0 \modu p$ the Fermat's equation $x^p+y^p+z^p=0$ leads to the three Barlow-Abel relations
\begin{displaymath}
\begin{split}
&\frac{x^p+y^p}{x+y}=t_1^p,\quad x+y=t^p,\quad t_1\in\Z-\{0\},\\
&\frac{y^p+z^p}{x+y}=r_1^p,\quad y+z=r^p,\quad r_1\in\Z-\{0\},\\
&\frac{z^p+x^p}{x+y}=s_1^p,\quad z+x=s^p,\quad s_1\in\Z-\{0\}.
\end{split}
\end{displaymath}
At first, if we assume that all primes $c_i$ dividing $c$ verify $c_i\not\equiv 1\modu p,$, we shall show that the Generalized Fermat-Wiles equation,
\begin{displaymath}
x^p+y^p+c\times z^p =0,
\quad x y (x^2-y^2)\not\equiv 0 \modu p,
\end{displaymath}
implies the partial Barlow-Abel equation
\begin{displaymath}
\frac{x^p+y^p}{x+y}=t_1^p,\quad x+y=t^p,\quad t_1\in\Z-\{0\}.
\end{displaymath}
which is thus more  general than equation (\ref{e05038a}) p.\pageref{e05038a} and a fortiori more general that the first case of Fermat,  $x^p+y^p+z^p=0,\quad xyz\not\equiv 0 \modu p$.
\item
We shall use the fact that several criteria on $p$  obtained for the first case of FLT in the literature and almost all the new criteria on $p$ that we have obtained rest on the factorization properties of {\bf one and only one} of the  Barlow-Abel equations $\frac{x^p+y^p}{x+y}=t_1^p$ where $t_1\in \Z-\{0\}$.
\item The first example of criterion on $p$ for the first case of FLT is Cauchy's criterion, see Ribenboim \cite{rib} p 120. As an example, we detail the proof of generalization of Cauchy's criterion for the Generalized Fermat-Wiles equation.
\item
For a comparison with results obtained for the Generalized Fermat-Wiles equation
$x^p+y^p+c z^p=0$, see Finch \cite{fin}, survey for this topic,
and Barlow-Abel equation $\frac{x^p+y^p}{x+y}=t_1^p$, see Terjanian \cite{ter}
\end{itemize}
%

\begin{subsection} {Some results}\label{ss08081}
\begin{lem}
Let $x,y,z\in\Z-\{0\}$ mutually co-prime. Let $c\in\Z-\{0\}$. Suppose that
$x\times y\times z\times c\not\equiv 0 \modu p$.
Let $c=\pm \prod_{i=1}^m c_i^{\alpha_i}$ be the factorization in primes of $c$.
Suppose that $c_i\not\equiv 1 \modu p,\quad i=1,\dots, m$.
Then $x^p+y^p+c\times z^p=0$ implies the generalization of Barlow-Abel relation $\frac{x^p+y^p}{x+y}=t_1^p$
and the relation $\frac{x+y}{c}=t^p$.
\begin{proof}
From $x^p+y^p+c\times z^p=0$, we get $(x+y)\times(\frac{x^p+y^p}{x+y})+c\times z^p=0$.
Classically, all the primes divisors $q$  of $\frac{x^p+y^p}{x+y}$ verify $q\equiv 1 \modu p$.
From hypothesis on $c$, we deduce that $c$ and $\frac{x^p+y^p}{x+y}$ are coprime and so $c$ divides $(x+y)$.
Therefore $(\frac{x+y}{c})\times (\frac{x^p+y^p}{x+y})+z^p=0$. Classically, $(x+y)$ and
$(\frac{x^p+y^p}{x+y})$ are coprime and so $(\frac{x+y}{c})$ and $(\frac{x^p+y^p}{x+y})$ are coprime and so $\frac{x+y}{c}=t^p$ and $\frac{x^p+y^p}{x+y}=t_1^p$.
\end{proof}
\end{lem}
%
\begin{thm}{Generalization of Cauchy's criterion}\label{t14081}

Let $p>5$ be a prime.  If $\frac{x^p+y^p}{x+y}=t_1^p,\quad x,y,t_1\in\Z, \quad xy(x^2-y^2)\not\equiv 0 \modu p$ then $p$ verifies the congruence
\begin{displaymath}
\sum_{j=1}^{(p-1)/2} j^{p-4}\equiv 0 \modu p.
\end{displaymath}
\begin{proof} $ $
\begin{itemize}
\item
Suppose, at first, that $\Q(\zeta)/\Q$ is a principal field.
Then, similarly to proposition \ref{p4} p.\pageref{p4}, it is possible to write
\begin{displaymath}
(x+\zeta y)=\zeta^v \times\eta\times\gamma^p ,\quad v\in \Z,\quad \eta\in\Z[\zeta+\zeta^{-1}]^*,
\quad \gamma\in \Z[\zeta].
\end{displaymath}
Then, by conjugation,
\begin{displaymath}
(x+\zeta^{-1} y)=\zeta^{-v} \times\overline{\eta}\times\overline{\gamma^p}.
\end{displaymath}
Then, observing that $\eta=\overline{\eta}$ we get, dividing these two relations
\begin{displaymath}
\frac{(x+\zeta y)}{(x+\zeta^{-1} y)}=\zeta^{2v}\frac{\gamma^p}{\overline{\gamma^p}}.
\end{displaymath}
Then, applying proposition \ref{p2} p.\pageref{p2}, we get $\gamma^p\equiv\overline{\gamma^p}\modu \pi^{p+1}$ and so
\begin{displaymath}
\frac{(x+\zeta y)}{(x+\zeta^{-1} y)}\equiv \zeta^{2v}\modu \pi^{p+1}.
\end{displaymath}
which leads to
\begin{displaymath}
x+\zeta y \equiv \zeta^{2v} x+\zeta^{2v-1} y\modu \pi^{p+1}
\end{displaymath}
Note, that this relation implies
\begin{displaymath}
x+\zeta y\equiv \zeta^{2v} x+ \zeta^{2v-1} y \modu p
\end{displaymath}
wich implies $2v\equiv 1 \modu p$
and then
\begin{displaymath}
(x-y)+\zeta (y-x)\equiv 0 \modu \pi^{p+1}
\end{displaymath}
so
\begin{displaymath}
(x-y)(1-\zeta)\equiv 0\modu \pi^{p+1},
\end{displaymath}
and finally $x-y\equiv 0 \modu p^2$, which contradicts hypothesis.
\item
Suppose that $\Q(\zeta)/\Q$ is not principal :
then, the proof of Cauchy's criterion needs only one Barlow-Abel congruence $\frac{x^p+y^p}{x+y}=t_1^p$ as it is possible to verify in classical proof, see for instance Ribenboim, \cite{rib} p 120-125.
\end{itemize}
\end{proof}

\end{thm}
%
\begin{thm}\label{t14082}
Let $p>5$ be a prime. Let $x,y,z\in\Z-\{0\}$ mutually coprime with $xy(x^2-y^2)\not\equiv 0 \modu p$.
\begin{itemize}
\item
If $\frac{x^p+y^p}{x+y}=t_1^p,\quad t_1\in\Z$, then $p$ verifies all the  criteria, congruences on $p$ or on $t=-\frac{x}{y}\modu p$, quoted {\bf (gfw)} in the table of results of the comparative survey of the previous section \ref{s06101} p.\pageref{s06101}
\item
 Let $c\in\Z-\{0\}$ where

$c=\pm\prod_{i=1}^m c_i^{\alpha_i}$ is the factorization of $c$ in primes. If  $c_i\not\equiv 1 \modu p,\quad i=1,\dots,m$ and
if $x^p+y^p+ c\times z^p=0$ then $p$ verifies all the  criteria, congruences on $p$ or on $t=-\frac{x}{y}\modu p$, quoted {\bf (gfw)} in the table of results of the comparative survey of the previous section \ref{s06101} p.\pageref{s06101}.
\end{itemize}
\begin{proof}
\begin{itemize} $ $
\item
Proof similar to theorem \ref{t14081} p.\pageref{t14081}.
\item
The proofs of all these criteria needs only to assume one Barlow-Abel congruence as hypothesis,
except all the case  with $\mu\geq 5$ which occurs in
corollary \ref{c21091} p.\pageref{c21091},
corollary \ref{c22091} p.\pageref{c22091},
corollary \ref{c09081} p.\pageref{c09081},
corollary \ref{c03104} p.\pageref{c03104},
corollary \ref{c2009}  p.\pageref{c2009},
corollary \ref{c05101} p.\pageref{c05101},
theorem   \ref{t05101} p.\pageref{t05101}
 and corollary \ref{c09101} p.\pageref{c09101},
where several Barlow-Abel equations are needed simultaneously for proof and cannot be used in this Fermat-Wiles Generalized equation.
\end{itemize}
\end{proof}
\end{thm}
%
{\bf Remark:} This result goes toward the {\bf Terjanian conjecture} which asserts that if
$xy(x^2-y^2)\not\equiv 0 \modu p$ and $\frac{x^p+y^p}{x+y}=t_1^p,\quad t_1\in\Z$, then the Kummer system of congruences of this Barlow relation has only the trivial solution, see \cite{ter}.
\end{subsection}
%

\clearpage
\section{FLT second case : miscelleanous }\label{s110022}
Contrary to a largely spread idea, it seems it is possible to use methods of this note to partial results on the second case of Fermat's Last Theorem. In this direction, we prove these, in our opinion,  new results, first on an elementary approach, then following some results in Washington, \cite{was}:
%
\subsection{ An elementary approach}
\begin{thm}{ *** } \label{t14102}
Let $p\in\N,\quad p>2$ be a prime.
Let $x,y,z\in\Z-\{0\},\quad p\parallel y $.
Then $x^p+y^p+z^p\not=0$.
\begin{proof} $ $
\begin{itemize}
\item
From hypothesis,  $(x+\zeta^u y )\not\equiv 0 \modu \pi$. Then all the prime ideals
$\mathbf q$ dividing $(x+\zeta^u y)$ are of degree 1, thus $N_{Q(\zeta)/Q}(\mathbf q)= q$ where $q\in\N$ is a prime with $q\equiv 1 \modu p$. Therefore, though we deal with second case of Fermat's Last Theorem, the relation (2.5) p 119 of \cite{rib}, and also proposition (3B) p 119 can be used.
In the same way the relation (4.7) p122 of \cite{rib} can be used to obtain
\begin{displaymath}
\prod_{i\in I(u,d)} (x+\zeta^{u_i} y)=\zeta^m \gamma^p,\quad m\in \N,\quad \gamma\in\Z[\zeta].
\end{displaymath}
\item
From $y\equiv 0 \modu p$ and
from $\prod_{i\in I(u,d)} (x+\zeta^{u_i} y)=\zeta^m\gamma^p$,
we obtain $x^{(p-1)/2}\equiv \zeta^m c \modu p, \quad c\in\Z$.
Then
$m=0$,
therefore
$\prod_{i\in I(u,d)} (x+\zeta^{u_i} y)=\gamma^p$,
so

$\prod_{i\in I(u,d)} ((x+y)+(\zeta^{u_i}-1)y)=\gamma^p$,
hence
\begin{displaymath}
(x+y)^{(p-1)/2}+(x+y)^{(p-3)/2}y\sum_{i\in I(u,d)} (\zeta^{u_i}-1)\equiv \gamma^p
\modu \pi^{2(p-1)}.
\end{displaymath}
\item
From Barlow-Abel formula, see \cite{rib}, formula (1.9) p 54, we have $x+y=t^p,\quad t\in \Z$. Let $t_1=t^{(p-1)/2}$. Then,
$(x+y)^{(p-1)/2}-\gamma^p=t_1^p-\gamma^p$. We have $t_1^p-\gamma^p\equiv 0 \modu \pi$, then  from proposition (\ref {p1}), we have
$t_1^p-\gamma^p\equiv 0 \modu \pi^{p+1}$,
therefore
$\sum_{i\in I(u,d)} (\zeta^{u_i}-1)x^{(p-3)/2}y\equiv 0 \modu \pi^{p+1}$
so
$\sum_{i\in I(u,d)} (\zeta^{u_i}-1)\equiv 0 \modu \pi^2$,

so
$\sum_{i\in I(u,d)} u_i(\zeta-1)\equiv 0 \modu \pi^2$,
by simplification
$\sum_{i\in I(u,d)} u_i\equiv 0 \modu \pi$,
hence
$\sum_{i\in I(u,d)} u_i\equiv 0 \modu p.$
\item
Let $d=1,2,\dots,p-2$.
Let $I(u,d)=
\{i|\quad 1\leq i \leq p-1,\quad u_{\nu-i}+u_{\nu-i+ind_u(d)}>p\}$.
For $i$ given and for $d$ varying from $1$ to $p-2$, the quantity
$u_{\nu-i+ind_u(d)}\modu p$ takes all integer values from $1$ to $p-1$ except the value $p-u_{\nu-i}$
which would correspond to $d=p-1$.
Therefore the quantity
$Q(i,d)= u_{\nu-i}+u_{\nu-i+ind_u(d)}$ takes all values
from $u_{\nu-i}+1$ to $u_{\nu-i}+p-1$
except $u_{\nu-i}+p-u_{\nu-i}=p$.
Therefore, when $d$ goes through $1$ to $p-2$,  the expression   $Q(i,d)$ takes all values from $p+1$ to $u_{\nu-i}+p-1$ with $Q(i,d)>p$, therefore $Q(i,d)$ takes $u_{\nu-i}-1$ different values with $Q(i,d)>p$.
\item
From $\sum_{i\in I(u,d)} u_i \equiv 0 \modu p$ for $d=1,\dots,p-2$, we have
$\sum_{d=1}^{p-2}\sum_{i\in I(u,d)} u_i\equiv 0\modu p$.
\item
In an other part, we obtain from previous  estimate of the cardinal number of $Q(i,d)>p$,
the congruence $\sum_{d=1}^{p-2}\sum_{i\in I(u,d)} u_i = \sum_{i=1}^{p-1} (u_{\nu-i}-1)u_i
\equiv (p-1)^2-\sum_{i=1}^{p-1} u_i\equiv 1\modu p$, a contradiction, which achieves the proof.
\end{itemize}
\end{proof}
\end{thm}
%
{\bf Remark :} Here, we have used the property $\alpha\equiv\beta,\quad \alpha\not\equiv 0 \modu \pi \Longrightarrow
\alpha^p\equiv \beta^p \modu \pi^{p+1}$ : it is important to observe  that the classical result
$\alpha^p\equiv \beta^p \modu \pi^{p-1}$ would not allow to conclude.
%

The next result is another proof of the same result interesting because it is strictly elementary in $\Z$.
%
\begin{thm}{ *** } \label{t14111}
Let $p\in\N,\quad p>3$ be a prime.
Let $x,y,z\in\Z-\{0\},\quad p\parallel y $.
Then $x^p+y^p+z^p\not=0$.
\begin{proof}
We shall show that if $x^p+y^p+z^p=0$ and $p|y$ then $p^2|y$.
\begin{itemize}
\item
From Barlow-Abel, see (1.9) p 54, we have $\frac{x^p+y^p}{x+y}=t_1^p$ with $t_1\in \N$. From little Fermat theorem, $x^p+y^p\equiv x+y \modu p$, so $t_1^p\equiv 1\modu p$, hence $t_1^p\equiv 1 \modu p^2$ and
\begin{equation}\label{e14111}
x^p+y^p \equiv x+y \modu p^2.
\end{equation}
\item
From Barlow-Abel, we have $\frac{x^p+z^p}{x+z}=p s_1^p$ where $s_1\in \N$. All the prime $q$ dividing $s_1$ verify $q\equiv 1 \modu p$,
so $s_1\equiv 1\modu p$ and  $s_1^p\equiv 1 \modu  p^2$ which leads to
\begin{equation}
\frac{x^p+z^p}{x+z}\equiv p\modu p^3.
\end{equation}
\item
We have $\frac{x^p+z^p}{x+z}=x^{p-1}-x^{p-2}z+\dots-xz^{p-2}+z^{p-1}$.
From $y|p$, we obtain $x+z\equiv 0 \modu p^{p-1}$ and thus $x+z\equiv 0 \modu p^3$ from hypothesis $p>3$. Then we obtain $x^{p-1}-x^{p-2}(-x)+\dots-x(-x)^{p-2}+(-x)^{p-1}\equiv p\modu p^3$, so $px^{p-1}\equiv p\modu p^3$, hence
\begin{equation}\label{e14112}
x^{p-1}\equiv 1 \modu p^2.
\end{equation}
\item
From relations (\ref{e14111}) p.\pageref{e14111} and (\ref{e14112}) p.\pageref{e14112}, we obtain $y\equiv 0 \modu p^2$, which completes the proof.
\end{itemize}
\end{proof}
\end{thm}
{\bf Remark:}
With an application of Furtwangler theorem in class field theory, Vandiver \cite{va7}, see also Ribenboim \cite{rib} (3E) p. 170, proved the  better result: Let $p\in\N,\quad p>3$ be a prime.
Let $x,y,z\in\Z-\{0\},\quad p^2\parallel y $.
Then $x^p+y^p+z^p\not=0$.

%
\subsection{ Cyclotomic number fields approach}
We apply similar technics to those of  Washington,\cite{was}, chapter 9 p 167,  for  FLT second case.
\begin{itemize}
\item
Let $x,y,z$ be a solution of second case with $y\equiv 0 \modu p$. Suppose that  $h^+\not\equiv 0 \modu p$.
We have
$\frac{x+\zeta y}{x+\zeta^{-1} y}=1+\frac{(\zeta-\zeta^{-1})y}{x+\zeta^{-1}y}$ and so
$\frac{x+\zeta y}{x+\zeta^{-1} y}\equiv 1 \modu \pi^p$. Then, we are in the situation to apply Washington \cite{was}, lemma 9.1 p 169 and lemma 9.2 p 170, to get
\begin{displaymath}
\frac{x+\zeta y}{x+\zeta^{-1}y} = \gamma^p,\quad \gamma \in \Q(\zeta).
\end{displaymath}
Therefore, from proposition \ref{p2} p.\pageref{p2},
$\frac{(\zeta-\zeta^{-1})y}{x+\zeta^{-1}y}=\gamma^p-1\equiv 0\modu \pi^{p+1}$, which implies that $y\equiv 0 \modu p^2$, result that we shall prove also, without assuming Vandiver conjecture, in theorems \ref{t14102} p.\pageref{t14102} and \ref{t14111} p.\pageref{t14111}.
\item
We derive also
\begin{displaymath}
\frac{(x+\zeta y)^2}{(x+\zeta y)(x+\zeta^{-1}y)}=\gamma^p,\quad\gamma\in\Q(\zeta).
\end{displaymath}
From $h^+\not\equiv 0\modu p$ we get
\begin{displaymath}
(x+\zeta y)(x+\zeta^{-1} y)=\eta_1\times \alpha^p,
\quad\alpha\in\Z[\zeta+\zeta^{-1}].
\end{displaymath}
and so
\begin{displaymath}
(x+\zeta y)^2=\eta_1\times\gamma^p\times\alpha^p,
\end{displaymath}
so, raising to $\frac{p+1}{2}$ power,
\begin{displaymath}
(x+\zeta y)(x+\zeta y)^{p}=\eta_1^{(p+1)/2}\times\gamma^{p(p+1)/2}\times \alpha^{p(p+1)/2},
\end{displaymath}
also
\begin{displaymath}
(x+\zeta y) =\eta_1^{(p+1)/2}\times (\frac{\gamma^{(p+1)/2}\times\alpha^{(p+1)/2}}{(x+\zeta y)})^p,
\end{displaymath}
and finally
\begin{displaymath}
x+\zeta y =\eta \times \beta^p,\quad \beta\in\Z[\zeta],\quad \eta\in\Z[\zeta+\zeta^{-1}]^*.
\end{displaymath}
\item
We have seen that $y\equiv 0 \modu p^2$ and, classically, $x+y+z\equiv 0 \modu p^2$, see for instance Ribenboim \cite{rib} (3B)-2 p 58,
therefore $x+z\equiv 0 \modu p^2$.
Thus
\begin{displaymath}
\frac{x+\zeta z}{x+\zeta^{-1} z}\equiv \frac{x(1-\zeta)}{x(1-\zeta^{-1})} \modu p^2,
\end{displaymath}
which leads to
\begin{displaymath}
-\zeta^{-1}\frac{x+\zeta z}{x+\zeta^{-1} z}\equiv 1 \modu \pi^{2(p-1)}.
\end{displaymath}
Let $\omega= -\zeta^{-1}\times\frac{x+\zeta z}{x+\zeta^{-1} z}$. Then
$\overline{\omega}=-\zeta\times\frac{x+\zeta^{-1}z}{x+\zeta z}$.
Therefore $\omega\times\overline{\omega}=1$ and $\omega\equiv 1 \modu \pi^{2(p-1)}$, and
also $\omega\Z[\zeta]=\mathbf b^p$, where $\mathbf b$ is a fractional ideal of $\Q(\zeta)$, because $(x+\zeta z)=\pi \times\s_y^p$ where $\s_y$ is an integral ideal of $\Z[\zeta]$. Thus, we are again in situation to apply  Washington lemma 9.1 p 169 and lemma 9.2 p 170, to get
\begin{displaymath}
-\zeta^{-1} \times\frac{x+\zeta z}{x+\zeta^{-1} z} =\gamma^p,\quad \gamma\in\Q(\zeta).
\end{displaymath}
\item
We get
\begin{equation}\label{e012092}
\zeta^{-1}\times (\frac{(x+\zeta z)^2}{(x+\zeta z)(x+\zeta^{-1} z)}=\gamma_y^p,
\quad \gamma_y\in\Q(\zeta).
\end{equation}
Let us note $\lambda=\zeta-1$ and $\mu=(\zeta-1)(\zeta^{-1}-1)$. From $h^+\not\equiv 0\modu p$, we derive
\begin{displaymath}
(x+\zeta z)(x+\zeta^{-1} z)=\mu\times\eta_1\times\alpha^p,
\quad \eta_1\in\Z[\zeta+\zeta^{-1}]^*,\quad \alpha\in\Z[\zeta+\zeta^{-1}].
\end{displaymath}
We have also
\begin{displaymath}
x+\zeta z=\lambda\times \delta,\quad \delta\in\Z[\zeta].
\end{displaymath}
Raising relation (\ref{e012092}) p.\pageref{e012092} to the $\frac{p+1}{2}$ power, we obtain
\begin{displaymath}
(x+\zeta z)(x+\zeta z)^p=-\zeta^{(p+1)/2}\times \mu^{(p+1)/2}\times
\eta_1^{(p+1)/2}\times \alpha^{p(p+1)/2}\times\gamma_y^{p(p+1)/2}
\end{displaymath}
and so
\begin{displaymath}
(x+\zeta z)=-\zeta^{(p+1)/2}\times \frac{\mu^{(p+1)/2}}{\lambda^p}\times
\eta_1^{(p+1)/2}\times (\frac{\alpha^{(p+1)/2}\times\gamma_y^{(p+1)/2}}{\delta})^p
\end{displaymath}
which gives
\begin{displaymath}
(x+\zeta z)=-\zeta^{(p+1)/2}\times \frac{\mu^{(p+1)/2}}{\lambda^p}\times
\eta\times \beta^p,\quad \beta\in\Z[\zeta],\quad \eta\in\Z[\zeta+\zeta^{-1}],
\end{displaymath}
so
\begin{displaymath}
(x+\zeta z)=-\zeta^{(p+1)/2}\times\frac{(\zeta-1)^{(p+1)/2}\times\zeta^{-(p+1)/2}
\times (1-\zeta)^{(p+1)/2}\times \eta\times\beta^p}{(\zeta-1)^p},
\end{displaymath}
and finally
\begin{displaymath}
(x+\zeta z)=(-1)^{(p-1)/2}\times (\zeta-1)\times\eta\times\beta^p,\quad \beta\in\Z[\zeta],\quad
\eta\in\Z[\zeta+\zeta^{-1}]^*.
\end{displaymath}
and, by symmetry on $x,z$,
\begin{displaymath}
\begin{split}
&(x+\zeta z)=(-1)^{(p-1)/2}\times (\zeta-1)\times\eta_{y,1}\times\beta_{y,1}^p,\quad \beta_{y,1}\in\Z[\zeta],\quad
\eta_{y,1}\in\Z\zeta+\zeta^{-1}]^*,\\
&(z+\zeta x)=(-1)^{(p-1)/2}\times (\zeta-1)\times\eta_{y,2}\times\beta_{y,2}^p,\quad \beta_{y,2}\in\Z[\zeta],\quad
\eta_{y,2}\in\Z\zeta+\zeta^{-1}]^*.
\end{split}
\end{displaymath}
\item
We get the two relations
\begin{displaymath}
\begin{split}
& \zeta\times (x+\zeta^{-1}z)
=(-1)^{(p-1)/2}\times (\zeta-1)\times \eta_{y,2}\times \beta_{y,2}^p,\\
& x+\zeta^{-1} z=(-1)^{(p-1)/2}\times(\zeta^{-1}-1)\times \eta_{y,1}
\times (\overline{\beta_{y,1})^{p}}),
\end{split}
\end{displaymath}
and finally
\begin{displaymath}
\frac {\eta_{y,2}}{\eta_{y,1}}=(-\frac{\overline{\beta_{y,1}}}{\beta_{y,2}})^p.
\end{displaymath}
\item
Finally we summarize all that in the set of simultaneous relations
%
\begin{thm} \label{t012091}
Suppose that $h^+\not\equiv 0 \modu p$
and $x^p+y^p+z^p=0,\quad x,y,z\in\Z-\{0\},\quad y\equiv 0 \modu p$;
then we have the set of relations
\begin{equation}\label{e01203}
\begin{split}
& y\equiv 0 \modu p^3,\\
& \frac{x+\zeta y}{x+\zeta^{-1} y}=\gamma_x^p,\quad \gamma_x\in\Q(\zeta),\\
& \frac{z+\zeta y}{z+\zeta^{-1} y}=\gamma_z^p,\quad \gamma_z\in\Q(\zeta),\\
& -\zeta^{-1} \times\frac{x+\zeta z}{x+\zeta^{-1} z}
=\gamma_y^p,\quad \gamma_y\in\Q(\zeta),\\
& (x+\zeta y)=\eta_x \times\beta_x^p,\quad \beta_x\in\Z[\zeta],
\quad \eta_x\in\Z[\zeta+\zeta^{-1}]^*,\\
& (z+\zeta y)=\eta_z \times\beta_z^p,\quad \beta_z\in\Z[\zeta],
\quad \eta_z\in\Z[\zeta+\zeta^{-1}]^*,\\
&(x+\zeta z)=(-1)^{(p-1)/2}\times (\zeta-1)\times\eta_{y,1}\times\beta_{y,1}^p,\quad \beta_{y,1}\in\Z[\zeta],\quad
\eta_{y,1}\in\Z\zeta+\zeta^{-1}]^*,\\
&(z+\zeta x)=(-1)^{(p-1)/2}\times (\zeta-1)\times\eta_{y,2}\times\beta_{y,2}^p,\quad \beta_{y,2}\in\Z[\zeta],\quad
\eta_{y,2}\in\Z\zeta+\zeta^{-1}]^*,\\
&\frac {\eta_{y,2}}{\eta_{y,1}}=(-\frac{\overline{\beta_{y,1}}}{\beta_{y,2}})^p.
\end{split}
\end{equation}
\end{thm}
\end{itemize}
%
%
\subsection{ Hilbert class field and representations approach}
Recall at first some definitions and properties obtained in representation and class field theory:
\begin{itemize}
\item
$C_p$ is the $p$-class group of $\Q(\zeta)$.
\item
$C_p^1$ is the subgroup of $C_p$ whose elements are of order $1$ or $p$. $r_p$ is the rank  of $C_p$ and $C_p^1$ is of order $p^{r_p}$.
\item
$r_p^-$ is the $p$-rank of the relative class group $C_p^-$ and $r_p^+$ is the $p$-rank of   $C_p^+$, class group of the field $\Q(\zeta+\zeta^{-1})$.
\item
There exists  a class of  ideals $\mathbf b$ of $\Z[\zeta]$ such that:
\begin{displaymath}
\begin{split}
&\mathbf b\simeq \mathbf b_1\dots\mathbf b_{r_p},\\
&\mathbf b\not\simeq \Z[\zeta],\quad \mathbf b^p\simeq \Z[\zeta],\\
&\mathbf b_i\not\simeq \Z[\zeta],\quad \mathbf b_i^p\simeq \Z[\zeta],
\quad i=1,\dots,r_p,\\
& \sigma(\mathbf b_i)\simeq \mathbf b_i^{\mu_i},
\quad \mu_i\in\N,\quad 1<\mu_i\leq p-1,\quad i=1,\dots,r_p,\\
& C_p^1=\oplus_{i=1}^{r_p} <Cl(\mathbf b_i)>.
\end{split}
\end{displaymath}
\item
There exists $C_i\in\Q(\zeta),\quad i=1,\dots,r_p^-$, such that
\begin{displaymath}
\begin{split}
& C_i\Z[\zeta]\simeq \mathbf b_i^{p}, \quad Cl(\mathbf b_i)\in C_p^-,\\
& C_i\equiv 1\modu \pi^{2m_i+1},\quad 1\leq m_i\leq\frac{p-3}{2}.
\end{split}
\end{displaymath}
\item
From Hilbert class field structure theorem \ref{t201013} p. \pageref{t201013},
$C_i\in\Q(\zeta)-\Z[\zeta]^*,\quad i=1,\dots,r_p^+$, are singular primary and correspond to Hecke component of the relative class group $C_p^-$ and $C_i,\quad i=r_p^++1,\dots,r_p^-$, are {\bf not} singular primary.
\end{itemize}
%
This theorem gives the general form of the class field theory structure of the Fermat's equation in the FLT second case:  $z\equiv 0\modu p$.
\begin{thm}{ *** }\label{t202051}
Suppose that $x^p+y^p+z^p=0$ with $z\equiv 0\modu p,\quad xy\not\equiv 0\modu p$. Then $x, y$ verify the relation
\begin{displaymath}
\begin{split}
& x+y\equiv 0\modu p^2,\\
& (\frac{x+\zeta y}{1-\zeta})/(\frac{x+\zeta^{-1} y}{1-\zeta^{-1}})
=  C_1^{\nu_1}\times\dots\times C_{r_p^+}^{\nu_{r_p^+}}\times
(\frac{\gamma_1}{\overline{\gamma}_1})^p,\quad
 \gamma_1\in\Q(\zeta),\\
& C_i\Z[\zeta]=\mathbf b_i^{p},\quad Cl(\mathbf b_i)\in C_p^-,\quad
\mbox{Hecke components},\quad i=1,\dots,r_p^+,\\
& C_i\equiv c_i^p\modu \pi^p,\quad i=1,\dots,r_p^+.
\end{split}
\end{displaymath}
\begin{proof}
Let $x^p+y^p+z^p=0$ with $z\equiv 0\modu p,\quad xy\not\equiv 0\modu p$.
We have proved in theorem \ref{t14111} p. \pageref{t14111} that $z\equiv 0\modu p^2$.
From Ribenboim \cite{rib} (3B).2 p 58, $x+y+z\equiv 0\modu p^2$, and so
$x+y\equiv 0\modu p^2$.
Therefore
\begin{displaymath}
\frac{x+\zeta y}{1-\zeta}\equiv x\modu p^2.
\end{displaymath}
By conjugation, we get
\begin{displaymath}
\frac{x+\zeta^{-1} y}{1-\zeta^{-1}}\equiv x\modu p^2.
\end{displaymath}
From these two relations we get
\begin{displaymath}
\frac{\frac{x+\zeta y}{1-\zeta}}{\frac{x+\zeta^{-1} y}{1-\zeta^{-1}}}
\equiv 1\modu p^2.
\end{displaymath}
From factorization of ideals in Dedekind ring $\Z[\zeta]$ we get
\begin{displaymath}
\frac{x+\zeta y}{1-\zeta}\Z[\zeta]=\s^p,
\end{displaymath}
where $\s$ is an integral ideal of $\Z[\zeta]$.
Let us consider the relative class group $C_p^-$.
Let $C_i,\quad i=1,\dots,r_p^-$, defined in lemma \ref{l108171} p. \pageref{l108171}.
Strictly similarly to relation obtained in FLT first case, we can write in the second case,
\begin{displaymath}
\begin{split}
& \frac{\frac{x+\zeta y}{1-\zeta}}{\frac{x+\zeta^{-1} y}{1-\zeta^{-1}}}
= \zeta^v\times C_1^{\nu_1}\times\dots\times  C_{r_p^-}^{\nu_{r_p^-}}\times \gamma_1^p,\\
& \gamma_1\in\Q(\zeta),\quad C_i\equiv 1\modu\pi^{2m_i+1},
\quad 1\leq m_i \leq\frac{p-3}{2},
\quad i=1,\dots,r_p^-\\
& v\in\N,\quad \nu_i\in \N,\quad i=1,\dots,r_p^-,
\end{split}
\end{displaymath}
where $C_i,\quad i=1,\dots, r_p^+$, are singular primary and
$C_i,\quad i=r_p^++1,\dots, r_p^-$, are not singular primary (see Hilbert class field structure theorem \ref{t201013} p. \pageref{t201013}).
These relations imply that $v\equiv 0\modu p$.
The number $(\frac{x+\zeta y}{1-\zeta})/(\frac{x+\zeta^{-1}y}{1-\zeta^{-1}})$ is singular primary, so
from theorem \ref{t201061} p. \pageref{t201061}, application of class field theory, we get
$\nu_{r_p^++1}\equiv\dots\equiv \nu_{r_p^-}\equiv 0\modu p$ and finally we get
\begin{displaymath}
(\frac{x+\zeta y}{1-\zeta})/(\frac{x+\zeta^{-1} y}{1-\zeta^{-1}})
=  C_1^{\nu_1}\times\dots\times C_{r_p^+}^{\nu_{r_p^+}}\times \gamma_2^p,\quad
 \gamma_2\in\Q(\zeta).
\end{displaymath}
\end{proof}
\end{thm}
%
\begin{cor}\label{c202031}
Suppose that $h^+\not\equiv 0\modu p$.
Suppose that $x^p+y^p+z^p=0$ with $z\equiv 0\modu p,\quad xy\not\equiv 0\modu p$. Then $x, y$ verify the relation
\begin{displaymath}
\begin{split}
& x+y\equiv 0\modu p^2,\\
& \frac{\frac{x+\zeta y}{1-\zeta}}{\frac{x+\zeta^{-1} y}{1-\zeta^{-1}}}
=  (\frac{\gamma_1}{\overline{\gamma}_1})^p,\quad
 \gamma_1\in\Z(\zeta),\\
& \frac{x+\zeta y}{1-\zeta}=  \eta\times\gamma_1^p,\\
& \eta\equiv d^p\modu \pi^p,\quad d\in\Z.
\end{split}
\end{displaymath}
\begin{proof}
Apply theorem \ref{t202051} p. \pageref{t202051} with the relation deduced of
$h^+\not\equiv 0\modu p$:
\begin{displaymath}
\begin{split}
& (\frac{x+\zeta y}{1-\zeta})\times(\frac{x+\zeta^{-1}y}{1-\zeta^{-1}})=
\eta^2\times\gamma_1^p\times\overline{\gamma}_1^p,\\
& \eta\in\Z[\zeta+\zeta^{-1}]^*,
\end{split}
\end{displaymath}
to get
\begin{displaymath}
\frac{x+\zeta y}{1-\zeta}=\eta\times \gamma_1^p,\quad \eta\in\Z[\zeta+\zeta^{-1}]^*.
\end{displaymath}
From $\frac{x+\zeta y}{1-\zeta}\equiv 1\modu p^2$ we derive
that $\eta\equiv d^p\modu \pi^p,\quad d\in\Z,\quad d\not\equiv 0\modu p$.
\end{proof}
\end{cor}
{\bf Remark:} observe that this corollary is consistent with Washington formulation \cite{was} in page 171
and with theorem \ref{t012091} p. \pageref{t012091} of this monograph.
%
%
The next result is an independant application of several results given in Washington,\cite{was}: the validation of this theorem implies to follow the proof with the Washington book. The notations are those of Washington.
%
\begin{thm}\label{t202052}
Let $I$ be the set of irregular indices.
Let $E_{2i},\quad i=1,\dots,\frac{p-3}{2}$, be the units defined in Washington \cite{was} p 155.
If $h^+\not\equiv 0\modu p$ then
\begin{displaymath}
\frac{x+\zeta y}{1-\zeta}=(\prod_{i\in I} E_{2i}^{\alpha_{j_i}})\times \gamma_2^p,
\quad \gamma_2\in\Z[\zeta].
\end{displaymath}
\begin{proof}
Suppose that $x^p+y^p+z^p=0,\quad z\equiv 0\modu p,\quad xy\not\equiv 0\modu p$.
We start of the relation
\begin{displaymath}
\frac{x+\zeta y}{1-\zeta}=\eta\times\gamma^p,\quad
\eta\in\Z[\zeta+\zeta^{-1}]^*,\quad \gamma\in\Z[\zeta+\zeta^{-1}].
\end{displaymath}
Let us consider the $E_{2i},\quad i=1,\dots,\frac{p-1}{2}$, defined in relation
\ref{e201011} p. \pageref{e201011} and also in Washington \cite{was} p. 155:
\begin{displaymath}
\begin{split}
& E_{2i}=
\prod_{a=1}^{p-1}
(\zeta^{(1-u)/2}\frac{1-\zeta^{au}}{1-\zeta^a})^{a^{2i}\tau_a^{-1}},\\
& \tau_a(\zeta)=\zeta^a,\quad a=1,\dots,p-1.
\end{split}
\end{displaymath}
From Washington \cite{was}, corollary 8.15 p 156, the
$E_{2i},\quad i=1,\dots,\frac{p-3}{2}$, generate the group
$\Z[\zeta+\zeta^{-1}]^*/(\Z[\zeta+\zeta^{-1}]^*)^p$.
Therefore we can write from corollary \ref{c202031} p. \pageref{c202031}
\begin{displaymath}
\frac{x+\zeta y}{1-\zeta}=(\prod_{i=1}^{(p-3)/2} E_{2i}^{\alpha_i})\times\gamma_1^p.
\end{displaymath}
Let $B_{2i},\quad i=1,\dots,\frac{p-3}{2}$, be the even Bernoulli Numbers.
Let $\rho$ be the index of iregularity of $p$. Let $i_j,\quad i=1,\dots,\rho$, be the irregular indices, so of the $B_{2i}\equiv0\modu p$.
From Washington \cite{was} exercise 8.9 p 165, the extension $H/\Q(\zeta)$
defined by
\begin{displaymath}
H=\Q(\zeta,E_{2i_1}^{1/p},\dots,E_{2i_{\rho}}^{1/p})/\Q(\zeta)
\end{displaymath}
is unramified.
Therefore, the extensions $\Q(\zeta,E_{2i_j}^{1/p})/\Q(\zeta),
\quad i=1,\dots,\rho$, are unramified and so
$E_{2i_j},\quad j=1,\dots,\rho$, are  primary units
and so
\begin{displaymath}
E_{i_j}\equiv c_{i_j}^p\modu\pi^p,\quad c_{i_j}\in\N,\quad  j=1,\dots, \rho.
\end{displaymath}
Let $I=\{i_1,\dots i_\rho\}$ be the set of irregular indices.
Let $J=\{i \ |\ i=1,\dots,\frac{p-3}{2},\quad i\not\in I\}$.
We deduce that
\begin{displaymath}
\frac{x+\zeta y}{1-\zeta}\equiv c_1\times \prod_{i\in J} E_{2i}^{\alpha_i}\modu \pi^p,
\quad c_1\in\Z,\quad c_1\not\equiv 0\modu p.
\end{displaymath}
We have seen in theorem \ref{t14111} p. \pageref{t14111} that $x+y\equiv 0\modu p^2$, so
\begin{displaymath}
\prod_{i\in J} E_{2i}^{\alpha_i}\equiv c_2
\modu \pi^p, \quad c_2\in\Z,\quad c_2\not\equiv 0\modu p.
\end{displaymath}
Let $\beta=(\zeta-1)(\zeta^{-1}-1)$.
From Washington \cite{was}, exercise 8.11 p 166, we get
\begin{displaymath}
\begin{split}
& E_{2i}\equiv a_i+b_i\times \beta^i\modu \pi^{2(i+1)},\quad i\in J,\\
& a_i\in \Z,\quad b_i\in\Z,\quad a_i\not\equiv 0\modu p,\quad b_i\not\equiv 0\modu p.
\end{split}
\end{displaymath}
Therefore
\begin{displaymath}
\prod_{i\in J} (a_i+b_i\times \beta^i)^{\alpha_i}\equiv c_2
\modu \pi^p, \quad c_2\in\Z,\quad c_2\not\equiv 0\modu p.
\end{displaymath}
Suppose that there exists a smallest  $i_1$ with $i\in J$ and with $\alpha_{i_1}\not\equiv 0\modu p$.
We deduce that
\begin{displaymath}
(a_{i_1}+\beta^i\times b_{i_1})^{\alpha_{i_1}}\equiv c\modu \pi^{2(i_1+1)}
\end{displaymath}
and so
\begin{displaymath}
(a_{i_1}+\alpha_{i_1}\times\beta^i\times b_{i_1})\equiv c\modu \pi^{2(i_1+1)}
\end{displaymath}
which implies that $\alpha_{i_1}\equiv 0\modu p$, contradiction.
Therefore
$\alpha_i\equiv 0\modu p$ for $i\in J$
and finally
\begin{displaymath}
\frac{x+\zeta y}{1-\zeta}=(\prod_{i\in I} E_{2i}^{\alpha_{j_i}})\times \gamma_2^p,
\quad \gamma_2\in\Z[\zeta].
\end{displaymath}
\end{proof}
\end{thm}
{\bf Remark:} This formulation is consistent with structure theorem on Hilbert-class field of $\Q(\zeta)$ \ref{t201013} p. \pageref{t201013}. Here $r_p^+=0$.
%
%
%
\clearpage
{\bf Acknowledgements}
I would like to thank Professors Guy Terjanian in Toulouse and Bruno Angl\`es in Caen, Yves Hellegouarch in Caen, Antoine Chambert-Loir in Paris, Vincent Fleckinger in Besanon,  Franz Lemmermeyer in Heidelberg and San Marcos, Gerhard Niklasch in Munchen and Preda Mihailescu in Paderborn Germany, to give me helpful advice during the writing of this article and to point out to me a very large number of errors in previous versions,  and also to Professor Ladislav Skula of Mazaryk University, Brno, Czech Republic for advice and help on Mirimanoff and Kummer congruences bibliography.
%

Roland Qu\^eme

13 avenue du ch\^ateau d'eau

31490 Brax

France

2002 october 31

mailto: roland.queme@free.fr

home page: http://roland.queme.free.fr/index.html
\end{document}